\documentclass[11pt]{article}
\usepackage{amsmath,amsbsy,amsfonts}
\usepackage[frenchb]{babel}
\begin{document}

\title{Sur le transfert des traces d'un groupe classique p-adique \`a un groupe lin\'eaire tordu.}
\author{C. M{\oe}glin et J.-L. Waldspurger\\ Institut de Math\'ematiques de Jussieu
\\ CNRS
}
\date{ }
\maketitle

Soit $F$ un corps local non archim\'edien et soit $N$ un entier. On note $\theta$ l'automorphisme de $GL(N,F)$ d\'efini pour tout \'el\'ement $g$ de $GL(N,F)$ par $\theta(g)=J\, ^tg^{-1}J$ o\`u $J$ est la matrice antidiagonale multipli\'e \`a gauche par la matrice diagonale $(1,-1, \cdots, (-1)^{N-2},(-1)^{N-1})$. On note $\tilde{G}(N)$ le produit semi-direct de $GL(N,F)$ par le groupe \`a 2 \'el\'ements engendr\'e par $\theta$. Pour justifier notre travail, admettons momentan\'ement que la conjecture de Ramanujan soit v\'erifi\'ee pour les repr\'esentations automorphes cuspidales de $GL(N,{\mathbb A})$ o\`u ${\mathbb A}$ est l'anneau des ad\`eles d'un corps global admettant $F$ comme localis\'e en une place. Les composantes locales en cette place des formes automorphes de carr\'e int\'egrable (et m\^eme des s\'eries d'Eisenstein sur l'axe unitaire) sont des induites irr\'eductibles de la forme $\sigma:=\times_{(\rho,a,b)}Sp(b,St(a,\rho))$ o\`u $(\rho,a,b)$ parcourt un ensemble de triplets o\`u $\rho$ est une repr\'esentation cuspidale irr\'eductible et unitaire d'un groupe $GL(d_{\rho},F)$ (ce qui d\'efiniti $d_{\rho}$) et $a,b$ sont des entiers. La notation $Sp(b,St(a,\rho))$ est importante pour nous: pour la d\'efinir on utilise des induites ce qui suppose que nous avons fix\'e un parabolique minimal, pour nous les matrices triangulaires sup\'erieures (pour fixer les id\'ees), alors $St(a,\rho)$ est la repr\'esentation de Steinberg de $GL(ad_{\rho},F)$ bas\'ee sur $\rho$, c'est-\`a-dire l'unique sous-module irr\'eductible de l'induite $\rho\vert\,\vert^{(a-1)/2}\times \cdots \times \rho\vert\,\vert^{-(a-1)/2}$ et $Sp(b,St(a,\rho))$ est l'unique quotient irr\'eductible de l'induite pour $GL(abd_{\rho},F)$:
$$
St(a,\rho)\vert\,\vert^{(b-1)/2}\times \cdots \times St(a,\rho)\vert\,\vert^{-(b-1)/2}.
$$C'est aussi l'unique sous-module irr\'eductible de l'induite
$$
St(a,\rho)\vert\,\vert^{-(b-1)/2}\times \cdots \times St(a,\rho)\vert\,\vert^{(b-1)/2}.
$$
On remarque aussi que $\rho$ est $\theta$-invariante si et seulement si $\rho$ est autoduale.
 Si l'ensemble des repr\'esentations $\rho$ qui interviennent dans ce produit est form\'e de repr\'esentations  dont la classe d'isomorphie est invariante sous l'action de $\theta$ (on dira repr\'esentations$\theta$-invariantes), l'induite $\sigma$ est elle aussi $\theta$-invariante et on peut donc la prolonger en une action de $\tilde{G}(N)$; il y a \'evidemment un choix, qu'il faudra pr\'eciser. D'un point de vue global, un choix pour l'action de $\theta$ s'impose assez naturellement: celui qui vient de la classification de Langlands des repr\'esentations. C'est  celui qui est sugg\'er\'e par Arthur. D'un point de vue local, ce n'est pas le plus simple et il n'a rien de canonique. Dans ce travail, on \'etudiera les diff\'erents choix possible, on calculera les signes qui interviennent quand on change de choix et on montrera l'influence de ces signes sur le transfert.

Mais le but principal de ce travail est de ramener le calcul de la trace $tr \sigma(g,\theta)$ pour tout $g\in GL(N,F)$ au calcul de la trace pour des repr\'esentations temp\'er\'ees. Pour cela, nous n'avons besoin d'aucune hypoth\`ese.

Ce cas a \'et\'e trait\'e  sous certaines hypoth\`eses (en particulier un lemme fondamental) en \cite{transferttempere} et sous ces m\^emes hypoth\`eses on obtient donc un transfert assez g\'en\'eral. Le but ultime est de d\'ecrire les paquets d'Arthur pour les groupes classiques; ici on r\`egle donc la combinatoire li\'ee \`a la distribution stable associ\'ee au paquet. Il restera \`a r\'egler la combinatoire li\'ee aux propri\'et\'es endoscopiques, et on voit clairement comment faire. Soyons plus pr\'ecis.

\

Arthur \cite{arthur} \cite{arthurrecent} a remarqu\'e que si $N$ est pair le groupe $H_{N}:=SO(N+1,F)$ (forme quasid\'eploy\'ee) est un groupe endoscopique ''stable'' pour la forme tordue $\tilde{G}(N)$ alors que si $N$ est impair, la m\^eme assertion est vraie pour $H_{N}:=Sp(2N,F)$. Ainsi si $tr(\sigma (g,\theta))$ est stable, cette distribution doit \^etre le tranfert d'une distribution stable sur l'un de ces groupes; quand cette distribution n'est pas stable d'atures groupes endoscopiques elliptiques remplace $H_{N}$ (cf \cite{arthurrecent} par. 30. L'id\'ee serait de calculer $tr(\sigma(g,\theta))$ dans le groupe de Grothendieck des repr\'esentations lisses de $\tilde{G}(N)$ avec comme base des induites de repr\'esentations temp\'er\'ees, en fait ce n'est pas tout \`a fait comme cela que l'on proc\`ede. On interpr\'etera la repr\'esentation $\sigma$ \`a l'aide des morphismes de $W_{F}\times SL(2,{\mathbb C}) \times SL(2,{\mathbb C})$ dans $GL(N,{\mathbb C})$ dont la restriction \`a $W_{F}$ est continue et born\'ee. Et les repr\'esentations temp\'er\'ees qui doivent alors intervenir sont celles qui sont associ\'ees au morphisme $\psi\circ \Delta$ o\`u $\Delta $ est le plongement diagonal de $SL(2,{\mathbb C})$ pour tous les sous-groupes de Levi de $GL(N,F)$ tel que le groupe dual contienne un conjugu\'e de $\psi\circ \Delta$. On oublie la notation $\sigma$ et on note $\pi(\psi)$ la repr\'esentation associ\'ee \`a un tel morphisme $\psi$. On veut \'etudier les $tr\pi(\psi)(g,\theta)$ quand l'image de $\psi$ est simplement $\theta$-invariante.

Le cas difficile est celui o\`u $\psi\circ \Delta$ est discret c'est-\`a-dire d\'efini une repr\'esentation semi-simple sans multiplicit\'es de $W_{F}\times SL(2,{\mathbb C})$; on dit alors que $\psi$ est de restriction discr\`ete \`a la diagonale. On  ram\`ene le cas g\'en\'eral (cf \ref{general})
 \`a ce cas, en utilisant le Jacquet ce qui est une op\'eration qui n\'ecessite de commencer par augmenter le rang du groupe mais cette op\'eration est tr\`es compatible au transfert (cf. \ref{proprietesgeneralesdutransfertstable}). 
 
  Supposons donc que $\psi\circ \Delta$ soit discret; dans ce cas, le centralisateur de $\psi$ et celui de $\psi\circ \Delta$ sont des produits de groupes $GL(1,{\mathbb C})$ et celui de $\psi$ est inclus dans celui de $\psi\circ \Delta$. On dit que $\psi$ est \'el\'ementaire si cette inclusion est une bijection; en termes simples cela veut dire que $\psi$ en tant que repr\'esentation de $W_{F}\times SL(2,{\mathbb C}) \times SL(2,{\mathbb C})$ est somme de repr\'esentations irr\'eductibles dont chacune est triviale sur au moins l'une des copies de $SL(2,{\mathbb C})$; la repr\'esentation de $GL(n,F)$ alors associ\'ee \`a $\psi$ est une induite irr\'eductible de repr\'esentations de Steinberg et de repr\'esentations de Speh (les repr\'esentations de Speh  sont les images par l'involution de Zelevinsky des repr\'esentations de  Steinberg). Ce cas se traite en r\'ealisant la repr\'esentation associ\'ee \`a $\psi$ comme cohomologie d'un complexe provenant de la repr\'esentation temp\'er\'ee associ\'ee \`a $\psi\circ \Delta$; c'est une r\'eminiscence des constructions d'Aubert \cite{aubert}. Ce complexe est muni d'une action de $\theta$ et il faut la comparer \`a l'action de $\theta$ mise  a priori. Ceci est fait en \ref{lecaselementaire}. Le travail combinatoire de toute la premi\`ere partie du papier est de ramener le cas des $\psi$ g\'en\'eraux de restriction discr\`ete \`a la diagonale au cas \'el\'ementaire; on r\'ealise la repr\'esentation associ\'ee \`a $\psi$ dans le groupe de Grothendieck avec comme base des induites de repr\'esentations associ\'ees \`a des morphismes \'el\'ementaires. Nous allons essayer de donner ici une id\'ee de la formule; on fixe une sous-repr\'esentation irr\'eductible de la repr\'esentation de $W_{F}\times SL(2,{\mathbb C}) \times SL(2,{\mathbb C})$. Une telle repr\'esentation est le produit tensorielle de 3 repr\'esentaions irr\'eductibles, l'une de $W_{F}$ que nous appelons $\rho$ et deux de $SL(2,{\mathbb C})$ qui sont connues d\`es que l'on a leurs dimensions, dimensions que l'on note $a,b$ respectivement. Ainsi $\psi$ est la somme d'une repr\'esentation associ\'ee \`a un morphisme du m\^eme type $\psi'$ et de cette repr\'esentation que nous \'ecrivons $\pi(\rho,a,b)$; comme $\psi$ n'est pas \'el\'ementaire, on suppose que $inf(a,b)>1$ et pour n'utiliser dans cette introduction que la repr\'esentation de Steinberg, on suppose ici que $a\geq b$. Gr\^ace \`a la correspondance de Langlands locale \cite{harris}, \cite {henniart}, on identifie $\rho$ \`a une repr\'esentation cuspidale d'un groupe $GL(d_{\rho},F)$ (ce qui d\'efinit $d_{\rho}$). Dans le groupe de Grothendieck on d\'efinit a priori la repr\'esentation suivante:
$$
\oplus_{c\in ]b,a]}(-1)^{a-c} St(c,\rho)\vert\,\vert^{(a-b-c+1)/2}\times \pi_{c}\times St(c,\rho)\vert\,\vert^{-(a-b-c+1)/2}$$
$$
\oplus (-1)^{[b/2]}St(a-b+1,\rho)\times \pi(\rho,a+1,b-1)\times \pi(\psi'),
$$
o\`u $\pi_{c}$, pour $c=a$ vaut $\pi(\rho,a,b-2)\times \pi(\psi')$ (si $b=2$, $\pi_{a}=\pi(\psi')$), est la repr\'esentation 0 si $b=2$ et $c<a$ et sinon est un module de Jacquet convenable de la repr\'esentation $\pi(\rho,a+2,b-2)\times \pi(\psi')$ qui d\'epend de $c$ \'evidemment; on renvoie au texte pour une 
d\'efinition pr\'ecise, le point important est que $inf(a+2,b-2)=b-2<b$ et de m\^eme $inf(a+1,b-1)=b-1<b$.   
Chacune des repr\'esentations \'ecrites dans cette combinaison lin\'eaire a une action de $\theta$ que l'on pr\'ecise de fa\c{c}on ad hoc de fa\c{c}on \`a obtenir le r\'esultat: cet \'el\'ement du groupe de Grothendieck de $\tilde{G}(N)$ v\'erifie pour tout \'el\'ement $g\in GL(n,F)$ suffisamment r\'egulier $tr\tilde{\pi}(\psi)(g,\theta)=tr \pi(\psi)(g,\theta)$. On conclut en faisant une r\'ecurrence sur $I_{\psi}:=\sum_{(\rho',a',b')}inf(a',b')$ o\`u $(\rho',a',b')$ parcourt l'ensemble des repr\'esentations irr\'eductibles incluses dans $\psi$. Revenons sur le probl\`eme de l'action de $\theta$; on  d\'efinit a priori une action de $\theta$ sur toutes les repr\'esentations associ\'ees \`a des morphismes du type de $\psi$. Avec les notations ci-dessus, on aura canoniquement une action sur chaque repr\'esentation $\pi_{c}$ qui se prolonge canoniquement en une action sur l'induite correspondant \`a $c$ ci-dessus. Toutes ces actions provenant canoniquement d'une action mise sur  $\pi(\rho,a+2,b-2)\times \pi(\psi')$, il est naturel qu'elles se combinent bien. Mais comme il y a encore le terme $St(a-b+1,\rho)\times \pi(\rho,a+1,b-1)\times \pi(\psi')$ dans la combinaison lin\'eaire, il faut que l'action mise sur $\pi(\rho,a+2,b-2)\times \pi(\psi')$ soit aussi li\'ee \`a celle mise sur $St(a-b+1,\rho )\times \pi(\rho,a+1,b-1)\times \pi(\psi')$. Si on prend l'action qui vient des classifications de Langlands (cf. \ref{definitiondelanormalisationdewhittaker}), ici, cela fonctionnerait bien car on a suppos\'e que $a\geq b$; mais dans le cas inverse, o\`u $b\geq a$ et o\`u les Steinberg sont remplac\'es par les ''Speh'' , cela ne fonctionnerait plus, il faudrait introduire un signe sur l'action de $\theta$. Plut\^ot que de mettre des signes, on a choisi de d\'efinir sur les repr\'esentations $\pi(\psi)$ une action de $\theta$ qui est ''ind\'ependante'' des signes $a'-b'$ (o\`u $(\rho',a',b')$ parcourt comme ci-dessus l'ensemble des sous-repr\'esentations irr\'eductibles incluses dans $\psi$). Et dans la partie \ref{normalisationdelactiondetheta}, on \'etudie ces signes. D'un point de vue local il n'y a de toute fa\c{c}on aucune raison de privil\'egier la classification de Langlands \`a celle de Zelevinsky. Appelons normalisation \`a la Whittaker, l'action de $\theta$  obtenue en utilisant la normalisation de Langlands; un peu plus pr\'ecis\'ement la classification de Langlands fournit une induite (n\'ecessairement $\theta$-invariante) tel que $\pi(\psi)$ soit l'unique quotient irr\'eductible (dans tout le travail on utilise sous-module et non quotient en dualisant convenablement). L'induite a un mod\`ele de Whittaker pour tout choix de caract\`ere additif de $F$. On fait un tel choix de caract\`ere additif et avec ce choix on construit un caract\`ere $\chi$ du sous-groupe unipotent sup\'erieur de $GL(n,F)$ invariant sous l'action de $\theta$; pour tout repr\'esentation $I$ de $GL(n,F)$, munie d'une action de $\theta$, l'espace des coinvariants pour ce sous-groupe unipotent et le caract\`ere $\chi$ h\'erite d'une action de $\theta$. Si l'espace des coinvariants est de dimension 1, ce qui est le cas ci-dessus, on a ainsi un moyen de normaliser l'action de $\theta$ en demandant que l'action induite soit l'action triviale. C'est ce que l'on appelle la  normalisation \`a la Whittaker.
On v\'erifiera que ce choix est ind\'ependant du caract\`ere additif fix\'e quand $\psi$ se factorise par le $L$-groupe de $H_{N}$. On note  $\theta_{W}(\psi)$ l'action ainsi d\'efinie de $\theta$. Quand on utilise la classification de Zelevinsky, on r\'ealise  $\pi(\psi)$ comme unique sous-module irr\'eductible d'une induite $\theta$-invariante; ici l'induite n'a pas de mod\`ele de Whittaker en g\'en\'eral, il faut donc commencer par mettre une action de $\theta$ sur l'induite puis la restreindre (c'est expliqu\'e en \ref{definitiondelanormalisationunipotente}). On note $\theta_{U}(\psi)$ cette normalisation et on note  $\theta(\psi)$ la normalisation ''ad hoc'', celle qui fait ''marcher'' le r\'esultat ci-dessus. Il n'y a plus besoin de supposer ici que $\psi$ est de restriction discr\`ete \`a la diagonale. On calcule  explicitement tous les signes $\sigma_{W}(\psi):= \theta(\psi)/\theta_{W}(\psi)$, $\sigma_{U}(\psi):=\theta(\psi)/\theta_{U}(\psi)$ $\sigma_{\emptyset}(\psi):=\theta_{W}(\psi)/\theta_{U}(\psi)$; ces signes ont une jolie interpr\'etation relativement au transfert endoscopique. Ici on se limite au cas o\`u $\psi$ se factorise par le $L$-groupe de $H_{N }$, cas o\`u $\pi(\psi)\circ \theta$ est stable d'apr\`es \cite{arthurrecent}.
 On montre alors  qu'il existe un caract\`ere $\epsilon_{W}(\psi),\epsilon_{U}(\psi)
$ et  $\epsilon_{\emptyset}(\psi)$ du centralisateur de $\psi$ dans le $L$-groupe de $H_{N}$ dont la restriction au centre de ce $L$-groupe est triviale et tel que en notant $z_{2}$ l'image par $\psi$ de l'\'el\'ement non trivial de la 2e copie de $SL(2,{\mathbb C})$, le signe $\sigma_{?}(\psi)$ soit la valeur du caract\`ere $\epsilon_{?}(\psi)$ sur $z_{2}$. 

Supposons encore que $\psi$ se factorise par le $L$-groupe de $H_{N}$;
pour motiver notre article on rappelle en \ref{constructiondesrepresentations} comment on a associ\'e \`a tout caract\`ere $\epsilon$ de $Cent (\psi)$ dans ce $L$-groupe une repr\'esentation $\pi(\psi,\epsilon)$ ou (parfois) 0 que l'on note $\pi^H(\psi,\epsilon)$; on a \'etudi\'e ces repr\'esentations sous certaines hypoth\`eses mais ici seule la construction nous importe. Le transfert que l'on obtient dans les cas o\`u on le conna\^{\i}t pour les repr\'esentations temp\'er\'ees (cf. \ref{hypothese} pour des hypoth\`eses 
pr\'ecises)  est de la forme:
$$ 
\sum_{\epsilon}\epsilon(z_{2})tr \pi^H(\psi,\epsilon)=tr\pi(\psi)\circ \theta,
$$
o\`u $\epsilon$ parcourt l'ensemble des caract\`eres du centralisateur de $\psi$ dans le $L$-groupe de $H_{N}$ triviaux sur le centre de ce $L$-groupe. Et avec la propri\'et\'e donn\'ee une telle formule ne doit \^etre montr\'ee que pour l'une des actions de $\theta$, on l'obtient ensuite pour les autres en changement l'application $\epsilon \mapsto \pi^H(\psi,\epsilon)$ en $\epsilon \mapsto \pi(\psi,\epsilon_{?}\epsilon)$. 
Cette formule est heureusement analogue \`a celle de \cite{arthurrecent} 30.1 qui fait suite \`a  \cite{stabilisation1}, \cite{stabilisation2}, \cite{stabilisation3}, \cite{stabilisation4}. Notre point de vue est  une description explicite des repr\'esentations qui interviennent. On explique d'ailleurs en \ref{remarquesurlhypothese} comment le transfert obtenu par Arthur est un point d\'epart pour nos constructions et ce qu'il entra\^{\i}ne des hypoth\`eses dont nous avons besoin pour aller au bout des descriptions; on pousse ce point de vue en \cite{classification} pour avoir toutes les hypoth\`eses dont nous avons besoin pour une description combinatoire mais cela sort du cadre de ce travail limit\'e \`a $GL(N,F)$ tordu.

\section{D\'efinitions et notations g\'en\'erales\label{definitionsgenerales}}
\subsection{Notations pour les intervalles}
A de nombreuses reprises dans cet article des sommes ou des produits sont index\'es par des \'el\'ement $\ell$ dans un segment; on a pris quelques conventions sur cette notations. Pour nous un segment est la donn\'ee de 2 nombres r\'eels $x,y$ tel que $x-y\in {\mathbb Z}$ (en g\'en\'eral $x,y$ seront des demi-entiers non n\'ecessairement entiers) et la notation $\ell\in [x,y]$ est l'ensemble des \'el\'ements $\ell$ de la forme $x,x-1, \cdots ,y+1,y$  si $x\geq y$ et $x,x+1, \cdots, y-1,y$ si $x\leq y$.
\subsection{Action de $\theta$ sur une repr\'esentation cuspidale autoduale}
 Soit $\rho$ une repr\'esentation cuspidale irr\'eductible de $GL(d_{\rho},F)$ (ce qui d\'efinit $d_{\rho}$). On suppose que $\rho$ est invariante par conjugaison sous $\theta$ (action d\'efinie ci-dessus quand on fait $n=d_{\rho}$). On fixe un espace de la repr\'esentation $\rho$ muni d'un automorphisme $N_{\rho,\theta}$ tel que $N_{\rho,\theta}$ entrelace $\rho$ et $\rho \circ \theta$ et on suppose que $N_{\rho,\theta}^2=1$. On consid\`ere $\rho$ comme une repr\'esentation du produit semi-direct $\tilde{G}(d_{\rho})$ o\`u $\theta$ agit par l'op\'erateur $N_{\rho,\theta}$. On a fait l\`a un choix qui dure dans tout le texte; une fa\c{c}on de normaliser est de normaliser \`a la  Whittaker (comme expliqu\'e dans l'inroduction), ce qui suppose que l'on a fix\'e un caract\`ere additif du corps de base; on v\'erifiera qu'un tel choix n'influe pas sur le transfert (cf. \ref{lesgroupes}).
\subsection{Action de $\theta$ sur des induites}
De fa\c{c}on g\'en\'erale, soit $\sigma$ une repr\'esentation de $\tilde{G}(d_{\sigma})$ (ce qui d\'efinit $d_{\sigma}$) et soit $\tau$ une repr\'esentation de $GL(d_{\tau},F)$ (ce qui d\'efinit $d_{\tau}$). On note $^\theta \tau$ l'image de $\tau$ par l'automorphisme induit par $\theta$ et soit $\tau'$ une repr\'esentation de $GL(d_{\tau},F)$ isomorphe \`a $^\theta \tau$. Soit $A$ un automorphisme de l'espace, $V_{\tau}$, de $\tau$ sur l'espace $V_{ \tau'}$, de $ \tau'$ tel que pour tout $v$ dans l'espace de $\tau$ et tout $g\in GL(d_{\tau},F)$, on ait $ \tau' (g) A(v)=A(\tau(\theta(g))v)$. On note $V_{\sigma}$ l'espace de $\sigma$ et $\theta_{\sigma}$ l'action de $\theta$ sur cet espace. On  pose pour tout $\phi$ dans l'induite de $\tau\otimes \sigma \otimes \tau'$  et tout $g\in GL(d_{\sigma}+2d_{\tau},F)$:
$$
\theta \phi(g):=\bigl(A^{-1} \otimes \theta_{\sigma} \otimes A\bigr)\iota (\phi(\theta(g))),
$$
o\`u $\iota$ est l'isomorphisme de $V_{\tau}\otimes V_{\sigma}\otimes V_{ \tau'}$ sur $V_{\tau'}\otimes V_{\sigma} \otimes V_{\tau}$ qui \'echange les termes extr\^emes. On v\'erifie que cela d\'efinit une action de $\theta$ sur l'induite $\tau \times \sigma \times  \tau'$; en effet soit $m_{1}\in Gl(d_{\tau},F)$, $m_{0}\in GL(d_{\sigma},F)$ et $m_{2}\in GL(d_{\tau},F)$. On consid\`ere $m:=(m_{1},m_{0},m_{2})$ comme un \'el\'ement de $GL(d_{\sigma}+2d_{\tau},F)$. Pour $g$ comme ci-dessus, on a $\theta(mg)=(\theta(m_{2}),\theta(m_{0}),\theta(m_{1}))\theta(g)$ et l'assertion s'en d\'eduit. Ainsi l'induite $\tau \times \sigma \times \tau'$ est une repr\'esentation de $\tilde{G}(d_{\sigma}+2 d_{\tau})$ et cette action  d\'epend  du choix de $A$ \`a ceci pr\`es qu'elle est ind\'ependante de la multiplication de $A$ par un scalaire; en particulier si $\tau$ est irr\'eductible, elle ne d\'epend  pas du choix de $A$. On g\'en\'eralise cette construction en rempla\c{c}ant $\tau$ par une induite $\tau_{1}\times\cdots \times \tau_{\ell}$, o\`u chaque $\tau_{i}$ est irr\'eductible; on voit clairement comment prolonger naturellement l'action de $\theta_{\sigma}$ \`a l'induite $\tau_{1}\times \cdots \times \tau_{\ell}\times \sigma \times \, ^\theta\tau_{\ell}\times \cdots \times ^\theta\tau_{1}$ de fa\c{c}on ind\'ependante de tout choix de morphisme, $A_{1}, \cdots, A_{\ell}$ entrela\c{c}ant, pour $i=1,\cdots, \ell$,  $^\theta\tau_{i}$ et $\tau_{i}\circ \theta$.

Pour $\rho$ fix\'e et $a$ un entier, on note $St(a,\rho)$ la repr\'esentation de Steinberg qui est l'unique sous-repr\'esentation irr\'eductible de l'induite $\sigma_{St,a,\rho}:=\rho\vert\,\vert^{(a-1)/2}\times \cdots \times \rho\vert\,\vert^{-(a-1)/2}$ et pour $b$ un entier, on note $Sp(b,\rho)$ la repr\'esentation de Speh qui est l'unique sous-repr\'esentation irr\'eductible de l'induite $\sigma_{Sp,b,\rho}:=\rho\vert\,\vert^{-(b-1)/2} \times \cdots \times \rho\vert\,\vert^{(b-1)/2}$.

\subsection{Notations relatives aux modules de Jacquet\label{notationdujac}} 
Soit $\pi$ une repr\'esentation de $\tilde{G}(N)$ de longueur finie et soit $\rho$ une repr\'esentation irr\'eductible cuspidale et autoduale de $GL(d_{\rho},F)$. On consid\`ere le parabolique standard de $GL(n,F)$ de Levi $GL(d_{\rho},F)\times GL(N-2d_{\rho},F) \times GL(d_{\rho},F)$; on remarque qu'il est $\theta$-invariante. On consid\`ere le module de Jacquet de $\pi$ le long du radical unipotent de ce  parabolique. On consid\`ere ce module de Jacquet dans le groupe de Grothendieck des repr\'esentations du Levi . On fixe $x\in {\mathbb R}$ et  on projette cet \'el\'ement semi-simplifi\'e sur le support cuspidal $\rho\vert\,\vert^x$ par l'action du premier facteur $GL(d_{\rho},F)$ et sur le support cuspidal $\rho\vert\,\vert^{-x}$ pour le deuxi\`eme facteur $GL(d_{\rho},F)$ (le dernier dans l'\'ecriture du Levi); cette projection semi-simplifi\'ee est une somme de repr\'esentations irr\'eductibles $\rho\vert\,\vert^x\otimes \sigma \otimes \rho\vert\,\vert^{-x}$ o\`u $\sigma$ parcourt un ensemble de repr\'esentations irr\'eductibles de $GL(N-2d_{\rho},F)$. On note $Jac^\theta_{x}(\pi)$ l'\'el\'ement du groupe de Grothendieck de $GL(N-2d_{\rho},F)$ qui est la somme de ces repr\'esentations $\sigma$. Comme $\pi$ a une action de $\theta$, on obtient une action de $\theta$ sur le module de Jacquet tout entier puis sur sa projection sur la composante correspond \`a l'action de $\rho\vert\,\vert^x$ pour le premier facteur $GL(d_{\rho},F)$ et $\rho\vert\,\vert^{-x}$ pour le dernier facteur. On obtient donc une action du produit semi-direct , $\tilde{M}$, de $GL(d_{\rho},F) \times GL(n-2d_{\rho},F) \times GL(d_{\rho},F)$ par le groupe $\{1,\theta\}$. On a d\'efini $\tilde{G}(n-2d_{\rho})$ d\'efinition que l'on prolonge au cas o\`u $n=2d_{\rho}$ en consid\'erant que dans ce cas le groupe est $\{1,\theta\}$. Et $\tilde{G}(n-2d_{\rho})$ se plonge naturellement dans $\tilde{M}$. Supposons d'abord que la projection consid\'er\'ee du module de Jacquet de $\pi$ est semi-simple en tant que repr\'esentation de $\tilde{M}$. Elle se r\'ealise alors dans un espace $V\otimes W \otimes V$ o\`u $V$ est un espace pour la repr\'esentation $\rho$ et $W$ est une repr\'esentation semi-simple, $\sigma_{W}$, de $GL(n-2d_{\rho},F)$. On fixe un op\'erateur $A$ de $Aut(V)$ qui entrelace $\rho$ et $^\theta \rho$ et on note $\iota$ l'\'echange des 2 facteurs $V$. Alors $D:=(A\otimes id _{W}\otimes A^{-1})\circ \iota \circ \theta$ est un op\'erateur de $V\otimes W\otimes V$ qui v\'erifie pour tout $m,m'\in GL(d_{\rho},F)$ et $h\in GL(n-2d_{\rho},F)$:
$$
D\circ(\rho\vert\,\vert^{x}\otimes \sigma_{W} \otimes \rho\vert\,\vert^{-x})(m,h,m')=
(\rho\vert\,\vert^{x}\otimes \sigma_{W} \otimes \rho\vert\,\vert^{-x})(m,\theta(h),m')\circ D.
$$
Par irr\'eductibilit\'e de $\rho$, $D$ est donc de la forme, $id\otimes D_{W}\otimes id$ et $D_{W}$ est ind\'ependant du choix de $A$ fait ci-dessus. Ainsi $W$ est muni d'une action de $\tilde{G}(n-2d_{\rho})$ en faisant agir $\theta$ via $D_{W}$. En filtrant convenablement le module de Jacquet de $\pi$, par cette m\'ethode, on munit $Jac^\theta(\pi)$ d'une action semi-simple de $\tilde{G}(n-2d_{\rho})$ qui ne d\'epend que de l'action de $\theta$ sur $\pi$. Ceci permet de consid\'erer $Jac^\theta_{x}(\pi)$ comme un \'el\'ement du groupe de Grothendieck de $\tilde{G}(N-2d_{\rho})$ et c'est ce que nous ferons. On remarque que cette d\'efinition n'utilise que l'image de $\pi$ dans le groupe de Grothendieck de $\tilde{G}(N)$. Quant au lieu d'avoir un seul r\'eel $x$, on en a une collection ordonn\'ee, $x_{1}, \cdots, x_{v}$, on it\`ere les d\'efinitions en posant:
$$
Jac^\theta_{x_{1}, \cdots, x_{v}}(\pi):=Jac^\theta_{x_{v}}\cdots Jac^\theta_{x_{1}}(\pi).
$$
Dans ce papier on reprend aussi la notation $Jac_{x}(\pi)$ o\`u ici $\rho,x$ sont comme ci-dessus, $\pi$ est simplement une repr\'esentation de longueur finie de $GL(N,F)$ et on consid\`ere le module de Jacquet de $\pi$ pour le parabolique standard de Levi $GL(d_{\rho},F)\times GL(N-d_{\rho},F)$ avec sa projection sur le support cuspidal $\rho\vert\,\vert^x$ pour l'action du premier facteur. La propri\'et\'e  que nous utiliserons est la suivante: soit $x,y$ tels que $x-y\neq \pm 1$ alors $Jac_{x,y}(\pi)=Jac_{y,x}(\pi)$ \'egalit\'e dans le groupe de Grothendieck. En effet pour calculer $Jac_{x,y}(\pi)$ on commence par calculer le module de Jacquet pour le parabolique $GL(2d_{\rho},F)\times GL(N-2d_{\rho},F)$; puis on projette sur le support cuspidal $\rho\vert\,\vert^x,\rho\vert\,\vert^y$ pour l'action du premier facteur. Toute repr\'esentation irr\'eductible de $GL(2d_{\rho},F)$ ayant ce support cuspidal est l'induite irr\'eductible $\rho\vert\,\vert^x\times \rho\vert\,\vert^y$. D'o\`u la sym\'etrie entre $x$ et $y$. Cette formule se g\'en\'eralise en $Jac^\theta_{x,y}=Jac^\theta_{y,x}$ sous les m\^emes hypoth\`eses (implicitement ici $N\geq 4d_{\rho}$).

\subsection{Notation pour le socle d'une induite\label{socle}}
Soient $m\in {\mathbb N}$ et $m_{1}, \cdots, m_{r}$ une collection d'entiers tels que $m=\sum_{i\in [1,r]}m_{i}$. Pour $i\in [1,r]$ soit $\sigma_{i}$ une repr\'esentation de $GL(m_{i},F)$. On consid\`ere $$\sigma:=\times_{i\in [1,r]}\sigma_{i}$$
la repr\'esentation induite de $GL(m,F)$. Le socle d'une repr\'esentation est la somme de ses sous-modules irr\'eductibles et on note ce socle
$$
<\sigma_{1}, \cdots, \sigma_{r}>.
$$
On n'utilisera cette d\'efinition que dans des cas o\`u le socle est irr\'eductible et essentiellement dans le cas suivant: $\rho$ est une repr\'esentation cuspidale fix\'e, $[x,y]$ est un intervalle au sens d\'efini ci-dessus et on prend pour $r:= \vert (x-y+1)\vert$ et pour $\sigma_{i}:=\rho\vert\,\vert^{x-i+1}$ si $x\geq y$ et $\rho\vert\,\vert^{x+i-1}$ si $y\geq x$. Le socle est alors not\'e 
$<\rho\vert\,\vert^x, \cdots, \rho\vert\,\vert^y>.$

\subsection{Notation pour les $\vert\,\vert^s$ et les op\'erateurs d'entrelacement standard\label{standard}}
Soit $M$ un sous-groupe de Levi de $GL(N,F)$ et $\sigma$ une repr\'esentation irr\'eductible de $M$; supposons que $M$ soit maximal, il s'\'ecrit alors comme produit $GL(N_{1},F)\times GL(N_{2},F)$ et $\sigma$ est un produit tensoriel $\sigma_{1}\otimes \sigma_{2}$.  Pour $s\in {\mathbb C}$, on d\'efinit l'induite pour le parabolique standard:
$$
\sigma_{1}\vert\,\vert^s\times \sigma_{2},
$$o\`u $\vert\,\vert^s$ est un raccourci pour $\vert det_{GL(N_{1},F)}\vert^s$.
On dispose alors des op\'erateurs d'entrelacement standard 
$$
M(s):=\sigma_{1}\vert\,\vert^s\times \sigma_{2} \rightarrow \sigma_{2}\times \sigma_{1}\vert\,\vert^s
$$
ces op\'erateurs sont obtenus par int\'egration pour $s$ tel que $Re\, s >>0$ puis par prolongement m\'eromorphe. Ils n'ont pas de raison d'\^etre holomorphe en $s=0$ mais on peut toujours trouver une fonction m\'eromorphe de $s$, $r(s)$ tel que l'op\'erateur $r(s)M(s)$ soit holomorphe en $s=0$ et que sa valeur en $s=0$ donne un op\'erateur non nul.

Cette d\'efinition d'op\'erateur d'entrelacement standard, se g\'en\'eralise \`a tout sous-groupe de Levi.

\subsection{Action de $\theta$ sur $Sp(b,\rho)$ et $St(a,\rho)$\label{definitionsimple}}

On suppose que $\rho\simeq \, ^\theta \rho$ et on fixe $N_{\theta,\rho}$ un isomorphisme de $\rho$ avec $^\theta \rho$. Pour tout r\'eel $x$, $N_{\theta,\rho}$ entrelace naturellement la repr\'esentation $\rho\vert\,\vert^x$ avec $^\theta \rho \vert\,\vert^{-x}$. Ainsi si $a$ est pair, la repr\'esentation induite  $\sigma_{St,a,\rho}$ est munie d'une structure de $\tilde{G}(ad_{\rho})$ module en it\'erant la construction faite ci-dessus pour une induite en commen\c{c}ant avec $\sigma$ la repr\'esentation triviale. Si $a$ est impair, on fait la m\^eme construction mais en partant de $\sigma=\rho$. Dans le cas $a$ pair, la structure ne d\'epend pas du choix de $N_{\theta,\rho}$ par contre elle en d\'epend dans le cas impair. On peut \'evidemment faire la m\^eme chose pour $\sigma_{Sp,b,\rho}$. Par restriction, on en d\'eduit une structure de $\tilde{G}(ad_{\rho})$-module sur $St(a,\rho)$ et une structure de $\tilde{G}(bd_{\rho})$-module sur $Sp(b,\rho)$.

\subsection{D\'efinition de l'action de $\theta$ sur $Sp(b,St(a,\rho))$\label{independancesimple}}

On g\'en\'eralise ces d\'efinitions, en d\'efinissant pour $a,b$ des entiers
 $Sp(b,St(a,\rho))$ comme l'unique sous-module irr\'eductible de l'induite:
$$
St(a,\rho)\vert\,\vert^{-(b-1)/2}\times \cdots \times St(a,\rho)\vert\,\vert^{(b-1)/2}.\eqno(1)
$$
On sait ou on v\'erifie facilement que cette d\'efinition a une forme de  sym\'etrie en $a,b$ ou plut\^ot entre $Sp$ et $St$ dans la mesure o\`u $Sp(b,St(a,\rho))$ est aussi l'unique sous-module irr\'eductible de l'induite:
$$
Sp(b,\rho)\vert\,\vert^{(a-1)/2}\times \cdots \times Sp(b,\rho)\vert\,\vert^{-(a-1)/2}.\eqno(2)
$$
Comme nous venons de munir $St(a,\rho)$ d'une action de $\theta$, on en d\'eduit une action de $\theta$ sur (1) et par restriction une action de $\theta$ sur $Sp(b,St(a,\rho))$ que l'on note momentan\'ement $\theta_{St}$. On peut faire la m\^eme chose pour (2) en utilisant l'action de $\theta$ sur $Sp(b,\rho)$ ce qui donne une action de $\theta$ sur $Sp(b,St(a,\rho))$ que l'on note $\theta_{Sp}$. Heureusement, ces actions sont les m\^emes comme nous allons le v\'erifier.

En effet, si $inf(a,b)=1$, c'est exactement la d\'efinition: pour clarifier supposons que $inf(a,b)=a$, c'est-\`a-dire $a=1$, alors (1) n'est autre que $\sigma_{Sp,b,\rho}$ et (2) est $Sp(b,\rho)$. On a bien d\'efini l'action de $\theta$ sur $Sp(b,\rho)$ par restriction de celle sur $\sigma_{Sp(b,\rho)}$.

On suppose donc maintenant que $inf(a,b)>1$ et on montre l'assertion par induction en l'admettant pour $Sp(b',St(a',\rho))$ pour tout couple $(a',b')$ tel que  $a+b>a'+b'$; dans ce qui suit si $a=2$ tout terme \'ecrit ci-dessous contenant $a-2$ est la repr\'esentation triviale de $GL(0,F)=\{1\}$ et on prend la m\^eme convention si $b=2$. L'induite (1) contient la repr\'esentation
$$
St(a,\rho)\vert\,\vert^{-(b-1)/2}\times Sp(b-2,St(a,\rho)) \times St(a,\rho)\vert\,\vert^{(b-1)/2}.\eqno(3)
$$
En rempla\c{c}ant $St(a,\rho)$ par sa d\'efinition, on v\'erifie que (3) est un sous-module de l'induite:
$$
\rho\vert\,\vert^{(a-b)/2}\times \cdots \times \rho\vert\,\vert^{-(a+b)/2+1}\times Sp(b-2,St(a,\rho)) \times \rho\vert\,\vert^{(a+b)/2-1}\times  \cdots \times \rho\vert\,\vert^{-(a-b)/2}.\eqno(4) $$
Sur (4), on a une action de $\theta$ en tant qu'induite mais il n'est pas clair a priori que l'action de $\theta$ sur (3) est la restriction de cette action sur (4). Toutefois c'est vrai car l'action sur (3) ne d\'epend pas du choix de $A$ l'op\'erateur qui identifie $^\theta (St(a,\rho)\vert\,\vert^{-(b-1)/2})$ avec $St(a,\rho)^{(b-1)/2}$. On peut donc prendre pour $A$ la restriction de l'op\'erateur qui se d\'eduit naturellement de $N_{\theta,\rho}$ et qui envoie $\rho\vert\,\vert^{(a-b)/2}\times \cdots \times \rho\vert\,\vert^{-(a+b)/2+1}$ sur $\rho\vert\,\vert^{(a+b)/2-1}\times  \cdots \times \rho\vert\,\vert^{-(a-b)/2}$.
Par r\'ecurrence on sait que l'action de $\theta$ sur $Sp(b-2,St(a,\rho))$ s'obtient en restreignant celle que l'on d\'efinit naturellement sur 
$$
Sp(b-2,\rho)\vert\,\vert^{(a-1)/2}\times Sp(b-2,St(a-2,\rho))\times Sp(b-2,\rho)\vert\,\vert^{-(a-1)/2}.\eqno(5)
$$
On repr\'esente $Sp(b-2,\rho)$ comme sous-module de $\sigma_{Sp,b-2,\rho}$ et comme ci-dessus on montre que l'action de $\theta$ sur (5) s'obtient en restreignant l'action de $\theta$ sur l'induite:
$$
\rho\vert\,\vert^{(a-b)/2+1}\times \cdots \times \rho\vert\,\vert^{(a+b)/2-2}\times Sp(b-2,St(a-2,\rho)) \times
\rho\vert\,\vert^{-(a+b)/2+2}\times  \cdots \times \rho\vert\,\vert^{-(a-b)/2-1}.
$$En remettant avec (4), on vient de montrer que l'action $\theta_{St}$ sur $Sp(b,St(a,\rho))$ est celle que l'on obtient en restreignant l'action naturelle de $\theta$ sur l'induite:
$$
\rho\vert\,\vert^{(a-b)/2} \times \cdots \times \rho\vert\,\vert^{-(a+b)/2+1}\times\rho\vert\,\vert^{(a-b)/2+1}\times \cdots \times \rho\vert\,\vert^{(a+b)/2-2}$$
$$\times Sp(b-2,St(a-2,\rho)) $$
$$\times \rho\vert\,\vert^{-(a+b)/2+2} \times \cdots \times \rho\vert\,\vert^{-(a-b)/2-1}\times\rho\vert\,\vert^{(a+b)/2-1}\times \cdots \times \rho\vert\,\vert^{-(a-b)/2}.\eqno(6)
$$
On fait les m\^emes manipulations en partant de (2) pour trouver que l'action $\theta_{Sp}$ est la restriction de l'action naturelle de $\theta$ sur l'induite:
$$
\rho\vert\,\vert^{(a-b)/2} \times \cdots \times \rho\vert\,\vert^{(a+b)/2-1}\times\rho\vert\,\vert^{(a-b)/2-1}\times \cdots \times \rho\vert\,\vert^{-(a+b)/2+2}$$
$$\times Sp(b-2,St(a-2,\rho)) $$
$$\times \rho\vert\,\vert^{(a+b)/2-2} \times \cdots \times \rho\vert\,\vert^{-(a-b)/2+1}\times\rho\vert\,\vert^{-(a+b)/2+1}\times \cdots \times \rho\vert\,\vert^{-(a-b)/2}.\eqno(7)
$$
Les induites (6) et (7) sont isomorphes (on le v\'erifiera) et il faut donc montrer que l'isomorphisme entrelace l'action de $\theta$, ainsi les restrictions seront les m\^emes sur leur unique sous-module $Sp(b,St(a,\rho))$.

Construisons un isomorphisme de (6) sur (7); d'abord on remarque que l'induite $$
\rho\vert\,\vert^{(a-b)/2} \times \rho\vert\,\vert^{(a-b)/2-1} \cdots \times \rho\vert\,\vert^{-(a+b)/2+1}\times\rho\vert\,\vert^{(a-b)/2+1}\times \cdots \times \rho\vert\,\vert^{(a+b)/2-2}$$  est isomorphe \`a
$$
\rho\vert\,\vert^{(a-b)/2} \times \rho\vert\,\vert^{(a-b)/2+1}\times \cdots \times \rho\vert\,\vert^{(a+b)/2-2}
 \times \rho\vert\,\vert^{(a-b)/2-1}\times  \cdots \times \rho\vert\,\vert^{-(a+b)/2+1},$$
 o\`u on utilise simplement le fait que pour $y\in [(a-b)/2-1,-(a+b)/2+1]$ et $y'\in [(a-b)/2+1,(a+b)/2-2]$, $y-y'<-1$. On peut construire un isomorphisme \`a l'aide des op\'erateurs d'entrelacement standard par exemple et en transportant cet isomorphisme par $\theta$, on construit un isomorphisme de (6) sur
 $$
 \rho\vert\,\vert^{(a-b)/2} \times \rho\vert\,\vert^{(a-b)/2+1}\times \cdots \times \rho\vert\,\vert^{(a+b)/2-2}
 \times \rho\vert\,\vert^{(a-b)/2-1}\times \cdots \times \rho\vert\,\vert^{-(a+b)/2+1} \times Sp(b-2,St(a-2,\rho)) 
 $$
 $$
 \times \rho\vert\,\vert^{(a+b)/2-1}\times \cdots\times \rho\vert\,\vert^{-(a-b)/2+1} \times \rho\vert\,\vert^{-(a+b)/2+2}
 \times \cdots \times  \rho\vert\,\vert^{-(a-b)/2-1)}\times  \rho\vert\,\vert^{-(a-b)/2}.\eqno(8)
 $$
 Cet isomorphisme est compatible aux actions de $\theta$ d\'efinies essentiellement car il n'y a qu'une action de $\theta$ naturelle sur ces induites une fois que l'on a fix\'e l'action sur $Sp(b-2,St(a-2,\rho)) $, mais on renvoie \`a \ref{entrelacement} pour une d\'emonstration pr\'ecise.
 De m\^eme on entrelace l'action de $\theta$ sur (7) avec l'action de $\theta$ sur l'induite (9) ci-dessous qui lui est isomorphe:
 $$
  \rho\vert\,\vert^{(a-b)/2} \times \rho\vert\,\vert^{(a-b)/2+1}\times \cdots \times \rho\vert\,\vert^{(a+b)/2-2}
 \times \rho\vert\,\vert^{(a-b)/2-1}\times \cdots \times \rho\vert\,\vert^{-(a+b)/2+2}\times  \rho\vert\,\vert^{(a+b)/2-1} $$ $$\times Sp(b-2,St(a-2,\rho)) \times\eqno(9)$$
 $$\rho\vert\,\vert^{-(a+b)/2+1}
 \times \rho\vert\,\vert^{(a+b)/2-2}\times \cdots\times \rho\vert\,\vert^{-(a-b)/2+1} \times \rho\vert\,\vert^{-(a+b)/2+2}
 \times \cdots \times  \rho\vert\,\vert^{-(a-b)/2-1}\times  \rho\vert\,\vert^{-(a-b)/2}.
 $$
Ici on a fait commuter $\rho\vert\,\vert^{(a+b)/2-1}$ avec les induites $\rho\vert\,\vert^y$ pour $y \in [(a-b)/2-1,-(a+b)/2+2]$ et  $\rho\vert\,\vert^{-(a+b)/2+1}$ avec les induites 
$\rho\vert\,\vert^{-y}$ pour les m\^emes valeurs de $y$. Il reste \`a construire un isomorphisme:
$$
\rho\vert\,\vert^{(a+b)/2-1}\times Sp(b-2,St(a-2,\rho)) \times \rho\vert\,\vert^{-(a+b)/2+1}\simeq
\rho\vert\,\vert^{-(a+b)/2+1}\times Sp(b-2,St(a-2,\rho))\times \rho\vert\,\vert^{(a+b)/2-1}
$$
entrela\c{c}ant les actions de $\theta$. On glisse un param\`etre $s\in {\mathbb C}$ et on regarde l'entrelacement:
$$
\rho\vert\,\vert^{(a+b)/2-1+s}\times Sp(b-2,St(a-2,\rho)) \times \rho\vert\,\vert^{-(a+b)/2+1 -s}\rightarrow$$
$$
\rho\vert\,\vert^{-(a+b)/2+1-s}\times Sp(b-2,St(a-2,\rho))\times \rho\vert\,\vert^{(a+b)/2-1+s}.\eqno(10)
$$
En $s=0$ les induites \'ecrites sont irr\'eductibles car le support cuspidal de $Sp(b-2,St(a-2,\rho))$ est form\'e de repr\'esentation de la forme $\rho\vert\,\vert^x$ avec $x\in [-(a+b)/2+3,(a+b)/2-3]$. Pour obtenir l'assertion cherch\'ee on va d\'emontrer plus g\'en\'eralement le lemme ci-dessous (\ref{entrelacement}).
\subsection{Op\'erateurs d'entrelacement et action de $\theta$\label{entrelacement}}
Fixons d'abord quelques notations:
soit $\pi_{i}$ pour $i=1,2$ des repr\'esentations de groupes lin\'eaires et soit $s_{i}\in {\mathbb C}$ pour $i=1,2$. On note $M(\pi_{1},\pi_{2},s)$ l'op\'erateur d'entrelacement standard de l'induite $$M(\pi_{1},\pi_{2},s): \qquad \pi_{1}\vert\,\vert^{s_{1}} \times \pi_{2}\vert\,\vert^{s_{2}} \rightarrow \pi_{2}\vert\,\vert^{s_{2}} \times \pi_{1}\vert\,\vert^{s_{1}}.$$ On d\'efinit de m\^eme $$M(^\theta\pi_{1},^\theta \pi_{2},s):\qquad ^\theta \pi_{2}\vert\,\vert^{-s_{2}} \times \, ^\theta\pi_{1}\vert\,\vert^{-s_{1}} \rightarrow \, ^\theta\pi_{1}\vert\,\vert^{-s_{1}}\times \, ^\theta\pi_{2}\vert\,\vert^{-s_{2}}.
$$Ces op\'erateurs ne d\'ependent que de $s=s_{1}-s_{2}$ dans les 2 cas. On identifie les espaces de $V_{\pi_{i}}$ et $V_{^\theta \pi_{i}}$ pour $i=1,2$.
On note $i$ l'isomorphisme de $V_{\pi_{2}}\otimes V_{\pi_{1}}$ sur $V_{\pi_{1}}\otimes V_{\pi_{2}}$ qui \'echange les facteurs . On remarque que l'application qui \`a $\phi \in \pi_{2}\vert\,\vert^{s_{2}}\times \pi_{1}\vert\,\vert^{s_{1}}$ associe $i\,\circ\,  \theta\, \phi$ est \`a valeurs dans $^\theta \pi_{1}\vert\,\vert^{-s_{1}}\times ^\theta \pi_{2}\vert\,\vert^{-s_{2}}$ et on note $i\circ \Theta$ cet op\'erateur. On d\'efinit $i_{\theta}$ l'isomorphisme inverse de $i$ et on a aussi l'op\'erateur $$i_{\theta} \circ \Theta: \pi_{1}\vert\,\vert^{s_{1}}\times \pi_{2}\vert\,\vert^{s_{2}}\rightarrow ^\theta\pi_{2}\vert\,\vert^{-s_{2}}\times \, ^\theta \pi_{1}\vert\,\vert^{-s_{1}}.$$

\bf Lemme: \sl on a l'\'egalit\'e des op\'erateurs: $i\circ \Theta \circ M(\pi_{1},\pi_{2},s)= M(^\theta\pi_{2},^\theta\pi_{1},s) \circ i_{\theta}\circ \Theta.$ en tant qu'op\'erateurs de $\pi_{1}\vert\,\vert^{s_{1}}\times \pi_{2}\vert\,\vert^{s_{2}}$ dans $^\theta \pi_{1}\vert\,\vert^{-s_{1}} \times ^\theta\pi_{2}\vert\,\vert^{-s_{2}}$.
\rm

On calcule pour $s$ de partie r\'eelle grande en utilisant la r\'ealisation des op\'erateurs d'entrelacement comme int\'egrales et l'assertion est imm\'ediate.

\subsection{Fin de la preuve de \ref{independancesimple}}
Revenons \`a (10); l'action de $\theta$ sur l'induite de gauche est, avec une g\'en\'eralisation des notations, $i\circ \Theta$ et sur l'induite de droite $i_{\theta}\circ \Theta$. L'entrelacement  cherch\'e r\'esulte ais\'ement du lemme. Cela termine la preuve de la proposition \ref{independancesimple}

\subsection{D\'efinition de $Jord(\psi)$,  $Jord_{\rho}(\psi)$ et convention d'\'ecriture\label{definitiondejord}}
On d\'ecompose $\psi$ en somme de repr\'esentations irr\'eductibles de $W_{F}\times SL(2,{\mathbb C}) \times SL(2,{\mathbb C})$; pour $c$ un entier on note, ici et ici seulement, $[c]$ la repr\'esentation de $SL(2,{\mathbb C})$ irr\'eductible de dimension $c$ une sous-repr\'esentation irr\'eductible incluse dans $\psi$ est un produit tensoriel, $\rho\otimes [a]\otimes [b]$ o\`u $\rho$ s'identifie par la correspondance de Langlands \`a une repr\'esentation irr\'eductible de $W_{F}$ et $a,b$ sont des entiers. Cette d\'ecomposition en sous-repr\'esentations irr\'eductibles donne donc un ensemble, avec multiplicit\'e, de triplets $\{(\rho,a,b)\}$ tel que $$\oplus_{(\rho,a,b)}\rho\otimes [a]\otimes [b]=\psi.$$ On transforme cet ensemble en un ensemble de quadruplets $(\rho,A,B,\zeta)$ en posant pour chaque $\rho,a,b$ comme ci-dessus, $A=(a+b)/2-1$, $B=\vert a-b\vert/2$ et $\zeta$ le signe de $a-b$ ou indiff\'eremment $\pm 1$ si $a=b$. La restriction de la repr\'esentation $\rho\otimes [a]\otimes [b]$ de $W_{F}\times SL(2,{\mathbb C})\times SL(2,{\mathbb C})$ \`a la diagonale est la somme des repr\'esentations $\rho\otimes [c]$ o\`u $c$ parcourt l'intervalle $[2A+1,2B+1]_{2}$, o\`u le $2$ en indice indique que l'on ne prend que les entiers ayant m\^eme parit\'e que les bornes; d'o\`u la signification de $A$ et $B$. On pose alors $\pi(\rho,A,B,\zeta):=Sp(b,St(a,\rho))$. On note  l'ensemble des quadruplets qui d\'ecrivent la repr\'esentation d\'efinie par $\psi$, $Jord(\psi)$; pour $\rho$ fix\'e, on note $Jord_{\rho}(\psi)$ l'ensemble des quadruplets $(\rho',A',B',\zeta')$ de $Jord(\psi)$ v\'erifiant $\rho\simeq \rho'$. 

\

\bf Attention \`a la convention commode sur $Jord_{\rho}(\psi)$.\rm

\

Cet ensemble, $Jord_{\rho}(\psi)$ est encore la r\'eunion de 2 ensembles qui n'ont rien \`a voir entre eux, $Jord_{\rho,ent}(\psi)$ et $Jord_{\rho,1/2 ent}(\psi)$ o\`u $(\rho,A,B,\zeta)\in Jord_{\rho,ent}(\psi)$ pr\'ecis\'ement quand $A,B$ sont des entiers. Dans la suite,  quand on aura un quadruplet $(\rho,A,B,\zeta)\in Jord(\psi)$ fix\'e on abr\' egera en $Jord_{\rho}(\psi)$ ce qui devrait \^etre ''ent'' ou ''1/2ent'' ceci \'etant fix\'e par $A,B$.

\subsection {D\'efinition de $\pi(\psi)$ pour $\psi$ de restriction discr\`ete \`a la diagonale \label{independance}}

On r\'eutilise la notation $\Delta$ de l'introduction: $\Delta$ est le plongement diagonal de $SL(2,{\mathbb C}) $ dans $SL(2,{\mathbb C})\times SL(2,{\mathbb C})$ que l'on prolonge par l'identit\'e sur $W_{F}$ en un morphisme de $W_{F}\times SL(2,{\mathbb C})$ dans $W_{F}\times SL(2,{\mathbb C}) \times SL(2,{\mathbb C})$.
On suppose ici que la restriction de $\psi \circ \Delta$ \`a  $W_{F}\times SL(2,{\mathbb C})$ d\'efinit une repr\'esentation sans multiplicit\'es. 
L'hypoth\`ese faite sur la restriction de $\psi$ \`a la diagonale se traduit par le fait que si $(\rho,A,B,\zeta)$ et $(\rho',A',B',\zeta')$ sont dans $Jord(\psi)$ alors soit $\rho\not \simeq \rho'$ soit $[B,A] \cap [B',A']=\emptyset$; cela entra\^{\i}ne que la repr\'esentation $<\rho\vert\,\vert^{\zeta B}, \cdots, \rho\vert\,\vert^{-\zeta A}> \times \pi(\rho',A',B',\zeta')$ est irr\'eductible et qu'il en est de m\^eme de la repr\'esentation $\pi(\rho',A',B',\zeta')  \times <\rho\vert\,\vert^{-\zeta A}, \cdots, \rho\vert\,\vert^{\zeta B}>$. Dans tout ce qui suit, quand on a fix\'e $(\rho,A,B,\zeta)\in Jord(\psi)$ on note $\psi'$ le morphisme qui se d\'eduit de $\psi$ en ''enlevant'' la repr\'esentation de $W_{F}\times SL(2,{\mathbb C}),\times SL(2,{\mathbb C})$ correspondant \`a $(\rho,A,B,\zeta)$; $\psi'$ est donc \`a valeurs dans $GL(N-d_{\rho}ab,{\mathbb C})$ (o\`u $(\rho,a,b)$ correspond \`a $(\rho,A,B,\zeta)$). Une \'ecriture $(\psi',(\rho,A',B',\zeta'))$ repr\'esente le morphisme $\psi''$ tel que $Jord(\psi'')=Jord(\psi')\cup (\rho,A',B',\zeta))$.

La d\'efinition de l'action de $\theta$ se fait en plusieurs  temps, c'est un point vraiment d\'elicat. D'abord on suppose que pour tout $(\rho,A,B,\zeta)\in Jord(\psi)$, $A=B=0$; dans ce cas $\pi(\psi)$ est une induite de cuspidales $\theta$-invariantes toutes distinctes; avec les op\'erateurs d'entrelacement, on d\'efinit une action de $\theta$. On  la normalise \`a la  Whittaker comme expliqu\'e dans l'introduction.

Ensuite on fait les constructions par r\'ecurrence \`a partir de ce choix. 
Supposons pour le moment, en plus, que pour tout $(\rho,A',B',\zeta')\in Jord_{\rho}(\psi)$, $A'=B'$; dans les notations initiales, c'est la condition $inf(a,b)=1$. Ici on note $(\rho,A_{0},B_{0},\zeta_{0})$ l'\'el\'ement de $Jord_{\rho}(\psi)$ tel que $B_{0}$ soit minimum. Il y a encore 2 cas \`a distinguer.

Le cas o\`u $B_{0}\neq 0$;
en tant que repr\'esentation de $GL(N)$, on v\'erifie que $\pi(\psi)$ est l'unique sous-module 
irr\'eductible de l'induite:
$$
\tau:=\rho\vert\,\vert^{\zeta B_{0}} \times \pi(\psi',(\rho,B_{0}-1,B_{0}-1,\zeta))\times \pi(\psi') \times \rho\vert\,\vert^{-\zeta B_{0}}
$$
$(\rho,B_{0}-1,B_{0}-1,\zeta)$ n'intervient pas si $B_{0}=1/2$; en effet que $\pi(\psi)$ soit un sous-module de $\tau$ r\'esulte imm\'ediatement des d\'efinitions et le point est de v\'erifier que l'induite de droite a un unique sous-module irr\'eductible. On v\'erifie cela par r\'eciprocit\'e de Frobenius. On utilise la notation $Jac^\theta_{x}$ de \ref{notationdujac} relative \`a notre $\rho$ fix\'e. Et il suffit de v\'erifier que $$Jac^\theta_{\zeta B_{0}}\tau=\pi(\psi',(\rho,B_{0}-1,B_{0}-1,\zeta)).$$ Or  on a $Jac_{x} \pi(\psi')\neq 0$ seulement si $x$ est de la forme $\zeta' B'$ avec $(\rho,A',B',\zeta')\in Jord(\psi')$; or pour un tel $\zeta' B'$, on a $\zeta' B' \notin [\zeta_{0} B_{0},-\zeta _{0}B_{0}]$ par l'hypoth\`ese sur la restriction de $\psi$ \`a la diagonale. Donc $Jac_{\zeta B_{0}}\pi(\psi',(\rho,B_{0}-1,B_{0}-1,\zeta))=0$. On a une assertion analogue pour calculer 
$Jac_{-\zeta B_{0}}$ ''\`a droite'', c'est \`a dire que l'on projette le module de Jacquet sur l'espace des repr\'esentations de la forme $\sigma'\otimes \rho\vert\,\vert^{-\zeta B_{0}}$.  Et les formules standard montrent alors l'assertion cherch\'ee. 

D'o\`u, en \'ecrivant les d\'efinitions:
$$\pi(\psi)\hookrightarrow
\rho\vert\,\vert^{\zeta B_{0}} \times \pi(\rho,B_{0}-1,B_{0}-1)\times \pi(\psi')  \times \rho\vert\,\vert^{-\zeta B_{0}}.
$$
L'induite \'ecrite a une action de $\theta$ d\`es que $\pi(\rho,B_{0}-1,B_{0}-1)\times \pi(\psi') $ en a une et l'action sur l'induite ne d\'epend que de l'action sur $\pi(\rho,B_{0}-1,B_{0}-1)\times \pi(\psi')$. Par unicit\'e du sous-module, on en d\'eduit une action sur $\pi(\psi)$; c'est ainsi que l'on met une action sur $\pi(\psi)$ par induction. 

On suppose maintenant que $B_{0}=0$ et  qu'il existe $(\rho,B_{1},B_{1},\zeta_{1}) \in Jord(\psi)$ avec $B_{1}\neq 0$. On fixe $B_{1}$ minimum avec cette propri\'et\'e et on note $\psi_{1}$ le morphisme qui se d\'eduit de $\psi$ en enlevant les 2 blocs, celui relatif \`a $B_{0}$ et celui relatif \`a $B_{1}$. Montrons d'abord l'inclusion (cf \ref{socle} pour la notation):
$$\pi(\psi) \hookrightarrow <\rho\vert\,\vert^{\zeta_{1}B_{1}}, \cdots, \rho> \times \pi(\psi_{1}) \times <\rho, \cdots, \rho\vert\,\vert^{-\zeta_{1}B_{1}}> \eqno(*)
. $$Or par minimalit\'e de $B_{1}$, l'entrelacement:
$$
\pi(\psi_{1}) \times <\rho, \cdots, \rho\vert\,\vert^{-\zeta_{1}B_{1}}> \rightarrow   <\rho, \cdots, \rho\vert\,\vert^{-\zeta_{1}B_{1}}> \times \pi(\psi_{1})
$$
est un isomorphisme. Comme
$$
\pi(\psi)\simeq \rho \times <\rho\vert\,\vert^{\zeta_{1}B_{1}}, \cdots, \rho\vert\,\vert^{-\zeta _{1}B_{1}}> \times \pi(\psi_{1})
$$
le seul point \`a v\'erifier est que $$ \rho \times <\rho\vert\,\vert^{\zeta_{1}B_{1}}, \cdots, \rho\vert\,\vert^{-\zeta _{1}B_{1}}> \hookrightarrow <\rho\vert\,\vert^{\zeta_{1}B_{1}}, \cdots, \rho>\times <\rho, \cdots, \rho\vert\,\vert^{-\zeta_{1}B_{1}}>$$
ce qui est un cas particulier des r\'esultats de Zelevinsky.

 L'induite (*) contient m\^eme $\pi(\psi)$ avec multiplicit\'e 1 en tant que sous-quotient par un calcul \'el\'ementaire de module de Jacquet (ici on utilise le fait que $\psi\circ \Delta$ est sans multiplicit\'es). Ainsi toute action de $\theta$ sur $\pi(\psi_{1})$ donne une action de $\theta$ sur l'induite et par restriction une action de $\theta$ sur $\pi(\psi)$. Cela termine les d\'efinitions dans le cas \'el\'ementaire.

Pour la suite, on a besoin de remarquer que $$
<\rho, \cdots, \rho\vert\,\vert^{-\zeta_{1}B_{1}}>\times  \pi(\psi_{1})\times <\rho\vert\,\vert^{\zeta_{1}B_{1}}, \cdots, \rho>. \eqno(**)
$$admet $\pi(\psi)$ comme unique quotient irr\'eductible: ce r\'esultat est dual du pr\'ec\'edent et r\'esulte donc du fait que $\pi(\psi)$ est autoduale. Toute action de $\theta$ sur $\pi(\psi_{1})$ se prolonge donc canoniquement \`a (**) et par restriction donne une action de $\theta$ sur $\pi(\psi)$. On montre que les 2 actions ainsi d\'efinies co\"{\i}ncident: pour cela on construit un entrelacement pour $s\in {\mathbb C}$ voisin de $0$:
$$
<\rho, \cdots, \rho\vert\,\vert^{-\zeta_{1}B_{1}}>\vert\,\vert^s\times  \pi(\psi_{1})\times <\rho\vert\,\vert^{\zeta_{1}B_{1}}, \cdots, \rho>\vert\,\vert^{-s}
$$
dans 
$$
<\rho\vert\,\vert^{\zeta_{1}B_{1}}, \cdots, \rho>\vert\,\vert^{-s}\times  \pi(\psi_{1})\times<\rho, \cdots, \rho\vert\,\vert^{-\zeta_{1}B_{1}}>\vert\,\vert^s.$$
D'apr\`es \ref{entrelacement}, les op\'erateurs d'entrelacement standard permettent de construire un tel entrelacement compatible aux prolongements canonique de l'action de $\theta$ sur $\pi(\psi_{1})$; en $s=0$, l'entrelacement n'a pas de raison d'\^etre holormorphe non nul mais on peut toujours obtenir cela \`a condition de multiplier par $s^z$ avec $z\in {\mathbb Z}$ convenable. On peut alors \'evaluer en $s=0$; l'image est  n\'ecessairement \'egale \`a $\pi(\psi)$ et on obtient l'entrelacement des actions de $\theta$ comme annonc\'e.

\

Supposons maintenant que $Jord(\psi)$ contient un quadruplet $(\rho,A,B,\zeta)$ avec $A>B$. Ici on v\'erifie que $\pi(\psi)$, comme repr\'esentation de $GL(N,F)$, est l'unique sous-module irr\'eductible de l'induite:
$$
\sigma:=<\rho\vert\,\vert^{\zeta B}, \cdots, \rho\vert\,\vert^{-\zeta A}>\times  \pi(\rho,A-1,B+1,\zeta) \times
<\rho\vert\,\vert^{\zeta A}, \cdots, \rho\vert\,\vert^{-\zeta B}>\times \pi(\psi'), \eqno(1)
$$o\`u $\pi(\rho,A-1,B+1,\zeta)$ est la repr\'esentation triviale de $GL(0,F)$ si $A=B+1$. L'assertion est un peu plus difficile que ci-dessus. En effet, il faut d'abord montrer avec les notations de \ref{notationdujac} que  $$Jac_{\zeta B, \cdots,  -\zeta A}\sigma= \pi(\rho,A-1,B+1,\zeta) \times
<\rho\vert\,\vert^{\zeta A}, \cdots, \rho\vert\,\vert^{-\zeta B}>\times \pi(\psi'). \eqno(2)$$
Pour cela, on v\'erifie que si ce n'est pas vrai il existe un sous-segment $[x,-\zeta A]$ de $[\zeta B, -\zeta A]$ tel que $$Jac_{x, \cdots, -\zeta A}\biggl(\pi(\rho,A-1,B+1,\zeta) \times
<\rho\vert\,\vert^{\zeta A}, \cdots, \rho\vert\,\vert^{-\zeta B}>\times \pi(\psi')\biggr)\neq 0.$$ Or $\rho\vert\,\vert^{-\zeta A}$ ne se trouve pas dans le support cuspidal de $\pi(\rho,A-1,B+1,\zeta) \times
<\rho\vert\,\vert^{\zeta A}, \cdots, \rho\vert\,\vert^{-\zeta B}>$, et il faudrait donc encore qu'il existe un sous-segment $[x',-\zeta A]$ de $[\zeta B,-\zeta A]$ tel que $Jac_{x',\cdots, -\zeta A}\pi(\psi')\neq 0$. On r\'ecrit $\pi(\psi')$ comme une induite de repr\'esentations de la forme $\pi(\rho',A',B',\zeta')$. Et quitte \`a reduire encore $[x',-\zeta A]$ il devrait exister $(\rho',A',B',\zeta') \in Jord(\psi')$ tel que $Jac_{x',\cdots,-\zeta A}\pi(\rho',A',B',\zeta')\neq 0$; cela entra\^{\i}ne $\rho\simeq \rho'$ et $x'=\zeta' B'$, en particulier $B'\leq A$; par l'hypoth\`ese de restriction discr\`ete \`a la diagonale, cela assure que  $A'$  v\'erifie $A'<B$ et $\rho\vert\,\vert^{-\zeta A}$ n'est pas dans le support cuspidal de $\pi(\rho,A',B',\zeta')$ ce qui donne la contradiction. La d\'emonstration n'est pas finie; on utilise le fait que dans (1) et dans (2) on peut mettre $<\rho\vert\,\vert^{\zeta A}, \cdots, \rho\vert\,\vert^{-\zeta B}>$ en derni\`ere position et par exemple en conjuguant par $\theta$, il suffit de d\'emontrer que l'induite $$Jac_{\zeta B, \cdots, -\zeta A}<\rho\vert\,\vert^{\zeta B}, \cdots, \rho\vert\,\vert^{-\zeta A}> \times \pi(\rho,A-1,B+1,\zeta) \times \pi(\psi')= \pi(\rho,A-1,B+1,\zeta) \times \pi(\psi').$$
C'est clairement la m\^eme d\'emonstration qui prouve l'assertion.

On a l'isomorphisme $$<\rho\vert\,\vert^{\zeta A}, \cdots, \rho\vert\,\vert^{-\zeta B}>\times \pi(\psi')\simeq \pi(\psi') \times <\rho\vert\,\vert^{\zeta A}, \cdots, \rho\vert\,\vert^{-\zeta B}>,$$un tel isomorphisme se construit avec des op\'erateurs d'entrelacement standard rendus holomorphes en $0$. En effet on plonge $\pi(\psi')$ dans une induite de la forme $$\times_{(\rho',A',B',\zeta')\in Jord(\psi')}\times _{k \in [0, A'-B']}<\rho'\vert\,\vert^{ \zeta' (B'+k)}, \cdots, \rho\vert\,\vert^{-\zeta' (A'-k)}>$$ et il suffit de d\'emontrer que pour tout $(\rho',A',B',\zeta'),k$ comme ci-dessus:
$$
<\rho\vert\,\vert^{\zeta A}, \cdots, \rho\vert\,\vert^{-\zeta B}>\times <\rho'\vert\,\vert^{ \zeta' (B'+k)}, \cdots, \rho\vert\,\vert^{-\zeta' (A'-k)}> \simeq $$
$$<\rho'\vert\,\vert^{ \zeta' (B'+k)}, \cdots, \rho\vert\,\vert^{-\zeta' (A'-k)}>\times <\rho\vert\,\vert^{\zeta A}, \cdots, \rho\vert\,\vert^{-\zeta B}>.
$$
L'hypoth\`ese que $\psi$ est de restriction discr\`ete \`a la diagonale assure que soit $A'<B$ soit $A<B'$ d'o\`u aussi soit $A'-k<B$ soit $A< B'+k$. Ainsi le segment $[\zeta A,-\zeta B]$ soit contient le segment $[\zeta' (B'+k),-\zeta'(A'-k)]$ soit est inclus dans ce segment m\^eme si les propri\'et\'es de croissance des segments ne sont pas les m\^emes. On construit alors l'isomorphisme comme restriction d'un isomorphisme o\`u on garde la repr\'esentation correspondant au segment le plus long et on remplace l'autre par l'induite \'evidente; par exemple supposons que $A<B'$,  pour tout $x\in [\zeta A,-\zeta B]$, $x \in [\zeta'(B'+k), -\zeta' (A'-k)]$ et l'induite $\rho\vert\,\vert^x\times <\rho'\vert\,\vert^{ \zeta' (B'+k)}, \cdots, \rho\vert\,\vert^{-\zeta' (A'-k)}>$ est irr\'eductibile d'apr\`es les r\'esultats de Zelevinsky. D'o\`u$$
\times _{x\in [\zeta A,-\zeta B]} \rho \vert\,\vert^x \times <\rho'\vert\,\vert^{ \zeta' (B'+k)}, \cdots, \rho\vert\,\vert^{-\zeta' (A'-k)}>\simeq <\rho'\vert\,\vert^{ \zeta' (B'+k)}, \cdots, \rho\vert\,\vert^{-\zeta' (A'-k)}>\times 
 _{x\in [\zeta A,-\zeta B]} \rho \vert\,\vert^x,
 $$isomorphisme que l'on peut construire avec les op\'erateurs d'entrelacements standard rendu holomorphes en $s=0$ (on a remplac\'e $<\rho'\vert\,\vert^{ \zeta' (B'+k)}, \cdots, \rho\vert\,\vert^{-\zeta' (A'-k)}>$ par $<\rho'\vert\,\vert^{ \zeta' (B'+k)}, \cdots, \rho\vert\,\vert^{-\zeta' (A'-k)}>\vert\,\vert^s$). D'o\`u l'assertion.

 On note $\psi_{-}$ le morphisme qui se d\'eduit de $\psi$ en rempl\c{c}ant $(\rho,A,B,\zeta)$ par $(\rho,A-1,B+1,\zeta)$ et on obtient donc que $\pi(\psi)$ est l'unique sous-module irr\'eductible de l'induite:
$$
<\rho\vert\,\vert^{\zeta B}, \cdots, \rho\vert\,\vert^{-\zeta A}> \times \pi(\psi_{-}) \times <\rho\vert\,\vert^{\zeta A}, \cdots, \rho\vert\,\vert^{-\zeta B}>.
$$
On remarque que le morphisme, $\psi_{-}$ a la m\^eme propri\'et\'e de restriction discr\`ete \`a la diagonale que $\psi$. Ainsi on peut munir par induction l'induite  que l'on vient d'\'ecrire d'une action de $\theta$ qui ne d\'epend que de l'action de $\theta$ sur $\pi(\psi_{-})$. Par unicit\'e du sous-module, cela donne une action de $\theta$ sur $\pi(\psi)$. 

\

\bf Lemme. \sl L'action que l'on vient de mettre sur $\pi(\psi)$ est ind\'ependante du choix de $(\rho,A,B,\zeta)$ tel que $A>B$.\rm

\

On fixe $(\rho,A,B,\zeta)$ et $(\rho',A',B',\zeta')$ dans $Jord(\psi)$ tel que $A>B$ et $A'>B'$. On note $\psi''$ le morphisme qui se d\'eduit de $\psi$ en enlevant ces 2 blocs.
On repr\'esente $\pi(\psi)$ comme sous-module irr\'eductible de l'induite:
$$
<\rho\vert\,\vert^{\zeta B}, \cdots, \rho\vert\,\vert^{-\zeta A}>\times <\rho'\vert\,\vert^{\zeta' B'}, \cdots, \rho'\vert\,\vert^{-\zeta' A'}>\times \pi(\psi'',(\rho,A-1,B+1,\zeta),(\rho',A'-1,B'+1,\zeta'))$$
$$
\times 
<\rho'\vert\,\vert^{\zeta' A'}, \cdots, \rho'\vert\,\vert^{-\zeta' B'}>\times <\rho\vert\,\vert^{\zeta B}, \cdots, \rho\vert\,\vert^{-\zeta A}>\eqno(3)
$$
et comme sous-module de l'induite:
$$
<\rho'\vert\,\vert^{\zeta' B'}, \cdots, \rho'\vert\,\vert^{-\zeta' A'}>\times <\rho\vert\,\vert^{\zeta B}, \cdots, \rho\vert\,\vert^{-\zeta A}>\times \pi(\psi'',(\rho,A-1,B+1,\zeta),(\rho',A'-1,B'+1,\zeta'))$$
$$
<\rho\vert\,\vert^{\zeta B}, \cdots, \rho\vert\,\vert^{-\zeta A}>\times 
<\rho'\vert\,\vert^{\zeta' A'}, \cdots, \rho'\vert\,\vert^{-\zeta' B'}>\eqno(4).
$$
Et le point est de montrer que les prolongements canoniques de l'action de $\theta$ sur $\pi(\psi'',(\rho,A-1,B+1,\zeta),(\rho',A'-1,B'+1,\zeta'))$ \`a (3) et \`a (4) ont la m\^eme restriction \`a $\pi(\psi)$. On remplace par exemple dans (3), $<\rho\vert\,\vert^{\zeta B}, \cdots, \rho\vert\,\vert^{-\zeta A}>$ par $<\rho\vert\,\vert^{\zeta B}, \cdots, \rho\vert\,\vert^{-\zeta A}>\vert\,\vert^s$ et $<\rho\vert\,\vert^{\zeta B}, \cdots, \rho\vert\,\vert^{-\zeta A}>$ en $<\rho\vert\,\vert^{\zeta B}, \cdots, \rho\vert\,\vert^{-\zeta A}>\vert\,\vert^{-s}$. On peut alors utiliser l'op\'erateur d'entrelacement standard 
$$
<\rho\vert\,\vert^{\zeta B}, \cdots, \rho\vert\,\vert^{-\zeta A}>\vert\,\vert^s \times <\rho'\vert\,\vert^{\zeta' B'}, \cdots, \rho'\vert\,\vert^{-\zeta' A'}>\rightarrow
$$
$$
<\rho'\vert\,\vert^{\zeta' B'}, \cdots, \rho'\vert\,\vert^{-\zeta' A'}>\times <\rho\vert\,\vert^{\zeta B}, \cdots, \rho\vert\,\vert^{-\zeta A}>\vert\,\vert^s
$$
Compos\'e avec son analogue pour \'echanger les 2 derniers facteurs. On vient de v\'erifier ci-dessus, qu'en le normalisant correctement, il devient un isomorphisme en $s=0$ (la normalisation consiste simplement \`a le rendre holomorphe en $s=0$). Il est clair comme en \ref{entrelacement} que l'on a ainsi construit un entrelacement entre le prolongement canonique de l'action de $\theta$ sur $\pi(\psi'',(\rho,A-1,B+1,\zeta),(\rho',A'-1,B'+1,\zeta'))$ \`a l'induite:
$$
<\rho\vert\,\vert^{\zeta B}, \cdots, \rho\vert\,\vert^{-\zeta A}>\vert\,\vert^s \times \times <\rho'\vert\,\vert^{\zeta' B'}, \cdots, \rho'\vert\,\vert^{-\zeta' A'}>\times \pi(\psi'',(\rho,A-1,B+1,\zeta),(\rho',A'-1,B'+1,\zeta'))$$
$$
\times 
<\rho'\vert\,\vert^{\zeta' A'}, \cdots, \rho'\vert\,\vert^{-\zeta' B'}>\times <\rho\vert\,\vert^{\zeta B}, \cdots, \rho\vert\,\vert^{-\zeta A}>\vert\,\vert^{-s}
$$
et le prolongement canonique de cette m\^eme action \`a l'induite:
$$
<\rho'\vert\,\vert^{\zeta' B'}, \cdots, \rho'\vert\,\vert^{-\zeta' A'}>\times <\rho\vert\,\vert^{\zeta B}, \cdots, \rho\vert\,\vert^{-\zeta A}>\vert\,\vert^s\times \pi(\psi'',(\rho,A-1,B+1,\zeta),(\rho',A'-1,B'+1,\zeta'))$$
$$
<\rho\vert\,\vert^{\zeta B}, \cdots, \rho\vert\,\vert^{-\zeta A}>\vert\,\vert^{-s}\times 
<\rho'\vert\,\vert^{\zeta' A'}, \cdots, \rho'\vert\,\vert^{-\zeta' B'}>.
$$En faisant $s=0$ on obtient le r\'esultat cherch\'e.

\

\bf Lemme. \sl L'action de $\theta$ mise sur $\pi(\psi)$ est ind\'ependante du choix de $\rho$.\rm

\
En effet fixons, $\rho,\rho'$ tel que $Jord_{\rho}(\psi)\neq \emptyset$ et $Jord_{\rho'}(\psi)\neq \emptyset$. On suppose que $\rho\not\simeq \rho'$.  Pour tout $x,y\in {\mathbb R}$, l'induite $\rho\vert\,\vert^x\times \rho'\vert\,\vert^y$ est irr\'eductible. Il suffit de faire commuter comme ci-dessus les d\'efinitions donn\'ees en utilisant $\rho$ et celles donn\'ees en utilisant $\rho'$ pour obtenir l'asssertion.

\subsection{Quelques propri\'et\'es de l'action de $\theta$\label{proprietes}}
Dans ce paragraphe on suppose encore que la restriction de $\psi \circ \Delta$ est sans multiplicit\'es vue comme repr\'esentation de $W_{F}\times SL(2,{\mathbb C})$. Soit $(\rho,A,B,\zeta) \in Jord(\psi)$ avec $B\geq 1$; on a d\'efini $\psi'$ dans \ref{independance}. On suppose ici que pour tout $(\rho',A',B',\zeta')\in Jord(\psi')$ avec $\rho'\simeq \rho$, $A' \neq B-1$. Le morphisme $\psi_{-}$ qui ici est d\'efini en rempla\c{c}ant $(\rho,A,B,\zeta)$ par $(\rho,A-1,B-1,\zeta)$  v\'erifie encore que $\psi_{-}\circ \Delta$ est sans multiplicit\'es. 

\

\bf Lemme. \sl Sous l'hypoth\`ese ci-dessus, l'induite $<\rho\vert\,\vert^{\zeta B},  \cdots, \rho\vert\,\vert^{\zeta A}> \times \pi(\psi_{-}) \times < \rho\vert\,\vert^{-\zeta A}, \cdots, \rho\vert\,\vert^{-\zeta B}>$ a un unique sous-module irr\'eductible et ce sous-module est $\pi(\psi)$. L'action de $\theta$ sur cette induite qui se d\'eduit naturellement de celle sur $\pi(\psi_{-})$ induit par restriction l'action de $\theta$ d\'ej\`a mise sur $\pi(\psi)$.\rm

\

L'hypoth\`ese est indispensable pour avoir un tel lemme.
On suppose d'abord que $A=B$;  dans ce cas c'est facile, on applique la d\'efinition et les commutations sont \'evidentes car on a \`a chaque fois irr\'eductibilit\'e gr\^ace \`a l'hypoth\`ese. On suppose donc que $A>B$.

L'unicit\'e du sous-module est un calcul de module de Jacquet plus simple que celui fait en \ref{independance}, on ne le refait pas. On \'ecrit la suite d'induites (avec la convention faite dans les d\'efinitions si $A-1=B$ pour interpr\'eter la d\'efinition de $\pi(\rho,A-2,B,\zeta)$ comme la repr\'esentation triviale de $GL(0)$):
$$
<\rho\vert\,\vert^{\zeta B},  \cdots, \rho\vert\,\vert^{\zeta A}> \times \pi(\psi_{-}) \times < \rho\vert\,\vert^{-\zeta A}, \cdots, \rho\vert\,\vert^{-\zeta B}>\hookrightarrow\eqno(1)
$$
$$
<\rho\vert\,\vert^{\zeta B},  \cdots, \rho\vert\,\vert^{\zeta A}> $$
$$ \times <\rho\vert\,\vert^{\zeta(B-1), \cdots, -\zeta(A-1)}>\times \pi(\rho,A-2,B,\zeta) \times \pi(\psi') \times <\rho\vert\,\vert^{\zeta(A-1)}, \cdots, \rho \vert\,\vert^{-\zeta(B-1)}>$$
$$ \times <\rho\vert\,\vert^{-\zeta A}, \cdots, \rho\vert\,\vert^{-\zeta B}>$$
$$
\hookrightarrow
$$
$$
\rho\vert\,\vert^{\zeta B}\times <\rho\vert\,\vert^{\zeta(B-1), \cdots, -\zeta(A-1)} >\times
<\rho\vert\,\vert^{\zeta (B+1)},  \cdots, \rho\vert\,\vert^{\zeta A}> $$
$$ \times  \pi(\rho,A-2,B,\zeta) \times \pi(\psi')$$
$$ \times 
<\rho\vert\,\vert^{-\zeta A}, \cdots, \rho\vert\,\vert^{-\zeta (B+1)}>\times <\rho\vert\,\vert^{\zeta(A-1)}, \cdots, \rho \vert\,\vert^{-\zeta(B-1)}> \times  \rho\vert\,\vert^{-\zeta B}.
$$
$$
\hookrightarrow
$$
$$
\rho\vert\,\vert^{\zeta B}\times <\rho\vert\,\vert^{\zeta(B-1), \cdots, -\zeta(A-1)} >\times
<\rho\vert\,\vert^{\zeta (B+1)},  \cdots, \rho\vert\,\vert^{\zeta (A-1)}>\times \rho\vert\,\vert^{\zeta A} $$
$$ \times  \pi(\rho,A-2,B,\zeta) \times \pi(\psi')$$
$$ \times \rho\vert\,\vert^{-\zeta A}\times
<\rho\vert\,\vert^{-\zeta (A-1)}, \cdots, \rho\vert\,\vert^{-\zeta (B+1)}>\times <\rho\vert\,\vert^{\zeta(A-1)}, \cdots, \rho \vert\,\vert^{-\zeta(B-1)}> \times  \rho\vert\,\vert^{-\zeta B}. \eqno(2)
$$
Il y a des actions de $\theta$ sur toutes ces induites qui ne d\'ependent que de l'action de $\theta$ sur $\pi(\rho,A-2,B,\zeta) \times \pi(\psi')$ et les inclusions sont compatibles aux actions de $\theta$. On sait que $A\geq B+1 \geq 2$ d'o\`u a fortiori $A>1/2$donc les induites ci-dessous sont irr\'eductibles et isomorphes: $$\rho\vert\,\vert^{-\zeta A}\times \pi(\rho,A-2,B,\zeta)\times \pi(\psi') \times \rho\vert\,\vert^{\zeta A}\simeq 
\rho\vert\,\vert^{\zeta A}\times \pi(\rho,A-2,B,\zeta)\times \pi(\psi') \times \rho\vert\,\vert^{-\zeta A}.$$ On peut continuer (2) par les isomorphismes:
$$
(2)\simeq 
\rho\vert\,\vert^{\zeta B}\times <\rho\vert\,\vert^{\zeta(B-1), \cdots, -\zeta(A-1)} >\times \rho\vert\,\vert^{-\zeta A}
$$
$$
\times 
<\rho\vert\,\vert^{\zeta (B+1)},  \cdots, \rho\vert\,\vert^{\zeta (A-1)}>\times \pi(\rho,A-2,B,\zeta) \times \pi(\psi')\times <\rho\vert\,\vert^{-\zeta (A-1)}, \cdots, \rho\vert\,\vert^{-\zeta (B+1)}>
$$
$$
\times \rho\vert\,\vert^{-\zeta A}\times <\rho\vert\,\vert^{\zeta(A-1)}, \cdots, \rho \vert\,\vert^{-\zeta(B-1)}> \times  \rho\vert\,\vert^{-\zeta B}.
$$
En admettant le lemme par recurrence, on voit que (2) contient $$
 \rho\vert\,\vert^{\zeta B}\times <\rho\vert\,\vert^{\zeta(B-1), \cdots, -\zeta(A-1)} >\times \rho\vert\,\vert^{-\zeta A}\times $$
 $$
  \pi(\rho,A-1,B+1,\zeta) \times \pi(\psi')$$
  $$
   \times \rho\vert\,\vert^{\zeta A}\times <\rho\vert\,\vert^{\zeta(A-1)}, \cdots, \rho \vert\,\vert^{-\zeta(B-1)}> \times  \rho\vert\,\vert^{-\zeta B}.
  $$L'action sur l'induite du milieu est celle que nous avons mise sur $\pi(\psi'')$ o\`u $\psi''$ est le morphisme qui se d\'eduit de $\psi'$ en ajoutant le bloc $(\rho,A-1,B+1,\zeta)$; d'o\`u $\pi(\psi'')=  \pi(\rho,A-1,B+1,\zeta) \times \pi(\psi')$.
  En tant que repr\'esentation de $\tilde{G}(N)$ cette induite contient encore:
  $$
   <\rho\vert\,\vert^{\zeta B}, \cdots, \rho{-\zeta A} >\times
\pi(\psi'') \times <\rho\vert\,\vert^{\zeta A}, \cdots, \rho \vert\,\vert^{-\zeta B}> .\eqno(3)
  $$
  Or (3) a encore comme unique sous-module irr\'eductible $\pi(\psi)$ avec la ''bonne'' action de $\theta$, par d\'efinition. L'action de $\theta$ sur $\pi(\psi)$ en tant que sous-module de (1) est donc la m\^eme que celle d\'ej\`a mise sur $\pi(\psi)$. Cela prouve le lemme.
  \vskip 0.5cm
  \noindent
  \bf Corollaire. \sl Avec les hypoth\`eses du lemme pr\'ec\'edent, soit $C\in [B,A]$ alors $Jac^\theta_{\zeta B, \cdots ,\zeta C}\pi(\psi)$ est l'unique sous-module irr\'eductible de l'induite
  $$
  <\rho\vert\,\vert^{\zeta(C+1)}, \cdots, \rho\vert\,\vert^{\zeta A}>\times \pi(\psi_{-}) \times <\rho\vert\,\vert^{-\zeta A}, \cdots, \rho\vert\,\vert^{-\zeta (C+1)}>.
  $$
  \rm
  Le lemme pr\'ec\'edent permet d'\'ecrire $\pi(\psi)$ comme sous-module de l'induite:
  $$
  <\rho\vert\,\vert^{\zeta B}, \cdots, \rho\vert\,\vert^{\zeta C}> \times$$
  $$ <<\rho\vert\,\vert^{\zeta (C+1)}, \cdots, \rho\vert\,\vert^{\zeta A}> \times \pi(\psi_{-}) \times <\rho\vert\,\vert^{-\zeta A},\cdots,\rho\vert\,\vert^{-\zeta(C+1)}>> $$
  $$\times <\rho\vert\,\vert^{-\zeta C}, \cdots,\rho\vert\,\vert^{-\zeta B}>,
  $$les $<$ $>$ ext\'erieurs dans l'induite du milieu disent que l'on prend l'unique sous-module irr\'eductible. Et  l'action de $\theta$ prolonge naturellement celle de $\theta$ sur l'induite du milieu.
  Et le r\'esultat s'en d\'eduit.

\subsection{Remarque \label{remarque}}
Soient $B,A$ des demi-entiers tel que $A-B\in {\mathbb N}_{>0}$ et soit $\zeta=\pm 1$. L'induite
$$
<\rho\vert\,\vert^{\zeta B}, \cdots, \rho\vert\,\vert^{-\zeta A}> \times <\rho\vert\,\vert^{\zeta A}, \cdots, \rho\vert\,\vert^{-\zeta B}>\eqno(1)$$ est de longueur 2, sans multiplicit\'es et contient comme quotient l'induite irr\'eductible:
$$
<\rho\vert\,\vert^{\zeta B}, \cdots, \rho\vert\,\vert^{-\zeta B}> \times <\rho\vert\,\vert^{\zeta A}, \cdots, \rho\vert\,\vert^{-\zeta A}>.\eqno(2)
$$
Sur (1), on a une action de $\theta$ d\`es que $\rho\simeq \, ^\theta\rho$ et cette action est ind\'ependante du choix d'un op\'erateur entrela\c{c}ant $\rho$ et $^\theta\rho$. Sur (2), on vient de mettre une action de $\theta$ sous la m\^eme hypoth\`ese, en consid\'erant que (2) correspond au morphisme $\psi$ tel que $Jord(\psi)=\{(\rho,A,A,\zeta); (\rho,B,B,\zeta)\}$.

\

\bf Remarque: \sl l'action de $\theta$ mise sur (2) est la restriction de celle mise sur (1).\rm

\

C'est la d\'efinition si $B=0$ et si $B$ est entier, avec $B>0$, on se ram\`ene \`a ce cas en faisant une r\'ecurrence: on \'ecrit
$$
<\rho\vert\,\vert^{\zeta B}, \cdots, \rho\vert\,\vert^{-\zeta A}>\times <\rho\vert\,\vert^{\zeta A}, \cdots, \rho\vert\,\vert^{-\zeta B}>\hookrightarrow
$$
$$
\rho\vert\,\vert^{\zeta B}\times <\rho\vert\,\vert^{\zeta (B-1)}, \cdots, \rho\vert\,\vert^{-\zeta A}>\times <\rho\vert\,\vert^{\zeta A},\cdots,\rho\vert\,\vert^{-\zeta (B-1)}>\times \rho\vert\,\vert^{-\zeta B}.\eqno(3)
$$
Par r\'ecurrence on admet que la projection de (3) sur
$$
\rho\vert\,\vert^{\zeta B}\times <\rho\vert\,\vert^{\zeta (B-1)}, \cdots, \rho\vert\,\vert^{-\zeta (B-1)}>\times <\rho\vert\,\vert^{\zeta A}, \cdots, \rho\vert\,\vert^{-\zeta A}>\times \rho\vert\,\vert^{-\zeta B} \eqno(4)
$$
entrelace le prolongement canonique de l'action de $\theta$ mise sur l'induite du milieu. L'induite (4) contient avec multiplicit\'e 1 la repr\'esentation $\tau:=$ $<\rho\vert\,\vert^{\zeta B}, \cdots, \rho\vert\,\vert^{-\zeta B}>$ $\times $ $<\rho\vert\,\vert^{\zeta A}, \cdots, \rho\vert\,\vert^{-\zeta A}>$. Cette repr\'esentation est un sous-module de (4)  et l'action de $\theta$ que nous y avons mise est la restriction de l'action de $\theta$ sur (4). Par un calcul de module de Jacquet, on v\'erifie que $\tau$ est l'unique sous-module irr\'eductible de (4). On v\'erifie aussi que $\tau$ intervient avec multiplicit\'e 1 comme sous-quotient de (3). On a donc construit une application $\theta$-\'equivariante de (1) dans (4);  si elle \'etait nulle, $\tau$ qui intervient dans (1) et dans l'image de (3) dans (4), interviendrait avec multiplicit\'e au moins 2 dans (3). Ceci est exclu et  $\tau$ est donc l'image de l'application construite. La $\theta$-\'equivariance prouve l'assertion.

Il reste \`a voir le cas o\`u $B$ est demi-entier; la d\'emonstration ci-dessus ram\`ene au cas o\`u $B=1/2$, cas que nous allons traiter de fa\c{c}on totalement similaire: on construit comme ci-dessus
$$
<\rho\vert\,\vert^{\zeta 1/2}, \cdots, \rho\vert\,\vert^{-\zeta A}>\times <\rho\vert\,\vert^{\zeta A}, \cdots, \rho\vert\,\vert^{-\zeta 1/2}>
\hookrightarrow
$$
$$
\rho\vert\,\vert^{\zeta 1/2}\times <\rho\vert\,\vert^{-\zeta 1/2}, \cdots, \rho\vert\,\vert^{-\zeta A}> \times <\rho\vert\,\vert^{\zeta A}, \cdots, \rho\vert\,\vert^{\zeta 1/2}> \times \rho\vert\,\vert^{-\zeta 1/2}\rightarrow
$$
$$
\rho\vert\,\vert^{\zeta 1/2}\times <\rho\vert\,\vert^{\zeta A}, \cdots, \rho\vert\,\vert^{-\zeta A}> \times \rho\vert\,\vert^{-\zeta 1/2}.
$$
Le point est de d\'emontrer que cette application est $\theta$-\'equivariante. Par d\'efinition, il n'y a plus qu'\`a d\'emontrer la $\theta$-\'equivariance de l'application 
$$
<\rho\vert\,\vert^{-\zeta 1/2}, \cdots, \rho\vert\,\vert^{-\zeta A}> \times <\rho\vert\,\vert^{\zeta A}, \cdots, \rho\vert\,\vert^{\zeta 1/2}> \rightarrow <\rho\vert\,\vert^{\zeta A}, \cdots, \rho\vert\,\vert^{-\zeta A}>.
$$
Or l'action de $\theta$ mise sur $<\rho\vert\,\vert^{\zeta A}, \cdots, \rho\vert\,\vert^{-\zeta A}>$ vient de l'inclusion
$$
<\rho\vert\,\vert^{\zeta A}, \cdots, \rho\vert\,\vert^{-\zeta A}>\hookrightarrow <\rho\vert\,\vert^{\zeta A}, \cdots, \rho\vert\,\vert^{\zeta 1/2}>\times <\rho\vert\,\vert^{-\zeta 1/2}, \cdots, \rho\vert\,\vert^{-\zeta A}>.
$$
Mais on a d\'ej\`a vu que la $\theta$-\'equivariance r\'esulte de \ref{entrelacement} en glissant un param\`etre $s\in {\mathbb C}$:
$$
<\rho\vert\,\vert^{-\zeta 1/2}, \cdots, \rho\vert\,\vert^{-\zeta A}>\vert\,\vert^s
\times 
<\rho\vert\,\vert^{\zeta A}, \cdots, \rho\vert\,\vert^{-\zeta A}>\vert\,\vert^{-s}\rightarrow
$$
$$
<\rho\vert\,\vert^{\zeta A}, \cdots, \rho\vert\,\vert^{-\zeta A}>\vert\,\vert^{-s} \times <\rho\vert\,\vert^{-\zeta 1/2}, \cdots, \rho\vert\,\vert^{-\zeta A}>\vert\,\vert^s.
$$
C'est une application $\theta$-\'equivariante qui, rendue holomorphe en $s=0$, donne l'\'egalit\'e des 2 actions de $\theta$ sur $<\rho\vert\,\vert^{\zeta A}, \cdots,\rho\vert\,\vert^{-\zeta A}>$. La fin de la d\'emonstration est comme ci-dessus.
\vskip 0.5cm
\noindent
\bf Remarque: \sl supposons que $B=0$; alors pour tout entier $A$,  $\pi(\rho, A,0,+)\simeq \pi(\rho, A,0,-)$.\rm
\vskip 0.5cm
\noindent
En tant que repr\'esentation de $GL(abd_{\rho},F)$ les 2 repr\'esentations \'ecrites sont isomorphes. Le probl\`eme est l'action de $\theta$.
On fait d'abord la d\'emonstration pour $A=1$; sur $\pi(\rho,1,0,+)$, $\theta$ agit par $\theta_{+}$ obtenu en restreignant l'action naturelle de $\theta$ sur l'induite:
$$
<\rho,\rho\vert\,\vert^{-1}>\times <\rho\vert\,\vert^{+1},\rho>.
$$
On inclut cette induite dans l'induite:
$$
\rho \times \rho\vert\,\vert^{-1}\times \rho\vert\,\vert^{+1}\times \rho\simeq \rho\times \rho'\vert\,\vert^{+1}\times \rho\vert\,\vert^{-1}\times \rho,
$$
l'isomorphisme pr\'eservant les actions naturelles de $\theta$. L'induite de droite contient encore comme sous-module l'induite $<\rho,\rho\vert\,\vert^{+1}>\times <\rho\vert\,\vert^{-1},\rho>$ munie de l'action naturelle de $\theta$ et cette induite contient $\pi(\rho,1,0,-)$ comme unique sous-module irr\'eductible. Pour conclure, il ne reste plus qu'\`a remarquer que l'isomorphisme du milieu a identifi\'e $\pi(\rho,1,0,+)$ et $\pi(\rho,1,0,-)$.

Le cas o\`u $A>1$ se traite par r\'ecurrence pour se ramener au cas ci-dessus.

 \section{``R\'esolution'' dans le groupe de Grothendieck}
\subsection{Le cas des repr\'esentations $\psi$ irr\'eductibles avec $inf(a,b)=2$ \label{simple}}
On fixe un triplet $(\rho,a,b)$ et son quadruplet associ\'e $(\rho,A,B,\zeta)$ cf. \ref{definitiondejord}. Dans ce paragraphe on a pos\'e $Sp(b,St(a,\rho))=\pi(\rho,A,B,\zeta)$ notation que nous allons utiliser. 
\vskip 0.5cm
\noindent 
Traitons d'abord le cas o\`u $A=B+1$, c'est-\`a-dire $inf(a,b)=2$. Alors $\pi(\rho,A,B,\zeta)$ est l'unique sous-module irr\'eductible de l'induite (on remplace syst\'ematiquement $A$ par $B+1$):
$$
<\rho\vert\,\vert^{\zeta B}, \cdots, \rho\vert\,\vert^{-\zeta (B+1)}> \times <\rho\vert\,\vert^{\zeta (B+1)}, \cdots, \rho\vert\,\vert^{-\zeta B}> .\eqno(1)
$$
Cette induite est en fait de longueur 2 avec comme quotient irr\'eductible la repr\'esentation induite
$$
<\rho\vert\,\vert^{\zeta (B+1)}, \cdots, \rho\vert\,\vert^{-\zeta(B+1)}>\times <\rho\vert\,\vert^{\zeta B}, \cdots, \rho\vert\,\vert^{-\zeta B}>. \eqno(2)
$$
On remarque d'ailleurs que (2) n'est autre que l'induite $\pi(\rho,A=B+1,B+1,\zeta) \times \pi(\rho,B,B,\zeta)$. Il faut comparer l'action de $\theta$ que nous avons mis sur (2) avec celle que l'on obtient en passant au quotient l'action de $\theta$ mise sur (1).  On a d\'emontr\'e en \ref{remarque} que ces 2 actions co\"{\i}ncidaient . Ainsi ici, $\pi(\rho,A,B,\zeta)$ est la diff\'erence des 2 induites celle \'ecrite en (1) moins celle \'ecrite en (2), diff\'erence dans le groupe de Grothendieck.

\subsection{Le cas des repr\'esentations $\psi$ irr\'eductibles\label{casquasigeneral}}
Dans ce paragraphe on va g\'en\'eraliser la formule de \ref{simple}. On reprend les notations $(\rho,A,B,\zeta)$ mais on suppose ici $A>B+1$.
On pose:
$$
\tilde{\pi}(\rho,A,B,\zeta):=$$
$$
\oplus_{C\in ]B,A]}(-1)^{A-C}<\rho\vert\,\vert^{\zeta B}, \cdots, \rho\vert\,\vert^{-\zeta C}>
\times
Jac^\theta_{\zeta (B+2), \cdots, \zeta C}\pi(\rho,A,B+2,\zeta)\times
<\rho\vert\,\vert^{\zeta C}, \cdots, \rho\vert\,\vert^{-\zeta B}>$$
$$
 \oplus (-1)^{[(A-B+1)/2]}\pi(\rho,A,B+1,\zeta)\times \pi(\rho, B,B,\zeta),$$
 o\`u par convention pour $C=B+1$, $Jac^\theta_{\zeta (B+2), \cdots, \zeta C}\pi(\rho,A,B+2,\zeta)$ est simplement $\pi(\rho,A,B+2,\zeta)$. Toutes les repr\'esentations intervenant ont \'et\'e munies d'une action de $\theta$ et $\tilde{\pi}(\rho,A,B,\zeta)$ est donc un \'el\'ement du groupe de Grothendieck des repr\'esentations de $\tilde{G}(abd_{\rho})$. Il n'est pas vrai en g\'en\'eral que $\pi(\rho,A,B,\zeta)=\tilde{\pi}(\rho,A,B,\zeta)$ mais on va montrer:
 \vskip 0.5cm
\noindent
\bf Proposition: \sl Dans sa description comme combinaison lin\'eaire de repr\'esentations irr\'eductibles de $\tilde{G}(abd_{\rho})$, $\tilde{\pi}(\rho,A,B,\zeta)$ contient exactement une repr\'esentation irr\'eductible qui reste irr\'eductible quand on la restreint \`a $GL(abd_{\rho},F)$; cette repr\'esentation est $\pi(\rho,A,B,\zeta)$ et elle intervient avec le coefficient $+1$. Plus \'el\'egamment, pour tout $g\in GL(n,F)$ $\theta$-semi-simple et r\'egulier, $$
tr(\tilde{\pi}(\rho,A,B,\zeta))(g,\theta)= tr \pi(\rho,A,B,\zeta)(g,\theta).
$$  
\rm
\subsubsection{Description des repr\'esentations dans la classification de Zelevinsky}
Avant de commencer la d\'emonstration, d\'ecrivons, pour $C\in ]B,A]$ la repr\'esentation $
Jac^\theta_{\zeta (B+2), \cdots, \zeta C}\pi(\rho,$ $A,B+2,\zeta)$ dans la classification de Zelevinsky.
Notons ${\cal T}_{A,B,\zeta}$ le tableau suivant dont les lignes sont des segments d\'ecroissants si $\zeta=+$ et croissants si $\zeta=-$ et les colonnes sont des segments ayant les propri\'et\'es de croissance oppos\'ees, en les consid\'erant du haut vers le bas.
$$
\begin{matrix}
\zeta B &\zeta (B-1)& \cdots &-\zeta A\cr
\vdots &\vdots&\vdots&\vdots \cr
\zeta A &\zeta (A-1)&\cdots &-\zeta B
\end{matrix}
$$Ce tableau a $sup(a,b)$ colonnes et $inf(a,b)$ lignes.

On peut consid\'erer le tableau analogue en rempla\c{c}ant $B$ par $B+2$ et en gardant $A$ et $\zeta$, c'est-\`a-dire ${\cal T}_{(A,B+2,\zeta)}$:
$$
\begin{matrix}
\zeta (B+2) &\zeta (B+1)& \cdots &-\zeta A\cr
\vdots &\vdots&\vdots&\vdots \cr
\zeta A &\zeta (A-1)&\cdots &-\zeta (B+2)
\end{matrix}
$$Ce tableau a $sup(a,b)+2$ colonnes et $inf(a,b)-2$ lignes.
Soit $C\in ]B,A]$; on suppose m\^eme que $C\geq B+2$ pour qu'il y ait quelque chose \`a d\'ecrire. Et notons ${\cal T}_{C}$ le tableau qui s'obtient \`a partir de celui \'ecrit ci-dessus mais en enlevant les $(C-B-1)$ premiers \'el\'ements de la premi\`ere colonne et les $(C-B-1)$ derniers \'el\'ements de la derni\`ere colonne; en particulier, on n'enl\`eve rien si $C=B+1$. Ecrivons ce que l'on obtient dans le cas o\`u $A+B\geq 2C$
$$
\begin{matrix}
&\zeta (B+1)& \cdots &\cdots &-\zeta A\cr
 &\vdots&\vdots&\vdots&\vdots \cr
\zeta (C+1) &\zeta (C)& \cdots&\cdots &-\zeta (A-C+B+1)\cr
\vdots &\vdots&\vdots&\vdots&\vdots \cr
\zeta(A-C+B+1) &\cdots &\cdots &-\zeta C&-\zeta(C+1)\cr
 \vdots&\vdots&\vdots&\vdots \cr
\zeta A &\cdots &\cdots &-\zeta (B+1)
\end{matrix}
$$
Dans le cas oppos\'e, c'est-\`a-dire $A+B< 2C$, on obtient:
$$
\begin{matrix}
&\zeta (B+1)& \cdots &\cdots &-\zeta A\cr
 &\vdots&\vdots&\vdots&\vdots \cr
 &\zeta(A+B-C)&\cdots&-\zeta C &-\zeta(C+1)\cr
 &\vdots &\vdots &\vdots\cr
\zeta (C+1) &\zeta (C)& \cdots &-\zeta ((A+B-C)\cr
\vdots &\vdots&\vdots&\vdots \cr
\zeta A &\cdots &\cdots &-\zeta (B+1)
\end{matrix}
$$
Quel que soit le tableau \'ecrit ci-dessus, la classification de Zelevinsky montre qu'il existe exactement une unique repr\'esentation irr\'eductible sous-module de l'induite associ\'ee aux multi-segments form\'es par l'ensemble des lignes mis dans l'ordre des lignes (du haut vers le bas). On la note $\sigma_{\cal T}$ o\`u ${\cal T}$ est le tableau correspondant. Avec un tableau \'ecrit comme ci-dessus avec les lignes form\'ees de segments et les colonnes form\'ees aussi de segments mais avec la propri\'et\'e de croissance oppos\'ee \`a celle des lignes, on v\'erifie que $\sigma_{\cal T}$  est aussi l'unique sous-repr\'esentation irr\'eductible de l'induite associ\'ee aux multi-segments form\'es par l'ensemble des colonnes mis dans l'ordre des colonnes (de gauche vers la droite). On a en tant que repr\'esentation du $GL$ convenable
$$
\pi(\rho,A,B,\zeta)=\sigma_{{\cal T}_{A,B,\zeta}},\qquad \qquad \pi(\rho,A,B+2,\zeta)=\sigma_{{\cal T}_{A,B+2,\zeta}}
$$
$$
Jac^\theta_{\zeta (B+2), \cdots, \zeta C}\pi(\rho,A,B+2,\zeta)= \sigma_{{\cal T}_{C}}.
$$
De plus, avec la description ci-dessus, on a $Jac_{x}Jac^\theta_{\zeta (B+2), \cdots, \zeta C}\pi(\rho,A,B+2,\zeta)\neq 0$ entra\^{\i}ne que $x$ est \`a la fois un d\'ebut de ligne et un d\'ebut de colonne, c'est-\`a-dire ici, n\'ecessairement $x=\zeta(B+1)$ ou $x=\zeta (C+1)$. On peut aussi avoir l'action de $\theta$: pour $\pi(\rho,A,B,\zeta)$ la d\'efinition dit qu'il faut privil\'egier la repr\'esentation comme sous-module en utilisant les segments form\'es par les lignes. Mais \ref{proprietes} et \ref{remarque} montrent que l'on obtient la m\^eme action de $\theta$ en utilisant les colonnes. Avec les m\^emes r\'ef\'erences, on d\'emontre le m\^eme r\'esultat pour $Jac^\theta_{\zeta (B+2), \cdots, \zeta C}\pi(\rho,A,B+2,\zeta)$.

\subsubsection{Propri\'et\'e des modules de Jacquet et preuve\label{jacsimple}}
Dans tout ce qui suit, on fixe $\tilde{\pi}$ une repr\'esentation irr\'eductible de $\tilde{G}(abd_{\rho})$ qui intervient avec un coefficient non nul dans la description de $\tilde{\pi}(\rho,A,B,\zeta)$.
\vskip 0.5cm
\noindent
\bf Lemme: \sl (i) Soit $x\in {\mathbb R}$. On a $Jac_{x}\tilde{\pi}=0$ pour tout $x\notin [\zeta B, \zeta A]$ et on a $Jac_{x,x}\tilde{\pi}=0$ pour tout $x$.

(ii) On suppose ici que la restriction de $\tilde{\pi}$ \`a $GL(abd_{\rho},F)$ est irr\'eductible; alors pour tout $x\in {\mathbb R}$, $Jac^\theta_{x}\tilde{\pi}=0$ sauf \'eventuellement pour $x=\zeta B$.

(iii) La proposition \ref{casquasigeneral} est vraie c'est-\`a-dire que $\tilde{\pi}(\rho,A,B,\zeta)$ contient une unique repr\'esentation irr\'eductible dont la restriction \`a $GL(abd_{\rho},F)$ est irr\'eductible, elle intervient avec coefficient 1 et est exactement $\pi(\rho,A,B,\zeta)$.
\vskip 0.5cm
\noindent
\rm
Comme $\tilde{\pi}$ est irr\'eductible, soit $\tilde{\pi}$ est l'induite $\pi(\rho,A,B+1,\zeta)\times \pi(\rho,B,B,\zeta)$ soit il existe $C\in ]B,A]$ tel que $\tilde{\pi}$ est un sous-quotient irr\'eductible de l'induite 
$$
\sigma_{C}:=<\rho\vert\,\vert^{\zeta B}, \cdots, \rho\vert\,\vert^{-\zeta C}>\times Jac^\theta_{\zeta (B+2), \cdots, \zeta C}\pi(\rho,A,B+2,\zeta)\times <\rho\vert\,\vert^{\zeta C}, \cdots, \rho\vert\,\vert^{-\zeta B}>.
$$Les 2 propri\'et\'es ne sont pas  exclusives l'une de l'autre.

On v\'erifie que $Jac_{x}\bigl(\pi(\rho,A,B+1,\zeta)\times \pi(\rho,B,B,\zeta)\bigr)=0$ sauf exactement pour $x=\zeta B$ et $x=\zeta(B+1)$ et on v\'erifie aussi ais\'ement que $Jac_{x,x}\bigl(\pi(\rho,A,B+1,\zeta)\times \pi(\rho,B,B,\zeta)\bigr)=0$ pour tout $x$ r\'eel. Supposons que $\tilde{\pi}$ est un sous-quotient irr\'eductible de $\sigma_{C}$ pour $C\in ]B,A]$. Si $Jac_{x}\tilde{\pi}\neq 0$, n\'ecessairement soit $x=\zeta B$, soit $x= \zeta C$ soit $Jac_{x}Jac^\theta_{\zeta (B+2), \cdots, \zeta C}\pi(\rho,A,B+2,\zeta)\neq 0$. Pour \'etudier cette derni\`ere possibilit\'e on utilise la r\'ealisation de $Jac^\theta_{\zeta (B+2), \cdots, \zeta C}\pi(\rho,A,B+2,\zeta)$ comme $\sigma_{{\cal T}_{C}}$;  on a vu que la non nullit\'e n\'ecessite que soit  $x=\zeta(B+1)$ soit $x=\zeta (C+1)$ avec les 2 cas particuliers suivant: si $C=B+1$ seul la valeur $x=\zeta (B+2)$ donne quelque chose de non nul et si $C=A$ seul $x=\zeta(B+1)$ donne quelque chose de non nul. On v\'erifie alors ais\'ement, dans tous les cas, que $Jac_{x,x}\sigma_{C}=0$ pour toute valeur de $x$. Cela termine la preuve de (i).

Prouvons (ii); la seule difficult\'e peut venir d'un r\'eel $x$ v\'erifiant $Jac_{x}\tilde{\pi}\neq 0$ et d'apr\`es (i) cela entra\^{\i}ne que $x\in [\zeta B,\zeta A]$. Il faut donc montrer que pour tout $C\in ]B,A]$, $Jac^\theta_{\zeta C}\tilde{\pi}=0$. 

On va d'abord traiter le cas o\`u $C\in [B+2,A]$ o\`u on va en fait montrer que $Jac^\theta_{\zeta C}\tilde{\pi}(\rho,A,B,\zeta)=0$. On pose $x=\zeta C$.
Montrons d'abord comment cela permettra de conclure; les arguments ci-dessous ne suppose que $\zeta x \in ]B,A]$, c'est-\`a-dire $C\in ]B,A]$.

Faisons remarquer tout de suite que, pour les valeurs que nous consid\'erons $x\neq 0$. Soit $\pi'$ une repr\'esentation irr\'eductible de $\tilde{G}(abd_{\rho})$ intervenant dans $\tilde{\pi}(\rho,A,B,\zeta)$ et telle que $Jac ^\theta_{x}\pi' \neq 0$; on consid\`ere la restriction de $\pi'$ \`a $GL(abd_{\rho},F)$ et son module de Jacquet relativement au parabolique de Levi $GL(d_{\rho}) \times GL((ab-2)d_{\rho},F) \times GL(d_{\rho})$. Par hypoth\`ese la projection de cette repr\'esentation sur la composante o\`u le premier facteur de $GL(d_{\rho})$ agit par $\rho\vert\,\vert^x$ et le troisi\`eme facteur agit  par $\rho\vert\,\vert^{-x}$ de fa\c{c}on isotypique est non nulle. Par r\'eciprocit\'e de Frobenius, il existe une repr\'esentation irr\'eductible $\pi_{0}$ de $GL((ab-2)d_{\rho},F)$ et une inclusion d'une des composantes irr\'eductibles de la restriction de $\pi'$ \`a $GL(abd_{\rho},F)$ dans l'induite $\rho\vert\,\vert^x \times \pi_{0} \times \rho\vert\,\vert^{-x}$. On note $\pi'_{0}$ une composante irr\'eductible de la restriction de $\pi'$ \`a $GL(abd_{\rho},F)$ incluse dans cette induite.

On calcule $Jac^\theta_{x}(\rho\vert\,\vert^x \times \pi_{0} \times \rho\vert\,\vert^{-x})=\pi_{0}$:  pour cela on remarque d'abord que si $Jac_{x}\pi_{0}\neq 0$, on aurait aussi $Jac_{x,x}\pi'_{0}\neq 0$ et certainement aussi $Jac_{x,x}\pi'\neq 0$ et une contractiction avec (i). Ainsi $Jac_{x}\pi_{0}=0$. Comme $x\neq 0$, on a aussi $Jac_{x}\bigl(\pi_{0}\times \rho\vert\,\vert^{-x}\bigr)=0$. D'o\`u l'assertion.

On d\'eduit de cette assertion que  l'induite $\rho\vert\,\vert^x\times \pi_{0}\times \rho\vert\,\vert^{-x}$ a un unique sous-$GL(abd_{\rho},F)$-module irr\'eductible qui ici est n\'ecessairement $\pi'_{0}$. L'unicit\'e assure aussi que $\pi'_{0}$ est $\theta$-invariante si et seulement si $\pi_{0}$ l'est. Supposons que $\pi_{0}$ est $\theta$-invariante; on a donc $\pi'_{0}=\pi'$ comme repr\'esentation de $GL(abd_{\rho},F)$. Mais $\pi'_{0}$ h\'erite alors d'une action de $\theta$ qui donne une action de $\theta$ sur $\pi_{0}$ uniquement d\'etermin\'ee. Par contre si $\pi_{0}$ n'est pas $\theta$-invariante $\pi'$ et l'action de $\theta$ sur $\pi'$ sont uniquement d\'etermin\'es par $\pi'_{0}$ et donc par $\pi_{0}$. Donc dans les 2 cas, la seule hypoth\`ese $Jac^\theta_{x} \pi'\neq 0$ (ajout\'ee \`a l'hypoth\`ese que $\pi'$ intervient dans $\tilde{\pi}(\rho,A,B,\zeta)$) entra\^{\i}ne que $Jac^\theta_{x}\pi'$ est irr\'eductible comme $\tilde{G}((ab-2)d_{\rho})$-module et d\'etermine uniquement $\pi'$. Cela prouve qu'il ne peut y avoir de simplification dans le groupe de Grothendieck pour $\tilde{G}((ab-2)d_{\rho})$ entre des $Jac^\theta_{x}\pi'$ faisant intervenir des $\pi'$ diff\'erents, ou encore que $Jac^\theta_{x}\tilde{\pi}(\rho,A,B,\zeta)=0$ entra\^{\i}ne que $Jac^\theta_{x}\pi'=0$ pour tout $\pi'$ repr\'esentation irr\'eductible de $\tilde{G}$ intervenant avec un coefficient non nul dans $\tilde{\pi(\rho,A,B,\zeta)}$.

Il nous reste  \`a montrer que pour tout $C\in ] (B+1),A]$, $Jac_{\zeta C}^\theta \tilde{\pi}(\rho,A,B,\zeta)=0$. Il y a exactement 2 termes dans la d\'efinition de $\tilde{\pi}(\rho,A,B,\zeta)$ qui v\'erifient simplement $Jac_{\zeta C}\neq 0$, il s'agit, d'apr\`es ce que l'on a vu dans la preuve de (i) de $\sigma_{C}$ et $\sigma_{C-1}$. On calcule $$Jac^\theta_{\zeta C}\sigma_{C}=<\rho\vert\,\vert^{\zeta B}, \cdots, \rho\vert\,\vert^{-\zeta(C-1)}>\times Jac^\theta_{\zeta(B+2), \cdots, \zeta C}\pi(\rho,A,B+2,\zeta) \times <\rho\vert\,\vert^{\zeta(C-1)}, \cdots, \rho\vert\,\vert^{-\zeta B}>;\eqno(1)
$$
Pour voir l'action de $\theta$ sur ce module de Jacquet, on consid\`ere d'abord l'inclusion:
$$
\sigma_{C} \hookrightarrow$$
$$ <\rho\vert\,\vert^{\zeta B}, \cdots, \rho\vert\,\vert^{-\zeta(C-1)}>\times \rho\vert\,\vert^{-\zeta C}\times Jac^\theta_{\zeta(B+2), \cdots, \zeta C}\pi(\rho,A,B+2,\zeta)\times \rho\vert\,\vert^{\zeta C}\times <\rho\vert\,\vert^{\zeta (C-1)}, \cdots , \rho\vert\,\vert^{-\zeta B}>.$$
C'est une inclusion qui respecte les actions naturelles de $\theta$. Le calcul de $Jac^\theta$ se fait avec la filtration de Bernstein Zelevinsky qui provient des doubles classes  du type $P_{0}w P_{0}$ pour $w$ parcourant un bon ensemble de matrices de permutations. L'action de $\theta$ respecte cette filtration \`a condition de regrouper $w$ avec $^\theta w$; le terme qui nous int\'eresse correspond \`a un unique $w$ dont la double classe est invariante sous l'action de $\theta$. L'action de $\theta$ obtenue sur $Jac^\theta_{\zeta C}$ est  l'action obtenue sur $Jac^\theta_{\zeta C}$ appliqu\'e \`a
$$
<\rho\vert\,\vert^{\zeta B}, \cdots, \rho\vert\,\vert^{-\zeta(C-1)}>\times \rho\vert\,\vert^{\zeta C}\times Jac^\theta_{\zeta(B+2), \cdots, \zeta C}\pi(\rho,A,B+2,\zeta)\times \rho\vert\,\vert^{-\zeta C}\times <\rho\vert\,\vert^{\zeta (C-1)}, \cdots , \rho\vert\,\vert^{-\zeta B}>,
$$
car on peut remplacer l'induite du milieu par son unique sous-module irr\'eductible. Mais cette derni\`ere induite est isomorphe \`a :
$$\rho\vert\,\vert^{\zeta C}\times 
<\rho\vert\,\vert^{\zeta B}, \cdots, \rho\vert\,\vert^{-\zeta(C-1)}>\times Jac^\theta_{\zeta(B+2), \cdots, \zeta C}\pi(\rho,A,B+2,\zeta)\times <\rho\vert\,\vert^{\zeta (C-1)}, \cdots , \rho\vert\,\vert^{-\zeta B}>\times \rho\vert\,\vert^{-\zeta C},
$$car $C\neq B+1$. L'action de $\theta$ sur $Jac^\theta_{\zeta C}$ de l'induite ci-dessus est l'action naturelle sans aucune ambigu\"{\i}t\'e.
On a aussi
$$
Jac^\theta_{\zeta C}\sigma_{C-1}=$$
$$<\rho\vert\,\vert^{\zeta B}, \cdots, \rho\vert\,\vert^{-\zeta(C-1)}>\times Jac^\theta_{\zeta C}Jac^\theta_{\zeta(B+2), \cdots, \zeta (C-1)}\pi(\rho,A,B+2,\zeta) \times <\rho\vert\,\vert^{\zeta(C-1)}, \cdots, \rho\vert\,\vert^{-\zeta B}>.\eqno(2)$$
Pour voir l'action de $\theta$ le plus simple est d'\'ecrire $Jac_{\zeta (B+2), \cdots, \zeta (C-1)}\pi(\rho,A,B+2,\zeta)$ comme l'unique sous-module irr\'eductible de l'induite:
$$
\rho\vert\,\vert^{\zeta C}\times Jac^\theta_{\zeta (B+2), \cdots, \zeta C}\pi(\rho,A,B+2,\zeta) \times \rho\vert\,\vert^{-\zeta C}$$
cette inclusion respecte les actions de $\theta$ par la d\'efinition m\^eme. Et l'action de $\theta$ sur (2) est l'action naturelle.

Ces  termes (1) et (2) sont donc exactement les m\^emes, action de $\theta$ comprise mais l'un intervient avec le signe $(-1)^{A-C}$ tandis que l'autre intervient avec le signe $(-1)^{A-C+1}$; ils s'\'eliminent donc, prouvant l'assertion cherch\'ee.

Consid\'erons pour finir le cas $x=\zeta(B+1)$; ici on d\'emontre moins que pour $C\in ]B+1,A]$ c'est-\`a-dire on d\'emontre seulement que pour $\tilde{\pi}$ une repr\'esentation irr\'eductible de $\tilde{G}(abd_{\rho})$ intervenant dans $\tilde{\pi}(\rho,A,B,\zeta)$ dont la restriction \`a $GL(abd_{\rho},F)$ est irr\'eductible, $Jac^\theta_{\zeta(B+1)}\tilde{\pi}=0$. Fixons un tel $\tilde{\pi}$ et pour \'eviter les confusions, on note $\pi$ la restriction de $\tilde{\pi}$ \`a $GL(abd_{\rho})$. Le support cuspidal de $\pi$ est un ensemble de repr\'esentations cuspidales de la forme $\rho\vert\,\vert^z$ avec $z$ r\'eel; on note ${\cal E}_{\pi}$ cet ensemble de nombres r\'eels compt\'es avec multiplicit\'e. Pour un ordre convenable sur ${\cal E}$, il existe une inclusion:
$$
\pi \hookrightarrow \times_{z\in {\cal E}}\rho\vert\,\vert^z.\eqno(1)
$$
Le r\'eel $-\zeta A$ est soit le plus petit \'el\'ement de ${\cal E}$ soit le plus grand (cela d\'epend de la valeur de $\zeta$). De fa\c{c}on formelle on montre qu'il existe $x_{0}\in {\cal E}$ tel que $[x_{0},-\zeta A]$ soit un segment et qu'il existe un ordre sur ${\cal E}$ tel que (1) soit encore vrai mais tel que les premiers \'el\'ements de ${\cal E}$ soient pr\'ecis\'ement les \'el\'ements du segment $[x_{0},-\zeta A]$. En effet fixons sur ${\cal E}$ un ordre tel que la place de $-\zeta A$ soit la plus petite possible. On note $x_{0}$ le premier \'el\'ement de ${\cal E}$ et on va v\'erifier qu'il convient, c'est-\`a-dire que $[x_{0},-\zeta A]$ est un segment.  On note $t$ le nombre d'\'el\'ements de ${\cal E}$ plus petit ou \'egaux (pour l'ordre fix\'e) \`a $-\zeta A$, c'est-\`a-dire que les $t$ premiers \'el\'ements de ${\cal E}$ sont $x_{i}$ pour $i\in [1,t]$ avec $x_{t}=-\zeta A$. On montre par r\'ecurrence d\'ecroissante sur $i$ que $[x_{i},-\zeta A]$ est un segment et que (1) se pr\'ecise en une inclusion:
$$
\pi\hookrightarrow \times_{j<i}\rho\vert\,\vert^{x_{j}} \times <\rho\vert\,\vert^{x_{i}}, \cdots, \rho\vert\,\vert^{-\zeta A}> \times_{z\in {\cal E}-\{x_{1}, \cdots, x_{t}\}}\rho\vert\,\vert^{z}.\eqno(2)
$$
Ceci est trivialement vrai si $i=t$. Admettons le pour $i$ et montrons le pour $i-1$; la repr\'esentation induite:
$$
\rho\vert\,\vert^{x_{i-1}}\times <\rho\vert\,\vert^{x_{i}}, \cdots, \rho\vert\,\vert^{-\zeta A}>
$$
est soit irr\'eductible soit de longueur 2; elle ne peut \^etre irr\'eductible, sinon on pourrait remplacer $t$ par $t-1$ en faisant commuter $x_{i-1}$ au dessus de $<\rho\vert\,\vert^{x_{i}}, \cdots, \rho\vert\,\vert^{-A}>$. Elle est donc de longueur 2, c'est-\`a-dire soit $x_{i-1}=x_{i}+\zeta$ soit $x_{i-1}=-\zeta A-\zeta$; la derni\`ere \'eventualit\'e est impossible par extr\'emalit\'e de $-\zeta A$. Donc $[x_{i-1},-\zeta A]$ est un segment et l'inclusion (2) se factorise soit par
$$
\pi\hookrightarrow \times_{j<i-1}\rho\vert\,\vert^{x_{j}} \times <\rho\vert\,\vert^{x_{i-1}}, \cdots, \rho\vert\,\vert^{-\zeta A}> \times_{z\in {\cal E}-\{x_{1}, \cdots, x_{t}\}}\rho\vert\,\vert^{z}\eqno(3)
$$
soit par
$$
\pi\hookrightarrow \times_{j<i-1}\rho\vert\,\vert^{x_{j}} \times <\rho\vert\,\vert^{x_{i}}, \cdots, \rho\vert\,\vert^{-\zeta A}>\times \rho\vert\,\vert^{x_{i}} \times_{z\in {\cal E}-\{x_{1}, \cdots, x_{t}\}}\rho\vert\,\vert^{z}.\eqno(4)
$$
Mais (4) est contraire \`a la minimalit\'e de $t$ et c'est donc (3) qui pr\'evaut, montrant ainsi notre assertion. Ainsi, il existe $x_{0}\in {\cal E}$ et une repr\'esentation irr\'eductible $\pi'_{0}$ convenable telle que $${\pi}\hookrightarrow <\rho\vert\,\vert^{x_{0}}, \cdots, \rho\vert\,\vert^{-\zeta A}>\times \pi'_{0}.$$Gr\^ace \`a (i), on sait que $x_{0}\in [\zeta B,\zeta A]$.

Supposons d'abord que $x_{0}=\zeta B$; dans ce cas, on va tr\`es facilement montrer que $\tilde{\pi}=\pi(\rho,A,B,\zeta)$ et donc en particulier $Jac_{\zeta (B+1)}\tilde{\pi}=0$, ce que l'on cherche. En effet, supposons d'abord qu'il existe $C\in ]B,A]$ tel que $\tilde{\pi}$ soit un sous-quotient de 
$$
\sigma_{C}:=<\rho\vert\,\vert^{\zeta B}, \cdots, \rho\vert\,\vert^{-\zeta C}> \times Jac^\theta_{\zeta (B+2), \cdots, \zeta C}(\pi(\rho,A,B+2,\zeta))\times <\rho\vert\,\vert^{\zeta C}, \cdots, \rho\vert\,\vert^{-\zeta B}>.
$$
On calcule $Jac_{\zeta B, \cdots -\zeta A}\sigma_{C}$ par les formules standard en se souvenant que $Jac_{x}Jac^\theta_{\zeta (B+2), \cdots, \zeta C}(\pi(A,B+2,\zeta))=0$ sauf \'eventuellement si $x=\zeta (B+1)$ ou $x=\zeta (C+1)$. Or ni $\zeta (B+1)$ ni $\zeta (C+1)$ ni $\zeta C$ ne sont dans le segment $[\zeta B,-\zeta A]$. Ainsi $$Jac_{\zeta B, \cdots -\zeta A}\sigma_{C}=(Jac_{\zeta B, \cdots, -\zeta A}<\rho\vert\,\vert^{\zeta B}, \cdots, \rho\vert\,\vert^{-\zeta C}>)\times $$
$$Jac^\theta_{\zeta (B+2), \cdots, \zeta C}(\pi(A,B+2,\zeta))\times <\rho\vert\,\vert^{\zeta C}, \cdots, \rho\vert\,\vert^{-\zeta B}>$$
Et cela est nul sauf si $C=A$ o\`u cela vaut $Jac^\theta_{\zeta (B+2), \cdots, \zeta A}(\pi(\rho,A,B+2,\zeta))\times <\rho\vert\,\vert^{\zeta A}, \cdots, \rho\vert\,\vert^{-\zeta B}>$. On v\'erifie que $Jac^\theta_{\zeta (B+2), \cdots, \zeta A}(\pi(\rho,A,B+2,\zeta))=\pi(\rho,A-1,B+1,\zeta)$ en particulier est irr\'eductible et on obtient l'inclusion:
$$
\tilde{\pi} \hookrightarrow <\rho\vert\,\vert^{\zeta B}, \cdots, \rho\vert\,\vert^{-\zeta A}>\times \pi(\rho,A-1,B+1,\zeta)\times <\rho\vert\,\vert^{\zeta A}, \cdots, \rho\vert\,\vert^{- \zeta B}>.
$$
Ainsi $\tilde{\pi}$ est l'unique sous-module irr\'eductible de cette induite avec l'action de $\theta$ qui provient de l'action sur l'induite. Par d\'efinition, $\tilde{\pi}=\pi(\rho,A,B,\zeta)$. Il reste \`a voir le cas o\`u $\tilde{\pi}$ est un sous-quotient de $\pi(\rho,A,B+1,\zeta)\times <\rho\vert\,\vert^{\zeta B}, \cdots, \rho\vert\,\vert^{-\zeta B}>$. Ici on sait que $Jac_{x}\pi(\rho,A,B+1,\zeta)=0$ sauf pour $x=\zeta(B+1)$ et il est facile de voir que cela entra\^{\i}ne que $$Jac_{\zeta B, \cdots, -\zeta A}\biggl(\pi(\rho,A,B+1,\zeta)\times <\rho\vert\,\vert^{\zeta B}, \cdots, \rho\vert\,\vert^{-\zeta B}>\biggr)=$$
$$\pi(\rho,A,B+1,\zeta) \times 
Jac_{\zeta B, \cdots, -\zeta A}<\rho\vert\,\vert^{\zeta B}, \cdots, \rho\vert\,\vert^{-\zeta B}>=0.
$$
Cela termine le cas o\`u $x_{0}=\zeta B$.

Supposons donc maintenant que $x_{0}\in ]\zeta B,\zeta A]$. On \'ecrit donc:
$$
\pi \hookrightarrow <\rho\vert\,\vert^{x_{0}}, \cdots, \rho\vert\,\vert^{-\zeta A}>\times \pi'_{0},
$$
pour $\pi'_{0}$ irr\'eductible convenable. On utilise le fait que $\pi$ est $\theta$-invariante et on obtient une autre inclusion:
$$
\pi \hookrightarrow \, ^\theta \pi'_{0}\times <\rho\vert\,\vert^{\zeta A}, \cdots, \rho\vert\,\vert^{-x_{0}}>.
$$
On applique $Jac_{x_{0}, \cdots, -\zeta A}$ \`a l'induite de droite et comme $\zeta A \notin [x_{0}, -\zeta A]$ cela vaut $$
Jac_{x_{0}, \cdots -\zeta A}\, ^\theta\pi'_{0}\times <\rho\vert\,\vert^{\zeta A}, \cdots, \rho\vert\,\vert^{-x_{0}}>.
$$
Par exactitude du foncteur de Jacquet, on obtient une inclusion:
$$Jac_{x_{0}, \cdots -\zeta A}\pi=
\pi'_{0}\hookrightarrow Jac_{x_{0}, \cdots -\zeta A} \, ^\theta\pi'_{0}\times <\rho\vert\,\vert^{\zeta A}, \cdots, \rho\vert\,\vert^{-x_{0}}>.
$$
Il existe donc un sous-quotient irr\'eductible $\pi''_{1}$ de $Jac_{x_{0}, \cdots , -\zeta A}\pi'_{0}$ tel que en posant $\pi'_{1}:= \, ^\theta \pi''_{1}$
$$
\pi'_{0}\hookrightarrow \pi'_{1}\times <\rho\vert\,\vert^{\zeta A}, \cdots, \rho\vert\,\vert^{-x_{0}}>.
$$
En remontant, on obtient:
$$
\pi\hookrightarrow <\rho\vert\,\vert^{x_{0}}, \cdots, \rho\vert\,\vert^{-\zeta A}>\times \pi'_{1}\times <\rho\vert\,\vert^{\zeta A}, \cdots, \rho\vert\,\vert^{-x_{0}}>.
$$
On v\'erifie que $\pi'_{1}=Jac^\theta_{x_{0}, \cdots, -\zeta A} \pi$ et h\'erite d'une action de $\theta$ provenant de l'action de $\theta$ sur $\tilde{\pi}$, on note alors $\pi_{1}$ la repr\'esentation ainsi \'etendue. Avec cette action, on obtient l'inclusion:
$$
\tilde{\pi} \hookrightarrow  <\rho\vert\,\vert^{x_{0}}, \cdots, \rho\vert\,\vert^{-\zeta A}>\times \pi_{1}\times <\rho\vert\,\vert^{\zeta A}, \cdots, \rho\vert\,\vert^{-x_{0}}>.\eqno(5)
$$
Gr\^ace \`a ce que l'on a d\'ej\`a d\'emontr\'e, $x_{0}=\zeta (B+1)$. De plus, $\tilde{\pi}$ est uniquement d\'etermin\'ee par $\pi_{1}$ et l'inclusion (5). Pour d\'emontrer qu'il n'existe pas de tel $\tilde{\pi}$ intervenant dans $\tilde{\pi}(\rho,A,B,\zeta)$, il suffit donc de d\'emontrer que $Jac^\theta_{\zeta(B+1), \cdots, -\zeta A}\tilde{\pi}(\rho,A,B,\zeta)$ ne contient pas de repr\'esentation irr\'eductible de $\tilde{G}(ab-2(sup(a,b)+1))d_{\rho})$ dont la restriction au GL correspondant reste irr\'eductible. 

On reprend la notation $\sigma_{C}$ introduite ci-dessus. On va avoir besoin des formules ci-dessous:
soit $C \in ]B+1, A[$ (pour que l'ensemble soit non vide, il faut $A>B+2$) et $x\in [\zeta (B+1),-\zeta A[$, alors 
$$
Jac_{y}Jac^\theta_{\zeta(B+1), \cdots, x}Jac^\theta_{\zeta (B+2), \cdots, \zeta C}\pi(\rho,A,B+2,\zeta)=0,
$$
pour tout $y \in ]x-\zeta,-\zeta A]$ et
$$Jac^\theta_{\zeta(B+1), \cdots, -\zeta A}Jac^\theta_{\zeta (B+2), \cdots, \zeta C}\pi(\rho,A,B+2,\zeta)=
Jac^\theta_{\zeta(B+3), \cdots, \zeta C}\pi(\rho,A-1,B+3,\zeta),
$$cette \'egalit\'e \'etant compatible aux actions de $\theta$ et si $A=B+3$ le terme de droite ci-dessus est la repr\'esentation triviale de $\tilde{G}(0)$.
Pour les v\'erifier on utilise la repr\'esentation de $Jac^\theta_{\zeta (B+2), \cdots, \zeta C}\pi(\rho,A,B+2,\zeta)$ donn\'ee \`a la fin de \ref{simple}; en calculant $Jac^\theta_{\zeta (B+1), \cdots x}$ on enl\`eve du tableau une partie (tout si $x=-\zeta A$) de la premi\`ere et sym\'etriquement de la derni\`ere ligne. Et ensuite on voit sur le tableau les valeurs de $y$ pour lesquelles $Jac_{y}$ du r\'esultat peut \^etre non nul, $y$ doit au moins \^etre un d\'ebut de ligne. L'action de $\theta$ a \'et\'e pr\'ecis\'ee en loc.cit, d'o\`u nos assertions.
Cela entra\^{\i}ne que pour $C\in ]B+1,A[$
$$
Jac^\theta_{\zeta (B+1), \cdots, -\zeta A}\sigma_{C}=
$$
$$
<\rho\vert\,\vert^{\zeta B}, \cdots,\rho\vert\,\vert^{-\zeta C}>\times Jac^\theta_{\zeta (B+3), \cdots, \zeta C}\pi(\rho,A-1,B+3,\zeta) \times <\rho\vert\,\vert^{\zeta C}, \cdots, \rho\vert\,\vert^{-\zeta B}>.
$$
Pour $C=B+1$, on a tout simplement $Jac_{x} \pi(\rho,A,B+2,\zeta)=0$ pour tout $x\in [\zeta B,-\zeta(B+1)]$ et 
$$
Jac^\theta_{\zeta (B+1), \cdots, -\zeta A}\sigma_{B+1}= Jac_{-\zeta(B+2),\cdots, -\zeta A}\pi(\rho,A,B+2,\zeta)=0.
$$
Pour $C=A$, $Jac^\theta_{\zeta (B+2), \cdots, \zeta A} \pi(\rho,A,B+2,\zeta)=\pi(\rho,A-1,B+1,\zeta)$ et 
$$
Jac^\theta_{\zeta (B+1), \cdots, -\zeta A}\sigma_{A}=
$$
$$
Jac^\theta_{\zeta (B+1), \cdots, -\zeta A} (<\rho\vert\,\vert^{\zeta B}, \cdots, \rho\vert\,\vert^{-\zeta A}> \times \pi(\rho,A-1,B+1,\zeta) \times <\rho\vert\,\vert^{\zeta A}, \cdots, \rho^{-\zeta B}>)=
$$
$$
Jac^\theta_{\zeta(B+1)}\pi(\rho,A-1,B+1,\zeta).
$$
On a aussi, sans l'action de $\theta$:
$$
Jac^\theta_{\zeta (B+1), \cdots, -\zeta A}(\pi(\rho,A,B+1,\zeta)\times \pi(\rho,B,B,\zeta))=\pi(\rho,A-1,B+2,\zeta)\times \pi(\rho,B,B,\zeta);
$$
pour avoir l'action de $\theta$, il faut revenir \`a la d\'efinition de $\pi(\rho,A,B+1,\zeta) \times \pi(\rho,B,B,\zeta)$. On veut montrer que l'action de $\theta$ sur le membre de gauche est l'action que nous avons mise sur le membre de droite. Par d\'efinition $\theta$ agit sur $\pi(\rho,A,B+1,\zeta)\times \pi(\rho,B,B,\zeta)$ par restriction de l'action naturelle de $\theta$ sur l'induite:
$$
<\rho\vert\,\vert^{\zeta (B+1)}, \cdots, \rho\vert\,\vert^{-\zeta A}> \times \pi(\rho,A-1,B+2,\zeta) \times \pi(\rho,B,B,\zeta) \times <\rho\vert\,\vert^{\zeta A}, \cdots, \rho\vert\,\vert^{-\zeta (B+1)}>.
$$
Et le r\'esultat est clair.

En regroupant tous les termes ensembles, on obtient:
$$
Jac^\theta_{\zeta (B+1), \cdots, -\zeta A}\tilde{\pi}(\rho,A,B,\zeta)=Jac^\theta_{\zeta (B+1)}\pi(\rho,A-1,B+1,\zeta)\oplus_{C\in ]B+1,A-1]}(-1)^{A-C} 
$$
$$
<\rho\vert\,\vert^{\zeta B}, \cdots, \rho\vert\,\vert^{-\zeta C}> \times Jac^\theta_{\zeta (B+3), \cdots, \zeta C}\pi(\rho, A-1,B+3,\zeta) \times <\rho\vert\,\vert^{\zeta C}, \cdots, \rho\vert\,\vert^{-\zeta B}>\eqno(6)
$$
$$
\oplus (-1)^{[(A-B+1)/2]} \pi(\rho,A-1,B+2,\zeta)\times <\rho\vert\,\vert^{\zeta B}, \cdots, \rho\vert\,\vert^{-\zeta B}>.\eqno(7)
$$
On remarque que la somme de (6) et (7) n'est autre que $(-1) Jac^\theta_{\zeta (B+1)}\tilde{\pi}(\rho,A-1,B+1,\zeta)$. On ne peut pas conclure tout de suite \`a cause de la diff\'erence entre $\pi$ et $\tilde{\pi}$. 

 Prenons  $\tau$ une repr\'esentation irr\'eductible de $Jac^\theta_{\zeta(B+1)}\tilde{\pi}(\rho,A-1,B+1,\zeta)$ qui reste irr\'eductible apr\`es restriction au $GL$ correspondant. On reprend un argument d\'ej\`a donn\'e; on fixe $\tilde{\tau}$ une repr\'esentation irr\'eductible intervenant dans $\tilde{\pi}(\rho,A-1,B+1,\zeta)$ et telle que $Jac^\theta_{\zeta(B+1)}\tilde{\tau}$ contienne $\tau$ comme sous-quotient. Comme on l'a fait ci-dessus, on d\'emontre que $Jac^\theta_{\zeta(B+1)}\tilde{\tau}=\tau$ et que $\tilde{\tau}$ est uniquement d\'etermin\'e par $\tau$ en particulier est de restriction irr\'eductible au $GL$ correspondant. Par exemple par r\'ecurrence on sait que la seule repr\'esentation de $\tilde{\pi}(\rho,A-1,B+1,\zeta)$ dont la restriction au GL correspondant est irr\'eductible est $\pi(\rho,A-1,B+1,\zeta)$. On vient donc de montrer que $Jac_{\zeta(B+1)}^\theta(\tilde{\pi}(\rho,A-1,B+1,\zeta)-\pi(\rho,A-1,B+1,\zeta))$ ne contient aucune repr\'esentation dont la restriction au GL est irr\'eductible.

Revenons \`a $\tilde{\pi}$ comme dans l'\'enonc\'e; on v\'erifie comme ci-dessus  que $Jac^\theta_{\zeta (B+1)}\pi$ est irr\'eductible et d\'etermine uniquement $\tilde{\pi}$. Il ne peut donc pas  y avoir de simplification dans le module de Jacquet et on vient donc de montrer que $Jac^\theta_{\zeta(B+1)}\pi=0$.

On d\'eduit  (iii): on vient de d\'emontrer que pour $\tilde{\pi}$ intervenant dans $\tilde{\pi}(\rho,A,B,\zeta)$ de restriction irr\'eductible \`a $GL(abd_{\rho},F)$ on a n\'ecessairement $x_{0}=\zeta B$ avec les notations ci-dessus. On a fait le calcul correspondant ci-dessus, il n'y a qu'un terme qui donne ce module de Jacquet, $\pi(\rho,A,B,\zeta)$. Cela termine la preuve de ce cas simple.


\subsection{R\'esolution dans le groupe de Grothendieck, cas de restriction discr\`ete \`a la diagonale\label{quasigeneral}}
Fixons $\psi$ et $(\rho,a,b)\in Jord(\psi)$ tel que $inf(a,b)>1$; on pose encore $A=(a+b)/2-1$, $B=\vert (a-b)\vert/2$ et $\zeta$ le signe de $a-b$ si $a\neq b$ et sinon $\zeta$ est un signe quelconque. On note $\psi'$ le morphisme qui se d\'eduit de $\psi$ en enlevant le bloc de Jordan $(\rho,a,b)$; on a donc d\'efini $\pi(\psi')$ comme repr\'esentation de $\tilde{G}(N_{\psi}-abd_{\rho})$, ce n'est autre que $\times_{(\rho',a',b')\in Jord(\psi)-(\rho,a,b)}Sp(b',St(a',\rho'))$. On pose encore $Sp(b,St(a,\rho))=\pi(\rho,A,B,\zeta)$ ce qui permet de d\'efinir $\pi(\rho,A,B+2,\zeta)$ et $\pi(\rho,A,B+1,\zeta)$ comme on l'a fait dans \ref{simple}. On suppose que $\psi \circ \Delta$ est sans multiplicit\'es et on a donc d\'efini $\pi(\psi)$ comme repr\'esentation de $\tilde{G}(n)$ dans \ref{independance}
\vskip 0.5cm
\noindent
On pose $$\tilde{\pi}(\psi)_{\rho,A,B,\zeta}:=
\oplus_{C\in ]B,A]}$$
$$(-1)^{A-C}<\rho\vert\,\vert^{\zeta B}, \cdots, \rho\vert\,\vert^{-\zeta C}> \times Jac^\theta_{\zeta(B+2), \cdots, \zeta C}(\pi(\psi')\times \pi(\rho,A,B+2,\zeta)) \times <\rho\vert\,\vert^{\zeta C}, \cdots, \rho\vert\,\vert^{-\zeta A}>
$$
$$
\oplus (-1)^{[(A-B+1)/2]} \pi(\psi')\times \pi(\rho,A,B+1,\zeta)\times <\rho\vert\,\vert^{\zeta B}, \cdots, \rho\vert\,\vert^{-\zeta B}>.
$$
\bf Th\'eor\`eme. \sl L'\'el\'ement, $\tilde{\pi}(\psi)_{\rho,A,B,\zeta}$ du groupe de Grothendieck associ\'e aux repr\'esentations irr\'eductibles de $\tilde{G}$ fait intervenir dans sa d\'ecomposition en repr\'esentations irr\'eductibles exactement une repr\'esentation irr\'eductible de $\tilde{G}$ dont la restriction au GL correspondant est encore irr\'eductible. C'est la repr\'esentation $\pi(\psi)$ et elle intervient avec le coefficient $+1$. Ou encore, pour tout $g\in GL(n,F)$ $\theta$-semi-simple et r\'egulier, $$tr\tilde{\pi}(\psi)(g,\theta)_{\rho,A,B,\zeta}=tr \pi(\psi)(g,\theta).$$
\rm
\vskip 0.5cm
\noindent
On traite d'abord le cas o\`u $A=B+1$. On v\'erifie que par d\'efinition
$
\tilde{\pi}(\psi)_{{\rho,A,B,\zeta}}$ $$
= <\rho\vert\,\vert^{\zeta B}, \cdots \rho\vert\,\vert^{-\zeta A}> \times \pi(\psi') \times <\rho\vert\,\vert^{\zeta A}, \cdots, \rho\vert\,\vert^{-\zeta B}> \ominus \pi(\psi') \times \pi(\rho,B+1,B+1,\zeta) \times \pi(\rho,B,B,\zeta).
$$L'action de $\theta$ sur le premier terme du membre de droite est l'action prolongeant naturellement celle de $\theta$ sur $\pi(\psi')$ par contre cette action sur $\pi(\psi')\times \pi(\rho,B+1,B+1,\zeta) \times \pi(\rho,B,B,\zeta)$ est plus difficile \`a d\'ecrire mais c'est celle que l'on a mise sur $\pi(\psi'')$ o\`u $\psi''$ se d\'eduit de $\psi'$ en ajoutant les quadruplets $(\rho,B+1,B+1,\zeta)$ et $(\rho,B,B,\zeta)$.

Il est imm\'ediat de v\'erifier qu'en tant que repr\'esentation de $GL(n,F)$ l'induite $<\rho\vert\,\vert^{\zeta B}, \cdots \rho\vert\,\vert^{-\zeta A}> \times \pi(\psi') \times <\rho\vert\,\vert^{\zeta A}, \cdots, \rho\vert\,\vert^{-\zeta B}>$ est de longueur exactement 2 avec comme sous-module $\pi(\psi')\times \pi(\rho,A,B,\zeta)$ (c'est-\`a-dire $ \pi(\psi)$) et comme quotient $\pi(\psi') \times \pi(\rho,B+1,B+1,\zeta) \times \pi(\rho,B,B,\zeta)$. Le probl\`eme est donc de v\'erifier que cette d\'ecomposition donne les ''bonnes'' actions de $\theta$. Pour le sous-module, il suffit de suivre les d\'efinitions et il n'y a pas de probl\`eme. Pour le quotient, si $\psi'$ est \'el\'ementaire et $B=0$, c'est exactement la d\'efinition. Si $B\geq 1/2$, c'est la d\'emonstration de \ref{remarque}. Reste le cas o\`u $\psi'$ n'est pas \'el\'ementaire. On fixe alors $(\rho',A',B',\zeta') \in Jord(\psi')$ tel que $A'>B'$ et on note $\psi''$ le morphisme qui se d\'eduit de $\psi'$ en rempla\c{c}ant $(\rho',A',B',\zeta')$ par $(\rho',A'-1,B'+1,\zeta')$ (ou en supprimant ce bloc si $A'=B'+1$). Toutes les actions de $\theta$ s'obtiennent alors comme on l'a vu en \ref{independance} en rempla\c{c}ant dans chaque repr\'esentation $\psi'$ par $\psi''$ d'o\`u une nouvelle repr\'esentation, $\sigma$, qu'il faut encore induire$$
<\rho'\vert\,\vert^{\zeta' B'}, \cdots, \rho'\vert\,\vert^{-\zeta' A'}>\times \sigma \times <\rho'\vert\,\vert^{\zeta' A'}, \cdots, \rho'\vert\,\vert^{-\zeta' B'}>,
$$
l'action de $\theta$ \'etant le prolongement canonique de l'action sur $\sigma$. Il est alors facile d'obtenir le r\'esultat par r\'ecurrence.
Cela termine la preuve du cas o\`u $A=B+1$.

On suppose donc maintenant que $A>B+1$.

La m\'ethode est tr\`es voisine de celle de \ref{casquasigeneral}. 

Soit $C\in ]B,A]$, notons $$\sigma_{C}:= <\rho\vert\,\vert^{\zeta B}, \cdots, \rho\vert\,\vert^{-\zeta C}> \times Jac^\theta_{\zeta (B+2), \cdots, \zeta C}\biggl(\pi(\psi')\times \pi(\rho,A,B+2,\zeta)\biggr) \times <\rho\vert\,\vert^{\zeta C}, \cdots, \rho\vert\,\vert^{-\zeta B}>.
$$
On montre que $\pi(\psi')\times \pi(\rho,A,B+2,\zeta)$ est l'unique sous-module irr\'eductible de l'induite:
$$
<\rho\vert\,\vert^{\zeta (B+2)}, \cdots, \rho\vert\,\vert^{\zeta A}>\times \pi(\psi')\times \pi(\rho,A-1,B+1,\zeta) \times <\rho\vert\,\vert^{-\zeta A}, \cdots, \rho\vert\,\vert^{-\zeta(B+2)}>
$$
l'inclusion \'etant compatible avec l'action de $\theta$. Cela et le fait que $Jac_{x}\pi(\psi')=0$ pour tout $x\in [\zeta B,\zeta A]$ entra\^{\i}nent ais\'ement que:
$$
Jac^\theta_{\zeta (B+2), \cdots, \zeta C}\biggl(\pi(\psi')\times \pi(\rho,A,B+2,\zeta)\biggr) \hookrightarrow $$
$$<\rho\vert\,\vert^{\zeta (C+1)}, \cdots,\rho\vert\,\vert^{\zeta A}> \times \pi(\psi')\times \pi(\rho,A-1,B+1) \times <\rho\vert\,\vert^{-\zeta A}, \cdots, \rho\vert\,\vert^{-\zeta(C+1)}>.
$$
On peut m\^eme pr\'eciser pour la suite que le terme de gauche est l'unique sous-module irr\'eductible du terme de droite, ce qui donne une construction pr\'ecise de $\sigma_{C}$. Evidemment ici, il est important que $\psi\circ \Delta$ soit sans multiplicit\'es.

On note $\psi_{>}$ le morphisme qui se d\'eduit de $\psi$ en enlevant tous les blocs de Jordan $(\rho,A',B',\zeta')$ tels que $B'\leq B$; \'eventuellement $\psi_{>}$ est trivial. 

Montrons encore que tout sous-quotient irr\'eductible de $\sigma_{C}$, irr\'eductible en tant que repr\'esentation de $\tilde{G}(n)$ et de restriction irr\'eductible \`a $GL(n,F)$ est de la forme $\pi(\psi_{>})\times \sigma'$ o\`u $\sigma'$ est une repr\'esentation convenable d'un groupe lin\'eaire dont le support cuspidal
est une collection de repr\'esentations $\rho\vert\,\vert^x$ avec $\vert x\vert \leq A$.  Montrons cette propri\'et\'e par r\'ecurrence sur $\sum_{(\rho,A',B',\zeta') \in Jord_{\rho}(\psi_{>})}(A'-B'+1)$; comme il y a une difficult\'e pour suivre l'action de $\theta$ alors que la conclusion ne fait pas intervenir cette action, on transforme la propri\'et\'e en oubliant $\tilde{G}$ et en parlant simplement de repr\'esentation irr\'eductible $\theta$-invariante de $GL(N,F)$. Soit $(\rho, A_{0},B_{0},\zeta_{0})$ l'\'el\'ement de $Jord_{\rho }(\psi_{>})$ avec $A_{0}$ maximum et soit $\tau$ un sous-quotient irr\'eductible de $\sigma_{C}$ pour l'action de $GL(N,F)$; on suppose que $\tau$ est  $\theta$-invariant. Suivant ce que l'on a d\'ej\`a vu, il existe $x_{0}$ un demi-entier tel que $\vert x_{0}\vert \leq A_{0}$ et $[x_{0},-\zeta A_{0}]$ est un segment et v\'erifiant:
$$
Jac_{x_{0}, \cdots, -\zeta_{0}A_{0}}\tau \neq 0.
$$
Un calcul simple montre que $x_{0}=\zeta_{0}B_{0}$; si $B_{0}=A_{0}$, on montre tout de suite que $\tau$ s'\'ecrit sous la forme $<\rho\vert\,\vert^{\zeta_{0}A_{0}}, \cdots, \rho\vert\,\vert^{-\zeta_{0}A_{0}}> \times \tau'$ o\`u $\tau'$ est irr\'eductible et $\theta$-invariante et v\'erifie l'analogue de $\tau$ quand on a enlev\'e le bloc $(\rho,A_{0},B_{0},\zeta_{0})$. Par contre si $B_{0}<A_{0}$, on montre encore que
$Jac^\theta_{\zeta_{0}B_{0}, \cdots, -\zeta_{0}A_{0}}\tau\neq 0$. Puis on montre que $\tau$, en tant que repr\'esentation de $GL$  est de la forme $<\rho\vert\,\vert^{\zeta_{0}B_{0}}, \cdots, \rho\vert\,\vert^{-\zeta_{0}A_{0}}> \times \tau' \times <\rho\vert\,\vert^{\zeta_{0}A_{0}}, \cdots, \rho\vert\,\vert^{-\zeta_{0}B_{0}}>$, o\`u $\tau'$ est un sous-quotient irr\'eductible, $\theta$-invariante de l'analogue de $\sigma_{C}$ obtenu en rempla\c{c}ant le bloc $(\rho,A_{0},B_{0},\zeta_{0})$ de $Jord(\psi)$ par $(\rho,A_{0}-1,B_{0}+1,\zeta_{0})$. Cela permet de conclure.

 Soit donc $\tilde{\sigma}$ un sous-quotient irr\'eductible de $\sigma_{C}$ dont on suppose que la restriction \`a $GL(n,F)$ est irr\'eductible; on note $\sigma$ cette restriction; on \'ecrit $\sigma$ sous la forme $\pi(\psi_{>})\times \sigma'$ conform\'ement \`a ce que l'on vient de voir; on v\'erifie ais\'ement que $\sigma'$ est n\'ecessairement $\theta$ invariante.

Comme dans \ref{casquasigeneral}, il existe $x_{0}$ un demi-entier tel que $\vert x_{0}\vert  \leq A$, $[x_{0},-\zeta A]$ soit un segment et tel que $Jac^\theta_{x_{0}, \cdots, -\zeta A}\sigma' \neq 0$. On en d\'eduit que $Jac^\theta_{x_{0}, \cdots, -\zeta A}\sigma\neq 0$ car $\rho\vert\,\vert^x \times \pi(\psi_{>})$ et $\pi(\psi_{>})\times \rho\vert\,\vert^x$ sont des induites irr\'eductibles pour tout $x$ verifiant $\vert x\vert \leq A$. On fixe donc un tel $x_{0}$.

Montrons  que n\'ecessairement $x_{0}\in [\zeta B, \zeta A]$:

ce que l'on obtient imm\'ediatement est que soit $x_{0}$ a la propri\'et\'e annonc\'ee soit $x_{0}=\zeta'B'$ pour un bon choix de $\zeta'$ et $B'$ tel qu'il existe $A'$ avec $(\rho,A',B',\zeta')\in Jord(\psi')$. Si cette derni\`ere \'eventualit\'e se produisait, on remarque que $\pi(\psi')=\pi(\psi_{>})\times \pi(\psi_{<})$ o\`u $\psi_{<}$ correspondant \`a tous les blocs de Jordan $(\rho,A'',B'',\zeta'')$ avec $B''<B$; mais le support cuspidal de $\pi(\psi_{<})$ ne contient aucun \'el\'ement du type $\rho\vert\,\vert^{\pm x}$ avec $x\in [B,A]$; donc il faudrait n\'ecessairement $$Jac_{-\zeta B, \cdots, -\zeta A}\biggl(<\rho\vert\,\vert^{\zeta B}, \cdots, \rho\vert\,\vert^{-\zeta C}>  \times <\rho\vert\,\vert^{\zeta (C+1)}, \cdots, \rho\vert\,\vert^{\zeta A}> 
$$
$$\times \pi(\rho,A-1,B+1,\zeta) \times <\rho\vert\,\vert^{-\zeta A}, \cdots, \rho\vert\,\vert^{-\zeta (C+1)}> \times <\rho\vert\,\vert^{\zeta C}, \cdots, \rho\vert\,\vert^{-\zeta B}>\biggr) \neq 0.
$$
Ce qui est exclu. D'o\`u l'assertion.

On montre la m\^eme assertion si $\tilde{\sigma}$ est $\pi(\psi')\times \pi(\rho,A,B+1,\zeta) \times \pi(\rho,B,B,\zeta)$; dans ce cas, on a m\^eme n\'ecessairement $x_{0}=\zeta (B+1)$.

Dans toute la suite de la d\'emonstration, on fixe $\tilde{\sigma}$ une repr\'esentation irr\'eductible intervenant dans la d\'ecomposition de $\tilde{\pi}(\psi)_{\rho,A,B,\zeta}$ et $x_{0}$ v\'erifiant la propri\'et\'e pr\'ec\'edente pour $\tilde{\sigma}$. Comme dans \ref{casquasigeneral} le point est de d\'emontrer que $x_{0}=\zeta B$. On suppose donc d'abord que $x_{0}\in ]\zeta B,\zeta A]$ et on veut alors montrer simplement que $Jac^\theta_{x_{0}}\sigma=0$. Exactement comme dans loc.cit., on v\'erifie qu'il suffit de montrer que $Jac^\theta_{x_{0}}\tilde{\pi}(\psi)_{\rho,A,B,\zeta}=0$. Le cas o\`u $x_{0}\in ]\zeta(B+1),\zeta A]$ est imm\'ediat avec la description pr\'ecise de $\sigma_{C}$; les termes qui contribuent au calcul est $\sigma_{\vert x_{0}\vert}$ et $\sigma_{\vert x_{0}\vert -1}$ et ils s'\'eliminent l'un l'autre. Il reste donc \`a voir le cas plus difficile o\`u $x_{0}=\zeta(B+1)$.  On  a vu ci-dessus que, pour tout $C\in ]B,A]$:
$$
\sigma_{C}=<\rho\vert\,\vert^{\zeta B}, \cdots, \rho\vert\,\vert^{-\zeta C}>\times $$
$$< <\rho\vert\,\vert^{\zeta (C+1)}, \cdots, \rho\vert\,\vert^{\zeta A}>\times \pi(\psi')\times \pi(\rho,A-1,B+1,\zeta) \times <\rho\vert\,\vert^{-\zeta A}, \cdots, \rho\vert\,\vert^{-\zeta(C+1)}>> $$
$$\times <\rho\vert\,\vert^{\zeta C}, \cdots, \rho\vert\,\vert^{-\zeta B}>,
$$
o\`u pour $C=A$ le terme du milieu est simplement $\pi(\psi')\times \pi(\rho,A-1,B+1,\zeta)$.  Pour $C\in ]B+1,A[$, on r\'ecrit $\sigma_{C}$. Pour qu'il y ait une valeur de $C$, on peut supposer que $A-1>B+1$, d'o\`u $\pi(\psi')\times \pi(\rho,A-1,B+1,\zeta)$ est par d\'efinition un sous-module (avec action de $\theta$) de $$
<\rho\vert\,\vert^{\zeta(B+1)}, \cdots, \rho\vert\,\vert^{-\zeta (A-1)}> \times \pi(\psi')\times \pi(\rho,A-2,B+2,\zeta) \times <\rho\vert\,\vert^{\zeta(A-1)}, \cdots, \rho\vert\,\vert^{-\zeta(B+1)}>.
$$En revenant \`a $\sigma_{C}$ pour ces valeurs de $C$,
$$
\sigma_{C}\hookrightarrow <\rho\vert\,\vert^{\zeta B}, \cdots, \rho\vert\,\vert^{-\zeta C}>\times
$$
$$
 <\rho\vert\,\vert^{\zeta (C+1)}, \cdots, \rho\vert\,\vert^{\zeta A}>\times<\rho\vert\,\vert^{\zeta(B+1)}, \cdots, \rho\vert\,\vert^{-\zeta (A-1)}> \times 
 $$
 $$
 \pi(\psi')\times \pi(\rho,A-2,B+2,\zeta)\times
 $$
 $$
 <\rho\vert\,\vert^{\zeta(A-1)}, \cdots, \rho\vert\,\vert^{-\zeta(B+1)}>\times <\rho\vert\,\vert^{-\zeta A}, \cdots, \rho\vert\,\vert^{-\zeta (C+1)}>
 $$
 $$
 \times <\rho\vert\,\vert^{\zeta C}, \cdots, \rho\vert\,\vert^{-\zeta B}>.
 $$
Comme l'action de $\theta$ sur l'induite est l'action qui se d\'eduit naturellement de celle mise sur $\pi(\psi')\times \pi(\rho,A-2,B+2,\zeta)$, on obtient:
 $$
 Jac^\theta_{\zeta (B+1), \cdots, -\zeta A}\sigma_{C}\hookrightarrow
 <\rho\vert\,\vert^{\zeta B}, \cdots, \rho\vert\,\vert^{-\zeta C}>\times $$
 $$
  <\rho\vert\,\vert^{\zeta (C+1)}, \cdots, \rho\vert\,\vert^{\zeta (A-1)}>\times \pi(\psi')\times \pi(\rho,A-2,B+2,\zeta) \times <\rho\vert\,\vert^{-\zeta (A-1)}, \cdots, \rho\vert\,\vert^{-\zeta C}> \eqno(3)
  $$
  $$
  \times <\rho\vert\,\vert^{\zeta C}, \cdots, \rho\vert\,\vert^{-\zeta B}>,
  $$
  avec l'action naturelle de $\theta$. Un calcul de module de Jacquet standard montre que cette inclusion se factorise par l'unique sous-module de (3) qui est en fait $Jac^\theta_{\zeta (B+3), \cdots, \zeta C}(\pi(\psi')\times \pi(\rho,A-1,B+3,\zeta))$. On a donc montr\'e 
  $$
  Jac^\theta_{\zeta (B+1), \cdots, -\zeta A}\sigma_{C}=$$
  $$
  <\rho\vert\,\vert^{\zeta B}, \cdots, \rho\vert\,\vert^{-\zeta C}>$$
  $$\times Jac^\theta_{\zeta (B+3), \cdots, \zeta C}(\pi(\psi')\times \pi(\rho,A-1,B+3,\zeta))
  \times$$
  $$ <\rho\vert\,\vert^{\zeta C}, \cdots, \rho\vert\,\vert^{-\zeta B}>.\eqno(3)$$Rappelons que ce terme n'appara\^{\i}t que si $(A-1)>(B+1)$, car $C\in ]B+1,A-1]$ et il 
 peut se r\'ecrire $Jac^\theta_{\zeta(B+1)}\sigma'_{C}$ o\`u $\sigma'_{C}$ est l'analogue de $\sigma_{C}$ pour $\psi$ remplac\'e par le morphisme o\`u le bloc $(\rho,A,B,\zeta)$ est remplac\'e par $(\rho,A-1,B+1,\zeta)$.   
  
  On v\'erifie par un calcul de module de Jacquet standard que pour $C=B+1$, $Jac_{\zeta (B+1), \cdots, -\zeta A}\sigma_{B+1}=0$ et donc a fortiori que $Jac^\theta_{\zeta (B+1), \cdots, -\zeta A}\sigma_{B+1}=0$.

  On suppose maintenant que $C=A$, et on calcule:
  $$
  Jac^\theta_{\zeta (B+1), \cdots, -\zeta A} \biggl(<\rho\vert\,\vert^{\zeta B}, \cdots, \rho\vert\,\vert^{-\zeta A}>\times \pi(\psi')\times \pi(\rho,A-1,B+1,\zeta) \times <\rho\vert\,\vert^{\zeta A}, \cdots, \rho\vert\,\vert^{-\zeta (B+1)}>\biggr)=
  $$
  $$
  Jac^\theta_{\zeta(B+1)}\biggl(\pi(\psi')\times \pi(\rho,A-1,B+1,\zeta)\biggl).\eqno(4)
  $$

  Il reste donc \`a calculer
  $$
  Jac^\theta_{\zeta (B+1), \cdots, -\zeta A}\biggl(\pi(\psi')\times \pi(\rho, A,B+1,\zeta) \times \pi(\rho,B,B,\zeta)\biggr).
  $$
  Comme on a suppos\'e que $A>B+1$ (le cas d'\'egalit\'e ayant d\'ej\`a \'et\'e trait\'e), on r\'ecrit $\pi(\psi')\times \pi(\rho, A,B+1,\zeta) \times \pi(\rho,B,B,\zeta)$ comme l'unique sous-module irr\'eductible de
  $$
 <\rho\vert\,\vert^{\zeta(B+1)}, \cdots, \rho\vert\,\vert^{-\zeta A}>\times  \pi(\psi')\times \pi(\rho,A-1,B+2,\zeta) \times \pi(\rho,B,B,\zeta) \times <\rho\vert\,\vert^{\zeta A}, \cdots, \rho\vert\,\vert^{-\zeta(B+1)}>.
 $$
 Ainsi le $Jac^\theta$ cherch\'e est exactement l'induite du milieu, c'est-\`a-dire $\pi(\psi')\times \pi(\rho,A-1,B+2,\zeta)\times \pi(\rho,B,B,\zeta)$ avec l'action de $\theta$ mise sur cette induite. On va v\'erifier que l'on peut r\'ecrire le r\'esultat sous la forme:
 $$
 Jac^\theta_{\zeta(B+1), \cdots, -\zeta A}\biggl( 	\pi(\psi')\times \pi(\rho,A,B+1,\zeta) \times \pi(\rho,B,B,\zeta)\biggr)=
 $$
 $$Jac^\theta_{\zeta(B+1)}\biggl(\pi(\psi')\times \pi(\rho,A-1,B+2,\zeta)\times \pi(\rho,B+1,B+1,\zeta)\biggr); \eqno(5)
 $$
gr\^ace aux d\'efinitions on \'ecrit:
$$
\pi(\psi')\times \pi(\rho,A,B+1,\zeta) \times \pi(\rho,B,B,\zeta)\hookrightarrow$$
$$
<\rho\vert\,\vert^{\zeta(B+1)}, \cdots,\rho\vert\,\vert^{-\zeta A}>$$
$$
 \times \pi(\rho,A-1,B+2,\zeta)\times \pi(\psi')\times \pi(\rho,B,B,\zeta) \times$$
 $$
  <\rho\vert\,\vert^{\zeta A}, \cdots, \rho\vert\,\vert^{-\zeta (B+1)}>.
$$
D'o\`u ais\'ement l'identification du $Jac^\theta$ cherch\'e avec l'induite du milieu, c'est-\`a-dire
$\pi(\rho,A-1,B+2,\zeta)\times \pi(\psi')\times \pi(\rho,B,B,\zeta)$. Mais on a vu en \ref{proprietes} que cette induite est aussi l'unique sous-module irr\'eductible de l'induite:
$$
\rho\vert\,\vert^{\zeta(B+1)} \times \pi(\rho,A-1,B+2,\zeta)\times \pi(\psi')\times \pi(\rho,B+1,B+1,\zeta) \times \rho\vert\,\vert^{-\zeta(B+1)},
$$
inclusion compatible avec les actions de $\theta$. D'o\`u le r\'esultat annonc\'e.

 On remarque que si $A-1=B+1$, ce terme (5) est le m\^eme que (4); en tenant compte des signes dans la d\'efinition de $\tilde{\pi}(\psi)_{\rho,A,B,\zeta}$, (4) est accompagn\'e de $+1$ tandis que (5) est accompagn\'e de $(-1)^{[(A-B+1)/2]}=-1$. Donc si $A=B+2$, on vient de montrer que $Jac^\theta_{\zeta(B+1), \cdots, -\zeta A}\tilde{\pi}(\psi)_{\rho,A,B,\zeta}=0$ et on conclut comme en \ref{quasigeneral}.

 Supposons que $A>B+2$.
 On note $\psi_{-1}$, le morphisme qui se d\'eduit de $\psi$ en rempla\c{c}ant le bloc de Jordan $(\rho,A,B,\zeta)$ par $(\rho,A-1,B+1,\zeta)$.
En remettant les signes, on a donc aussi montr\'e que dans le groupe de Grothendieck de $\tilde{G}(N')$ pour $N'$ convenable:
$$
Jac^\theta_{\zeta(B+1), \cdots, -\zeta A}\tilde{\pi}(\psi)_{\rho,A,B,\zeta}=
$$
$$
Jac^\theta_{\zeta(B+1)}\pi(\psi_{-1})\ominus Jac_{\zeta(B+1)}^\theta\tilde{\pi}(\psi_{-1})_{\rho,A-1,B+1,\zeta}.$$
On conclut encore comme en \ref{quasigeneral}. Cela termine la preuve.

\subsection{Le cas g\'en\'eral \label{general}}
Le r\'esultat de \ref{quasigeneral} n'est sans doute pas vrai dans le cas g\'en\'eral, c'est-\`a-dire o\`u l'on ne fait pas l'hypoth\`ese que la restriction \`a la diagonale est discr\`ete, on rappelle d'ailleurs que pour le moment nous n'avons pas mis d'action de $\theta$ sur $\pi(\psi)$ sans l'hypoth\`ese que $\psi$ est de restriction discr\`ete \`a la diagonale. Pour ce cas, on r\'ealise $\pi(\psi)$ comme un module de Jacquet convenable. Pr\'ecis\'ement, pour tout $\rho$; 
 on met un ordre total sur $Jord_{\rho}(\psi)$ en posant, pour $(\rho,A,B,\zeta), (\rho,A',B',\zeta')\in Jord(\psi)$,

$(\rho, A,B,\zeta) \geq  (\rho,A',B',\zeta')$ si

soit $A>A'$, soit $A=A'$ mais $B>B'$ soit $A=A', B=B'$ mais $B>0$ et $\zeta=+$.

\noindent
Il n'y a aucune justification \`a avoir priviligi\'e $A$ sur $B$, on aurait tout aussi bien pu faire l'inverse. Il n'y a pas non plus de justification pour avoir privil\'egi\'e le signe $+$ sur le signe $-$. On \'etudiera pr\'ecis\'ement l'influence de nos choix en \ref{comparaisondesnormalisations}

Ayant fait ces choix, on dit qu'un morphisme $\tilde{\psi}$ de $W_{F}\times SL(2,{\mathbb C}) \times  SL(2,{\mathbb C})$ dans $GL(N,F)$ avec $N\geq n$ domine $\psi$ s'il est de restriction discr\`ete \`a la diogonale et s'il  existe une bijection de $Jord(\tilde{\psi})$ sur $Jord(\psi)$ not\'ee $bij$ v\'erifiant:

pour tout $\rho$, $bij$ envoie $Jord_{\rho}(\tilde{\psi})$ sur $Jord_{\rho}(\psi)$ et $bij$ pr\'eserve l'ordre. De plus pour tout $(\rho,\tilde{A},\tilde{B},\tilde{\zeta})\in Jord_{\rho}(\tilde{\psi})$, en notant $(\rho,A,B,\zeta)$, l'image de cet \'el\'ement par $bij$, on a:
$$
\zeta=\tilde{\zeta}, A-B=\tilde{A}-\tilde{B}.
$$La bijection est clairement uniquement d\'etermin\'ee par le seul fait qu'elle respecte l'ordre.

Comme on part plut\^ot de $\psi$ que de $\tilde{\psi}$, on notera plut\^ot $(\rho,\tilde{A}_{\rho,A,B,\zeta},\tilde{B}_{\rho,A,B,\zeta},\zeta)$ l'image inverse de $(\rho,A,B,\zeta) \in Jord_{\rho}(\psi)$.

On pose alors pour tout $(\rho,A,B,\zeta) \in Jord(\psi)$,  ${\cal E}_{\rho,A,B,\zeta}:=\cup_{D\in [\tilde{B}_{\rho,A,B,\zeta},B[} [\zeta D, \zeta (D+A-B)]$, union d'ensembles avec multiplicit\'es. On met sur ${\cal E}_{\rho,A,B,\zeta}$ un ordre total
 en consid\'erant d'abord l'ordre d\'ecroissant sur les $D$ puis l'ordre sur le segment $[\zeta D, \zeta (D+A-B)]$.Pour $\rho$ fix\'e, l'ensemble ${\cal E}_{\rho,\tilde{\psi}}\cup_{(A,B,\zeta) \in Jord_{\rho}(\psi)}{\cal E}_{\rho,A,B,\zeta}$ est lui aussi totalement ordonn\'e, en mettant d'abord l'ordre croissant sur les \'el\'ements de $Jord_{\rho}(\psi)$ puis l'ordre sur les ensembles ${\cal E}_{\rho,A,B,\zeta}$ que l'on vient de d\'efinir. Quand on fait encore l'union sur toues les repr\'esentation $\rho$, on obtient un ensemble ${\cal E}_{\tilde{\psi}}$ partiellement ordonn\'e. On pose alors, en respectant dans le terme de droite ci-dessous, l'ordre partiel:
$$
\tilde{\pi}(\psi):= Jac^\theta_{x \in {\cal E}_{\tilde{\psi}}}\pi(\tilde{\psi}).
$$Supposons que 
Cela donne l'action de $\theta$ sur la repr\'esentation $\tilde{\pi}(\psi)$.

\vskip 0.5cm
\noindent
\bf Th\'eor\`eme, d\'efinition: \sl $\tilde{\pi}(\psi)$ est une repr\'esentation irr\'eductible de $\tilde{G}(n)$ dont la restriction \`a $GL(n,F)$ reste irr\'eductible. Pr\'ecis\'ement cette restriction est
$\times_{(\rho,a,b)\in Jord(\psi)}Sp(b,St(a,\rho))$. De plus l'action de $\theta$ sur $\tilde{\pi}(\psi)$ est ind\'ependante du choix fait de $\tilde{\psi}$ dominant $\psi$. On posera donc $\pi(\psi)=\tilde{\pi}(\psi)$.

\rm
\vskip 0.5cm
\noindent
On montre ce th\'eor\`eme par r\'ecurrence sur $\vert Jord(\psi'')\vert $ pour tous les morphismes $\psi''$ de $W_{F}\times SL(2,{\mathbb C}) \times SL(2,{\mathbb C})$ dans un groupe $GL(N'',{\mathbb C})$; en effet si ce nombre est r\'eduit \`a $1$, $\psi''\circ \Delta$ est n\'ecessairement sans multiplicit\'es et le th\'eor\`eme est le corollaire de  \ref{proprietes} appliqu\'e de proche en proche de $C=\tilde{A}$ jusqu' \`a $C=A+1$.

En fait il y a 2 assertions dans ce th\'eor\`eme et on va les d\'emontrer successivement. On montre d'abord que la restriction de $\tilde{\pi}(\psi)$ \`a $GL(n,F)$ est ce que l'on attend \`a savoir $\times_{(\rho,a,b)\in Jord(\psi)}Sp(b,St(a,\rho))$.

On fixe $(\rho,A,B,\zeta) \in Jord(\psi)$, on note $\psi_{<}$ le morphisme qui correspond \`a la repr\'esentation de $W_{F}\times SL(2,{\mathbb C}) \times SL(2,{\mathbb C})$ dont les blocs de Jordan sont ceux de $$Jord(\psi)-Jord_{\rho}(\psi) \cup \{(\rho,A',B',\zeta')\in Jord_{\rho}(\psi), (\rho,A',B',\zeta')<(\rho,A,B,\zeta)\} $$
$$\cup \{(\rho, \tilde{A}',\tilde{B'},\tilde{\zeta}') \in Jord_{\rho}(\tilde{\psi}); (\rho, \tilde{A}',\tilde{B}',\tilde{\zeta}')> (\rho, \tilde{A}_{\rho,A,B,\zeta}, \tilde{B}_{\rho,A,B,\zeta},\zeta)\}.$$

 On admet par r\'ecurrence que la restriction de  $\tilde{\pi}(\psi_{<})$ au $GL$ correspondant est ce que l'on attend (on note $\pi(\psi_{<})$ cette restriction). On va montrer que:
$$
Jac^\theta_{x\in {\cal E}_{A,B,\zeta}}(\pi(\psi_{<})\times \pi(\rho,\tilde{A}_{\rho,A,B,\zeta}, \tilde{B}_{\rho,A,B,\zeta},\zeta))=\pi(\psi_{<})\times \pi(\rho,A,B,\zeta). \eqno(*)
$$
De proche en proche cela donnera la premi\`ere assertion du th\'eor\`eme.

Pour cela, on  montre que pour tout $D\in [\tilde{B}_{\rho,A,B,\zeta},B]$:
$$
Jac^\theta_{x \in[ \zeta D, \cdots, \zeta (D+A-B)]}\biggl(\pi(\psi_{<})\times \pi(\rho,D+A-B,D,\zeta) \biggr)=\eqno(1)$$
$$
\pi(\psi_{<})\times Jac^\theta_{x \in[ \zeta D, \cdots, \zeta (D+A-B)]} \pi(\rho,D+A-B,D,\zeta)= \pi(\psi_{<})\times \pi(\rho, D+A-B-1,D-1,\zeta).\eqno(2)
$$
Soit  ${\cal I}$ un sous-ensemble de $[\zeta D, \zeta(D+A-B)]$ tel que $Jac_{x\in {\cal I}}\pi(\rho,D+A-B,D,\zeta)\neq 0$. On se rappelle que $\pi(\rho,D+A-B,D,\zeta)$, comme repr\'esentation du $GL$ correspondant, est l'unique sous-module irr\'eductible de l'induite:
$$
<\rho\vert\,\vert^{\zeta D}, \cdots, \rho\vert\,\vert^{\zeta(D+A-B)}> \times \pi(\rho, D+A-B-1,D-1,\zeta)\times <\rho\vert\,\vert^{-\zeta (D+A-B}, \cdots, \rho\vert\,\vert^{-\zeta D}>. \eqno(3)
$$
Comme pour tout $x\in [\zeta D, \zeta (D+A-B)]$, $Jac_{x}\pi(\rho, D+A-B-1,D-1,\zeta)=0$ et qu'il en est de m\^eme pour $Jac_{x}<\rho\vert\,\vert^{-\zeta (D+A-B)}, \cdots, \rho\vert\,\vert^{-\zeta D}>$, la non nullit\'e suppos\'ee entra\^{\i}ne que $Jac_{x\in {\cal I}}<\rho\vert\,\vert^{\zeta D}, \cdots, \rho\vert\,\vert^{\zeta(D+A-B)}> \neq 0$. En particulier ${\cal I}$ est un sous-segment de $[\zeta D, \zeta (D+A-B)]$ commen\c{c}ant par $\zeta D$. Par les formules standard, pour avoir l'\'egalit\'e de (1) avec le premier membre de (2),  la seule chose \`a d\'emontrer est donc que pour tout ${\cal I}'$ sous-segment de $[\zeta D, \zeta (D+A-B)]$ contenant l'extr\^emit\'e $\zeta (D+A-B)$, $Jac_{x\in {\cal I}'}\pi(\psi_{<})=0$

On d\'ecompose encore $\pi(\psi_{<})$ en $\pi(\psi') \times \pi(\psi_{1})\times \pi(\psi_{2})$ o\`u $Jord (\psi')=Jord(\psi)-Jord_{\rho}(\psi)$, $Jord(\psi_{1})=\{(\rho,A',B',\zeta') \in Jord(\psi); (\rho,A',B',\zeta')< (\rho,A,B,\zeta)\}$, $Jord(\psi_{2})=\{(\rho,\tilde{A}',\tilde{B}',\zeta')\in Jord(\tilde{\psi});
(\rho,\tilde{A}',\tilde{B}',\zeta')> (\rho,\tilde{A}_{\rho,A,B,\zeta},\tilde{B}_{\rho,A,B,\zeta},\zeta)\}.$ Pour ${\cal I}'$ un sous-ensemble de $[\zeta D, \zeta(D+A-B)]$, on a $Jac_{x\in {\cal I}}'\pi(\psi')=0$ pour des raisons de support cuspidal $Jord_{\rho}(\psi')=\emptyset$. On a aussi $Jac_{x \in {\cal I}'}\pi(\psi_{2})=0$ car sinon  il faudrait que le premier \'el\'ement de ${\cal I}$, not\'e $z$, v\'erifie $Jac_{z}\pi(\psi_{2})\neq 0$; il serait alors de la forme $\tilde{\zeta}'\tilde{B'}$ avec $(\rho, \tilde{A}',\tilde{B}',\tilde{\zeta}')\in Jord(\psi_{2})$. Or ceci est exclu puisque l'on a certainement $\tilde{B'}>\tilde{A}_{\rho,A,B,\zeta}$ par l'hypoth\`ese  $\tilde{\psi}$ de restriction discr\`ete \`a la diagonale, d'o\`u $\tilde{B}'>\tilde{A}_{\rho,A,B,\zeta}\geq D+\tilde{A}_{\rho,A,B,\zeta}-\tilde{B}_{\rho,A,B,\zeta}=D+A-B$. Il reste \`a voir que $Jac_{x\in [\zeta D, \zeta(D+A-B)]}\pi(\psi_{1})=0$; or $D+A-B >A$ d'o\`u $D+A-B>A'$ pour tout $(\rho,A',B',\zeta')\in Jord(\psi_{1})$ et $\rho\vert\,\vert^{\zeta (D+A-B)}$ n'est pas dans le support cuspidal de $\pi(\psi_{1})$. On a donc prouv\'e l'\'egalit\'e de (1) avec le premier terme \'ecrit dans (2).

Et l'\'egalit\'e de (2) r\'esulte de la repr\'esentation de $\pi(\rho,D+A-B,D,\zeta)$ donn\'ee en (3). Ceci termine la preuve de (*).

Il reste \`a montrer que l'action de $\theta$ sur $\tilde{\pi}(\psi)$ est ind\'ependante du choix fait de $\tilde{\psi}$. Fixons 2 choix $\tilde{\psi}_{i}$ pour $i=1,2$ dominant $\psi$; il est clair qu'il existe $\tilde{\psi}$ qui domine les trois morphismes $\tilde{\psi}_{i}$ pour $i=1,2$ et $\psi$ lui-m\^eme. On a d\'efini les ensembles ordonn\'es ${\cal E}_{\tilde{\psi}}$, ${\cal E}_{\tilde{\psi}_{i}}$ pour $i=1,2$; ces 2 derniers ensembles sont inclus dans ${\cal E}_{\tilde{\psi}}$. Fixons $i$ et posons ${\cal F}_{i}:={\cal E}_{\tilde{\psi}}-{\cal E}_{\tilde{\psi}_{i}}$. Cet ensemble h\'erite d'un ordre. Il r\'esulte du corollaire de \ref{proprietes} appliqu\'e de nombreuses fois dans l'ordre impos\'e par l'ordre de ${\cal F}_{i}$ que:
$$
Jac^\theta_{{\cal F}_{i}}\pi(\tilde{\psi})=\pi(\tilde{\psi}_{i}).
$$
Le point est donc de d\'emontrer que l'on a:
$$
Jac^\theta_{{\cal E}_{\tilde{\psi}}}= Jac^\theta_{{\cal E}_{\tilde{\psi}_{i}}}Jac^\theta_{{\cal F}_{i}}\pi(\tilde{\psi}),
$$
c'est-\`a-dire que l'on peut perturber l'ordre de ${\cal E}_{\tilde{\psi}}$; on ne perturbe pas n'importe comment (sinon ce serait faux) on utilise simplement le fait que $Jac^\theta_{x}Jac^\theta_{y}=Jac^\theta_{y}Jac^\theta_{x}$ si $\vert y-x\vert >1$ ce qui ne pose pas de probl\`eme  (cf \ref{notationdujac}). D'o\`u le th\'eor\`eme.

\subsubsection{Exemples et commentaires\label{exemple}}
Dans le paragraphe  \ref{general}, on a fait des choix. Ici, on va donner un exemple qui prouve que le r\'esultat n'est pas indiff\'erent aux choix. On consid\`ere l'exemple o\`u $G=GL(4,F)$ et o\`u $\psi$ est trivial sur $W_{F}$ et sur $SL(2,{\mathbb C})\times SL(2,{\mathbb C})$ est $(Rep_{2}\otimes 1) \oplus (1\otimes Rep_{2})$ o\`u $Rep_{2}$ est la repr\'esentation irr\'eductible de dimension 2 de $SL(2,{\mathbb C})$. La repr\'esentation $\pi(\psi)$ consid\'er\'ee uniquement comme repr\'esentation de $GL(4,F)$ est l'induite irr\'eductible $St(2)\times Sp(2)$ o\`u $St(2)$ est la repr\'esentation de Steinberg de $GL(2,F)$ et $Sp(2)$ est la repr\'esentation triviale de ce groupe.

Pour d\'ecrire l'action de $\theta$, on consid\`ere $GL(6,F)$ et le morphisme $\tilde{\psi}$ trivial sur $W_{F}$ et valant $Rep_{4}\otimes 1 \oplus 1\otimes Rep_{2}$ \`a valeurs dans $GL(6,{\mathbb C})$; ici $Rep_{4}$ est la repr\'esentation irr\'eductible de dimension 4. Alors:
$$
\pi(\tilde{\psi}) \hookrightarrow \rho\vert\,\vert^{-1/2}\times St(4) \times \rho\vert\,\vert^{1/2}
$$
avec l'action de $\theta$ qui prolonge naturellement l'action de $\theta$ mise sur $St(4)$. Et on a par d\'efintion, $\pi(\psi)=Jac^\theta_{3/2}\pi(\tilde{\psi})$ et cela entra\^{\i}ne que l'action de $\theta$ sur $\pi(\psi)$ est celle qui se d\'eduit naturellement de l'action de $\theta$ mise sur $St(2)$ en utilisant l'inclusion dans l'induite:
$$
\pi(\psi) \hookrightarrow \rho\vert\,\vert^{-1/2}\times St(2) \times \rho\vert\,\vert^{1/2}.\eqno(1)
$$
Pour distinguer, on va noter $\pi_{st}(\psi)$ cette repr\'esentation de $\tilde{G}(4)$.
L'autre choix raisonnable aurait \'et\'e d'inverser les r\^oles de $St$ et $Sp$ c'est-\`a-dire de mettre sur $\pi(\psi)$ l'action de $\theta$ qui prolonge naturellement celle mise sur $Sp(2)$ en utilisant l'inclusion:
$$
\pi(\psi) \hookrightarrow \rho\vert\,\vert^{1/2}\times Sp(2) \times \rho\vert\,\vert^{-1/2}.\eqno(2)
$$
On note $\pi_{sp}(\psi)$ cette repr\'esentation de $\tilde{G}(4)$. Les 2 repr\'esentations $\pi_{st}(\psi)$ et $\pi_{sp}(\psi)$ ne sont pas isomorphes comme nous allons le v\'erifier.

On consid\`ere la s\'erie principale:
$$
\sigma:=\rho\vert\,\vert^{-1/2}\times \rho\vert\,\vert^{1/2}\times \rho\vert\,\vert^{-1/2}\times \rho\vert\,\vert^{1/2}.
$$
Cette induite est munie d'une action naturelle de $\theta$ et l'induite (1) en est un sous $\tilde{G}(4)$ module. Ainsi $\pi_{st}(\psi)$ est un sous- $\tilde{G}(4)$-module. On v\'erifie par un calcul de module de Jacquet que $\sigma$, en tant que repr\'esentation de $GL(4,F)$ est de longueur 4: elle contient avec multiplicit\'e 1 la repr\'esentation $Sp(2)\times Sp(2)$ et la repr\'esentation $St(2)\times St(2)$ et avec multiplicit\'e 2 la repr\'esentation $Sp(2)\times St(2)$. On v\'erifie aussi que $\pi_{sp}(\psi)$ est un quotient de $\sigma$; en effet $Sp(2)$ muni de l'action de $\theta$ que nous avons mis est quotient de l'induite ''du milieu'', c'est-\`a-dire $\rho\vert\,\vert^{1/2}\times \rho\vert\,\vert^{-1/2}$ et $\pi_{sp}(\psi)$ est quotient de l'induite $\rho\vert\,\vert^{-1/2}\times Sp(2) \times \rho\vert\,\vert^{1/2}$. On note $\pi_{temp}$ et $\pi_{unip}$ les repr\'esentations de $\tilde{G}(4)$ intervenant dans $\sigma$ dont la restriction \`a $GL(4,F)$ sont respectivement la repr\'esentation temp\'er\'ee $St(2)\times St (2)$  et sa duale $Sp(2)\times Sp(2)$. Il suffit maintenant de d\'emontrer que pour tout $g\in GL(4,F)$, 
$$
tr \sigma (g,\theta)= tr \pi_{temp}(g,\theta)+ tr \pi_{unip}(g,\theta).
$$
Pour cela on r\'ealise diff\'eremment la repr\'esentation de $GL(4,F)$, $Sp(2)\times St(2)$ comme sous-quotients de $\sigma$; on remarque que $Sp(2)$ est un sous-module de $\rho\vert\,\vert^{-1/2}\times \rho\vert\,\vert^{1/2}$ et que $St(2)$ en est un quotient; donc $Sp(2)\times St(2)$ est sous-quotient de $$\biggl(\rho\vert\,\vert^{-1/2}\times \rho\vert\,\vert^{1/2}\biggr) \times \biggl(\rho\vert\,\vert^{-1/2}\times \rho\vert\,\vert^{1/2}\biggr). $$
Mais $\theta$ ne laisse pas stable cette r\'ealisation et au contraire l'\'echange avec l'autre r\'ealisation naturelle comme sous-quotient o\`u on consid\`ere $St(2)\times Sp(2)$. Donc la trace de $(g,\theta)$ est nulle sur l'espace sous-jacent \`a la somme de ces 2 r\'ealisations. Ceci prouve notre assertion.

\section{R\'esolution dans le cas \'el\'ementaire\label{lecaselementaire}}
On fixe encore $\psi$ un morphisme $\theta$-invariante, en supposant $\psi$ \'el\'ementaire; par contre,  il n'est pas utile qu'il soit de restriction discr\`ete \`a la diagonale.  Il faut faire pour GL(n,F) ce qui a \'et\'e fait pour les groupes classiques en \cite{elementaire}, c'est-\`a-dire repr\'esenter $\pi(\psi)$ \`a l'aide de la cohomologie d'un complexe exact en tout degr\'e sauf 1. C'est une r\'eminiscence de \cite{aubert}.  On abr\`ege la notation $GL(?,F)$ en $GL(?)$.

\subsection{D\'efinition du complexe sans action de $\theta$\label{definitionducomplexe}}
Soient $\rho$ une repr\'esentation irr\'eductible cuspidale telle que $\theta(\rho)=\rho$ (ou, ce qui est \'equivalent $\rho\simeq \rho^*$) et $d\in {\mathbb N}$, $d>0$. Quel que soit l'entier $m$, on note $\Pi(\rho,d)$ les repr\'esentation irr\'eductibles de $GL(m)$ dont le support cuspidal est form\'e de repr\'esentations de la forme $\rho\vert\,\vert^z$ avec $\vert z\vert \leq (d-1)/2$ et $z-(d-1)/2 \in {\mathbb Z}$.

Soient $n_{1},...,n_{s}$ des entiers $\geq 1$, posons:
$$M=GL(n_{1})\times...\times GL(n_{s}).\eqno(1)$$
Soit $\sigma=\sigma_{1}\otimes...\otimes \sigma_{s}$ une repr\'esentation irr\'eductible de $M$. Notons ${\cal S}_{\rho,\leq d}(\sigma)$ l'ensemble des $t\in \{1,...,s\}$ tels que $\sigma_{t}\in \Pi(\rho,\leq d)$.
Soit $\pi$ une repr\'esentation de $M$. On peut d\'ecomposer $\pi$ en somme directe $\pi^1\oplus \pi^2$ de sorte que, pour tout sous-quotient irr\'eductible $\sigma$ de $\pi^1$, resp. $\pi^2$, on ait $\vert{\cal S}_{\rho,\leq d}(\sigma)\vert\geq s-1$, resp. $\vert{\cal S}_{\rho,\leq d}(\sigma)\vert\leq s-2$. On  note $proj_{\rho,\leq d}$ la projection de $\pi$ sur $\pi^1$ de noyau $\pi^2$. C'est un endomorphisme de $\pi$. On note \'egalement  $proj_{\rho,\leq d}$ le foncteur qui, \`a la repr\'esentation $\pi$ de $M$, associe la repr\'esentation $\pi^1$. Remarquons que cet endomorphisme, resp. ce foncteur, est l'identit\'e si $s=1$.

Soit $n\in {\mathbb N}$.  Posons $\Delta=\{1,...,n-1\}$ et fixons un espace vectoriel complexe $E$ de dimension $n-1$, muni d'une base $(e_{j})_{j\in \Delta}$. Pour tout entier $j\in \{0,...,n-1\}$, notons ${\cal M}_{j}(GL(n))$ l'ensemble des $M\in {\cal M}(GL(n))$ tels que la dimension du centre de $M$ soit \'egale \`a $j+1$, autrement dit tels que $M$ s'\'ecrive sous la forme (1), avec  $s=j+1$ et des entiers $n_{1},...,n_{s}>0$. Pour $M\in {\cal M}_{j}(GL(n))$ \'ecrit sous cette forme, on pose:
$$\Delta^M=\{1,...,n_{1}-1\}\cup \{n_{1}+1,...,n_{1}+n_{2}-1\} \cup...\cup\{n_{1}+...+n_{j}+1,...,n-1\},$$
\noindent et on d\'efinit l'\'el\'ement $e^M=e_{n_{1}}\wedge e_{n_{1}+n_{2}}\wedge...\wedge e_{n_{1}+...+n_{j}}$ de $\bigwedge^{j}E$. Supposons $j\leq n-2$, soient $M\in {\cal M}_{j}(GL(n))$ et $M'\in {\cal M}_{j+1}(GL(n))$, supposons $M'\subset M$. Il y a un unique \'el\'ement $m\in \Delta$ tel que $\Delta^M=\Delta^{M'}\cup\{m\}$. On a une \'egalit\'e $e^{M'}=\pm e^M\wedge e_{m}$ et on note $\xi_{M',M}$ le signe qui y intervient, c'est-\`a-dire tel que $e^{M'}=\xi_{M',M}e^M\wedge e_{m}$.

Soit $\pi$ une repr\'esentation de $GL(n)$. On va d\'efinir un complexe ${\cal E}(\pi)=({\cal E}_{j}(\pi))_{j\in {\mathbb Z}}$. Pour $j<0$ ou $j\geq n$, ${\cal E}_{j}(\pi)=\{0\}$. Pour $j\in\{0,...,n-1\}$,
$${\cal E}_{j}(\pi)=\oplus_{M\in {\cal M}_{j}(GL(n))}Ind_{P}^{GL(n)}(proj_{\rho,\leq d}(\pi_{P})).$$
Pour $j<0$ ou $j\geq n-1$, on note $\epsilon_{j}:{\cal E}_{j}(\pi)\to {\cal E}_{j+1}(\pi)$ l'application nulle (c'est d'ailleurs l'unique application possible puisqu'au moins un des deux espaces est nul). Soit $j\in \{0,...,n-2\}$. On d\'efinit $\epsilon_{j}:{\cal E}_{j}(\pi)\to {\cal E}_{j+1}(\pi)$ comme la somme d'applications:
 $$\epsilon_{j,M',M}:Ind_{P}^{GL(n)}(proj_{\rho,\leq d}(\pi_{P}))\to Ind_{P'}^{GL(n)}(proj_{\rho,\leq d}(\pi_{P'})),$$ 
 \noindent la somme \'etant prise sur les couples $(M',M)\in {\cal M}_{j+1}(GL(n))\times {\cal M}_{j}(GL(n))$. Si  $M'\not\subset M$, $\epsilon_{j,M',M}=0$. Supposons $M'\subset M$. On a une application canonique:
 $$\pi_{P}\to Ind_{M\cap P'}^M(\pi_{P'}).$$
 \noindent Consid\'erons l'application compos\'ee:
 $$proj_{\rho,\leq d}(\pi_{P})\to \pi_{P}\to Ind_{M\cap P'}^M(\pi_{P'})\to Ind_{M\cap P'}^M(proj_{\rho,\leq d}(\pi_{P'})).$$
 \noindent Par fonctorialit\'e, on en d\'eduit une application:
 $$\bar{\epsilon}_{M',M}:Ind_{P}^{GL(n)}(proj_{\rho,\leq d}(\pi_{P}))\to Ind_{P}^{GL(n)}Ind_{M\cap P'}^M(proj_{\rho,\leq d}(\pi_{P'}))=Ind_{P'}^{GL(n)}(proj_{\rho,\leq d}(\pi_{P'})).$$
 \noindent Alors $\epsilon_{j,M',M}=\xi_{M',M}\bar{\epsilon}_{M',M}$.
 
 Pour que ${\cal E}(\pi)$ soit un complexe, on doit v\'erifier:
 
 (2) pour tout $j\in {\mathbb Z}$, on a l'\'egalit\'e $\epsilon_{j+1}\circ\epsilon_{j}=0$.
 
 Pour cela, on a besoin de la propri\'et\'e auxiliaire suivante. Pour deux L\'evi $M,M'\in {\cal M}(GL(n))$ tels que $M'\subset M$, notons $\tau_{M',M}:\pi_{P}\to \pi_{P'}$ la projection canonique. Dans $Hom(\pi_{P},\pi_{P'})$, on a l'\'egalit\'e:
 $$proj_{\rho,\leq d}\circ\tau_{M',M}\circ proj_{\rho,\leq d}=proj_{\rho,\leq d}\circ\tau_{M',M}.\eqno(3)$$
 
 Preuve de (3). D\'ecomposons $\pi_{P}$ en $\pi_{P}^1\oplus \pi_{P}^2$ comme au d\'ebut du paragraphe, avec $\pi_{P}^1=proj_{\rho,\leq d}(\pi_{P})$. Sur $\pi_{P}^1$, $proj_{\rho,\leq d}$ est l'identit\'e et l'\'egalit\'e (2) est claire. Sur $\pi_{P}^2$, $proj_{\rho,\leq d}$ est nul et on doit montrer que $proj_{\rho,\leq d}\circ\tau_{M',M}$ l'est aussi. Ecrivons $M$ sous la forme (1) et $M'=GL(n'_{1})\times...\times GL(n'_{s'})$. Il suffit de prouver l'assertion suivante: soient $\sigma=\sigma_{1}\otimes...\otimes\sigma_{s}$ une repr\'esentation irr\'eductible de $M$ et $\sigma'=\sigma'_{1}\otimes...\otimes\sigma'_{s'}$ une repr\'esentation irr\'eductible de $M'$; supposons que $\sigma'$ intervienne comme sous-quotient de $\sigma_{M\cap P'}$ et que ${\cal S}_{\rho,\leq d}(\sigma)$ ait au plus $s-2$ \'el\'ements; alors ${\cal S}_{\rho,\leq d}(\sigma')$ a au plus $s'-2$ \'el\'ements. D\'emontrons cela. Par hypoth\`ese, on peut fixer deux \'el\'ements distincts $u,v\in \{1,...,s\}$ tels que $\sigma_{u}$ et $\sigma_{v}$ n'appartiennent pas \`a $\Pi(\rho,\leq d)$. L'intersection $GL(n_{u})\cap M'$ est \'egale \`a $GL(n'_{a})\times...\times GL(n'_{b})$ pour un sous-intervalle $\{a,...,b\}$ de $\{1,...,s'\}$. Puisque $\sigma'$ intervient dans $\sigma_{M\cap P'}$, le support cuspidal de $\sigma_{u}$ est r\'eunion, en un sens facile \`a pr\'eciser, des supports cuspidaux des $\sigma'_{k}$ pour $k\in\{a,...,b\}$. Puisque $\sigma_{u}\not\in \Pi(\rho,\leq d)$, il existe $u'\in \{a,...,b\}$ tel que $\sigma'_{u'}\not\in \Pi(\rho,\leq d)$, autrement dit $u'\not\in {\cal S}_{\rho,\leq d}(\sigma')$. En rempla\c{c}ant $u$ par $v$, on obtient un deuxi\`eme indice $v'$ tel que $v'\not\in {\cal S}_{\rho,\leq d}(\sigma')$. Cela prouve (3). 
 
 \
 
 Preuve de (2). On peut supposer $j\in \{0,...,n-3\}$. Alors $\epsilon_{j+1}\circ\epsilon_{j}$ est somme d'applications:
 $$\gamma_{M'',M}:Ind_{P}^{GL(n)}(proj_{\rho,\leq d}(\pi_{P}))\to Ind_{P''}^{GL(n)}(proj_{\rho,\leq d}(\pi_{P''})),$$
 \noindent la somme \'etant prise sur les couples $(M'',M)\in {\cal M}_{j+2}(GL(n))\times {\cal M}_{j}(GL(n))$. Ces applications sont d\'efinies de la fa\c{c}on suivante. Si $M''\not\subset M$, $\gamma_{M,M''}=0$. Supposons $M''\subset M$. Notons $m_{1}, m_{2}$ les deux \'el\'ements de $\Delta$ tels que $\Delta^M=\Delta^{M''}\cup \{m_{1},m_{2}\}$. Notons $M'_{1}$ et $M'_{2}$ les \'el\'ements de ${\cal M}_{j+1}(GL(n))$ tels que $\Delta^M=\Delta^{M'_{k}}\cup \{m_{k}\}$ pour $k=1,2$. Alors:
$$\gamma_{M'',M}=\epsilon_{j+1,M'',M'_{1}}\circ\epsilon_{j,M'_{1},M}+\epsilon_{j+1,M'',M'_{2}}\circ\epsilon_{j,M'_{2},M}.$$
 \noindent On a bien s\^ur l'\'egalit\'e $\xi_{M'',M'_{1}}\xi_{M'_{1},M}+\xi_{M'',M'_{2}}\xi_{M'_{2},M}=0$. Il suffit alors de v\'erifier que l'application:
  $$\bar{\epsilon}_{M'',M'_{k}}\circ\bar{\epsilon}_{M'_{k},M}:Ind_{P}^{GL(n)}(proj_{\rho,\leq d}(\pi_{P}))\to Ind_{P''}^{GL(n)}(proj_{\rho,\leq d}(\pi_{P''}))$$
 \noindent est ind\'ependante de $k\in \{1,2\}$.  Soit $f\in Ind_{P}^{GL(n)}(proj_{\rho,\leq d}(\pi_{P}))$. Alors $\bar{\epsilon}_{M'',M'_{k}}\circ\bar{\epsilon}_{M'_{k},M}(f)$ est la fonction qui, en un point $g\in GL(n)$, vaut:
 $$proj_{\rho,\leq d}\circ\tau_{M'',M'_{k}}\circ proj_{\rho,\leq d}\circ\tau_{M'_{k},M}(f(g)).$$
\noindent D'apr\`es (2), ceci est \'egal \`a:
$$proj_{\rho,\leq d}\circ\tau_{M'',M'_{k}}\circ\tau_{M'_{k},M}(f(g)),$$
\noindent ou encore:
$$proj_{\rho,\leq d}\circ\tau_{M'',M}(f(g)),$$
\noindent qui est bien ind\'ependant de $k$. 
\subsubsection{Un lemme d'irr\'eductibilit\'e\label{lemmedirreductibilite}}
Pour prouver l'exactitude du complexe sauf en 1 degr\'e, nous allons utiliser \`a plusieurs reprise le lemme d'irr\'eductibilit\'e suivant. On fixe $\rho$ une repr\'esentation irr\'eductible $\theta$-invariante et $d$ un entier $\geq 1$ comme dans le paragraphe pr\'ec\'edent; on reprend aussi la notation $\Pi(\rho,\leq d)$ de ce paragraphe.  Soient encore $\pi$ une repr\'esentation irr\'eductible de $GL(n(\pi)$ (ce qui d\'efinit $n(\pi)$) et  $\psi'$ un morphisme \'el\'ementaire tel que $Jord(\psi')$ soit un ensemble de triplets $(\rho_{i},a_{i},b_{i})$ pour $i$ dans un ensemble $I$, tel que pour tout $i\in I$ l'une au moins des hypoth\`eses suivantes est v\'erifi\'ee:

$\rho_{i}\not = \rho$, $\rho_{i}$ est une repr\'esentation cuspidale irr\'eductible autoduale d'un groupe $GL(d_{\rho_{i}})$;

$sup(a_{i},b_{i})-(d-1)/2 \notin {\mathbb Z}_{\leq 0}$.

On pose $n(\psi'):=\sum_{i\in I}a_{i}b_{i}d_{\rho_{i}}$.

\bf Lemme. \sl Pour tout $\pi\in \Pi(\rho,\leq d)$, les induites $\pi(\psi')\times \pi$ et $\pi\times \pi(\psi')$ sont irr\'eductibles et isomorphes.\rm

\

 Rappelons les r\'esultats d'adjonction suivants. Soient  $n,n_{1},n_{2}$ trois entiers $\geq 0$ tels que $n=n_{1}+n_{2}$. Posons $M=GL(n_{1})\times GL(n_{2})$, $M'=GL(n_{2})\times GL(n_{1})$. Ce sont des L\'evi de $GL(n)$. Soient $\sigma_{1}$ et $\sigma_{2}$ des repr\'esentations de $GL(n_{1})$ et $GL(n_{2})$ et $\sigma$ une repr\'esentation de $GL(n)$. Alors:
$$(1) \qquad Hom_{GL(n)}(\sigma,Ind_{P}^{GL(n)}(\sigma_{1}\otimes \sigma_{2}))=Hom_{M}(\sigma_{P},\sigma_{1}\otimes \sigma_{2}),$$
$$(2) \qquad Hom_{GL(n)}(Ind_{P'}^{GL(n)}(\sigma_{2}\otimes \sigma_{1}),\sigma)=Hom_{M}(\sigma_{1}\otimes \sigma_{2},\sigma_{P}).$$La premi\`ere est la r\'eciprocit\'e de Frobenius et la 2e est due \`a Bernstein.

On va utiliser ces propri\'et\'es pour $n_{1}=n(\psi')$, $n_{2}=n(\pi)$, $n=n_{1}+n_{2}$. Pour tout repr\'esentation $\delta$ de $M$, on d\'ecompose $\delta$ en $\delta=\delta^{\pi}\oplus \delta^{\not=\pi}$ o\`u tous les sous-quotients irr\'eductibles $\delta_{1}\otimes \delta_{2}$ de $\delta^{\pi}$, resp. $\delta^{\not=\pi}$, v\'erifient la propri\'et\'e: $\delta_{2}$ a m\^eme support cuspidal que $\pi$, resp. $\delta_{2}$ n'a pas m\^eme support cuspidal que $\pi$. Ceci s'applique en particulier quand $\delta$ est le module de Jacquet $\sigma_{P}$ d'une repr\'esentation $\sigma$ de $GL(n)$.

Posons $\Sigma=Ind_{P}^{GL(n)}(\pi({\psi'})\otimes \pi)$. Soit $\sigma$ un quotient de $\Sigma$. Montrons que:

$$\sigma_{P}^{\pi}\not=\{0\}.\eqno(3)$$
\noindent Notons ${\bf x}=(x_{1},...,x_{e})$ une famille de demi-entiers telle que, en posant $\Pi=\rho\vert.\vert^{x_{1}}\times...\times \rho\vert.\vert^{x_{e}}$, $\pi$ soit un quotient de $\Pi$. Alors $\sigma$ est un quotient de:
$$Ind_{P}^{GL(n)}(\pi({\psi'})\otimes \Pi)=(\times_{i\in I}<\rho_{i},a_{i},b_{i}>)\times \rho\vert.\vert^{x_{1}}\times...\times \rho\vert.\vert^{x_{e}}$$
\noindent (l'ensemble $I$ \'etant muni d'un ordre arbitraire). Cette repr\'esentation est isomorphe \`a:
$$Ind_{P'}^{GL(n)}(\Pi\otimes \pi({\psi'}))=\rho\vert.\vert^{x_{1}}\times...\times \rho\vert.\vert^{x_{e}}\times (\times_{i\in I}<\rho_{i},a_{i},b_{i}>).$$
\noindent En effet, on passe d'une induite \`a l'autre en permutant des termes $<\rho_{i},a_{i},b_{i}>$ et $\rho\vert.\vert^{x_{m}}$ pour $i\in I$ et $m=1,...,e$. D'apr\`es l'hypoth\`ese $\pi\in \Pi(\rho,\leq d)$, on a $x_{m}\equiv \frac{d-1}{2}\,\,mod\, {\mathbb Z}$ et $\vert x_{m}\vert\leq \frac{d-1}{2}$. L'hypoth\`ese sur $\psi'$ et les r\'esultats de Zelevinsky entra\^{\i}nent que les induites $<\rho_{i},a_{i},b_{i}>\times \rho\vert.\vert^{x_{m}}$ et $\rho\vert.\vert^{x_{m}}\times <\rho_{i},a_{i},b_{i}>$ sont irr\'eductibles et isomorphes. D'o\`u l'assertion. Alors $\sigma$ est quotient de $Ind_{P'}^{GL(n)}(\Pi\otimes \pi(\psi'))$. En appliquant (2), on obtient $Hom_{M}(\pi(\psi')\otimes \Pi,\sigma_{P})\not=\{0\}$, ce qui implique (3).

Montrons que l'on a:
$$ \Sigma_{P}^{\pi}=\pi(\psi')\otimes \pi.\eqno(4)$$
\noindent On utilise la filtration de Bernstein-Zelevinsky. Elle est index\'ee par $[W^M\backslash W/W^M]$. Le quotient correspondant \`a $w\in [W^M\backslash W/W^M]$ est:
$$Ind_{M\cap wPw^{-1} }^M(w((\pi(\psi')\otimes \pi)_{M\cap w^{-1} Pw})).$$
\noindent Le terme pour $w=1$ donne $\pi(\psi')\otimes \pi$. Pour un autre $w$, soit $\sigma_{1}\otimes \sigma_{2}$ un sous-quotient irr\'eductible de l'induite ci-dessus. Le support cuspidal de $\sigma_{2}$ contient forc\'ement un terme $\rho_{i}\vert.\vert^{\frac{b_{i}-1}{2}}$ ou $\rho_{i}\vert.\vert^{\frac{1-a_{i}}{2}}$ pour un $i\in I$. D'apr\`es les hypoth\`eses sur $\psi'$ et $\pi$, un tel terme ne peut pas intervenir dans le support cuspidal de $\pi$. D'o\`u (4).

Soit maintenant $\sigma$ un sous-module irr\'eductible de $\Sigma$ et $\bar{\Sigma}$ le module quotient. D'apr\`es (1), $\sigma_{P}^{\pi}$ contient $\pi({\psi'})\otimes \pi$. Par exactitude du foncteur $\delta \mapsto \delta_{P}^{\pi}$ et d'apr\`es (4), on a $\bar{\Sigma}_{P}^{\pi}=\{0\}$. Cela contredit (3) sauf si $\bar{\Sigma}=\{0\}$. Donc cette derni\`ere \'egalit\'e est v\'erifi\'ee et $\Sigma=\sigma$ est irr\'eductible.

Posons $\Sigma'=Ind_{P'}^{GL(n)}(\pi\otimes \pi(\psi'))$. On montre de m\^eme que $\Sigma'$ est irr\'eductible. D'apr\`es (2) et (4), $Hom_{GL(n)}(\Sigma',\Sigma)\not=\{0\}$. Par irr\'eductibilit\'e, cela entra\^{\i}ne $\Sigma=\Sigma'$.

\subsubsection{Exactitude du complexe \label{exactitudeducomplexe}} 
Soient $\rho$, $d$ et $n$ comme dans le paragraphe pr\'ec\'edent, $\psi$ un morphisme \'el\'ementaire; on pose $Jord(\psi):=\{(\rho_{i},a_{i},b_{i})_{i\in I}\}$ pour un ensemble convenable d'indices $I$.  On d\'efinit  les ensembles $J$, $J_{\leq d}$ par $J:=\{i\in I; \rho_{i}=\rho$ et $sup(a_{i},b_{i})\equiv d mod\, {\mathbb Z}\}$ et  $J_{\leq d}:=\{i\in J; sup(a_{i},b_{i})\leq d$. Posons:
$$j_{0}(\psi)=\sum_{i\in J_{\leq d}}sup(a_{i},b_{i}),$$
$$j(\psi)=\left\lbrace\begin{array}{cc}j_{0}(\psi),&\,\,si\,\,I\not=J_{\leq d}\\ j_{0}(\psi)-1,&\,\,si\,\,I=J_{\leq d}.\\ \end{array}\right.$$
Notons $\psi^{\sharp}$ l'\'el\'ement de $\Psi_{elem,n}$ form\'e des $(\rho_{i},a_{i},b_{i})$ pour $i\not\in J_{\leq d}$ et des $(\rho_{i},b_{i},a_{i})$ pour $i\in J_{\leq d}$.

{\it {\bf Proposition.} Le complexe ${\cal E}(\pi(\psi))$ est exact en degr\'e $<j(\psi)$. Il est nul en degr\'e $>j(\psi)$. Le conoyau de $\epsilon_{j(\psi)-1}$ est isomorphe \`a $\pi(\psi^{\sharp})$.}

Preuve. Si $I=J_{\leq d}$,  la repr\'esentation $\pi(\psi)$ appartient \`a $\Pi(\rho,\leq d)$ et tous les foncteurs $proj_{\rho,\leq d}$ intervenant dans la d\'efinition de ${\cal E}(\pi(\psi))$ agissent par l'identit\'e. Le complexe ${\cal E}(\pi(\psi))$ est celui qui d\'efinit l'involution de Zelevinsky et l'assertion est connue (\cite{aubert}, \cite{ss}). On suppose d\'esormais $I\not=J_{\leq d}$.Notons $\psi'$ la famille form\'ee des $(\rho_{i},a_{i},b_{i})$ pour $i\not\in J_{\leq d}$. Elle n'est pas vide. Notons ${\cal R}$ l'ensemble des repr\'esentations de la forme $\pi(\psi')\times \pi$, o\`u $\pi\in \Pi(\rho,\leq d)$ (il s'agit de repr\'esentations de groupes $GL(m)$, o\`u $m$ est quelconque). Soient $j\in \{0,...,n-1\}$, $M\in {\cal M}_{j}(GL(n))$ et $\sigma$ un sous-quotient irr\'eductible de $proj_{\rho,\leq d}(\pi(\psi)_{P})$. Ecrivons:
$$M=GL(n_{1})\times...\times GL(n_{j+1}),\,\,\sigma=\sigma_{1}\otimes..\otimes\sigma_{j+1}.$$
\noindent Montrons qu'il existe $t\in \{1,...,j+1\}$ tel que:

$$(1) \left\lbrace\begin{array}{cc}\sigma_{t}\in {\cal R};&\\ \,\,pour\,\,r\in\{1,...,j+1\},&\,\,r\not=t,\,\,\sigma_{r}\in \Pi(\rho,\leq d). \\ \end{array}\right.$$

Par d\'efinition de $proj_{\rho,\leq d}(\pi(\psi))$, on peut en tout cas trouver $t$ tel que la deuxi\`eme condition soit satisfaite. Fixons un tel $t$. Rappelons que, par d\'efinition, $\pi(\psi)$ est une induite; on \'ecrit $\sigma(\psi)$ la repr\'esentation induisante c'est-\`a-dire $\times_{i\in I}S(\rho_{i},a_{i},b_{i})$ o\`u pour $i\in I$, $S(\rho_{i},a_{i},b_{i})$ est $St(\rho_{i},a_{i})$ si $b_{i}=1$ et  est  $Sp(\rho_{i},b_{i})$ si $a_{i}=1$. On note $P_{\psi}$ le parbolique (standard) qui permet l'induction.

Filtrons $\pi(\psi)_{P}$ par la filtration de Bernstein-Zelevinsky.  Cette filtration est index\'ee par $[W^M\backslash W/W^{M_{\psi}}]$ et le terme du gradu\'e associ\'e \`a un \'el\'ement $w$ de cet ensemble est:
$$\delta_{w}=Ind_{M\cap wP_{\psi}w^{-1} }^M(w(\sigma_{\psi,M_{\psi}\cap w^{-1} Pw})).$$
   Soit $w\in [W^M\backslash W/W^{M_{\psi}}]$ tel que $\sigma$ intervienne dans $\delta_{w}$. On a:

(2) pour $i\in I\setminus J_{\leq d}$, $w$ envoie le bloc de $M_{\psi}$ correspondant \`a $(\rho_{i},a_{i},b_{i})$ dans le $t$-i\`eme bloc $GL(n_{t})$ de $M$. 

En effet, notons $M_{\psi,i}$ le bloc de $M_{\psi}$ en question.  Consid\'erons l'ensemble des \'el\'ements $r\in\{1,...,j+1\}$ tels que $wM_{\psi,i}w^{-1} \cap GL(n_{r})\not=\{1\}$. Notons $r_{1}<...<r_{s}$ ses \'el\'ements. Si la conclusion de (2) \'etait fausse, on aurait $r_{1}\not=t$ ou $r_{s}\not=t$. Alors $\sigma_{r_{1}}$, resp. $\sigma_{r_{s}}$, contiendrait dans son support cuspidal une repr\'esentation $\rho_{i}\vert.\vert^{\pm \frac{sup(a_{i},b_{i})-1}{2}}$, ce qui contredirait la propri\'et\'e $\sigma_{r_{1}}\in \Pi(\rho,\leq d)$, resp. $\sigma_{r_{s}}\in \Pi(\rho,\leq d)$.

Gr\^ace \`a (2), $\sigma_{t}$ est sous-quotient d'une induite de la forme $\pi_{1}\times...\times \pi_{u}$ o\`u la famille $(\pi_{1},...,\pi_{u})$ est r\'eunion des $<\rho_{i},a_{i},b_{i}>$ pour $i\in I\setminus J_{\leq d}$ et d'une famille d'\'el\'ements de $\Pi(\rho,\leq d)$. Le lemme \ref{lemmedirreductibilite} entra\^{\i}ne que $\sigma_{t}$ appartient \`a ${\cal R}$, ce qui prouve (1).

De (1) r\'esulte d'abord que, si $proj_{\rho,\leq d}(\pi(\psi)_{P})\not=0$, il existe $t\in \{1,...,j+1\}$ et des entiers $m_{r}\geq 0$ pour $r\in \{1,...,j+1\}$ de sorte que:
$$n_{t}=n(\psi')+m_{t}d_{\rho},\,\,n_{r}=m_{r}d_{\rho}\,\,pour\,\,r\not=t.$$
\noindent Cela entra\^{\i}ne $j\leq j(\psi)$. Donc ${\cal E}_{j}(\pi(\psi))$ est nul pour $j>j(\psi)$.

Introduisons le L\'evi:
$$\bar{M}=GL(n(\psi'))\times GL({d_{\rho}})\times...\times GL(d_{\rho})$$
\noindent de $GL(n)$. Soit ${\cal X}$ l'ensemble des familles ${\bf x}=(x_{k})_{k=1,...,j(\psi)}$ form\'ees de demi-entiers tels que $x_{k}\equiv \frac{d-1}{2}\,\,mod\,\,{\mathbb Z}$ et $\vert x_{k}\vert\leq \frac{d-1}{2}$. Pour toute telle famille ${\bf x}$, on d\'efinit la repr\'esentation de $\bar{M}$:
$$\bar{\sigma}_{{\bf x}}=\pi(\psi')\otimes \rho\vert.\vert^{x_{1}}\otimes...\otimes \rho\vert.\vert^{x_{j(\psi)}}.$$
\noindent Pour tout couple de repr\'esentations $\pi_{1},\sigma_{1}$ d'un L\'evi $M_{1}$, avec $\sigma_{1}$ irr\'eductible, notons $\pi_{1}[\sigma_{1}]$ le facteur direct de $\pi_{1}$ dont tous les sous-quotients ont m\^eme support cuspidal que $\sigma_{1}$. Soit $\pi$ une repr\'esentation de $GL(n)$ appartenant \`a ${\cal R}$. Il r\'esulte du lemme \ref{lemmedirreductibilite} et de sa preuve que:

- pour tout ${\bf x}\in {\cal X}$, $\pi_{\bar{P}}[\bar{\sigma}_{{\bf x}}]$ est isotypique, de type $\bar{\sigma}_{{\bf x}}$;

- il existe ${\bf x}\in {\cal X}$ tel que $\pi_{\bar{P}}[\bar{\sigma}_{{\bf x}}]\not=0$.

Soit $j\in \{1,...,j(\psi)\}$. D'apr\`es (1) et le lemme \ref{lemmedirreductibilite}, tout sous-quotient irr\'eductible de ${\cal E}_{j}(\pi(\psi))$ appartient \`a ${\cal R}$. D'apr\`es les propri\'et\'es ci-dessus et d'apr\`es l'exactitude du foncteur $\pi\mapsto \pi_{\bar{P}}[\bar{\sigma}_{{\bf x}}]$, l'exactitude du complexe ${\cal E}(\pi(\psi))$ en degr\'es $\not= j(\psi)$ r\'esulte de l'assertion suivante:

(3) pour tout ${\bf x}\in {\cal X}$, la suite:
$$0 \to {\cal E}_{0}(\pi(\psi))_{\bar{P}}[\bar{\sigma}_{{\bf x}}]\to {\cal E}_{1}(\pi({\psi}))_{\bar{P}}[\bar{\sigma}_{{\bf x}}]\to ...\to {\cal E}_{j(\psi)}(\pi(\psi))_{\bar{P}}[\bar{\sigma}_{{\bf x}}]$$
\noindent est exacte.

On reprend maintenant la preuve de \cite{aubert}. Munissons $W$ de l'ordre de Bruhat, \'ecrivons $W=\{w_{1},...,w_{N}\}$ de sorte que, si $w_{m}<w_{m'}$, alors $m<m'$. Pour $m\in \{1,...,N+1\}$ et $M\in {\cal M}(GL(n))$, notons $\Omega_{M,m}$ le plus grand sous-ensemble de $\{w_{m},...,w_{N}\}$ qui soit invariant par multiplication \`a droite par $W^M$ et \`a gauche par $W^{\bar{M}}$. Cet ensemble v\'erifie les propri\'et\'es suivantes:

(4) soient $m,h,h'\in \{1,...,N\}$, supposons $w_{h}\in \Omega_{M,m}$ et $w_{h}<w_{h'}$. Alors $w_{h'}\in \Omega_{M,m}$.

En effet, soient $\ell,\ell'$ les entiers tels que $w_{\ell},w_{\ell'}$ soit des \'el\'ements de $ [W^{\bar{M}}\backslash W/W^M]$  de longueur minimale dans leur double classe  et $w_{h}\in W^{\bar{M}}w_{\ell}W^M$, $w_{h'}\in W^{\bar{M}}w_{\ell'}W^M$. On montre ais\'ement que l'in\'egalit\'e $w_{h}<w_{h'}$ entra\^{\i}ne $w_{\ell}\leq w_{\ell'}$. Donc $\ell\leq \ell'$. Puisque $w_{h}\in \Omega_{M,m}$ et par d\'efinition de cet ensemble, on a aussi $w_{\ell}\in \Omega_{M,m}$, donc $\ell\geq m$, puis $\ell'\geq m$. Soit $w_{k}\in W^{\bar{M}}w_{\ell'}W^M$. Alors $w_{k}\geq w_{\ell'}$ donc $k\geq \ell'\geq m$. Cela montre que $W^{\bar{M}}w_{\ell'}W^M\subseteq \{w_{m},...,w_{N}\}$. Par maximalit\'e de $\Omega_{M,m}$, on a aussi $W^{\bar{M}}w_{\ell'}W^M\subseteq \Omega_{M,m}$. A fortiori $w_{h'}\in \Omega_{M,m}$. Cela prouve (4);

(5) pour $m\leq N$, on a:
$$\Omega_{M,m}=\left\lbrace\begin{array}{cc}\Omega_{M,m+1},&\,\,si\,\,w_{m}\not\in [W^{\bar{M}}\backslash W/W^M],\\ (W^{\bar{M}}w_{m}W^M)\cup \Omega_{M,m+1},&\,\,si\,\,w_{m}\in [W^{\bar{M}}\backslash W/W^M].\\ \end{array}\right.$$

Il est clair que $\Omega_{M,m}\not=\Omega_{M,m+1}$ si et seulement si $w_{m}\in \Omega_{M,m}$. Si cette condition est v\'erifi\'ee, on a $W^{\bar{M}}w_{m}W^M \subseteq \Omega_{M,m}\subseteq \{w_{m},...,w_{N}\}$. Donc $w_{m}$ est minimal pour l'ordre de Bruhat dans sa double classe $W^{\bar{M}}w_{m}W^M $ et cela entra\^{\i}ne $w_{m}\in [W^{\bar{M}}\backslash W/W^M]$. Inversement, si $w_{m}\in [W^{\bar{M}}\backslash W/W^M]$, on a pour la m\^eme raison $W^{\bar{M}}w_{m}W^M
\subseteq \{w_{m},...,w_{N}\}$, donc $W^{\bar{M}}w_{m}W^M
\subseteq \Omega_{M,m}$ par maximalit\'e de cet ensemble. Il est clair que $\Omega_{M,m+1}\subseteq \Omega_{M,m}\setminus W^{\bar{M}}w_{m}W^M$. Mais ce dernier ensemble est bien inclus dans $\{w_{m+1},...,w_{N}\}$ et v\'erifie les conditions d'invariance requises. La derni\`ere inclusion est donc une \'egalit\'e;

(6) si $M'\in {\cal M}(GL(n))$ est inclus dans $M$, alors $\Omega_{M,m}\subseteq \Omega_{M',m}$.

C'est \'evident.

Pour $M\in {\cal M}(GL(n))$, on d\'efinit une filtration sur $Ind_{P}^{GL(n)}(proj_{\rho,\leq d}(\pi(\psi)_{P}))$ index\'ee par $\{1,...,N+1\}$. Le $m$-i\`eme terme est form\'e des fonctions \`a support dans $\bigcup_{w\in \Omega_{M,m}}Bw^{-1} B$, ou encore $\bigcup_{w\in \Omega_{M,m}}Pw^{-1} \bar{P}$, ces deux ensembles \'etant \'egaux. Remarquons que, d'apr\`es (4), cet ensemble est ouvert dans $GL(n)$. Chaque terme de la filtration est stable par l'action de $\bar{P}$. De ces filtrations r\'esulte une filtration $({\cal E}_{j,m}(\pi(\psi)))_{m=1,...,N+1}$ de chaque espace ${\cal E}_{j}(\pi(\psi))$. Cette filtration est compatible avec les diff\'erentielles $\epsilon_{j}$ d'apr\`es (6). Par exactitude du foncteur $\pi\mapsto \pi_{\bar{P}}[\bar{\sigma}_{{\bf x}}]$ (qui est d\'efini pour toute repr\'esentation de $\bar{P}$), on a aussi des filtrations $({\cal E}_{j,m}(\pi(\psi))_{\bar{P}}[\bar{\sigma}_{{\bf x}}])_{m=1,...,N+1}$ sur les termes ${\cal E}_{j}(\pi(\psi))_{\bar{P}}[\bar{\sigma}_{{\bf x}}]$, compatibles aux diff\'erentielles. Pour d\'emontrer que la suite (3) est exacte, il suffit de d\'emontrer la m\^eme assertion pour la suite des $m$-i\`emes gradu\'es, pour tout $m=1,...,N$. Cette suite s'\'ecrit:

$$(7) \qquad ...\to\oplus_{M\in {\cal M}_{j-1}(GL(n))}\delta_{M,m}[\bar{\sigma}_{{\bf x}}]\to \oplus_{M\in {\cal M}_{j}(GL(n))}\delta_{M,m}[\bar{\sigma}_{{\bf x}}]\to...\to\oplus_{M\in {\cal M}_{j(\psi)}(GL(n))}\delta_{M,m}[\bar{\sigma}_{{\bf x}}]$$
\noindent o\`u $\delta_{M,m}$ est le $m$-i\`eme quotient de la filtration de $(Ind_{P}^{GL(n)}(proj_{\rho,\leq d}(\pi_{P})))_{\bar{P}}$.

 On va d\'ecrire cette suite. Soient $m\in \{1,...,N\}$ et $M\in {\cal M}(GL(n))$. Le $m$-i\`eme gradu\'e de la filtration de $Ind_{P}^{GL(n)}(proj_{\rho,\leq d}(\pi(\psi)_{P}))$ est nul, d'apr\`es (5), si $w_{m}\not\in [W^{\bar{M}}\backslash W/W^M]$. Si $w_{m}\not\in [W^{\bar{M}}\backslash W]$, cette condition est v\'erifi\'ee pour tout $M$ et la suite (7) est enti\`erement nulle donc exacte. On suppose d\'esormais $w_{m}\in [W^{\bar{M}}\backslash W]$. La condition $w_{m}\in [W^{\bar{M}}\backslash W/W^M]$ est alors \'equivalente \`a $w_{m}\in [W/W^{M}]$. Supposons-la v\'erifi\'ee. D'apr\`es (5), on a alors:
 $$\delta_{M,m}=Ind_{\bar{M}\cap w_{m}Pw_{m}^{-1} }^{\bar{M}}(w_{m}((proj_{\rho,\leq d}(\pi(\psi)_{P}))_{M\cap w_{m}^{-1} \bar{P}w_{m}})).$$
 \noindent Introduisons quelques notations. Si $w_{m}^{-1} \bar{M}w_{m}$ est un L\'evi standard, on le note $\bar{M}_{m}$ et, comme toujours, on note $\bar{P}_{m}$ le sous-groupe parabolique standard qui lui est associ\'e (qui n'a aucune raison d'\^etre \'egal \`a $w_{m}^{-1} \bar{P}w_{m}$). On pose:
 $$\bar{\delta}_{m}=w_{m}(\pi(\psi)_{\bar{P}_{m}}[w_{m}^{-1} (\bar{\sigma}_{{\bf x}})]).$$
 \noindent Montrons que:
 
 $$(8)\qquad \left\lbrace\begin{array}{cc}\,\,si\,\,w_{m}^{-1} \bar{M}w_{m}\not\subseteq M,&\,\,\delta_{M,m}[\bar{\sigma}_{{\bf x}}]=0,\\ \,\,si\,\,w_{m}^{-1} \bar{M}w_{m}\subseteq M,&\,\,\delta_{M,m}[\bar{\sigma}_{{\bf x}}]= \bar{\delta}_{m}.\\ \end{array}\right.$$
 Remarquons que, puisque $w_{m}\in [W^{\bar{M}}\backslash W/W^M]$, l'hypoth\`ese $w_{m}^{-1} \bar{M}w_{m}\subseteq M$ entra\^{\i}ne que $w_{m}^{-1} \bar{M}w_{m}$ est un groupe de L\'evi standard de $GL(n)$, ce qui donne un sens \`a $\bar{\delta}_{m}$.
 
  Ecrivons:
 $$M=GL(n_{1})\times...\times GL(n_{r}),\,\,\bar{M}=G_{\bar{n}_{1}}\times...\times G_{\bar{n}_{s}}.$$
 \noindent Pour tous $h\in \{1,...,r\}$, $k\in \{1,...,s\}$, d\'ecomposons le sous-groupe $GL(n_h)\cap w_{m}^{-1} \bar{M}w_{m}$ de $M\cap w_{m}^{-1} \bar{M}w_{m}$ en $GL(n_{h,1})\times ...\times GL({n_{h,s}})$ et le sous-groupe $G_{\bar{n}_{k}}\cap w_{m}Mw_{m}^{-1} $ de $\bar{M}\cap w_{m}Mw_{m}^{-1} $ en $G_{\bar{n}_{k,1}}\times...\times G_{\bar{n}_{k,r}}$. Les $n_{h,k}$ et $\bar{n}_{k,h}$ sont des entiers $\geq 0$ tels que $n_{h,k}=\bar{n}_{k,h}$. L'action de $w_{m}$ envoie le bloc $GL({n_{h,k}})$ de  $M\cap w_{m}^{-1} \bar{M}w_{m}$ sur le bloc $G_{\bar{n}_{k,h}}$ de $\bar{M}\cap w_{m}Mw_{m}^{-1} $. Supposons $\delta_{M,m}[\bar{\sigma}_{{\bf x}}]\not= 0$. On peut fixer un sous-quotient irr\'eductible $\sigma=\sigma_{1}\otimes...\otimes \sigma_{r}$ de $proj_{\rho,\leq d}(\pi(\psi)_{P})$ et une repr\'esentation irr\'eductible:
 $$\mu=\mu_{1,1}\otimes...\otimes\mu_{1,s}\otimes\mu_{2,1}\otimes...\otimes \mu_{2,s}\otimes...\otimes \mu_{r,1}\otimes...\otimes \mu_{r,s}$$
 \noindent de $M\cap w_{m}^{-1} {\bar{M}}w_{m}$ de sorte que les propri\'et\'es suivantes soient v\'erifi\'ees:
 
 - $\mu$ intervient dans $\sigma_{M\cap w_{m}^{-1} \bar{M}w_{m}}$;
 
 - une repr\'esentation de m\^eme support cuspidal que $\bar{\sigma}_{\bf x}$ intervient dans $Ind_{\bar{M}\cap w_{m}Pw_{m}^{-1} }^{\bar{M}}(w_{m}(\mu))$.
 
 Remarquons que la repr\'esentation $w_{m}(\mu)$ de $\bar{M}\cap w_{m}Mw_{m}^{-1} $ s'\'ecrit:
 $$w_{m}(\mu)=\mu_{1,1}\otimes...\otimes\mu_{r,1}\otimes\mu_{1,2}\otimes...\otimes \mu_{r,2}\otimes...\otimes \mu_{1,s}\otimes...\otimes \mu_{r,s}.$$
 A $\sigma$ associons l'entier $t$ v\'erifiant (1). On a donc $\sigma_{t}=\pi(\psi')\times \pi$ pour une repr\'esentation $\pi\in \Pi(\rho,\leq d)$. Le m\^eme raisonnement qui nous a permis de prouver (1) montre que l'une des deux propri\'et\'es suivantes est v\'erifi\'ee:
 
 - il existe $u\in \{1,...,s\}$ tel que le support cuspidal de $\mu_{t,u}$ contient enti\`erement celui de $\pi(\psi')$;
 
 - il existe deux \'el\'ements distincts $u,v\in \{1,...,s\}$ tels que $\mu_{t,u}$ et $\mu_{t,v}$ n'appartiennent pas \`a $\Pi(\rho,\leq d)$.
 
Rappelons que:
$$\bar{\sigma}_{{\bf x}}=\pi(\psi')\otimes \rho\vert.\vert^{x_{1}}\otimes...\otimes\rho\vert.\vert^{x_{j(\psi)}}.$$
\noindent La seconde propri\'et\'e ci-dessus entra\^{\i}ne que deux facteurs de cette expression n'appartiennent pas \`a $\Pi(\rho,\leq d)$. Or seul le premier facteur n'appartient pas \`a cet ensemble. La seconde propri\'et\'e ci-dessus est donc exclue. Consid\'erons la premi\`ere. Pour la m\^eme raison, on a n\'ecessairement $u=1$ et il y a en fait \'egalit\'e des supports cuspidaux de $\mu_{t,1}$ et $\pi(\psi')$. Donc $n(\psi')=n_{t,1}=\bar{n}_{1,t}$. Ces \'egalit\'es signifient que $w_{m}^{-1} $ envoie le premier bloc de $\bar{M}$ dans le $t$-i\`eme bloc de $M$. Pour le $k$-i\`eme bloc de $\bar{M}$, avec $k\not=1$, puisque le facteur correspondant de $\bar{\sigma}_{{\bf x}}$ est cuspidal, il est clair que la suite $\bar{n}_{k,1},...,\bar{n}_{k,r}$ est form\'ee de $r-1$ termes nuls et d'un terme \'egal \`a $d_{\rho}$. Cela entra\^{\i}ne encore que $w_{m}^{-1} $ envoie le $k$-i\`eme bloc de $\bar{M}$ dans un bloc de $M$. Cela prouve que $w_{m}^{-1} \bar{M}w_{m}\subseteq M$ et  la premi\`ere assertion de (8) s'ensuit.

Supposons maintenant $w_{m}^{-1} \bar{M}w_{m}\subseteq M$. Alors:
 $$\delta_{M,m}=w_{m}((proj_{\rho,\leq d}(\pi(\psi)_{P})_{M\cap w_{m}^{-1} \bar{P}w_{m} })$$ 
 \noindent et:
 $$\delta_{M,m}[\bar{\sigma}_{{\bf x}}]=w_{m}((proj_{\rho,\leq d}(\pi(\psi)_{P})_{M\cap w_{m}^{-1}\bar{P}w_{m} }[w_{m}^{-1} (\bar{\sigma}_{{\bf x}})]).$$
 Avec la notation introduite dans la preuve de \ref{definitionducomplexe}(3), on a:
 $$(proj_{\rho,\leq d}(\pi(\psi)_{P}))_{M\cap w_{m}^{-1} \bar{P}w_{m}}[w_{m}^{-1} (\bar{\sigma}_{{\bf x}})]=(\tau_{\bar{M}_{m},M}\circ proj_{\rho,\leq d}\circ\tau_{M,GL(n)}(\pi(\psi)))[w_{m}^{-1} (\bar{\sigma}_{{\bf x}})].$$
 \noindent Par construction, on a $proj_{\rho,\leq d}(w_{m}^{-1} (\bar{\sigma}_{{\bf x}}))=w_{m}^{-1} (\bar{\sigma}_{{\bf x}})$ et l'expression ci-dessus est \'egale \`a:
$$ (proj_{\rho,\leq d}\circ\tau_{\bar{M}_{m},M}\circ proj_{\rho,\leq d}\circ\tau_{M,GL(n)}(\pi(\psi)))[w_{m}^{-1} (\bar{\sigma}_{{\bf x}})].$$
\noindent D'apr\`es \ref{definitionducomplexe}(3), ceci est \'egal \`a:
$$ (proj_{\rho,\leq d}\circ\tau_{\bar{M}_{m},M}\circ \tau_{M,GL(n)}(\pi(\psi)))[w_{m}^{-1} (\bar{\sigma}_{{\bf x}})],$$
\noindent ou encore \`a:
$$ (proj_{\rho,\leq d}\circ\tau_{\bar{M}_{m},GL(n)}(\pi(\psi)))[w_{m}^{-1} (\bar{\sigma}_{{\bf x}})].$$
\noindent Toujours parce que  $proj_{\rho,\leq d}(w_{m}^{-1} (\bar{\sigma}_{{\bf x}}))=w_{m}^{-1} (\bar{\sigma}_{{\bf x}})$, on peut supprimer le $proj_{\rho,\leq d}$ de cette expression et il ne reste plus que $\tau_{\bar{M}_{m},GL(n)}(\pi(\psi)))[w_{m}^{-1} (\bar{\sigma}_{{\bf x}})]$. Alors $\delta_{M,m}[\bar{\sigma}_{{\bf x}}]=w_{m}(\tau_{\bar{M}_{m},GL(n)}(\pi(\psi)))[w_{m}^{-1} (\bar{\sigma}_{{\bf x}})])$, ce qui d\'emontre la seconde assertion de (8).

Supposons que $w_{m}^{-1} \bar{M}w_{m}$ ne soit pas un L\'evi standard. D'apr\`es la remarque qui suit (8), on a alors $\delta_{M,m}[\bar{\sigma}_{{\bf x}}]=0$ pour tout $M$. De nouveau, la suite (7) est enti\`erement nulle, donc exacte. Supposons maintenant que $w_{m}^{-1} \bar{M}w_{m}$ soit un L\'evi standard, notons $\Delta_{m}\subseteq \Delta$ son ensemble associ\'e. Identifions $w_{m}$ \`a une permutation de $\{1,...,n\}$, notons $\Delta'_{m}\subseteq \Delta$ l'ensemble des $k\in \Delta$ tels que $w_{m}(k+1)>w_{m}(k)$. Puisque $w_{m}\in[W^{\bar{M}}\backslash W]$, on a $\Delta_{m}\subseteq \Delta'_{m}$. La condition $w_{m}^{-1} \bar{M}w_{m}\subseteq M$ \'equivaut \`a $\Delta_{m}\subseteq \Delta^M$. La condition $w_{m}\in [W/W^M]$ \'equivaut \`a $\Delta^M
\subseteq \Delta'_{m}$. Le r\'esultat de nos calculs s'exprime par l'\'egalit\'e suivante:

$$\delta_{M,m}[\bar{\sigma}_{{\bf x}}]=\left\lbrace\begin{array}{cc}\bar{\delta}_{m},&\,\,si\,\, \Delta_{m}\subseteq \Delta^M
\subseteq \Delta'_{m},\\ 0,&\,\,sinon.\\ \end{array}\right.$$
La suite (7) a donc pour $j$-i\`eme terme:
$$(9) \qquad \oplus_{M\in {\cal M}_{j}(GL(n)), \Delta_{m}\subseteq \Delta^M
\subseteq \Delta'_{m}} \bar{\delta}_{m}.$$
\noindent Les diff\'erentielles de ce complexe sont \'evidentes: ce sont les sommes directes d'identit\'es multipli\'ees par les signes $\xi_{M',M}$. Il est bien connu que ce complexe est exact, sauf dans le cas o\`u $\Delta_{m}=\Delta'_{m}$, auquel cas il est exact sauf en degr\'e $j=\vert\Delta\vert-\vert\Delta_{m}\vert$. Mais $\vert\Delta\vert=n(\psi')+j(\psi)d_{\rho}-1$ et $\vert\Delta_{m}\vert=\vert\Delta^{\bar{M}}\vert=n(\psi')-1+j(\psi)(d_{\rho}-1)$. Le degr\'e $j$ ci-dessus est donc \'egal \`a $j(\psi)$. Cela prouve que le complexe ${\cal E}(\pi(\psi))$ est exact en degr\'e $\not= j(\psi)$.

Notons $\pi^{\sharp}$ le conoyau de $\epsilon_{j(\psi)-1}$. C'est une repr\'esentation dont on ne sait pas encore qu'elle est irr\'eductible. Toutefois, gr\^ace \`a (1) et au lemme \ref{lemmedirreductibilite}, tous ses sous-quotients irr\'eductibles appartiennent \`a ${\cal R}$. On peut donc l'\'ecrire $\pi^{\sharp}=\pi(\psi')\times \pi''$, o\`u $\pi''$ est une repr\'esentation de $GL({n-n(\psi')})$ dont tous les sous-quotients appartiennent \`a $\Pi(\rho,\leq d)$. Pour ${\bf x}=(x_{k})_{k=1,...,j(\psi)}\in {\cal X}$, on a introduit la repr\'esentation $\bar{\sigma}_{{\bf x}}$ de $\bar{M}$. Introduisons maintenant la repr\'esentation:
$$\sigma''_{{\bf x}}=\rho\vert.\vert^{x_{1}}\otimes...\otimes\rho\vert.\vert^{x_{j(\psi)}}$$
\noindent du L\'evi $ M''=GL({d_{\rho}})\times...\times GL({d_{\rho}})$ de $GL({n-n(\psi')})$. On a $\bar{\sigma}_{{\bf x}}=\pi(\psi')\otimes \sigma''_{{\bf x}}$. Une preuve analogue \`a celle du lemme \ref{lemmedirreductibilite} prouve que:
$$(10) \qquad \pi^{\sharp}_{\bar{P}}[\bar{\sigma}_{{\bf x}}]=\pi(\psi')\otimes \pi''_{P''}[\sigma''_{{\bf x}}].$$
On calcule $\pi^{\sharp}_{\bar{P}}[\bar{\sigma}_{{\bf x}}]$ en reprenant les calculs ci-dessus: il suffit de remplacer le $(j(\psi)+1)$-i\`eme terme du complexe ${\cal E}(\pi(\psi))$ par $\pi^{\sharp}$. On obtient que $\pi^{\sharp}_{\bar{P}}[\bar{\sigma}_{{\bf x}}]$ est filtr\'e par une filtration index\'ee par $\{1,...,N+1\}$. Le $m$-i\`eme terme du gradu\'e associ\'e est nul sauf si $w_{m}\in [W^{\bar{M}}\backslash W]$ et $w_{m}^{-1} \bar{M}w_{m}$ est un L\'evi standard. Si ces conditions sont v\'erifi\'ees, ce $m$-i\`eme terme du gradu\'e est le dernier conoyau de la suite dont les termes index\'es par $j\leq j(\psi)$ sont les termes (9). Comme on l'a vu cette suite est exacte et le conoyau est nul, sauf si $\Delta_{m}=\Delta'_{m}$. Dans ce dernier cas, la suite en question est r\'eduite \`a un seul terme $\bar{\delta}_{m}$ plac\'e en degr\'e $j(\psi)$ et le conoyau est \'evidemment \'egal \`a $\bar{\delta}_{m}$. Or on v\'erifie qu'il n'y a qu'un \'el\'ement $m$ v\'erifiant toutes ces conditions, i.e. $w_{m}\in [W^{\bar{M}}\backslash W]$, $w_{m}^{-1} \bar{M}w_{m}$ est un L\'evi standard, $\Delta_{m}=\Delta'_{m}$. Notons $m_{0}$ cet unique $m$.  On a $w_{m_{0}}=w_{max}^{\bar{M}}w_{max}$, o\`u $w_{max}$ est l'\'el\'ement de longueur maximale de $W$ et $w_{max}^{\bar{M}}$ l'\'el\'ement de longueur maximale de $W^{\bar{M}}$. Finalement, on obtient l'\'egalit\'e:
$$(11) \qquad \pi^{\sharp}_{\bar{P}}[\bar{\sigma}_{{\bf x}}]=\bar{\delta}_{m_{0}}.$$
\noindent Notons $\psi''$ la famille form\'ee des $(\rho_{i},a_{i},b_{i})_{i\in J_{\leq d}}$. Posons $\iota({\bf x})=(x_{j(\psi)+1-k})_{k=1,...,j(\psi)}$. On a $\pi(\psi)=\pi({\psi'})\times \pi({\psi''})$ et on prouve une \'egalit\'e analogue \`a (10):
$$\pi(\psi)_{\bar{P}_{m_{0}}}[w_{m_{0}}^{-1} (\bar{\sigma}_{{\bf x}})]=\pi(\psi'')_{P''}[\sigma''_{\iota({\bf x})}]\otimes \pi(\psi').$$ 
D'o\`u:
$$\bar{\delta}_{m_{0}}=w_{m_{0}}(\pi(\psi'')_{P''}[\sigma''_{\iota({\bf x})}]\otimes \pi(\psi')).$$
En utilisant (10) et (11), on en d\'eduit que la multiplicit\'e de $\sigma''_{{\bf x}}$ dans $\pi''_{P''}$ est \'egale \`a celle de $\sigma''_{\iota({\bf x})}$ dans $\pi(\psi'')_{P''}$. Notons $\psi^{_{''}\sharp}$ la famille form\'ee des $(\rho_{i},b_{i},a_{i})$ pour $i\in J_{\leq d}$. On sait que $\pi({\psi^{_{''}\sharp}})$ est l'image par l'involution de Zelevinsky de $\pi({\psi''})$. Gr\^ace \`a [A], on sait que la multiplicit\'e de $\sigma''_{\iota({\bf x})}$ dans $\pi(\psi'')_{P''}$ est \'egale \`a celle de $\sigma''_{{\bf x}}$ dans $\pi(\psi^{_{''}\sharp})_{P''}$. Finalement, on a deux repr\'esentations $\pi''$ et $\pi({\psi^{_{''}\sharp}})$ de $GL({n-n(\psi')})$ qui v\'erifient les propri\'et\'es suivantes:

- $\pi({\psi^{_{''}\sharp }})$ est irr\'eductible et appartient \`a $\Pi(\rho,\leq d)$;

- tout sous-quotient irr\'eductible de $\pi''$ appartient \`a $\Pi(\rho,\leq d)$;

- pour tout ${\bf x}\in {\cal X}$, les multiplicit\'es de $\sigma''_{{\bf x}}$ dans $\pi(\psi^{_{''}\sharp})_{P''}$ et dans $\pi''_{P''}$ sont \'egales.

Il en r\'esulte que $\pi''=\pi({\psi^{_{''}\sharp}})$. En effet, munissons ${\cal X}$ de l'ordre lexicographique. Pour toute repr\'esentation irr\'eductible $\mu$ de $GL({n-n(\psi')})$, appartenant \`a $\Pi(\rho,\leq d)$, notons ${\bf x}(\mu)$ le plus grand \'el\'ement de l'ensemble des ${\bf x}\in {\cal X}$ tels que $\mu_{P''}[\sigma''_{{\bf x}}]\not=0$. D'apr\`es [Z] proposition 6.9, l'application $\mu\mapsto {\bf x}(\mu)$ est injective. Consid\'erons l'\'el\'ement ${\bf x}^{\sharp}={\bf x}(\pi(\psi^{_{''}\sharp}))$. On a $\pi(\psi^{_{''}\sharp})_{P''}[\sigma''_{{\bf x}^{\sharp}}]\not=0$, donc aussi $\pi''_{P''}[\sigma''_{{\bf x}^{\sharp}}]\not=0$. Fixons un sous-quotient irr\'eductible $\pi''_{0}$ de $\pi''$ tel que $\pi''_{0,P''}[\sigma''_{{\bf x}^{\sharp}}]\not=0$. Alors ${\bf x}^{\sharp}$ est encore maximal parmi les ${\bf x}\in {\cal X}$ tels que $\pi''_{0,P''}[\sigma''_{{\bf x}}]\not=0$. Donc ${\bf x}(\pi''_{0})={\bf x}^{\sharp}$. D'apr\`es le r\'esultat de Zelevinsky que l'on vient de citer, $\pi''_{0}=\pi({\psi^{_{''}\sharp}})$. Mais alors on a aussi $\pi''=\pi''_{0}$ car les propri\'et\'es ci-dessus impliquent qu'il n'y a plus de place pour d'autres sous-quotients. Cela d\'emontre notre assertion.

On a maintenant:
$$\pi({\psi}^{\sharp})=\pi({\psi'})\times \pi''=\pi({\psi'})\times \pi({\psi^{_{''}\sharp}})=\pi({\psi^{\sharp}}),$$
\noindent ce qui ach\`eve la preuve.
\subsection{Action de $\theta$\label{actiondetheta}}
Les donn\'ees et les notations sont les m\^emes que dans \ref{exactitudeducomplexe}. Prolongeons $\pi({\psi})$ en une repr\'esentation $\pi(\psi)^+$ de $\tilde{G}({n})$. Alors on peut munir chaque terme  ${\cal E}_{j}(\pi(\psi))$  du complexe ${\cal E}(\pi(\psi))$ d'un prolongement \`a $\tilde{G}(n)$ que l'on note ${\cal E}_{j}(\pi(\psi)^+)$. Il se d\'efinit de la fa\c{c}on suivante. Soit $M\in {\cal M}_{j}(GL(n))$. De $\pi(\psi)^+(\theta)$ se d\'eduit naturellement un homomorphisme de $\pi(\psi)_{P}$ dans $\pi(\psi)_{\theta(P)}$, qui se restreint \'evidemment en un homomorphisme de $proj_{\rho,\leq d}(\pi(\psi)_{P})$ dans $proj_{\rho,\leq d}(\pi(\psi)_{\theta(P)})$. On induit cet homomorphisme en un homomorphisme:
 $$Ind_{P}^{GL(n)}(proj_{\rho,\leq d}(\pi(\psi)_{P})) \to Ind_{\theta(P)}^{GL(n)}(proj_{\rho,\leq d}(\pi(\psi)_{\theta(P)})).$$
 \noindent L'action de $\theta$ sur ${\cal E}_{j}(\pi(\psi))$ est la somme sur les $M\in {\cal M}_{j}(GL(n))$ de ces homomorphismes, multipli\'ee par $(-1)^{[\frac{j}{2}]}$. On v\'erifie sans peine que les diff\'erentielles $\epsilon_{j}$ sont \'equivariantes pour cette action de $\theta$.

{\bf Remarque.} C'est pour assurer cette compatibilit\'e que l'on a introduit la multiplication par $(-1)^{[\frac{j}{2}]}$. En effet, si on omet ces signes, on a \'equivariance non pas des diff\'erentielles $\epsilon_{j}$, mais des applications $\bar{\epsilon}_{j}$ d\'efinies comme $\epsilon_{j}$, mais en oubliant les signes $\xi_{M',M}$. Avec notre d\'efinition, on obtient la bonne \'equivariance en vertu de l'\'egalit\'e suivante: soient $M\in {\cal M}_{j}(GL(n))$, $M'\in {\cal M}_{j+1}(GL(n))$ deux L\'evi tels que $M'\subset M$; alors $(-1)^{[\frac{j}{2}]}\xi_{M',M}=(-1)^{[\frac{j+1}{2}]}\xi_{\theta(M'),\theta(M)}$.

Le conoyau de $\epsilon_{j(\psi)-1}$ se retrouve naturellement muni d'un prolongement \`a $\tilde{G}(n)$. Gr\^ace \`a la proposition 15, on obtient ainsi un prolongement $\pi(\psi^{\sharp})^+$ de $\pi(\psi^{\sharp})$. 

on pose $\beta(\psi,\rho,\leq d)=\begin{cases}(-1)^{\vert J_{\leq d}\vert (\vert J_{\leq d}\vert -1)/2}\prod_{i\in J_{\leq d}}(-1)^{(a_{i}b_{i}-1)/2}\hbox{, si $d$ est impair}\\\prod_{i\in J_{\leq d}}(-1)^{(a_{i}b_{i}-1)/2}\hbox{, si $d$ est pair}\end{cases}$.

On remarque que ce signe ne d\'epend que de $\psi\circ \Delta$, on le note donc $\beta(\psi\circ \Delta,\rho,\leq d)$. On v\'erifie encore ais\'ement que $\beta(\psi\circ \Delta,\rho,\leq d)$ vaut $(-1)^{[j(\psi)/2]}$, o\`u $j(\psi)$ a \'et\'e d\'efini au d\'ebut de \ref{exactitudeducomplexe}.

\

\bf Lemme. \sl L'action de $\theta$ sur $\pi(\psi^\sharp)$ est celle que nous y avons mise si et seulement si  $\beta(\psi,\rho,\leq d)=+1$.\rm 

\

Comme $\psi$ est \'el\'ementaire pour tout $(\rho,a,b)\in Jord(\psi)$, le quadruplet correspondant $(\rho,A,B,\zeta)$ v\'erifie $A=B$. On supprimera donc le $A$ de la notation. On va faire la d\'emonstration dans le cas o\`u $B$ est demi-entier pour tout $(\rho,B,\zeta)\in Jord(\psi)$. Le cas oppos\'e donne une \'ecriture plus compliqu\'ee mais n'introduit pas de signe. On note $\psi_{>d}$ le morphisme qui se d\'eduit de $\psi$ en enlevant tous les blocs $(\rho,B,\zeta)$ avec $B\leq (d-1)/2$.

Par d\'efinition $\pi(\psi)$ est un  sous-module irr\'eductible de l'induite:
$$\pi(\psi)\hookrightarrow
\biggl(\times_{(\rho,B,\zeta)\in Jord(\psi); B\leq (d-1)/2} <\rho\vert\,\vert^{\zeta B}, \cdots, \rho\vert\,\vert^{\zeta 1/2}>\biggr) \times \pi(\psi_{> d})$$
$$\times 
\biggl( \times_{(\rho,B,\zeta)\in Jord(\psi); B\leq (d-1)/2}<\rho\vert\,\vert^{- \zeta1/2}, \cdots, \rho\vert\,\vert^{-\zeta B}>\biggr),\eqno(1)
$$
le produit se fait en prenant les $B$ dans l'ordre croissant pour le terme de gauche et d\'ecroissant pour le terme de droite. On vient de d\'emontrer que l'image de $\pi(\psi)$ par $inv_{\leq d}$ est un sous-module irr\'eductible de l'induite:
$$\pi(\psi^\sharp)\hookrightarrow
\biggl(\times_{(\rho,B,\zeta)\in Jord(\psi); B\leq (d-1)/2} <\rho\vert\,\vert^{-\zeta B}, \cdots, \rho\vert\,\vert^{-\zeta 1/2}>\biggr) \times \pi(\psi_{> d})$$
$$
\biggl(\times_{(\rho,B,\zeta)\in Jord(\psi); B\leq (d-1)/2} \times <\rho\vert\,\vert^{\zeta 1/2}, \cdots, \rho\vert\,\vert^{\zeta B}>\biggr).\eqno(2)
$$
L'action $\theta(\psi)$ mise sur $\pi(\psi)$ est la restriction \`a $\pi(\psi)$ du prolongement canonique de $\theta(\psi_{>d})$ (cf. (1)) et on veut  montrer que l'action de $\theta$ que l'on obtient par la r\'esolution du complexe a la m\^eme propri\'et\'e en utilisant (2). Pour cela on note $\pi(\psi^\sharp)^+$ la repr\'esentation munie de l'action $\theta(\psi^\sharp)$ et $\pi(\psi^\sharp)^-$ la repr\'esentation munie de l'action $-\theta(\psi^\sharp)$.  On pose $$\sigma:=\biggl(\otimes_{(\rho,B,\zeta)\in Jord(\psi); B\leq (d-1)/2 }\otimes _{k\in [-\zeta B,-\zeta 1/2]} \rho\vert\,\vert^{-k}\biggr)\otimes \pi(\psi_{>d})$$
$$
\biggl(\otimes_{(\rho,B,\zeta)\in Jord(\psi); B\leq (d-1)/2}\otimes _{k\in [\zeta 1/2, \zeta B]}\rho\vert\,\vert^k\biggr).
$$On note $M_{0}$ le Levi sur lequel vit $\sigma$ et $P_{0}$ le parabolique standard de Levi $M_{0}$.
On note encore $\sigma^+$ la repr\'esentation $\sigma$ munie de l'action de $\theta$ qui \'echange la partie gauche avec la partie droite et qui agit par $\theta(\psi_{>d})$ sur la repr\'esentation du milieu. Et on note $\sigma^-$ la m\^eme repr\'esentation o\`u on remplace $\theta(\psi_{>d})$ par son oppos\'e. On v\'erifie ais\'ement que $\sigma^+$ intervient dans le module de Jacquet de $\pi(\psi^\sharp)^+_{P_{0}}$ alors que $\sigma^-$ n'y intervient pas. Pour $\pi(\psi^\sharp)^-$, les r\^oles sont invers\'es. Pour $\tau$ une repr\'esentation, avec action de $\theta$, du Levi sur lequel vit $\sigma$, on note $\tilde{m}(\tau)$ la multiplicit\'e de $\sigma^+$ moins la multiplicit\'e de $\sigma^-$. L'int\'er\^et d'une telle multiplicit\'e est qu'elle vaut $0$ pour toute repr\'esentation de la forme $\tau\oplus \tau'$ si $\theta$ \'echange les 2 facteurs. En effet la trace de $\theta$ y vaut 0 et donc la multiplicit\'e avec laquelle $\sigma^+$ intervient est la m\^eme que celle avec laquelle $\sigma^-$ intervient.

Soit $P$ un sous-groupe parabolique standard de Levi $M$;  pour toute repr\'esentation $\delta$ de $M$, on utilise la description du semi-simplifi\'e de $(Ind_{P}^{GL(n)}\delta)_{P_{0}}$ comme la somme des termes index\'es par les \'el\'ements du groupe de Weyl, $w$, de longueur minimale \`a gauche modulo le groupe de Weyl de $M_{0}$ et \`a droite modulo le groupe de Weyl de $M$, $Ind_{M_{0}\cap w P w^{-1}}^{M_{0}} w(\delta_{M \cap w^{-1}P_{0}w})$. On note $I_{P,w}$ ce terme pour $\delta$ $=$$proj_{\rho,\leq d}\pi(\psi)_{P}$.  Montrons que, la premi\`ere  somme ci-dessous porte sur les paraboliques standard de rang fix\'e:
$$
\tilde{m}(\sum_{P}(Ind_{P}^{GL(n)}(proj_{\rho,\leq d}\pi(\psi)_{P})_{P_{0}})=\sum_{P; \theta(M)=M}\tilde{m}(Ind_{P}^{GL(n)}(proj_{\rho,\leq d}\pi(\psi)_{P})_{P_{0}}).
$$
En effet, on regroupe dans la premi\`ere somme le parabolique de Levi $M$ avec celui de Levi $\theta(M)$; si ces 2 Levi sont diff\'erents, $\theta$ \'echange les 2 termes $\sigma$ isotypiques correspondants et agit donc par une repr\'esentation de trace 0. On fixe donc $P$ de Levi $\theta$-invariante; de la m\^eme fa\c{c}on, on montre que l'on peut limiter la somme aux \'el\'ements $w$ qui sont $\theta$-invariantes. Soient $M$ et $w$ $\theta$-invariantes; montrons que $\tilde{m}(I_{P,w})=0$ si  $ w^{-1}M_{0}w \not\subset M$. En effet, 
on \'ecrit $$M=GL(n_{1})\times \cdots \times GL(n_{r})\times GL(n_{0}) \times GL(n_{r})\times \cdots \times GL(n_{1}),$$
avec des $n_{i}>0$ pour $i\in [1,r]$ et $n_{0}\geq 0$. On a d\'efini $proj_{\rho,\leq d}$ pour les repr\'esentations de $M$; on note $proj^\theta_{\rho,\leq d}$ la somme des sous-modules irr\'eductibles de $proj_{\rho,\leq d}$ de la forme $\otimes_{i\in [1,r]}\sigma_{i}\otimes \sigma_{0}\otimes_{i\in [r,1]}\, ^\theta\sigma_{i}$. Par les arguments ci-dessus, on voit que l'on peut remplacer dans les calculs $proj_{\rho,\leq d}\pi(\psi)_{P}$ par $proj^\theta_{\rho,\leq d}\pi(\psi)_{P}$.  Et comme par d\'efinition une au plus des repr\'esentations $\delta'$ dans l'union des 3 ensembles  $$ 
\{\sigma_{i}; i\in [1,r]\} \cup \{^\theta\sigma_{i}; i\in [1,r]\}\cup \{ \sigma_{0}\}$$ pourrait contenir dans son support cuspidal $\rho\vert\,\vert^y$ avec $\vert y\vert >(d-1)/2$ ou $y-(d-1)/2)\notin {\mathbb Z}$, ce n'est pas le cas des repr\'esentations $\sigma_{i}$ pour $i\in [1,r]$. Le fait que $w$ est de longueur minimale dans sa classe modulo $M$ et la forme de $M_{0}$ entra\^{\i}ne alors l'assertion et entra\^{\i}ne m\^eme que $w^{-1}M_{0}w=M_{0}$.

Soit donc $P,M,w$ comme ci-dessus, tels que $M$ et $w$ soient $\theta$-invariantes et tels que $w^{-1}M_{0}w =M_{0}$. Alors $\tilde{m}I_{P,w}=\tilde{m}(w(\pi(\psi)_{P_{0}}))$. 

On note $W(M_{0})$ l'ensemble des \'el\'ements du groupe de Weyl $w$ tel que $wM_{0}w^{-1}=M_{0}$ et $W_{M_{0}}$ le groupe de Weyl de $M_{0}$.
Notons ici $\pi_{\sharp}$ la repr\'esentation $\pi(\psi^\sharp)$ munie de l'action de $\theta$ qui lui vient du complexe. Pour tout entier $j$, notons encore ${\cal M}^\theta_{j}$ l'ensemble des Levi standard de corang $j$ qui sont $\theta$ invariants et qui contiennent $M_{0}$. On a alors montr\'e que:
$$
\tilde{m}((\pi_{\sharp})_{P_{0}})=\sum_{j}(-1)^{j+[j/2]-j(\psi)}\sum_{M\in {\cal M}^\theta_{j}} \sum_{w\in W_{M_{0}}\backslash W(M_{0}); \theta(w)=w}\tilde{m}(w(\pi(\psi)_{P_{0}})).\eqno(3)
$$
On remarque que $(-1)^{j+[j/2]}=(-1)^{-j+[j/2]}=(-1)^{[(j+1)/2]}$. On peut inverser l'ordre des sommations, c'est-\`a-dire fixer $w$ et sommer sur les Levi standard $M$, $\theta$-invariantes tels que $w$ soit de longueur minimale dans sa classe \`a droite modulo $M$. Pour $w$ fix\'e, on pose $\Delta^w$ l'ensemble des racines simples positives hors de $M_{0}$ dont les images par $w$ soient encore positives. Cet ensemble est $-\theta$-invariante et on note $\Delta_{ \theta}^w$ les classes de conjugaison pour cette action $-\theta$. Un Levi standard $M$, $\theta$-invariante tel que $w$ soit minimal dans sa classe \`a droite modulo $M$ est pr\'ecis\'ement d\'efini pour un sous-ensemble de $\Delta_{\theta}^w$. Ainsi, si $\Delta_{\theta}^w$ est non vide, la somme altern\'ee est nulle. Il n'y a qu'un seul \'el\'ement $w$ pour lequel $\Delta^w$ est vide, celui qui envoie toutes les racines positives hors de $M_{0}$ sur des racines n\'egatives; notons le $w_{0}$. On a 
$$
w_{0}\sigma=\biggr(\otimes_{(\rho,B,\zeta)\in Jord(\psi); B\leq (d-1)/2 }\otimes _{k\in [\zeta B,\zeta 1/2]}\rho\vert\,\vert^{-k}\biggr) \otimes \pi(\psi_{>d}) $$
$$
\biggl(\otimes_{(\rho,B,\zeta)\in Jord(\psi); B\leq (d-1)/2}\otimes _{k\in [-\zeta 1/2, -\zeta B]}\rho\vert\,\vert^k\biggr).
$$
Si $\psi$ est de restriction discr\`ete \`a la diagonale (cas que nous utiliserons), la multiplicit\'e de cette repr\'esentation dans $\pi(\psi)_{P_{0}}$ est pr\'ecis\'ement $1$ avec comme action de $\theta$ l'action qui se d\'eduit naturellement de $\theta(\psi_{>d})$. Avec (3), cela donne:
$$
\tilde{m}((\pi_{\sharp})_{P_{0}})=(-1)^{j(\psi)+[j(\psi)/2]-j(\psi)}=(-1)^{[j(\psi)/2]}.
$$
C'est le r\'esultat annonc\'e.

\section{Transfert, contexte et propri\'et\'es g\'en\'erales}
\subsection{Les groupes\label{lesgroupes}}
Soit $n\in {\mathbb N}$. On reprend la notation $\tilde{G}(n)$ pour le produit semi-direct de $GL(n)$ avec le groupe engendr\'e par $\theta$ et on note $G^+_{n}$ la composante non neutre de ce groupe.

Si $n$ est pair, on note $H_{n}$ la forme d\'eploy\'ee du groupe sp\'ecial orthogonal $SO(n+1)$ sur $F$. Si $n$ est impair, on note $H_{n}$ le groupe symplectique $Sp(n-1)$ sur $F$.  Plus exactement,  $H_{n}$ est le groupe des points sur $F$ du groupe alg\'ebrique indiqu\'e. Le groupe $H_{n}$ est un groupe endoscopique de $\tilde{G}(n)$, c'est celui qui contr\^ole les distributions stables sur $G^+_{n}$.  Notons $G^+_{n,reg}$, resp. $H_{n,reg}$, l'ensemble des \'el\'ements semi-simples r\'eguliers de $G^+_{n}$, resp. $H_{n}$. On sait d\'efinir une application, d'ailleurs bijective, entre classes de conjugaison stable  dans $G^+_{n,reg}$ et classes de conjugaison stable dans $H_{n,reg}$; cf. \cite{kottwitz}, \cite{labesse} et pour une description concr\`ete \cite{waldspurger1} III.2.

 Pour $\tilde{g}\in G^+_{n,reg}$ et $h\in H_{n,reg}$, on note $\tilde{g}\sim h$ si les classes de conjugaison stable de $\tilde{g}$ et de $h$ se correspondent. Soient ${\bf D}$, resp. ${\bf D}^H$, deux distributions sur $\tilde{G}(n)$, resp. $H_{n}$, invariantes par conjugaison par $GL(n)$, resp. $H_{n}$, et localement int\'egrables, c'est-\`a-dire d\'efinies par des fonctions localement int\'egrables not\'ees $D$, resp. $D^H$. Supposons de plus $D$, resp. $D^H$, localement constantes sur $G^+_{n,reg}$, resp. $H_{n,reg}$. Alors ${\bf D}$, resp. ${\bf D}^H$ sont stablement invariantes si et seulement si $D$, resp. $D^H$, est constante sur les classes de conjugaison stable contenues dans $G^+_{n,reg}$, resp. $H_{n,reg}$. La distribution ${\bf D}$ est un transfert de ${\bf D}^H$ si et seulement si on a l'\'egalit\'e $D(\tilde{g})=D^H(h)$ pour tout couple $(\tilde{g},h)\in G^+_{n,reg}\times H_{n,reg}$ tel que $\tilde{g}\sim h$. En particulier, soient $\pi$, resp. $\pi^H$, des repr\'esentations virtuelles de $\tilde{G}(n)$, resp. $H_{n}$, c'est-\`a-dire des combinaisons lin\'eaires \`a coefficients dans ${\mathbb Z}$ de repr\'esentations irr\'eductibles. Notons $trace\,\pi^H$ le caract\`ere de $\pi^H$ (d\'efini par lin\'earit\'e \`a partir du cas o\`u $\pi^H$ est irr\'eductible). Notons $trace_{G^+_{n}}\pi$ la restriction \`a $G^+_{n}$ du caract\`ere de $\pi$. Pour simplifier les notations, on consid\`ere ces caract\`eres soit comme des distributions, soit comme des fonctions localement constantes sur $G^+_{n,reg}$ ou $H_{n,reg}$. On dit que $\pi$, resp. $\pi^H$, est stable si et seulement si $trace_{G^+_{n}}\pi$ l'est, resp. $trace\,\pi^H$ l'est. On dit que $\pi$ est un transfert de $\pi^H$ si et seulement si $trace_{G^+_{n}}\pi$ est un transfert de $trace\,\pi^H$.
\subsection{Propri\'et\'es g\'en\'erales du transfert stable\label{proprietesgeneralesdutransfertstable}}
Comme souvent quand on travail avec $GL(n),H_{n}$ on doit aussi travailler avec des Levi de ces groupes qui se correspondent; consid\'erons une partition de $n$ de la forme sym\'etrique suivante $$(n_{1}, \cdots, n_r,n_{0}, n_{r}, \cdots, n_{1});n=\sum_{i\in [1,r]}2n_{i}+n_{0}.$$  On note $M$ le Levi standard de $GL(n)$ correspondant \`a cette partition et $M^H$ celui de $H_{n}$; $M^H$ est isomorphe \`a $\times_{i\in [1,r]}GL(n_{i})\times H_{n_{0}}$ et $M$ est $\theta$ invariant. Ceci permet de d\'efinir $\tilde{M}$ et $M^+$ comme sous-groupe, resp. sous-ensemble de $\tilde{G}({n})$ et $G^+_{n}$. On fixe un ensemble $\cal C$ de repr\'esentations cuspidales (non autoduales en g\'en\'eral). Pour $\pi$ une repr\'esentation de $\tilde{G}$, on d\'efinit $Jac^\theta_{\cal C}\pi$ de la fa\c{c}on suivante: on consid\`ere d'abord la restriction de $\pi$ le long du parabolique $P$ de Levi $M$. On obtient donc une repr\'esentation de $M$; on projette cette repr\'esentation sur son facteur direct tel que l'action des $r$ premiers (resp. derniers) facteurs $GL(n_{i})$ se fasse via des repr\'esentations de support cuspidal les repr\'esentations de ${\cal C}$ (resp. $^\theta {\cal C}$); on note cette op\'eration $proj^\theta_{\cal C}$.  L'action de $\theta$ sur $\pi$ donne une action de $\theta$ sur cette projection et on obtient donc une repr\'esentation de $\tilde{M}$. Si $\pi^H$ est une repr\'esentation de $H$, on d\'efinit de fa\c{c}on analogue (et plus simple) $Jac_{{\cal C}}\pi^H$, ce qui utilise la projection $proj_{\cal C}$. On \'etend ces d\'efinitions lin\'eairement au groupe de Grothendieck des repr\'esentations de $\tilde{G}(n)$ et de $H$.

\

{\bf Lemme.} \sl Soient $\pi^H$ une repr\'esentation virtuelle de $H_{n}$, $\pi$ une repr\'esentation virtuelle de $\tilde{G}(n)$ et $M, M^H$ comme ci-dessus.  Supposons que $\pi^H$ et $\pi$ sont stables et que $\pi$ est un transfert de $\pi^H$.

Alors $Jac_{\cal C}\pi^H$ et  $Jac^\theta_{\cal C}\pi$ sont stables et $Jac^\theta_{\cal C}\pi$ est un transfert de $Jac_{\cal C}\pi^H$.\rm

\

On d\'emontre en 2 temps cette assertion; d'abord on d\'emontre que les modules de Jacquet ont les propri\'et\'es de stabilit\'e et de transfert annonc\'es puis on montre que leur projection suivant ${\cal C}$ ont aussi ces propri\'et\'es. C'est un r\'esultat qui ne surprend personne mais dont on redonne la d\'emonstration.

On note $\pi_{P}$ et $\pi^H_{P}$ les modules de Jacquet de $\pi$ et $\pi^H$ respectivement; on rappelle que $\theta$ agit sur $\pi_{P}$ et que $\pi_{P}$ est donc une repr\'esentation de $\tilde{M}$. Il suffit de montrer que, pour tous $\tilde{g}\in M^+_{reg}$ et $h\in M^H_{reg}$ tels que $\tilde{g}\sim h$, on a l'\'egalit\'e $trace_{\tilde{M}}\pi^+_{P}(\tilde{g})=trace\,\pi^H_{P^H}(h)$. Fixons de tels $\tilde{g}$ et $h$. Ecrivons:
$$M=GL({n_{1}})\times...\times GL({n_{r}})\times GL({n_{0}})\times GL({n_{r}})\times...\times GL({n_{1}})$$
$$M^H=GL({n_{1}})\times...\times GL({n_{r}})\times H_{n_{0}}.$$
Soient $z_{1},...,z_{r}\in F^{\times}$, posons:
$$z_{M}=z_{1}{\bf 1}\times...\times z_{r}{\bf 1}\times {\bf 1}\times z_{r}^{-1} {\bf 1}\times...\times z_{1}^{-1} {\bf 1},$$
$$z_{M^H}= z_{1}{\bf 1}\times...\times z_{r}{\bf 1}\times {\bf 1},$$
\noindent o\`u on note uniform\'ement ${\bf 1}$ les unit\'es des groupes en question. Supposons que les r\'eels suivants soient assez petits:

-  $\vert \frac{z_{i}}{z_{i+1}}\vert$ pour $i=1,...,r-1$ et $\vert z_{r}\vert$.

Alors un r\'esultat  de Casselman affirme que l'on a les \'egalit\'es:
$$trace_{M^+}\pi_{P}(z_{M}\tilde{g})=\delta_{P}(z_{M}\tilde{g})^{-1/2}trace_{G^+_{n}}\pi(z_{M}\tilde{g}),$$
$$trace\,\pi^H_{P^H}(z_{M^H}h)=\delta_{P^H}(z_{M^H}h)^{-1/2}trace\,\pi^H(z_{M^H}h),$$ 
o\`u $\delta_{P}$ et $\delta_{P^H}$ sont les modules usuels. De la d\'efinition de la correspondance entre classes de conjugaison stables et de l'hypoth\`ese $\tilde{g}\sim h$ r\'esulte que $z_{M}\tilde{g}\sim z_{M^H}h$. On v\'erifie par ailleurs l'\'egalit\'e des facteurs $\delta$ ci-dessus. Puisque $\pi$ est un transfert de $\pi^H$, les membres de droite des \'egalit\'es ci-dessus sont \'egaux. Donc les membres de gauche le sont aussi. Mais, parce que $z_{M}$ appartient au centre de $M^+$, la fonction:
 $$(z_{1},...,z_{r})\mapsto trace_{\tilde{M}}\pi_{P}(z_{M}\tilde{g})$$
 \noindent est une combinaison lin\'eaire finie de caract\`eres de $(F^{\times})^r$. Il en est de m\^eme de la fonction:
 $$(z_{1},...,z_{r})\mapsto trace\,\pi^H_{P^H}(z_{M^H}h).$$
 \noindent On vient de voir que ces fonctions co\"{\i}ncident dans un "c\^one ouvert" de l'espace de d\'epart. Il en r\'esulte qu'elles sont \'egales. En particulier, leurs valeurs en $(z_{1},...,z_{r})=(1,...,1)$ sont \'egales, c'est-\`a-dire $trace_{M^+}\pi_{P}(\tilde{g})=trace\,\pi^H_{P^H}(h)$, ce que l'on voulait d\'emontrer.
 
 On reprend les notations $proj^\theta_{\cal C}$ et $proj_{\cal C}$ introduites avant l'\'enonc\'e du lemme. Le point \`a d\'emontrer est maintenant le suivant: soit $\pi_{M}$ une repr\'esentation virtuelle de $\tilde{M}$ et $\pi^H_{M}$ une repr\'esentation virtuelle de $M^H$; on suppose qu'elles sont stables et que $\pi_{M}$ est un transfert de $\pi^H_{M}$, alors la m\^eme propri\'et\'e est vraie  pour $proj^\theta _{\cal C}\pi_{M}$ et $proj_{\cal C}\pi^H_{M}$. On \'ecrit $\pi_{M}$ dans la base des repr\'esentations irr\'eductibles de $\tilde{M}$ et on remarque tout de suite que toute repr\'esentation irr\'eductible de $\tilde{M}$ dont la restriction \`a $M$ n'est pas irr\'eductible donne une trace nulle sur $M^+$; on peut donc supposer que $$\pi_{M}=\sum_{\sigma,\sigma_{0}}c(\sigma,\sigma_{0}) \sigma \otimes \sigma_{0}\otimes \theta(\sigma),$$o\`u $\sigma$ parcourt l'ensemble des repr\'esentations irr\'eductibles de $\times_{i\in [1,r]}GL(n_{i})$ et $\sigma_{0}$ l'ensemble des repr\'esentations irr\'eductibles de $\tilde{G}(n_{0})$ et o\`u $c(\sigma,\sigma_{0})$ est un nombre complexe. On regroupe en fonction des $\sigma$:
 $$
 \pi_{M}=\sum_{\sigma}\sigma\otimes \sigma_{0}[\sigma] \otimes \theta(\sigma),
 $$
 o\`u $\sigma$ est comme ci-dessus mais o\`u $\sigma_{0}[\sigma]$ est une repr\'esentation virtuelle de $\tilde{G}(n_{0})$; remarquons que si $n_{0}=0$ il reste quand m\^eme une action de $\theta$, c'est-\`a-dire un signe. On \'ecrit de fa\c{c}on analogue
 $$
 \pi_{M}^H=\sum_{\sigma}\sigma\otimes \sigma_{0}^H[\sigma].
 $$
 On va d\'emontrer que sous l'hypoth\`ese que $\pi_{M}$ et $\pi^H_{M}$ sont  stables et des transferts l'une de l'autre alors il en est de m\^eme de $\sigma_{0}[\sigma]$ et $\sigma_{0}^H[\sigma]$ pour tout $\sigma$ comme ci-dessus. Ce sera largement suffisant.

  Soient $g=(g_{1},...,g_{r},g_{0},g'_{r},...,g'_{0})\in M_{reg }$ et $h=(h_{1},...,h_{r},h_{0})\in M^H_{reg}$. On v\'erifie les propri\'et\'es suivantes:
 
 - $g\theta\sim h$ si et seulement si $g_{0}\theta\sim h_{0}$ et, pour tout $i=1,...,r$, $g_{i}\theta(g'_{i})$ est conjugu\'e \`a $h_{i}$ dans $GL(n_{i})$;
 
 - pour toute repr\'esentation irr\'eductible $\sigma=\sigma_{1}\otimes...\otimes \sigma_{r}$ de $GL(n_{1})\times...\times GL(n_{r})$ et pour toute repr\'esentation virtuelle $\sigma_{0}^+$ de $GL(n_{0})^+$, on a l'\'egalit\'e:
 $$trace_{\tilde{M}}(\sigma\otimes \sigma_{0}^+)(g\theta)=trace_{G^+_{n_{0}}}\sigma_{0}^+(g_{0}\theta)
 \prod_{i=1,...,r}trace\,\sigma_{i}(g_{i}\theta(g'_{i}))).$$
 
 Dire que $\pi_{M}$ et $\pi_{M}^H$ sont stables et que $\tau_{M}$ est un transfert de $\pi_{M}^H$ signifie donc que pour tout $\gamma\in GL(n_{1}){reg}\times...\times GL(n_{r}){reg}$, tout $\tilde{g}_{0}\in G^+_{n_{0}}$, tout $h_{0}\in H_{n_{0},reg}$ tels que $\tilde{g}_{0}\sim h_{0}$, on a l'\'egalit\'e:
$$\sum_{\sigma}trace\,\sigma(\gamma)trace_{G^+_{n_{0}}}\sigma_{0}^+[\sigma](\tilde{g}_{0})=\sum_{\sigma}trace\,\sigma(\gamma)trace\,\sigma_{0}^H[\sigma](h_{0}).$$
 Consid\'erons $\tilde{g}_{0}$ et $h_{0}$ comme fix\'es et $\gamma$ comme variable. L'\'egalit\'e ci-dessus est une \'egalit\'e entre combinaisons lin\'eaires (\`a coefficients complexes) de traces de repr\'esentations irr\'eductibles de $GL({n_{1}})\times...\times GL({n_{r}})$. On sait que ces traces sont lin\'eairement ind\'ependantes. Ainsi cette \'egalit\'e entra\^{\i}ne que, pour tout $\sigma$, on a l'\'egalit\'e:
 $$trace_{G^+_{n_{0}}}\sigma_{0}^+[\sigma](\tilde{g}_{0})=trace\,\sigma_{0}^H[\sigma](h_{0}).$$
 Ces \'egalit\'es  ach\`event la d\'emonstration.

\subsection{Les repr\'esentations}
Jusqu'\`a pr\'esent, nous avions regard\'e un ensemble de morphismes $\psi$ $\theta$-invariant, il faut maintenant tenir compte des questions de parit\'e qui sont li\'es aux choix de $H_{n}$ si l'on veut que l'objet attach\'e soit un transfert d'un analogue pour $H_{n}$. Pour $\rho$ une repr\'esentation cuspidale irr\'eductible autoduale, on attache  un signe $\eta_{\rho}$ de la fa\c{c}on suivante. D'apr\`es la conjecture de Langlands (th\'eor\`eme de Harris-Taylor et Henniart), on sait associer \`a $\rho$ une repr\'esentation irr\'eductible de dimension $d_{\rho}$ du groupe de Weil de $F$. Parce que $\theta(\rho)=\rho$, l'image de cette repr\'esentation est incluse dans un sous-groupe de $GL_{d_{\rho}}({\mathbb C})$ qui est soit un groupe orthogonal, soit un groupe symplectique. On pose $\eta_{\rho}=1$ dans le premier cas, $\eta_{\rho}=-1$ dans le second. On rappelle d'autre part que, pour toute repr\'esentation $\pi$ irr\'eductible, on note $\chi_{\pi}$ son caract\`ere central. L'hypoth\`ese $\theta(\pi)=\pi$ entra\^{\i}ne $\chi_{\pi}^2=1$. Pour $\rho$ comme ci-dessus telle que $\eta_{\rho}=-1$, $\chi_{\rho}$ est n\'ecessairement trivial. On va maintenant noter $\Psi^H_{n}$ l'ensemble des morphismes $\psi$ de $W_{F}\times SL(2,{\mathbb C}) \times  SL(2,{\mathbb C})$ dans $GL(n,{\mathbb C})$ v\'erifiant les propri\'et\'es suivantes (les morphismes sont regard\'es \`a conjugaison pr\`es):

1- l'image de $\psi$ est \`a valeurs dans $Sp(n)$ si $n$ est pair et \`a valeurs dans $SO(n)$ si $n$ est impair. On note $H^*$ le groupe $Sp(n)$ ou $SO(n)$, comme ci-dessus;

2- le centralisateur de $\psi$ dans $H^*$ est un groupe \`a centre fini.

Cela se traduit encore de fa\c{c}on combinatoire ainsi: $\psi$ d\'efinit une repr\'esentation de dimension $n$ de $W_{F}\times SL(2,{\mathbb C}) \times SL(2,{\mathbb C})$ que l'on d\'ecompose en sous-repr\'esentations irr\'eductibles. Chaque sous-repr\'esentation irr\'eductible est associ\'ee \`a un triplet $(\rho,a,b)$ o\`u $\rho$ est une repr\'esentation irr\'eductible de $W_{F}$ et $a,b$ sont les dimensions des repr\'esentations irr\'eductibles de $SL(2,{\mathbb C})$. La condition 1 assure que la multiplicit\'e de la repr\'esentation associ\'ee au triplet $(\rho,a,b)$ dans cette d\'ecomposition est la m\^eme que celle de la repr\'esentation $(\rho^*,a,b)$. La condition 2 entra\^{\i}ne alors que $\rho\simeq \rho^*$ (sinon le centralisateur contient un groupe lin\'eaire). Et la condition 1 peut alors se r\'ecrire:
$$
\eta_{\rho}(-1)^{a+b}=(-1)^{n+1} \qquad \hbox{ et }\qquad \prod_{(\rho,a,b)}\chi_{\rho}^{ab}=1.
$$
{\bf Remarque.} \sl  Pour $\psi\in \Psi_{n}$, on a d\'efini $\pi(\psi)$ comme repr\'esentation de $\tilde{G}(n)$. Cette d\'efinition d\'epend du choix initial de l'action de $\theta$ sur chaque $\rho$ intervenant dans la d\'ecomposition. Un choix raisonnable, c'est celui que nous avons fait, est de fixer un caract\`ere additif $\tau$ est d'imposer que $\theta$ agisse trivialement sur l'espace des vecteurs de Whittaker associ\'e \`a $\pi(\psi)$ (l'espace des coinvariants comme expliqu\'e dans l'introduction). L'action est ind\'ependante du choix de $\tau$.\rm

\

Rempla\c{c}ons $\tau$ par le caract\`ere $x\mapsto \tau(zx)$ avec $z\in F^{\times}$. L'action de $\theta$ sur une repr\'esentation autoduale $\rho$ d'un groupe $GL(d_{\rho})$ est alors multipli\'ee par $\chi_{\rho}(z)^{d_{\rho}+1}$ (on applique les d\'efinitions). Pour d\'emontrer l'assertion d'ind\'ependance, il suffit de le faire dans le cas o\`u $\psi$ est de restriction discr\`ete \`a la diagonale. On reprend les d\'efinitions et le changement de $\zeta$ induit donc la multiplication par le signe:
$$
\prod_{(\rho,a,b)\in Jord(\psi)}\chi_{\rho}(z)^{(d_{\rho}+1)(a-2[a/2])(b-2[b/2])}=\prod_{(\rho,a,b)}\chi_{\rho}(z)^{ab}=1$$
par les conditions ci-dessus. D'o\`u l'assertion.

\

\bf D\'efinition. \sl On note $\Psi^{temp}_{n}$ le sous-ensemble de $\Psi$ form\'e des morphismes triviaux sur la 2e copie de $SL(2,{\mathbb C})$.\rm

\subsection{Hypoth\`ese\label{hypothese}}

On fixe un ensemble ${\cal E}$ de repr\'esentations irr\'eductibles de $W_{F}$  et pour tout $n$, on note $\Psi_{n,\cal E}$ l'ensemble des \'el\'ements de $\Psi_{n}$ dont la restriction \`a $W_{F}$ est une somme de repr\'esentations isotypiques de type des repr\'esentations dans ${\cal E}$. On d\'efinit de fa\c{c}on \'evidente $\Psi_{n,\cal E}^{temp}$. 
\

\bf Hypoth\`ese. \sl Nous allons supposer \`a partir de maintenant que ${\cal E}$ est fix\'e de telle sorte que pour tout $n$ et tout $\phi \in \Psi_{n,\cal E}^{temp}$, la repr\'esentation temp\'er\'ee, $\pi(\phi)$, associ\'ee \`a $\phi$ (gr\^ace aux travaux de Zelevinsky et \`a la preuve de la conjecture de Langlands)  est un transfert d'une repr\'esentation virtuelle $\pi^H(\phi)$ de $H_{n}$.\rm

\

Remarquons que l'on n'a pas pr\'ecis\'e l'action de $\theta$; il est clair que si l'on a un transfert pour un choix d'action de $\theta$, on en aura pour tout choix, la seule chose est que $\pi(\phi)$ n'est d\'efini qu'au signe pr\`es. Toutefois, les repr\'esentations irr\'eductibles de $H_{n}$ intervenant dans la d\'ecomposition de $\pi^H(\phi)$ sont uniquement d\'etermin\'ees par $\phi$ et l'\'egalit\'e des traces.
On notera $\Pi^H(\phi)$ cet ensemble de repr\'esentations.

\subsection {Remarques sur l'hypoth\`ese; points de r\'eductibilit\'e des induites de cuspidales\label{remarquesurlhypothese}}
On suppose que le lemme fondamental pour $H_{n},\tilde{G}(n)$ est connu, les travaux d'Arthur d\'emontrent alors cette hypoth\`ese. Ces travaux la d\'emontrent d'ailleurs sans supposer que $\phi$ est temp\'er\'e. Ce que nous faisons ici, ramener le cas g\'en\'eral au cas temp\'er\'e, n'a donc d'int\'er\^et que si l'on peut \^etre plus pr\'ecis et en particulier si on peut d\'ecomposer les repr\'esentations virtuelles obtenues.  Sans connaissance pr\'ecise de $\pi^H(\phi)$ dans le cas cas temp\'er\'e, nous n'ajoutons pas grand chose.

 Ici, nous allons voir rapidement ce que l'on peut en d\'eduire de cette hypoth\`ese tr\`es g\'en\'erale et reposer le probl\`eme pr\'ecis\'ement.

\bf Remarque 1. \sl Arthur annonce en \cite{arthurrecent} 30.1 que  pour un bon choix de l'action de $\theta$,  $\pi^H(\phi)$ est pour tout $\phi$ une combinaison lin\'eaire \`a coefficients positifs de repr\'esentations. Alors $\pi^H(\phi)$ est une repr\'esentation temp\'er\'ee.\rm

En effet, les modules de Jacquet non nul de $\pi^H(\phi)$ se transf\`erent en des modules de Jacquet non nuls de $\pi(\phi)$ d'apr\`es ce que l'on a vu en \ref{proprietesgeneralesdutransfertstable} et il est alors facile de v\'erifier la positivit\'e des exposants. On peut conclure que chaque repr\'esentation irr\'eductible constituant $\pi^H(\phi)$ est temp\'er\'ee car il ne peut y avoir de simplification dans le calcul des modules de Jacquet sous l'hypoth\`ese faite que $\pi^H(\phi)$ est une repr\'esentation et pas seulement une repr\'esentation virtuelle. Une d\'emonstration plus conceptuelle est dans \cite{waldspurger1} VI.1.

\

\bf Remarque 2. \sl Soit $\sigma$ une repr\'esentation cuspidale de $H_{n}$; supposons qu'il existe un morphisme $\phi\in \Psi_{n,\cal E}^{temp}$ tel que $\sigma\in \Pi^H(\phi)$.  Alors

(i) $Jord(\phi)$ est sans multiplicit\'es et sans trou, c'est-\`a-dire que si une repr\'esentation $\rho\otimes [a]$ de $W_{F}\times SL(2,{\mathbb C})$ ($[a]$ est la repr\'esentation de dimension $a$ de $SL(2,{\mathbb C})$) intervient dans la repr\'esentation d\'efinie par $\phi$, elle y intervient avec multiplicit\'e 1 et si $a>2$ la repr\'esentation $\rho\otimes [a-2]$ y intervient aussi avec multiplicit\'e 1.

(ii) Soit $\rho$ une repr\'esentation irr\'eductible cuspidale de $GL(d_{\rho},F)$ (ce qui d\'efinit $d_{\rho}$). On suppose que $\rho\simeq \rho^*$ et que l'ensemble $Jord_{\rho}(\phi)\neq 0$ et on note $a_{max}$ son \'el\'ement maximum. Soit $x_{0}\in {{\mathbb R}}$ tel que l'induite $\rho\vert\,\vert^{x_{0}}\times \sigma$ soit r\'eductible. Alors $x_{0}=\pm (a_{max}+1)/2$.

(iii)Soit $\rho$ comme en (ii) mais tel que $Jord_{\rho}(\phi)=\emptyset$; on suppose que $\eta_{\rho}=(-1)^{n+1}$, alors pour $x_{0}$ comme en (ii), $x_{0}=\pm 1/2$.

\rm

\

Avant de donner la preuve, remarquons que l'\'enonc\'e ne dit rien si $\eta_{\rho}=(-1)^n$ et l'ensemble de (ii) est vide. On traitera ce cas dans \cite{classification}. 
 Toutefois, les propri\'et\'es d\'emontr\'ees ici sont suffisantes pour construire avec \cite{europe}, \cite{ams}, pour tout morphisme $\phi\in \Psi^{temp}_{n,\cal E}$ discret (i.e. $Jord(\phi)$ est sans multiplicit\'es) une application (non triviale) de l'ensemble des caract\`eres du centralisateur de $\phi$ dans $^LH_{n}$, $\epsilon \mapsto \pi^H(\phi,\epsilon)$, avec $\pi^H(\phi,\epsilon)$ est soit 0 soit une s\'erie discr\`ete sans multiplicit\'es. Les questions sont alors:

1- a-t-on $\Pi^H(\phi)=\{\pi^H(\phi,\epsilon)\}$?  la r\'eponse est n\'ecessairement oui si l'on sait a priori que les $c_{\pi}$ sont des r\'eels positifs pour un bon choix de l'action de $\theta$, ce qui fait partie des r\'esultats d'Arthur d\'ej\`a cit\'es.

2- a-t-on plus pr\'ecis\'ement $tr\pi(\phi)(g,\theta)=\pm\sum_{\epsilon} tr\pi^H(\phi,\epsilon)(h)$ quand les classes de conjugaison stable de $h$et de $ (g,\theta)$ se correspondent?  C'est cette formule que nous  a en vue avec une description pr\'ecise des repr\'esentations $\pi^H(\phi,\epsilon)$.

\

Supposons d'abord que $Jord(\phi)$ ne soit pas sans multiplicit\'es; il existe donc $\rho\otimes [a]$ qui intervient dans $\phi$ avec multiplicit\'e au moins $2$ (notations de l'\'enonc\'e). On d\'ecompose:
$$
\phi= \phi_{0}\oplus \rho\otimes [a]\oplus \rho\otimes [a].
$$
Pour $n_{0}$ convenable, $\phi_{0}\in \Psi^{temp}_{n_{0},\cal E}$. On v\'erifie que pour tout $g \in \tilde{G}_{n,reg}^+$
$$
tr\, \pi(\phi)(g)=tr\, St(a,\rho)\times \pi(\phi_{0})\times St(a,\rho) (g)
$$
qui dans le transfert donne une \'egalit\'e en tout point $h\in H_{n,reg}$
$$
tr \pi^H(\phi)(h)=tr \, (ind(St(a,\rho)\otimes \pi^{H_{n_{0}}}(\phi_{0}))(h).
$$
Ainsi $\Pi^H(\phi)$ est inclus dans la r\'eunion des composantes irr\'eductibles des induites de la forme $St(a,\rho)\times \sigma_{0}$ pour $\sigma_{0}\in \Pi^{H_{n_{0}}}(\phi_{0})$. A priori, il pourrait y avoir des simplifications et c'est pour cela que l'on n'a qu'une inclusion. Mais comme $\sigma$ est cuspidale par hypoth\`ese, elle ne pourrait \^etre dans $\Pi^H(\phi)$. D'o\`u la premi\`ere assertion de (i).

Montrons le reste de la remarque: on fixe, en admettant que  c'est possible, $(\rho,a)$ avec $a\geq 0$ tel que $\rho\otimes [a+2]$ ne soit pas une sous-repr\'esentation de de $\phi$ (i.e. $(\rho,a+2)\notin Jord(\phi)$) alors que $\rho\otimes [a]$ en est une si  $a>0$; si $a=0$ on suppose que  $\eta_{\rho}=(-1)^{n+1}$. On note $\phi_{1}$ le morphisme qui se d\'eduit de $\phi$ de sorte que $$Jord(\phi_{1})=Jord(\phi)-(\rho,a)\cup (\rho,a+2).$$
Clairement $\phi_{1}\in \Psi^{temp}_{n+2d_{\rho},{\cal E}}$. On calcule (les $Jac$ sont comme en \ref{notationdujac} pour $\rho$ fix\'e):
$$
Jac^\theta_{{(a+1)/2}}\pi(\phi_{1})=\pi(\phi); \qquad Jac^\theta_{{-(a+1)/2}}\pi(\phi_{1})=0.
$$En particulier $Jac_{(a+1)/2}\pi^H(\phi_{1})=\pi^H(\phi)$ et $Jac_{-(a+1)/2}\pi^H(\phi_{1})=0$.
Ainsi, il existe $\sigma_{1}$ dans $\Pi^H(\phi_{1})$ tel que $Jac_{(a+1)/2}\sigma_{1}$ contient $\sigma$ comme sous-quotient. Ainsi $\sigma_{1}$ est un sous-module irr\'eductible de l'induite $\rho\vert\,\vert^{(a+1)/2}\times \sigma$. Si cette induite est irr\'eductible, $\sigma_{1}$ co\"{\i}ncide avec toute l'induite et est donc l'unique \'el\'ement de $\Pi^H(\phi_{1})$ sous-quotient de cette induite. Et on ne pourrait avoir de simplification de $Jac_{-(a+1)/2}\sigma_{1}$ quand on calcule $Jac_{-(a+1)/2}\pi^H(\phi_{1})$. On vient donc de montrer que $\rho\vert\,\vert^{(a+1)/2}\times \sigma$ est r\'eductible. D'apr\`es un r\'esultat de Silberger (\cite{silberger}) $(a+1)/2$ est uniquement d\'etermin\'e par $\rho$ et $\sigma$ (il n'y a qu'un point de r\'eductibilit\'e r\'eel positif ou nul).

Supposons maintenant que $\rho$ est tel que $Jord_{\rho}(\phi)\neq 0$. En prenant $a=a_{max}$ l'\'el\'ement maximum de $Jord_{\rho}(\phi)$, l'hypoth\`ese $(\rho,a+2)\notin Jord(\phi)$ est bien satisfaite; ainsi $(a_{max}+1)/2$ est le point de r\'eductibilit\'e cherch\'e.  L'in\'egalit\'e prouv\'ee en \cite{algebra} (cf. l'intorduction de cet article) montre que $\phi$ est sans trou et on a prouv\'e (ii). On d\'emontre (iii) en prenant $a=0$. Cela termine la preuve.

\subsection{Construction des repr\'esentations de $H_{n}$\label{constructiondesrepresentations}}
On fixe $\psi\in \Psi_{n,\cal E}$ avec l'hypoth\`ese faite en \ref{hypothese}; pour pouvoir traduire en termes combinatoires, on d\'ecompose la repr\'esentation $\psi$ de $W_{F}\times SL(2,{\mathbb C})\times SL(2,{\mathbb C})$ en repr\'esentations irr\'eductibles, d'o\`u un ensemble de triplets $(\rho,a,b)$ constituant par d\'efinition $Jord(\psi)$. On rappelle que $\psi\circ \Delta$ est la restriction de $\psi$ \`a $W_{F}$ fois la diagonale de $SL(2,{\mathbb C})$. Il est facile de v\'erifier que le centralisateur de $\psi$ dans $GL(n,{\mathbb C})$ est un sous-groupe du centralisateur de $\psi\circ \Delta$. On dit que $\psi$ est discret si le centralisateur de $\psi$ est un groupe fini et que $\psi$ est de restriction discr\`ete \`a la diagonale si le centralisateur de $\psi\circ \Delta$ est un groupe fini (notation d\'ej\`a utilis\'ee ici). On verra souvent $\psi\circ \Delta$ comme un \'el\'ement de $\Psi_{n,\cal E}^{temp}$.

\

On suppose d'abord que $\Psi$ est \'el\'ementaire, discret, c'est-\`a-dire par d\'efinition  que l'inclusion du centralisateur de $\psi$ dans le centralisateur de $\psi\circ \Delta$ est un isomorphisme et que ces groupes sont finis. En termes combinatoires cela se traduit par le fait que $Jord(\psi)$ est sans multiplicit\'es et que  pour tout $(\rho,a,b)\in Jord(\psi)$, $inf(a,b)=1$. On sait donc parfaitement d\'efinir $\pi^H(\psi\circ \Delta)$ gr\^ace \`a l'hypoth\`ese de \ref{hypothese}. Pour construire $\pi^H(\psi)$ en suivant \cite{elementaire}, on a besoin de la construction $inv_{\rho,\leq d}$, o\`u $\rho$ est une cuspidale autoduale irr\'eductible d'un groupe lin\'eaire et $d$ est un entier strictement sup\'erieur \`a 1 et o\`u pour toute repr\'esentation $\pi^H$ de $H_{n}$, on pose, ${\cal C}_{\rho,\leq d}$ l'ensemble des repr\'esentations cuspidales de la forme $\rho\vert\,\vert^x$ avec $x$ un r\'eel v\'erifiant $\vert x\vert \leq (d-1)/2$ et 
 $$
inv_{\rho,\leq d}:=\sum_{P^H} (-1)^{corang M^H} Ind_{P^H}^{H_{n}} Jac^{P^H}_{{\cal C}_{\rho,\leq d}} \pi^H.
$$
On donne une d\'efinition analogue pour $<d$ en rempla\c{c}ant partout $\leq $ par $<$.
On d\'efinit alors
$$
\pi^H(\psi):=\circ_{(\rho,a,b)\in Jord(\psi); sup(a,b)=b} \biggl(\beta(\psi,\rho, b)inv_{\rho,<b}\circ inv_{\rho,\leq b}\biggr)\pi^H(\psi\circ \Delta),
$$
o\`u $\beta(\psi,\rho, b)$ est un signe que l'on  pr\'ecisera ci-dessous  et o\`u on prend n'importe quel ordre pour faire les op\'erations; en \cite{elementaire}, on a bien montr\'e que l'ordre n'importait pas mais avec des hypoth\`eses un peu diff\'erentes (a priori au moins) de celles que l'on a ici. Comme on va montrer que quel que soit l'ordre $\pi^H(\psi)$ est stable et que  $\pi(\psi)$ en est un transfert, on aura bien l'ind\'ependance de l'ordre. Le signe n'est pas simple \`a interpr\'eter, il vaut $\beta(\psi,\rho,\leq b)\beta(\psi,\rho,\leq b')$ o\`u ces signes ont \'et\'e d\'efinis  en \ref{actiondetheta} avec o\`u $b'$ est le plus grand entier s'il existe (sinon le signe correspondant est +) tel qu'il existe $(\rho,a'',b'')\in Jord(\psi)$ avec $sup(a'',b'')=b'<b$.

Dans ce travail, nous n'avons besoin que de la d\'efinition mais pour motiver un peu rappelons le r\'esultat de \cite{elementaire}. En fait tel que cet article est r\'edig\'e on a suppos\'e un peu plus que seulement \ref{hypothese}; on a suppos\'e qu'il existe une bijection entre l'ensemble des caract\`eres cuspidaux du centralisateur d'un morphisme discret sans trou, $\phi$, (cf. \ref{remarquesurlhypothese}) et les repr\'esentations cuspidales incluses dans $\Pi^H(\phi)$. On pourrait facilement reformuler les r\'esultats pour \'eviter cette derni\`ere hypoth\`ese; il faut remplacer irr\'eductible par nul ou de longueur finie sans multiplicit\'es; ici on va continuer avec les hypoth\`eses de \cite{elementaire} qui sont pr\'ecis\'ement cette asssertion de bijectivit\'e, l'hypoth\`ese \ref{hypothese} et une r\'eponse positive \`a la question 2 de \ref{remarquesurlhypothese}.

En \cite{elementaire}, on a alors montr\'e que $\pi^H(\psi)$ est de la forme $\sum_{\epsilon}\epsilon(z_{2})\pi^H(\psi,\epsilon)$ o\`u $\epsilon$ parcourt l'ensemble des caract\`eres du centralisateur de $\psi$ de restriction triviale au centre de $^LH$, o\`u  $\pi^H(\psi,\epsilon)$ est une repr\'esentation irr\'eductible et $z_{2}$ est le centre de la 2e copie de $SL(2,{\mathbb C})$.

\

On suppose maintenant que $\psi$ est de restriction discr\`ete \`a la diagonale. On reprend les notations bien commodes $(\rho,A,B,\zeta)$ de tout ce travail pour d\'ecrire les \'el\'ements de $Jord(\psi)$ et on d\'efinit par r\'ecurrence $\pi^H(\psi)$, la r\'ecurrence porte sur $\sum_{(\rho,A,B,\zeta)}(A-B)$. On fixe $(\rho,A,B,\zeta)\in Jord(\psi)$ tel que $A>B$, a priori la construction d\'epend de ce choix mais comme ci-dessus \`a posteriori elle n'en d\'epend pas, donc on ne s'ennuie pas avec cela et on pose:$$
\pi^H(\psi):=\oplus_{C\in ]B,A]}(-1)^{A-C}<\rho\vert\,\vert^{\zeta C}, \cdots, \rho\vert\,\vert^{\zeta A}> \times Jac_{\rho\vert\,\vert^{\zeta (B+2)}, \cdots, \rho\vert\,\vert^{\zeta C}}\pi^H(\psi',(\rho,A,B+2,\zeta))$$
$$
\oplus (-1)^{[(A-B+1)/2]}\pi^H(\psi', (\rho,A,B+1,\zeta),(\rho,B,B,\zeta)).
$$
En \cite{paquet}, on a d\'ecrit cette repr\'esentation virtuelle comme combinaison lin\'eaire explicite de repr\'esentations index\'ees par les caract\`eres du centralisateur de $\psi$, sous certaines hypoth\`eses, le r\'esultat obtenu (qui n\'ecessite les m\^emes hypoth\`eses que dans le cas \'el\'ementaire) est que cette distribution s'\'ecrit elle aussi sous la forme $\sum_{\epsilon}\epsilon(z_{2})\pi^H(\psi,\epsilon)$ o\`u chaque $\pi^H(\psi,\epsilon)$ est une repr\'esentation en g\'en\'eral non irr\'eductible. Une d\'ecomposition pr\'ecise est donn\'ee en loc.cit. Ici comme ci-dessus, nous n'avons besoin que de la d\'efinition.

\

On supprime encore l'hypoth\`ese $\psi$ de restriction discr\`ete \`a la diagonale, en fixant $N>>n$ et un morphisme $\tilde{\psi}\in \Psi_{N,{\cal E}}$ dominant $\psi$ suivant la d\'efinition de \ref{general}. On a alors d\'efini en \ref{general} l'ensemble de repr\'esentations ${\cal E}_{\tilde{\psi}}$ et on d\'efinit $Jac_{{\cal E}_{\tilde{\psi}}}$ comme en loc. cit. (attention, il faut tenir compte de l'ordre comme indiqu\'e) et on pose:

$$
\pi^H(\psi):=Jac_{{\cal E}_{\tilde{\psi}}} \pi^H(\tilde{\psi}).
$$
En \cite{casgeneral}, on d\'ecrit un peu moins pr\'ecis\'ement que dans le cas de restriction discr\`ete \`a la diagonale cette repr\'esentation virtuelle, sous les m\^emes  hypoth\`eses que ci-dessus, on garde toujours le m\^eme r\'esultat que cette distribution s'\'ecrit sous la forme $\sum_{\epsilon}\epsilon(z_{2})\pi^H(\psi,\epsilon)$ o\`u chaque $\pi^H(\psi,\epsilon)$ est une repr\'esentation en g\'en\'eral non irr\'eductible qui peut m\^eme ici \^etre nulle.

\subsection {Preuve du transfert}
On fait ici l'hypoth\`ese de \ref{hypothese} et on reprend les d\'efinitions de \ref{constructiondesrepresentations}; on garde  l'action de $\theta$ mise dans tout ce travail sur $\pi(\psi)$ pour tout $\psi\in \Psi_{n}$. On a en vue le th\'eor\`eme suivant:

\

\bf Th\'eor\`eme. \sl Les repr\'esentations $\pi(\psi)$ et $\pi^H(\psi)$ sont stables et $\psi(\psi)$ est un transfert de $\pi^H(\psi)$.\rm

Remarquons que l'assertion de stabilit\'e est un corollaire imm\'ediat de \ref{proprietesgeneralesdutransfertstable} et du fait facile avec les formules explicites  que l'induction respecte la stabilit\'e.

Ici on admet le th\'eor\`eme dans le cas \'el\'ementaire, le cas g\'en\'eral se ram\`ene au cas de restriction discr\`ete \`a la diagonale gr\^ace \`a \ref{general} et \`a \ref{proprietesgeneralesdutransfertstable}. Le cas de restriction discr\`ete \`a la diagonale r\'esulte de \ref{quasigeneral} et \ref{proprietesgeneralesdutransfertstable}. Il n'y a donc \`a d\'emontrer le th\'eor\`eme que dans le cas \'el\'ementaire.

Supposons donc que $\psi$ est \'el\'ementaire.
On veut montrer que $\pi(\psi)$ est un transfert de $\pi^H(\psi)$, sous l'hypoth\`ese que $\pi(\psi\circ \Delta)$ est un transfert de $\pi^H(\psi\circ \Delta)$. Si $\psi$ est trivial sur la 2e copie de $SL(2,{\mathbb C})$, il n'y a donc rien \`a d\'emontrer. On fait donc la preuve par r\'ecurrence sur le nombre d'\'el\'ements de $Jord(\psi)$, $(\rho,a,b)$ tel que $sup(a,b)=b>1$. On suppose que $\psi$ est de restriction discr\`ete \`a la diagonale et on fixe $(\rho,a,b)\in Jord(\psi)$ tel que $sup(a,b)=b>1$. On a d\'efini $\psi^\sharp$ en \ref{exactitudeducomplexe} en prenant $d=b$. Puis on passe de $\psi^\sharp$ \`a $(\psi^\sharp)^\sharp$ en prenant maitenant $d$ le plus grand entier, not\'ee $d_{0}$, strictement inf\'erieur \`a $b$ et tel qu'il existe $d'$ avec soit $(\rho,d,d')\in Jord(\psi)$ soit $(\rho,d',d)\in Jord(\psi)$; si ce nombre n'existe pas, on pose $(\psi^\sharp)^\sharp=\psi^\sharp$. On a que $Jord((\psi^\sharp)^\sharp)$ co\"{\i}ncide avec $Jord(\psi)$ sauf que $(\rho,a,b)$ a \'et\'e chang\'e en $(\rho,b,a)$. Par r\'ecurrence, on peut donc admettre que $\pi((\psi^\sharp)^\sharp)$ est un transfert de $\pi^H((\psi^\sharp)^\sharp)$, attention avec l'action que nous avons mise sur ces repr\'esentations. En tenant compte de ce qu'\^etre un transfert commute \`a l'induction est \`a la restriction, on voit que $\pi(\psi^\sharp)$ est un transfert de $inv_{\rho,\leq d_{0}}\pi^H((\psi^\sharp)^\sharp)$ si $\beta(\psi\circ \Delta,\rho,\leq d_{0})=+1$ et de son oppos\'e sinon. Puis $\pi(\psi)$ est un transfert de $inv_{\rho,\leq b}\circ inv_{\rho,\leq d_{0}}\pi^H((\psi^\sharp)^\sharp)$ si $\beta(\psi\circ \Delta,\rho,\leq b)\beta(\psi\circ \Delta,\rho, \leq d_{0})=1$ et de son oppos\'e sinon. En tenant compte de ce signe, on a exactement $\pi^H(\psi)$ tel qu'on l'a d\'efini dans \ref{constructiondesrepresentations}. Cela termine donc la preuve.
\section{Normalisation de l'action de $\theta$\label{normalisationdelactiondetheta}}
\subsection{D\'efinition de la normalisation de Whittaker \label{definitiondelanormalisationdewhittaker}}
On fixe $\delta$ un caract\`ere  de $F$; si $\psi$ est temp\'er\'e, la repr\'esentation $\pi(\psi)$ a un 
mod\`ele de Whittaker et l'action de $\theta$ est celle qui est normalis\'ee  \`a la Whittaker. 
On consid\`ere le morphisme de $$W_{F}\times SL(2,{\mathbb C}) \rightarrow W_{F}\times SL(2,{\mathbb C}) \times SL(2,{\mathbb C})$$
$$
(w\in W_{F}, h\in SL(2,{\mathbb C}))\mapsto (w,h,
\biggl( \begin{matrix} &\vert w\vert^{1/2}& 0\\
&0&\vert\,\vert^{-1/2}\end{matrix}\biggr))$$ et en composant on obtient un morphisme, $\psi^2$, de $W_{F}\times SL(2,{\mathbb C})$ dans $GL(n,{\mathbb C})$. A ce morphisme correspond 
une repr\'esentation par la classification de Langlands, qui n'est autre que $\pi(\psi)$; mais dans la construction de Langlands, $\pi(\psi)$ est quotient d'une induite \`a partir d'une repr\'esentation 
temp\'er\'ee. Il est facile de voir que cette induite est $\theta$-invariante. On normalise l'action de $\theta$ sur l'induite en utilisant les fonctionnelles de Whittaker (cf. introduction).

Par restriction ce choix donne une action sur $\pi(\psi)$. Pour $z=\pm$, on note $\pi(\psi)^z$ la repr\'esentation $\pi(\psi)$ munie de cette action de $\theta$ si $z=+$ et de l'action oppos\'ee si $z=-$. Cette normalisation, facile \`a pr\'esenter n'est pas tr\`es simple non plus mais elle a le gros avantage par rapport \`a celle utilis\'ee ici d'\^etre d'origine globale. Localement il introduit des signes qu'il nous faut de toute fa\c{c}on calculer. On note $\pi(\psi)_{W}$ la repr\'esentation $\pi(\psi)$ avec l'action de $\theta$ telle que l'on vient de la d\'efinir et on note $\theta_{W}(\psi)$ cette action.

\subsection{Explicitation de la normalisation de Whittaker\label{explicitationdelanormalisationdewhittaker}}
On \'ecrit l'induite dont $\pi(\psi)$ est l'unique sous-module irr\'eductible (on transforme quotient de Langands en sous-module sans probl\`eme). Pour cela, pour tout $k$ demi-entier  positif ou nul, on 
d\'efinit $Jord(\psi)_{\geq k}:=\{(\rho,a,b)\in Jord(\psi)$; $(b-1)/2-k$ est un entier positif ou  nul$\}$. Ainsi $\pi(\psi)$ est l'unique sous-module irr\'eductible de l'induite:
$$
\biggl(\times_{k}\times_{(\rho,a,b)\in Jord(\psi)_{\geq k}}St(a,\rho)\vert\,\vert^{-k}\biggr)\times_{(\rho,a,b)\in Jord(\psi)_{\geq 0}}St(a,\rho) \biggl(\times_{k }\times _{(\rho,a,b)\in Jord(\psi)_{\geq k}}St(a,\rho)\vert\,\vert^{k}\biggr),\eqno(1)$$
o\`u les $k$ sont d'abord pris dans l'ordre d\'ecroissant puis dans l'ordre croissant. Et l'action  $\theta_{W}$ sur $\pi(\psi)$ normalis\'ee ''\`a la Whittaker'' se d\'eduit de l'action de $\theta$ sur cette induite normalis\'ee \`a la  Whittaker. On remarque que l'induite du milieu, correspondant \`a $k=0$, n'est autre que $\prod_{(\rho,a,b)\in Jord(\psi); (-1)^b=-1}St(a,\rho)$; on la note de fa\c{c}on coh\'erente avec ce qui pr\'ec\`ede $\pi(\psi^2_{imp})$. L'action de $\theta$ sur cette repr\'esentation est avec les notations d\'ej\`a introduites $\theta_{W}(\psi^2_{imp})$.

\

\bf Lemme. \sl L'action  $\theta_{W}(\psi)$ sur (1) est celle qui se d\'eduit canoniquement de l'action de $\theta$ sur $\pi(\psi^2_{imp})$.\rm

\

Notons $V'$ le sous-espace de l'induite (1) form\'e des fonctions \`a support dans la grosse cellule $PwB$ o\`u $w$ est l'\'el\'ement du groupe de Weyl de longueur maximal ($P$ est le parabolique qui sert dans l'induction et $B$ le Borel standard). La restriction d'une fonctionnelle de Whittaker \`a $V'$ se fait par l'int\'egration, o\`u il faut fixer une fonctionnelle de Whittaker, $\ell_{M}$ sur la repr\'esentation que l'on induit, $\theta$ invariante (ce qui est possible) $$f\in V' \mapsto \int_{U\cap M\backslash U}\ell^M(f(wu))\chi(u)^{-1}\, du.$$L'action canonique de $\theta$ laisse stable $V'$ et cette forme lin\'eaire, ce qui prouve le lemme.

\

On peut encore expliciter un peu plus; on note $\psi_{imp,imp}$ le morphisme dont les blocs de Jordan sont les $(\rho,a-2[a/2],b-2[b/2])$ pour $(\rho,a,b)$ parcourant $Jord(\psi)$ (avec les m\^emes multiplicit\'es). On a encore une inclusion:
$$
\pi(\psi^2_{imp}) \hookrightarrow $$
$$\biggl(\times_{k\geq 1/2}\times _{\begin{matrix}(\rho,a,b)\in Jord(\psi);\\ (-1)^b=-1;\\
(a-1)/2-k\in {\mathbb Z}_{\geq 0}\end{matrix}}\rho\vert\,\vert^k\biggr) \times \pi(\psi_{imp,imp})\times \biggl(\times _{\begin{matrix}(\rho,a,b)\in Jord(\psi);\\ (-1)^b=-1;\\
(a-1)/2-k\in {\mathbb Z}_{\geq 0}\end{matrix}}\rho\vert\,\vert^{-k}\biggr), \eqno(2)
$$
o\`u dans le premier produit les $k$ sont les demi-entiers, $>0$,  pris dans l'ordre d\'ecroissant et dans le deuxi\`eme produit les $k$ sont les demi-entiers, $>0$, pris dans l'ordre croissant. On obtient alors avec (1) une s\'erie d'inclusion, pour des repr\'esentations $\sigma, \sigma'$ convenables qui sont explicit\'ees dans (1) et (2)
$$
\pi(\psi)\hookrightarrow \sigma \times \pi(\psi^2_{imp}) \times \, ^\theta\sigma \hookrightarrow \sigma\times \sigma' \times \pi(\psi_{imp,imp}) \times \, ^\theta\sigma' \times \, ^\theta \sigma. \eqno(3)
$$
En suivant les d\'efinitions $\theta_{W}(\psi)$ est la restriction \`a $\pi(\psi)$ de l'action de $\theta$ sur le dernier membre de (3) qui se d\'eduit canoniquement de $\theta_{W}(\pi(\psi)_{imp,imp})$.

\subsection{D\'efinition de la normalisation unipotente\label{definitiondelanormalisationunipotente}}
Dans la d\'efinition de la normalisation pour l'action de $\theta$ \`a la Whittaker, le c\^ot\'e temp\'er\'e, c'est-\`a-dire la premi\`ere copie de $SL(2,{\mathbb C})$ a \'et\'e privil\'egi\'e. D'un point de vue local, il n'y a pas de raison \`a ce choix; pour faire le choix oppos\'e, il faut d'abord traiter le cas des repr\'esentations ''unipotentes'' c'est-\`a-dire les duales au sens de l'involution de Zelevinsky des repr\'esentations temp\'er\'ees. Pour celles-l\`a, on prend la normalisation \`a la Whittaker, c'est-\`a-dire qu'on repr\'esente cette repr\'esentation comme quotient de Langlands d'une s\'erie principale g\'en\'eralis\'ee qui elle a un mod\`ele de Whittaker. On fixe l'action de $\theta$ sur cette s\'erie principale en demandandant que $\theta$ induise l'action naturelle  sur le mod\`ele de Whittaker. Fixons maintenant $\psi$ g\'en\'eral et on consid\`ere le morphisme de $W_{F}$ dans la premi\`ere copie de $SL(2,{\mathbb C})$ analogue \`a celui d\'efini en \ref{definitiondelanormalisationdewhittaker} que l'on note $i_{1}$ et  on obtient un morphisme, $\psi_{unip}$ de $W_{F}\times SL(2,{\mathbb C})$ qui sur $W_{F}$ co\"{\i}ncide avec $\psi\circ i_{1}$ sur $W_{F}$ et avec la restriction de $\psi$ \`a la deuxi\`eme copie de $SL(2,{\mathbb C})$ sur $SL(2,{\mathbb C})$. Ce morphisme $\psi_{unip}$ d\'efinit une induite \`a partir d'une repr\'esentation de Speh d'un Levi, c'est la classification originelle de Zelevinsky. Explicitons; on pose ici pour distinguer de \ref{explicitationdelanormalisationdewhittaker} $Jord^a(\psi)_{\geq k}:=\{(\rho,a,b)\in Jord(\psi); (a-1)/2- k\in {\mathbb Z}_{\geq 0}\}$  pour tout $k$ demi-entier strictement positif et on pose $\psi_{imp}^1$ le morphisme dont les blocs de Jordan sont les $(\rho,a-2 [a/2],b)$ pour tout $(\rho,a,b)\in Jord(\psi)$ avec $a$ impair. Alors
$$
\pi(\psi)\hookrightarrow \biggl( \times_{k\geq 1/2}\times_{(\rho,a,b)\in Jord^a(\psi)_{\geq k}}Sp(b,\rho)\vert\,\vert^{-k}\biggr) \times \pi(\psi_{imp}^1)\times  \biggl(\times_{k\geq 1/2}\times_{(\rho,a,b)\in Jord^a(\psi)_{\geq k}}Sp(b,\rho)\vert\,\vert^{k}\biggr),\eqno(1)
$$
o\`u les $k$ sont d'abord pris dans l'ordre d\'ecroissant puis dans l'ordre croissant. Comme on a mis une action de $\theta$ sur $\pi(\psi^1_{imp})$, on prolonge canoniquement cette action de $\theta$ sur (1). Par restriction on obtient une action de $\theta$ sur $\pi(\psi)$ que l'on note $\theta_{U}(\psi)$.

Comme pour $\theta_{W}(\psi)$ on peut retrouver cette action par une suite d'inclusion analogue \`a (3) de \ref{explicitationdelanormalisationdewhittaker}; pour cela on consid\`ere l'inclusion
$$
\pi(\psi^1_{imp})\hookrightarrow$$
$$ \times_{k\geq 1/2} \times_{\begin{matrix}(\rho,a,b)\in Jord(\psi)\\(-1)^a=-1; \\(b-1)/2-k\in {\mathbb Z}_{\geq 0}
\end{matrix}}\rho\vert\,\vert^{-k} \times \pi(\psi_{imp,imp}) \times_{k\geq 1/2}\times_{\begin{matrix}(\rho,a,b)\in Jord(\psi)\\(-1)^a=-1; \\(b-1)/2-k\in {\mathbb Z}_{\geq 0}
\end{matrix}}\rho\vert\,\vert^{k}, \eqno(2)
$$D'o\`u pour des bons choix de $\tau,\tau'$
$$
\pi(\psi)\hookrightarrow \tau \times \pi(\psi^1_{imp}) \times \, ^\theta\tau \hookrightarrow \tau \times \tau' \times \pi(\psi_{imp,imp}) \times \, ^\theta \tau' \times ^\theta \tau.\eqno(3)
$$
Et $\theta_{U}(\psi)$ s'obtient en restreignant l'action  de $\theta$ sur le dernier membre de (3) obtenu en prolongeant canoniquement $\theta_{W}(\psi_{imp,imp})$. C'est bien $\theta_{W}(\psi_{imp,imp})$ car $\pi(\psi_{imp,imp})$ est une induite de cuspidales et les choix se voient sur l'espace de Whittaker relatif au caract\`ere additif que nous avons fix\'e. Pour garder la sym\'etrie on pose $$\theta(\psi_{imp,imp}):=\theta_{U}(\psi_{imp,imp}):=\theta_{W}(\psi_{imp,imp}).$$
\subsection{Une propri\'et\'e commune \`a ces actions\label{uneproprietecommuneacesactions}}
\bf Proposition. \sl 
Soit $\psi$ et soit $(\rho,a,b)\in Jord(\psi)$.

(i) Supposons que $b\geq 2$; on note $\psi'$ le morphisme qui se d\'eduit de $\psi$ en rempla\c{c}ant $(\rho,a,b)$ par $(\rho,a,b-2)$. On suppose qu'il existe une inclusion:
$$
\pi(\psi)\hookrightarrow St(a,\rho)\vert\,\vert^{-(b-1)/2}\times \pi(\psi') \times St(a,\rho)\vert\,\vert^{(b-1)/2}.\eqno(1)
$$
Alors $\theta_{W}(\psi)$ est par restriction l'action $\theta_{W}(\psi')$ \'etendue canoniquement \`a l'induite de droite.

(ii) Supposons que $a\geq 2$; on note $\psi'$ le morphisme qui se d\'eduite de $\psi$ en rempla\c{c}ant $(\rho,a,b)$ par $(\rho,a-2,b)$. On suppose qu'il existe une inclusion:
$$
\pi(\psi) \hookrightarrow Sp(b,\rho)\vert\,\vert^{(a-1)/2}\times \pi(\psi') \times Sp(b,\rho)\vert\,\vert^{-(a-1)/2}.\eqno(2)
$$
Alors $\theta_{U}(\psi)$ est par restriction l'action $\theta_{U}(\psi')$ \'etendue canoniquement \`a l'induite de droite.\rm

\

On montre par exemple (ii); on utilise l'inclusion de $\pi(\psi')$ dans \ref{definitiondelanormalisationunipotente} (1) (o\`u on remplace $\psi$ par $\psi')$); on remarque que $\psi^1_{imp}=(\psi')^1_{imp}$  d'o\`u avec l'hypoth\`ese de (2) une inclusion
$$
\pi(\psi)\hookrightarrow Sp(b,\rho)\vert\,\vert^{(a-1)/2}\times \biggl(\times_{k\geq 1/2}\times _{(\rho',a',b')\in Jord(\psi'); (a'-1)/2-k\in {\mathbb Z}_{\geq 0}}Sp(\rho',b')\vert\,\vert^{k}\biggr)\times \pi(\psi^1_{imp}) $$
$$\times \biggl(\times _{k\geq 1/2}\times_{(\rho',a',b')\in Jord(\psi'); (a'-1)/2-k\in {\mathbb Z}_{\geq 0}}Sp(b',\rho)\vert\,\vert^{-k}\biggr)\times Sp(b,\rho)\vert\,\vert^{-(a-1)/2}.\eqno(3)
$$
Cette inclusion ressemble \`a l'inclusion \ref{definitiondelanormalisationunipotente} (1) pour $\pi(\psi)$ \`a ceci pr\`es que $Sp(b,\rho)\vert\,\vert^{(a-1)/2}$ et $Sp(b,\rho)\vert\,\vert^{-(a-1)/2}$ ne sont pas \`a leur bonne place. On sait quand m\^eme que $\pi(\psi)$ intervient avec multiplicit\'e exactement 1 en tant que sous-quotient  \`a la fois dans (3) et dans \ref{definitiondelanormalisationunipotente} (1). On construit un entrelacement de (3) vers \ref{definitiondelanormalisationunipotente} (1), en rempla\c{c}ant d'abord $Sp(b,\rho)^{(a-1)/2}$ par $Sp(b,\rho)\vert\,\vert^{(a-1)/2+s}$ et $Sp(b,\rho)\vert\,\vert^{-(a-1)/2}$ par $Sp(b,\rho)\vert\,\vert^{-(a-1)/2-s}$; pour cela on prend un op\'erateur d'entrelacement standard qui 
d\'epend m\'ero\-morphi\-quement de $s$. On le note $M(s)$ et on le multiplie par une fonction 
m\'eromorphe de $s$ de sorte qu'il deviennent holomorphe non nul en $s=0$. On note $M_{0}$ la valeur de cet op\'erateur en $s=0$. Comme $\pi(\psi)$ est l'unique sous-module irr\'eductible de \ref{definitiondelanormalisationunipotente} (1), l'image de $M_{0}$ contient $\pi(\psi)$ et envoie donc $\pi(\psi)$ le sous-module de (3) sur $\pi(\psi)$ le sous-module de \ref{definitiondelanormalisationunipotente} (1). Il est facile de v\'erifier que $M_{0}$ entrelace l'action de $\theta$ sur (3) qui prolonge canoniquement une action fix\'ee de $\theta$ sur $\pi(\psi^1_{imp})$ en l'action de $\theta$ sur (1) qui prolonge canoniquement cette m\^eme action. D'o\`u l'assertion de (ii).

\subsection{Comparaison des normalisations \label{comparaisondesnormalisations}}
\subsubsection{Le cas o\`u $\vert Jord(\psi)\vert=1$}
\bf Lemme. \sl Supposons que  $Jord(\psi)$ a un seul \'el\'ement $(\rho,a,b)$; dans ce cas $\theta_{W}(\psi)=\theta_{U}(\psi)$.
\rm

\

Si $inf(a,b)=1$ cela r\'esulte des d\'efinitions. On d\'emontre l'assertion par r\'ecurrence. Supposons donc que $inf(a,b)>1$. On a une inclusion
$
\pi((\rho,a,b))\hookrightarrow$
$$ St(a,\rho)\vert\,\vert^{-(b-1)/2}\times \pi(\rho,a,b-2) \times St(a,\rho)\vert\,\vert^{(b-1)/2}
$$
qui d'apr\`es \ref{uneproprietecommuneacesactions} entrelace $\theta_{W}((\rho,a,b))$ et le prolongement canonique de $\theta_{W}((\rho,a,b-2))$. L'hypoth\`ese de r\'ecurrence dit que $\theta_{W}((\rho,a,b-2))=\theta_{U}((\rho,a,b-2))$. On a encore une inclusion:
$
\pi((\rho,a,b-2))\hookrightarrow $
$$Sp(b-2,\rho)\vert\,\vert^{(a-1)/2}\times \pi((\rho,a-2,b-2)) \times Sp(b-2,\rho)\vert\,\vert^{-(a-1)/2}.
$$
D'apr\`es \ref{uneproprietecommuneacesactions} cette inclusion entrelace $\theta_{U}((\rho,a,b-2))$ et le prolongement canonique de $\theta_{U}((\rho,a-2,b-2))$. D'o\`u une inclusion
$
\pi((\rho,a,b))\hookrightarrow $
$$St(a,\rho)\vert\,\vert^{-(b-1)/2}\times Sp(b-2,\rho)\vert\,\vert^{(a-1)/2}\times \pi((\rho,a-2,b-2)) \times Sp(b-2,\rho)\vert\,\vert^{-(a-1)/2}\times St(a,\rho)\vert\,\vert^{(b-1)/2},\eqno(4)
$$
qui entrelace $\theta_{W}((\rho,a,b))$ avec le prolongement canonique de $\theta_{U}((\rho,a-2,b-2))$. En travaillant sym\'etriquement, on montre une inclusion 
$
\pi((\rho,a,b))\hookrightarrow$
$$Sp(b,\rho)\vert\,\vert^{(a-1)/2}\times St(a-2,\rho)\vert\,\vert^{-(b-1)/2}\times \pi((\rho,a-2,b-2)) \times St(a-2,\rho)\vert\,\vert^{(b-1)/2}\times Sp(b,\rho)\vert\,\vert^{-(a-1)/2}\eqno(5)
$$
qui entrelace $\theta_{U}((\rho,a,b))$ avec le prolongement canonique de $\theta_{W}(\rho,a-2,b-2)$; si $(a-2)(b-2)=0$ cette action est l'action triviale. Par r\'ecurrence on sait encore que $\theta_{U}((\rho,a-2,b-2))=\theta_{W}((\rho,a-2,b-2))$. On remarque que 
$
St(a,\rho)^{-(b-1)/2}\hookrightarrow 
\rho\vert\,\vert^{(a-b)/2}\times St(a-2,\rho)^{-(b-1)/2}\times \rho\vert\,\vert^{-(a+b)/2+1}
$
et que 
$
Sp(b,\rho)\vert\,\vert^{(a-1)/2}\hookrightarrow \rho\vert\,\vert^{(a-b)/2}\times Sp(b-2,\rho)\vert\,\vert^{(a-1)/2}\times \rho\vert\,\vert^{(a+b)/2-1}.
$
L'induite $St(a-2,\rho)\vert\,\vert^{-(b-1)/2}$ est un sous-module de l'induite $\times_{\ell \in [(a-b)/2-1,-(a+b)/2+2}\rho\vert\,\vert^{\ell}$ tandis que $Sp(b-2,\rho)\vert\,\vert^{(a-1)/2}$ est un sous-module de l'induite $\times_{\ell \in [(a-b)/2+1,(a+b)/2-2]}\rho\vert\,\vert^\ell$. D'o\`u les isomorphismes (cf. la preuve de \ref{independance}):
$$
\rho\vert\,\vert^{(a+b)/2-1} \times St(a-2,\rho)\vert\,\vert^{-(b-1)/2}\simeq St(a-2,\rho)\vert\,\vert^{-(b-1)/2}\times \rho\vert\,\vert^{(a+b)/2-1};$$
$$   \rho\vert\,\vert^{-(a+b)/2+1}\times Sp(b-2,\rho)\vert\,\vert^{(a-1)/2}\simeq Sp(b-2,\rho)\vert\,\vert^{(a-1)/2}\times \rho\vert\,\vert^{-(a+b)/2+1};$$
$$St(a-2,\rho)^{-(b-1)/2}\times Sp(b-2,\rho)\vert\,\vert^{(a-1)/2}\simeq Sp(b-2,\rho)\vert\,\vert^{(a-1)/2}\times St(a-2,\rho)^{-(b-1)/2}.
$$

On peut donc  remplacer (5) par l'inclusion
$
\pi((\rho,a,b))\hookrightarrow$
$$ \rho\vert\,\vert^{(a-b)/2}\times Sp(b-2,\rho)\vert\,\vert^{(a-1)/2}\times St(a-2,\rho)\vert\,\vert^{-(b-1)/2}\times \rho\vert\,\vert^{(a+b)/2-1}\times \pi((\rho,a-2,b-2))$$
$$ \times  \rho\vert\,\vert^{-(a+b)/2+1}\times   St(a-2,\rho)\vert\,\vert^{(b-1)/2}\times Sp(b-2,\rho)\vert\,\vert^{-(a-1)/2}\times \rho\vert\,\vert^{-(a-b)/2} \eqno(5)'
$$
et remplacer (4) par l'inclusion
$
\pi((\rho,a,b))\hookrightarrow$
$$ \rho\vert\,\vert^{(a-b)/2}\times St(a-2,\rho)\vert\,\vert^{(b-1)/2}\times Sp(b-2,\rho)\vert\,\vert^{-(a-1)/2}\times \rho\vert\,\vert^{-(a+b)/2+1}\times \pi((\rho,a-2,b-2))$$
$$ \times  \rho\vert\,\vert^{(a+b)/2-1}\times  Sp(b-2,\rho)\vert\,\vert^{(b-1)/2}\times St(a-2,\rho)\vert\,\vert^{-(a-1)/2}\times \rho\vert\,\vert^{-(a-b)/2}. \eqno(4)'
$$
Et ces inclusions entrelacent pour (4)$'$, $\theta_{W}((\rho,a,b))$ avec le prolongement canonique de $\theta((\rho,a-2,b-2))$ et pour (5)$'$, $\theta_{U}((\rho,a,b))$ avec le prolongement canonique de $\theta((\rho,a-2,b-2))$. On peut encore \'echanger dans (4)$'$ $St(a-2,\rho)\vert\,\vert^{-(b-1)/2}\times Sp(b-2,\rho)\vert\,\vert^{(a-1)/2}$ en $Sp(b-2,\rho)\vert\,\vert^{(a-1)/2}\times St(a-2,\rho)\vert\,\vert^{-(b-1)/2}$ et dualement
$Sp(b,\rho)\vert\,\vert^{-(a-1)/2}\times St(a-2,\rho)\vert\,\vert^{(b-1)/2}$ en $St(a-2,\rho)\vert\,\vert^{(b-1)/2}\times Sp(b-2,\rho)\vert\,\vert^{-(a-1)/2}$ puis  remarquer que $\rho\vert\,\vert^{(a+b)/2-1}\times \pi((\rho,a-2,b-2)) \times \rho\vert\,\vert^{-(a+b)/2+1}$ est irr\'eductible; cela r\'esulte de ce que les repr\'esentations $\rho\vert\,\vert^\ell$ dans le support cuspidal de $\pi((\rho,a-2,b-2))$ v\'erifient $(a+b)/2-1 -\ell >1$. On peut donc encore \'echanger et transformer le deuxi\`eme membre de (4)$'$en celui de 
(5)$'$par des entrelacements qui entrelacent les prolongements canoniques de $\theta((\rho,a-2,b-2))$ (cf. \ref{entrelacement}). Pour pouvoir conclure que $\theta_{W}((\rho,a,b))=\theta_{U}((\rho,a,b))$, il faut encore remarquer que (5)$'$admet $\pi((\rho,a,b))$ comme sous-module irr\'eductible avec multiplicit\'e 1; supposons qu'il existe une inclusion de $\Pi:=\pi((\rho,a,b))\oplus \pi((\rho,a,b))$ dans le deuxi\`eme membre de (5)$'$. On reprend les arguments en sens inverse; on \'ecrit une suite exacte:
$$
0\rightarrow Sp(b,\rho)\vert\,\vert^{(a-1)/2}\times St(a-2,\rho)\vert\,\vert^{-(b-1)/2}\rightarrow$$
$$ \rho\vert\,\vert^{(a-b)/2}\times Sp(b-2,\rho)\vert\,\vert^{(a-1)/2}\times St(a-2,\rho)\vert\,\vert^{-(b-1)/2}\times \rho\vert\,\vert^{(a+b)/2-1}$$
$$\rightarrow \tau \rightarrow 0,
$$
ce qui d\'efinit $\tau$. On v\'erifie qu'il existe $x\neq (a-b)/2$ tel que $Jac_{x}\tau\neq 0$. On consid\`ere l'image de $\Pi$ dans $$\tau\times \pi((\rho,a-2,b-2)) \times   \rho\vert\,\vert^{-(a+b)/2+1}\times  St(a-2,\rho)\vert\,\vert^{(b-1)/2}\times Sp(b-2,\rho)\vert\,\vert^{-(a-1)/2}\times \rho\vert\,\vert^{-(a-b)/2} .$$ Si elle \'etait non nulle, pour un $x$ comme ci-dessus, on aurait $Jac_{x}\pi((\rho,a,b))\neq 0$, ce qui est impossible. On montre de fa\c{c}on analogue que $\Pi$ doit \^etre inclus dans 
$$
Sp(b,\rho)\vert\,\vert^{(a-1)/2}\times St(a-2,\rho)\vert\,\vert^{-(b-1)/2}\times \pi((\rho,a-2,b-2)) \times St(a-2,\rho)\vert\,\vert^{(b-1)/2}\times Sp(b,\rho)\vert\,\vert^{-(a-1)/2}.
$$
Un calcul de module de Jacquet maintenant imm\'ediat montre que cette induite n'a qu'un unique sous-module irr\'eductible. D'o\`u a fortiori notre assertion qui termine la preuve.

\subsubsection{Action de $\theta$ et op\'erateurs d'entrelacement normalis\'es 1.\label{operateursdentrelacement1}}
Pour pouvoir comparer les actions $\theta_{W}(\psi)$ et $\theta_{U}(\psi)$ on va les comparer \`a des actions obtenues avec des op\'erateurs d'entrelacement normalis\'es et on va montrer que le signe dont diff\`erent ces 2 actions est un quotient de facteurs de normalisation.

 D'une fa\c{c}on tr\`es g\'en\'erale, soit $\delta,\delta'$ des repr\'esentations $\theta$-invariantes; on consid\`ere un op\'erateur d'entrelacement, pour $s\in {\mathbb C}$ voisin de 0 :
$$
N(s): \delta \vert\,\vert^s \times \delta' \rightarrow \delta' \times \delta\vert\,\vert^s,
$$
d\'ependant m\'eromorphiquement de $s$; exactement on note $M(s)$ l'op\'erateur d'entrelacement standard (cf. \ref{standard}) et on suppose qu'il existe une fonction m\'eromorphe de $s$, $r(s)$ tel que $N(s)$ soit l'op\'erateur $r(s)M(s)$. On suppose que $N(s)$ est holomorphe en $s=0$. On suppose aussi que $\delta$ et $\delta'$ sont munis d'une action de $\theta$, not\'ee $\theta_{\delta}$ et $\theta_{\delta'}$. On peut alors d\'efinir une action de $\theta$ sur l'induite $\delta\times \delta'$, en posant pour toute section $f$ \`a valeurs dans $\delta \otimes \delta'$ et pour tout \'el\'ement $g$ du groupe:
$$
\theta. f(g):= (\theta_{\delta}\otimes \theta_{\delta'}\circ inv)   (N(0)f)(\theta(g))=: (\Theta\circ N(0))f \eqno(1)
$$
o\`u $inv$ est l'\'echange $\delta'\otimes \delta$ dans $\delta\otimes \delta'$.  Pour que $\theta^2.f$ soit \'egal  \`a $1$, il faut que l'on ait $$\biggl((\Theta\circ N(-s)\circ  \Theta) \circ N(s)\biggr)_{s=0}=1.$$ Pour $Re(-s)>>0$, on v\'erifie que pour l'op\'erateur d'entrelacement standard
$
\Theta \circ M(-s)\circ \Theta$ est l'op\'erateur d'entrelacement standard:
$$
M'(-s): \quad \delta'\times \delta\vert\,\vert^s \rightarrow \delta\vert\,\vert^s \times \delta';
$$
cela se fait en \'ecrivant explicitement l'int\'egrale. La condition devient donc $$\biggl(r(s)r(-s)M'(-s)M(s)\biggr)_{s=0}=1.$$ Quand on repr\'esente l'induite $\delta\vert\,\vert^s\times \delta'$ dans un espace 
ind\'ependant de $s$, le produit $M'(-s)\circ M(s)$ est une homoth\'etie avec un coefficient qui est une fonction m\'eromorphe de $s$, $c(s):=M'(-s)\circ M(s)$.Pour avoir une action de $\theta$ avec la d\'efinition ci-dessus, suffit donc que $r(s)r(-s)=c(s)^{-1}$. Des choix de $r(s)$, plusieurs sont possibles; on consid\`ere d'abord celui sugg\'er\'e par Langlands et calcul\'e explicitement par Shahidi \cite{shahidi}, \`a ceci pr\`es que l'on modifiera les facteurs $\epsilon$. On pose donc:
$$
N_{L}(s):=E(\delta,\delta',s) L(\delta\times \delta',s)/L(\delta \times \delta',1+s),
$$
o\`u $E(s):=\prod_{\tau,\tau'}\epsilon(\tau\times \tau',s)$; dans ce produit $\tau$ parcourt l'ensemble des repr\'esentations cuspidales formant le support cuspidale de $\delta$ et $\tau'$ parcourt l'ensemble analogue pour $\delta'$.

On applique cela \`a $\pi(\psi)$ de la fa\c{c}on suivante: on fixe $(\rho,a,b)\in Jord(\psi)$ tel que $b$ soit maximum. On note $\psi'$ le morphisme qui se d\'eduit de $\psi$ en enlevant ce triplet; on a d\'ej\`a d\'efini  $\pi((\rho,a,b))$ et on a l'isomorphisme:
$$
\pi(\psi)=\pi((\rho,a,b))\times \pi(\psi').
$$On a d\'emontr\'e en \cite{mw} un certain nombre de propri\'et\'es pour les op\'erateurs d'entrelacement normalis\'es \`a la Langlands-Shahidi:
$$
N_{L}(s): \pi(\rho,a,b)\vert\,\vert^s \times \pi(\psi') \rightarrow \pi(\psi')\times \pi(\rho,a,b)\vert\,\vert^s.
$$
La premi\`ere propri\'et\'e est une formule de produit; par exemple quand on inclut $\pi(\rho,a,b)$ dans $St(a,\rho)^{(b-1)/2}\times \pi(\rho,a-2,b-2) \times St(a,\rho)^{-(b-1)/2}$, $N_{L}(s)$ est un produit, ici, de 3 op\'erateurs d'entrelacements normalis\'es \`a la Langlands-Shahidi. C'est facile car cela r\'esulte uniquement la formule de produit \'evidente pour les facteurs de normalisations introduits. On utilisera aussi la formule de produit suivant la d\'ecomposition $\pi(\psi')=\times_{(\rho',a',b')\in Jord(\psi')}\pi(\rho',a',b')$.

La deuxi\`eme propri\'et\'e est nettement plus difficile; on dit en suivant \cite{mw} que deux repr\'esentations $\pi(\rho,a,b)\vert\,\vert^s$ ($s$ r\'eel) et $\pi(\rho',a',b')$ sont li\'ees si $\rho=\rho'$ et $$
\vert (a-a')/2\vert +\vert (b-b')/2\vert <\vert s\vert \leq (a+a'+b+b')/2-1.
$$ C'est une g\'en\'eralisation de la notion de segments li\'es introduite par Zelevinsky, condition qui est naturelle quand on travaille avec les composantes locales des formes automorphes de carr\'e int\'egrables r\'esiduelles des groupes $GL$.

Et la propri\'et\'e est que l'op\'erateur d'entrelacement:
$$
N_{L}((\rho,a,b),(\rho',a',b'),s): \pi(\rho,a,b)\vert\,\vert^s\times \pi(\rho',a',b') \rightarrow \pi(\rho',a',b')\times \pi(\rho,a,b)\vert\,\vert^s
$$
est holomorphe en tout point $s$ r\'eel  tel que $\pi(\rho,a,b)\vert\,\vert^s$ et $\pi(\rho',a',b')$ ne sont pas li\'es. 

On admet que l'on sait d\'efinir par r\'ecurrence une action $\theta_{L}(\psi')$ sur $\pi(\psi')$ et on pose $$\theta_{L}(\psi):=\bigl((\theta(\rho,a,b)\otimes \theta_{L}(\psi'))\circ inv\bigr) \circ N_{L}(0).$$
La propri\'et\'e de multiplicativit\'e des op\'erateurs d'entrelacement normalis\'es montre que $\theta_{L}(\psi)$ est ind\'ependant du choix de $(\rho,a,b)\in Jord(\psi)$.

On suppose ici que $b\geq 2$ et on  construit les morphismes suivants, les inclusions sont naturelles et les fl\`eches sont les op\'erateurs d'entrelacement normalis\'es \`a la Langlands-Shahidi:
$$
\pi(\rho,a,b)\vert\,\vert^s\times \pi(\psi') \hookrightarrow St(a,\rho)^{-(b-1)/2+s}\times \pi(\rho,a,b-1)\vert\,\vert^{1/2+s}\times \pi(\psi') \rightarrow\eqno(1)
$$
$$
St(a,\rho)\vert\,\vert^{-(b-1)/2+s} \times \pi(\psi')\times \pi(\rho,a,b-1)\vert\,\vert^{s+1/2} \hookrightarrow$$
$$
St(a,\rho)\vert\,\vert^{-(b-1)/2+s} \times \pi(\psi')\times \pi(\rho,a,b-2)\vert\,\vert^s \times St(a,\rho)\vert\,\vert^{(b-1)/2+s}\eqno(2)$$
$$
\rightarrow St(a,\rho)\vert\,\vert^{-(b-1)/2+s} \times \pi(\rho,a,b-2)\vert\,\vert^s \times \pi(\psi') \times St(a,\rho)\vert\,\vert^{(b-1)/2+s}.\eqno(3)
$$
L'op\'erateur d'entrelacement d\'efini entre (2) et (3) est holormophe en $s=0$ car pour tout $(\rho,a',b')\in Jord(\psi')$, $\vert (a-a')/2\vert +\vert (b-1-b')/2\vert \geq 1/2$ sauf a priori si $a=a'$ et $b'=b+1$ mais la diff\'erence de parit\'e suffit alors \`a l'holomorphie. L'op\'erateur qui passe de (2) \`a (3) est aussi holomorphe en $s=0$; on construit donc une application avec des op\'erateurs d'entrelacement  normalis\'es 
 $$n_{L}(0):
\pi(\rho,a,b)\times \pi(\psi') \rightarrow St(a,\rho)^{-(b-1)/2} \times \pi(\rho,a,b-2)\times \pi(\psi') \times St(a,\rho)^{(b-1)/2};
$$
cette application est certainement non nulle car tous les op\'erateurs d'entrelacement utilis\'es sont inversibles (condition de non liaison qui ne d\'epend que de $\vert s\vert$). 

\

\bf Lemme. \sl L'op\'erateur $n_{L}(0)$ entrelace $\theta_{L}(\psi)$ avec le prolongement canonique de $\theta_{L}(\psi',(\rho,a,b-2))$.
\rm

\

En tant qu'op\'erateur m\'eromorphe $n_{L}(s)$ co\"{\i}ncide avec l'op\'erateur d'entrelacement normalis\'e \`a la Langlands-Shahidi:
 $$
 \pi(\rho,a,b)\vert\,\vert^{s} \times \pi(\psi')\hookrightarrow St(a,\rho)\vert\,\vert^{(b-1)/2+s}\times \pi(\rho,a,b-2)\vert\,\vert^{s} \times St(a,\rho)\vert\,\vert^{(b-1)/2}\times \pi(\psi') \rightarrow\eqno(4)
 $$
 $$
 St(a,\rho)\vert\,\vert^{-(b-1)/2+s}\times \pi(\rho,a,b-2)\vert\,\vert^{s} \times \pi(\psi') \times St(a,\rho)\vert\,\vert^{(b-1)/2+s}\eqno(5)
 $$
 On note $\theta_{L}^+(\psi',(\rho,a,b-2))$ le prolongement canonique de $\theta_{L}(\psi',(\rho,a,b-2)$; il s'obtient comme la valeur en $s=0$ du compos\'e de l'op\'erateur d'entrelacement normalis\'e qui \'echange les 2 facteurs de  $\pi(\rho,a,b-2)\vert\,\vert^s\times \pi(\psi')$  compos\'e avec l'action naturelle de $\theta$ sur l'induite
$$
St(a,\rho)\vert\,\vert^{-(b-1)/2+s}\times \pi(\psi')\times \pi(\rho,a,b-2)\vert\,\vert^s\times St(a,\rho)\vert\,\vert^{(b-1)/2+s}\rightarrow$$
$$
St(a,\rho)\vert\,\vert^{-(b-1)/2-s}\times\pi(\rho,a,b-2)\vert\,\vert^s\times  \pi(\psi')\times St(a,\rho)\vert\,\vert^{(b-1)/2-s}.
$$
On note $\theta_{L}^+(\psi',(\rho,a,b-2),s)$ cet op\'erateur.

Pour appliquer d'abord $\theta_{L}(\psi)$, il faut r\'ealiser $\theta_{L}(\psi)$ comme la valeur en $s=0$ de l'action naturelle de $\theta: \pi(\rho,a,b)\vert\,\vert^s\times \pi(\psi')$ dans $\pi(\psi')\times \pi(\rho,a,b)^{-s}$ compos\'ee avec l'op\'erateur d'entrelacement normalis\'e de $\pi(\psi')\times \pi(\rho,a,b)^{-s}$ dans $\pi(\rho,a,b)\vert\,\vert^{-s}\times \pi(\psi')$; op\'erateur que l'on note $\theta_{L}(\psi,s)$. Ensuite on applique $n_{L}(-s)$ et on montre alors facilement l'\'egalit\'e d'op\'erateurs m\'eromorphes:
$$
\theta_{L}^+(\psi',(\rho,a,b-2),s)\circ n_{L}(s)=n_{L}(-s) \circ \theta_{L}(\psi,s).
$$
En faisant $s=0$, on obtient l'assertion.

\

\bf Corollaire. \sl L'inclusion de Langlands de \ref{explicitationdelanormalisationdewhittaker} (1) entrelace $\theta_{L}(\psi)$ et le prolongement canonique de $\theta_{L}(\psi^2_{imp})$.\rm

\

C'est un corollaire imm\'ediat. Remarquons qu'il n'y a pas d'espoir d'avoir une \'egalit\'e entre $\theta_{L}(\psi)$ et $\theta_{W}(\psi)$ car $\theta_{L}(\psi)$ d\'epend de $E(s)$ (cf  ci-dessus) qu'on a fix\'e sans pr\'ecaution. Pr\'ecis\'ement, il semble que l'on puisse toujours tr\`es raisonablement remplacer $N_{L}(s)$ par $-N_{L}(s)$ ce qui \'evidemment change $\theta_{L}$ en son oppos\'e. On cherche \`a calculer $\theta_{W}(\psi)/\theta_{U}(\psi)$ et on va le r\'ealiser \`a l'aide de quotients de facteurs de normalisation o\`u $E(s)$ ne jouera plus de r\^ole; les seuls termes qui compteront seront les fonctions $L$ qui ne d\'ependent d'aucun choix.

\subsubsection{Action de $\theta$ et op\'erateurs d'entrelacement normalis\'es 2. \label{operateursdentrelacement2}}
Ce que l'on a fait ci-dessus avec la normalisation de Langlands-Shahidi peut se faire en rempla\c{c}ant les repr\'esentations de Steinberg par les repr\'esentations de Speh. On va le faire explicitement. Rappelons la valeur du facteur de normalisation pour l'entrelacement entre 2 blocs:
$$
N_{L}(s): \qquad \pi(\rho,a,b)^s\times \pi(\rho',a',b') \rightarrow \pi(\rho',a',b')\times \pi(\rho,a,b)\vert\,\vert^{s}.
$$
On n'a pas besoin de la valeur de $E(s)$ introduite dans \ref{operateursdentrelacement1}, pour $r_{L}(s)$ tel que $r_{L}(s)N_{L}(s)$ soit l'op\'erateur d'entrelacement standard:
$$
r_{L}(s)=E(s)\prod_{\begin{matrix}j\in [-(b-1)/2,(b-1)/2],\\j'\in [-(b'-1)/2,(b'-1)/2]\end{matrix}}L(\rho\times \rho',j+j'+\vert (a-a')/2\vert +s)/L(\rho\times \rho',j+j'+(a+a')/2+s).
$$
On pose, a priori, pour le m\^eme $E(s)$:
$$r_{S}(s):=E(s)\prod_{\begin{matrix}i\in [(a-1)/2,-(a-1)/2],\\i'\in [(a'-1)/2,-(a'-1)/2]\end{matrix}}L(\rho\times \rho',i+i'+\vert (b-b')/2\vert -s)/L(\rho\times \rho',i+i'+(b+b')/2-s).
$$
Cela permet de poser $N_{S}(s):=r_{S}(s)^{-1}M(s)$, o\`u $M(s)$ est encore l'op\'erateur d'entrelacement normalis\'e. Et on g\'en\'eralise \`a plusieurs blocs en composant les op\'erateurs.

\

\bf Remarque. \sl Pour $\rho$ et $\rho'$ des repr\'esentations cuspidales autoduales irr\'eductibles, le quotient $r_{L}(s)/r_{S}(s)$ est holormophe en $s=0$; en ce point il vaut 1 si $\rho\not\simeq \rho'$  et vaut le signe $(-1)^{inf(a,a')inf(b,b')}$ si $\rho\simeq \rho'$. Plus g\'en\'eralement ce quotient est holomorphe en tout point $s$ tel que les repr\'esentations $\pi(\rho,a,b)\vert\,\vert^s$ et $\pi(\rho',a',b')$ ne sont pas li\'ees. \rm

\

On pose $a_{+}=sup(a,a')$, $a_{-}=inf(a,a')$, $b_{+}=sup(b,b')$, $b_{-}=inf(b,b')$. On r\'ecrit
$$
\prod_{j\in [-(b-1)/2,(b-1)/2],j'\in [-(b'-1)/2,(b'-1)/2]}\frac{L(\rho\times \rho',j+j'+\vert (a-a')/2\vert +s)}{L(\rho\times \rho',j+j'+(a+a')/2+s)}=$$
$$
\prod_{\begin{array}{l}k'\in [-(b_{+}-1)/2,(b_{+}-1)/2]\\k\in [-(b_{-}-1)/2,(b_{-}-1)/2]; \\ \ell \in [(a_{-}-1)/2,-(a_{-}-1)/2]\end{array}}\frac{L(\rho\times \rho',k+k'+(a_{+}-1)/2+\ell+s)}{L(\rho\times \rho',k+k'+(a_{+}-1)/2+\ell+1+s)}=
$$
$$
=\prod_{\begin{array}{l}k\in [(b_{-}-1)/2,-(b_{-}-1)/2] , \\\ell \in [(a_{-}-1)/2,-(a_{-}-1)/2]\end{array}}\frac{L(\rho\times \rho',k+\ell+(a_{+}-1)/2-(b_{+}-1)/2+s)}{L(\rho\times \rho',k+\ell+(a_{+}-1)/2+(b_{+}-1)/2+1+s)}.
$$
De fa\c{c}on analogue, on obtient:
$$
\prod_{\begin{array}{l}i\in [(a-1)/2,-(a-1)/2]\\i'\in [(a'-1)/2,-(a'-1)/2]\end{array}}\frac{L(\rho\times \rho',i+i'+\vert (b-b')/2\vert -s)}{L(\rho\times \rho',i+i'+(b+b')/2-s)}=
$$
$$\prod_{\begin{array}{l}k\in [(b_{-}-1)/2,-(b_{-}-1)/2] \\\ell \in [(a_{-}-1)/2,-(a_{-}-1)/2]\end{array}}\frac{ L(\rho\times \rho',k+\ell+(b_{+}-1)/2-(a_{+}-1)/2-s)}{L(\rho\times \rho',k+\ell+(a_{+}-1)/2+(b_{+}-1)/2+1-s)}.
$$
Le seul cas o\`u les fonctions $L$ ne valent pas identiquement 1 et celui  o\`u  $\rho'=\chi \otimes \rho$ avec $\chi$ est un caract\`ere non ramifi\' e de $F^*$, y compris, \'evidemment le carat\`ere trivial.  En particulier $\rho$ et $\rho'$ sont des repr\'esentations d'un m\^eme groupe $GL(d,F)$.

Pour ce cas, on revient \`a la formule explicite pour $L(\rho\times \rho',s)$ donn\'ee par Shahidi \ref{shahidi}; le produit porte sur tous les caract\`eres non ramifi\'es $\tau$ de $F^*$ vus comme caract\`eres de $GL(d,F)$ via le d\'eterminant tels que $\rho\simeq \tau\otimes \rho'$ et $\omega$ est une uniformisante:
$$
L(\rho\times \rho',s)=\prod_{\tau}(1-\tau(\omega)q^{-s})^{-1}.
$$
On note $M$ l'ordre du groupe des carac\`eres non ramifi\'es, $\tau'$ tel que $\rho\simeq \tau'\otimes \rho$ et on a encore:
$
L(\rho\times \rho',s)=(1-\zeta_{\rho,\rho'} q^{-Ms})$ o\`u $\zeta_{\rho,\rho'}=+1$ si $\rho=\rho'$ et $-1$ sinon.

Or pour tout $k,\ell$ $$L(\rho\times \rho',k+\ell+(b_{+}-1)/2-(a_{+}-1)/2-s)=-\zeta_{\rho,\rho'}q^{Ms}\chi(L(\rho\times \rho',-k-\ell+(a_{+}-1)/2-(b_{+}-1)/2+s)$$ pour $M$ convenable (ne d\'ependant d'ailleurs que de $\rho$). On remarque que le produit porte sur des ensembles sym\'etriques par rapport \`a $0$. Ainsi, pour $m$ convenable  $$r_{L}(s)/r_{S}(s)=(-\zeta_{\rho,\rho'})^{a_{-}b_{-}}q^{ms}\prod_{k,\ell}L(\rho\times \rho',k+\ell+(a_{+}+b_{+})/2-s)/L(\rho\times \rho',k+\ell+(a_{+}+b_{+})/2+s),
$$
o\`u $k,\ell$ varient comme ci-dessus. Il suffit de d\'emontrer qu'en $s=0$ et plus g\'en\'eralement en tout point $s$ satisfaisant:
$$
\vert s\vert \leq (a_{+}-a_{-})/2+(b_{+}-b_{-})/2 \hbox{ ou } (a_{+}+a_{-}+b_{+}+b_{-})/2-1 <\vert s\vert
$$
num\'erateur et d\'enominateur sont holomorphes. Or on a, pour tout $k,\ell$ comme ci-dessus:
$$
k+\ell+(a_{+}+b_{+})/2\leq (b_{-}-1)/2+(a_{-}-1)/2+(a_{+}+b_{+})/2=(a_{+}+a_{-}+b_{+}+b_{-})/2-1
$$
$$
 k+\ell+(a_{+}+b_{+})/2 \geq -(b_{-}-1)/2-(a_{-}-1)/2 +(a_{+}-1)/2+(b_{+}+1)/2+1=$$
 $$(a_{+}-a_{-})/2+(b_{+}-b_{-})/2+1>(a_{+}-a_{-})/2+(b_{+}-b_{-})/2.
 $$
 Cela prouve l'holomorphie pour les valeurs de $s$ fix\'ees, puisque $L(\rho\times \rho',s')$ est holomorphe non nul en tout $s'$ r\'eel tel que $s'\neq 0$.
 
 \
 
 \bf Corollaire. \sl Les op\'erateurs $N_{S}(s)$ ont les m\^emes propri\'et\'es de multiplicativit\'e et d'holomorphie que les op\'erateurs $N_{L}(s)$  quand on remplace Steinberg par Speh .\rm
 
 \
 
 La premi\`ere assertion r\'esulte de la construction m\^eme et la deuxi\`eme du lemme ci-dessus.
 
 \
 
 On d\'efinit alors $\theta_{S}(\psi)$ comme on a d\'efini $\theta_{L}(\psi)$ en rempla\c{c}ant $N_{L}(s)$ par $N_{S}(s)$ et on a le corollaire:

 \
 
 \bf Corollaire. \sl L'inclusion de $\pi(\psi)$ dans \ref{definitiondelanormalisationunipotente} (1) entrelace $\theta_{S}(\psi)$ avec le prolongement canonique de $\theta_{S}(\psi^1_{imp})$.\rm
\subsubsection{Calcul de $\theta_{W}(\psi)/\theta_{U}(\psi)$.\label{calculdusigne}}
On a d\'efini les facteurs de normalisations $r_{L}(s)$, $ r_{S}(s)$; ils d\'ependent de $\psi$ et ici on met $\psi$ dans la d\'efinition. Et on pose $a(\psi):=\biggl( r_{L}(\psi)(s)/r_{S}(\psi)(s)\biggr)_{s=0}$. Pour \'eviter toute ambigu\"{\i}t\'e, on pr\'ecise que quand on somme sur les couples de triplets $((\rho,a,b),(\rho',a',b'))\in Jord(\psi)$, on regarde ces couples avec leur ordre, c'est-\`a-dire qu'un tel couple est diff\'erent de $((\rho',a',b'),(\rho,a,b))$ mais on supposera toujours $(\rho,a,b)$ est diff\'erent de $(\rho,a',b')$;  c'est pour cela que beaucoup de somme sont affect\'ees du coefficient 1/2. Il y a encore une autre ambigu\"{\i}t\'e qui vient du fait que $Jord(\psi)$ peut avoir de la multiplicit\'e; il faut en tenir compte, c'est-\`a-dire qu'un couple $((\rho,a,b),(\rho',a',b'))$ vient dans les formules avec la multiplicit\'e produit de la multiplicit\'e de $(\rho,a,b)$ dans $Jord(\psi)$ avec la multiplicit\'e de $(\rho',a',b')$ dans $Jord(\psi)$. 

\

\bf Th\'eor\`eme. \sl L'action $\theta_{W}(\psi)$ diff\`ere de $\theta_{U}(\psi)$ par le signe $a(\psi)a(\psi^1_{imp})a(\psi^2_{imp})a(\psi_{imp,imp})$. Explicitement ce signe vaut:
$(-1)^{1/2\sum_{(\rho,a,b),(\rho,a',b')\in Jord(\psi)}inf(a,a')(1+sup(a,a'))inf(b,b')(1+sup(b,b'))}$.\rm

\

On reprend \ref{explicitationdelanormalisationdewhittaker} (5); pour des bonnes repr\'esentations $\sigma,\tau$, o\`u $\sigma$ est un produit de repr\'esentations de Steinberg tordues par des caract\`eres n\'egatifs  et $\tau$ un produit de cuspidales tordues par des caract\`eres positifs et que l'on voit donc comme des repr\'esentations de Speh tordues par des caract\`eres positifs,  on a des inclusions:
$$
\pi(\psi) \hookrightarrow \sigma \times \pi(\psi^2_{imp}) \times \, ^\theta \sigma \hookrightarrow
\sigma\times \tau \times \pi(\psi_{imp,imp}) \times \, ^\theta\tau \times \, ^\theta\sigma.
$$
L'action $\theta_{W}(\psi)$ est la restriction \`a $\psi$ du prolongement canonique de l'action de $\theta$ sur $\pi(\psi_{imp,imp})$ qui agit naturellement sur l'espace des fonctionnelles de Whittaker pour le caract\`ere additif fix\'e. On note $\theta_{f}$ cette action de $\theta$ sur $\pi(\psi_{imp,imp})$.

La derni\`ere inclusion entrelace le prolongement canonique de $\theta_{S}(\psi^2_{imp})$ avec le prolongement canonique de $\theta_{S}(\psi_{imp,imp})$ tandis que la premi\`ere inclusion entrelace le prolongement canonique de $\theta_{L}(\psi^2_{imp})$ avec $\theta_{L}(\psi)$.  Posons $c_{2}:=\theta_{f}/\theta_{S}(\psi_{imp,imp})$. Ainsi le prolongement canonique de $\theta_{f}$ se restreint en $c_{2}$ fois le prolongement canonique de $\theta_{S}(\psi_{imp}^2)$ sur l'induite du milieu et donc par d\'efinition de $a(\psi^2_{imp})$ en $c_{2}a(\psi^2_{imp})$ fois le prolongement canonique de $\theta_{L}(\psi^2_{imp})$. Ainsi $\theta_{W}(\psi)$ diff\`ere de $\theta_{L}(\psi)$ par le signe $c_{2}a(\psi^2_{imp})$. On proc\`ede de fa\c{c}on analogue avec les normalisations unipotentes. On pose $c_{1}$ le signe $\theta_{f}/\theta_{L}(\psi_{imp,imp})$ et on voit que $\theta_{U}(\psi)$ diff\`ere de $\theta_{S}(\psi)$ par le signe $c_{1}a(\psi^1_{imp})$. Par d\'efinition $c_{1}c_{2}=a(\psi_{imp,imp})$ et on obtient donc:
$$
\theta_{W}(\psi)/\theta_{U}(\psi)=c_{2}a(\psi^2_{imp})c_{1}a(\psi^1_{imp})\theta_{L}(\psi)/\theta_{S}(\psi)=a(\psi_{imp,imp})a(\psi ^1_{imp})a(\psi^2_{imp})a(\psi).
$$
C'est la premi\`ere assertion du th\'eor\`eme. Passons au calcul explicite. On rappelle que $$a(\psi)=(-1)^{\sum_{(\rho,a,b),(\rho,a',b')\in Jord(\psi)}inf(a,a')inf(b,b')}.$$ On a un calcul analogue pour les autres $a(?)$; pour $a(\psi^1_{imp})$ il faut consid\'erer $Jord(\psi^1_{imp})$ c'est-\`a-dire uniquement les couples $(\rho,a,b),(\rho,a',b')$ tel que $inf(a,a')$ et $sup(a,a')$ sont impairs et la contribution d'un tel couple est $(-1)^{inf(b,b')}$. Pour consid\'erer tous les couples, il suffit de prendre comme contribution $(-1)^{inf(a,a')sup(a,a')inf(b,b')}$. On proc\'ede de fa\c{c}on  analogue pour $\psi^2_{imp}$; $\psi_{imp,imp}$ la contribution ne vient que des couples $(\rho,a,b),(\rho,a',b')$ tel que $a,a',b,b'$ sont impairs et la contribution est $(-1)$, on peut donc l'\'ecrire, en consid\'erant tous les couples, sous la forme $(-1)^{inf(a,a')sup(a,a')inf(b,b')sup(b,b')}$, on trouve le signe $(-1)^x$ o\`u 
$$
x=\sum_{(\rho,a,b),(\rho,a',b')\in Jord(\psi)}
inf(a,a')inf(b,b')+inf(a,a')inf(b,b')sup(b,b')+$$
$$+inf(b,b')inf(a,a')sup(b,b')+inf(a,a')sup(a,a')inf(b,b')sup(b,b')=$$
$$=\sum_{(\rho,a,b),(\rho,a',b')\in Jord(\psi)}inf(a,a')(1+sup(a,a'))inf(b,b')(1+sup(b,b')).
$$
C'est le r\'esultat annonc\'e.
\subsection{Comparaison avec notre normalisation de l'action de $\theta$\label{comparaisonavecnotrenormalisation}}
Dans ce travail, nous avons d\'efini une normalisation adhoc pour l'action de $\theta$. Elle s'exprime comme un m\'elange de $\theta_{W}$ et $\theta_{U}$; au lieu de plonger $W_{F}$ dans la premi\`ere o\`u la deuxi\`eme copie de $SL(2,{\mathbb C})$, on choisit la copie en fonction  de $(\rho,a,b)$ et plus pr\'ecis\'ement en fonction du signe de $a-b$.  L'ordre dans lequel on effectue les op\'erations est important pour les $\psi$ \'el\'ementaires m\^eme de restriction discr\`ete \`a la diagonale (c'est -\`a-dire tel que pour tout $(\rho,a,b)\in Jord(\psi)$, $inf(a,b)=1$). Ensuite, nous sommes pass\'es des $\psi$ \'el\'ementaires aux $\psi$ de restriction discr\`ete \`a la diagonale sans rajouter de choix . Le choix redevient crucial quand on traite le cas des $\psi$ g\'en\'eraux. Il y a toutefois une astuce simple pour ne pas recommencer les calculs.

On pose ${\cal Z}(\psi)$ l'ensemble des couples ordonn\'es $((\rho,a,b),(\rho,a',b'))\in Jord(\psi)$ tels que  $sup(b,b')$ et $sup(a,a')$ sont pairs tandis que  $inf(b,b')$ et  $inf(a,a')$ sont impairs; on remarque imm\'ediatement que ces conditions excluent la possibilit\'e $(\rho,a,b)=(\rho',a',b')$. On pose aussi $${Z\cal }_{W}(\psi):= \{(\rho,a,b),(\rho,a',b')\in {\cal Z}(\psi); a+a'<b+b', \hbox{ ou } a-b=a'-b'\neq 0 \hbox{ et } \zeta(a+b)+\zeta'(a'+b')<0\},$$
o\`u $\zeta$ est le signe de $a-b$ et $\zeta'$ celui de $a'-b'$. On pose alors $${\cal Z}_{U}(\psi)={\cal Z}(\psi)\setminus {\cal Z}_{W}=$$
$$
\{(\rho,a,b),(\rho,a',b')\in {\cal Z}(\psi); a+a'>b+b', \hbox{ ou }a-b=b'-a'=0, $$
$$\hbox { ou } a-b=b'-a'\neq 0, \zeta (a+b)+\zeta' (a'+b') \geq 0\}.
$$
Ces 3 ensembles admettent donc une sym\'etrie (celle qui \'echange les 2 \'el\'ements des couples) qui n'a pas de points fixes. Les cardinaux sont donc pairs.

Comme nous en aurons aussi besoin, on peut r\'ecrire ces 2 ensembles en utilisant les quadruplets $(\rho,A,B,\zeta)$, $(\rho,A',B',\zeta')$ plut\^ot que les triplets et cela donne:
$$
{\cal Z}_{W}(\psi):=\{(\rho,A,B,\zeta),(\rho,A',B',\zeta'); \zeta B+\zeta' B'<0, \hbox{ ou } B=B'\neq 0, \zeta\zeta'=-1 , \zeta (A-A')<0\};
$$
$$
{\cal Z}_{U}(\psi):=\{(\rho,A,B,\zeta),(\rho,A',B',\zeta'); \zeta B+\zeta' B'>0, \hbox{ ou }B=B'=0, $$
$$\hbox{ ou } B=B'\neq 0, \zeta\zeta'=-1, \zeta (A-A')\geq 0\};
$$
en effet $a-b+a'-b'$ vaut 2 fois $\zeta B+\zeta' B$. Quand $\zeta B=-\zeta' B'$, on a soit $B=B'=0$ soit $BB'\neq 0$ et sous cette derni\`ere hypoth\`ese, n\'ecessairement $\zeta\zeta'=-1$. Ensuite, on \'ecrit $\zeta(a+b)+\zeta'(a'+b')= \zeta (2A+2)+\zeta'(2A'+2)=2(\zeta (A-A'))$ si $\zeta\zeta'=-1$. D'o\`u la r\'ecriture faite.

Pour $(\rho,a,b)\in Jord(\psi)$, on \'ecrit ${\cal Z}(\psi)_{(\rho,a,b)}$, ${\cal Z}_{W}(\psi)_{(\rho,a,b)}$, ${\cal Z}_{U}(\psi)_{(\rho,a,b)}$ l'ensemble des \'el\'ements $(\rho,a',b')\in Jord(\psi)$ tels que $((\rho,a,b),(\rho,a',b'))$ soit respectivement dans ${\cal Z}(\psi)$, ${\cal Z}_{W}(\psi)$ et ${\cal Z}_{U}(\psi)$. Il faut ici remarquer que les conditions de parit\'e qui d\'efinissent ${\cal Z}$, entra\^{\i}nent que si $((\rho,a,b),(\rho,a',b'))\in {\cal Z}$ alors $a'\neq a$ et $b'\neq b$. En particulier les ensembles que l'on vient de 
d\'efinir relativement \`a $(\rho,a,b)$ fix\'e sont des sous-ensembles de $Jord(\psi')$ o\`u $\psi'$ s'obtient \`a partir de $\psi$ en enlevant $(\rho,a,b)$.

\

\bf Th\'eor\`eme. \sl L'action de $\theta$ mise sur $\pi(\psi)$ diff\`ere de $\theta_{W}(\psi)$ par le signe $z_{W}(\psi):=(-1)^{1/2\vert {\cal Z}_{W}(\psi)\vert }$ et de $\theta_{U}(\psi)$ par le signe $z_{U}(\psi):=(-1)^{1/2\vert {\cal Z}_{U}(\psi)\vert }$. \rm

\

Tenant compte du th\'eor\`eme de \ref{calculdusigne}, il suffit de d\'emontrer l'une des assertions. On fait d'abord les remarques suivantes.
Soit $(\rho,a,b)\in Jord(\psi)$. On note $\psi'$ le morphisme qui se d\'eduit de $\psi$ en enlevant ce triplet.

\

Que ce soit pour les morphismes $\psi$ \'el\'ementaires ou pour ceux de restriction discr\`ete \`a la diagonale, on a d\'efini l'action  $\theta(\psi)$ sur $\pi(\psi)$, en pr\'ecisant un quadruplet $(\rho,A,B,\zeta)\in Jord(\psi)$ tel que, en notant $\psi'$ le morphisme qui se d\'eduit de $\psi$ en enlevant ce quadruplet, l'action $\theta(\psi)$ est le prolongement canonique, utilisant l'inclusion ci-dessous, de l'action $\theta(\psi',(\rho,A-1,B+1,\zeta)$ si $A>B$ et $\theta(\psi',(\rho,A-1,B-1,\zeta))$ si $A=B\neq 0$. Le premier cas est celui de $\psi$ de restriction discr\`ete \`a la diagonale o\`u l'on peut prendre n'importe quel quadruplet $(\rho,A,B,\zeta)$ du moment que $A>B$ et le deuxi\`eme cas est celui de $\psi$ \'el\'ementaire o\`u il y a la condition $B$ minimal parmi les $(\rho,A',B',\zeta')$ v\'erifiant $B'\neq 0$ et si $Jord(\psi)$ contient \`a la fois $(\rho,A,B,\zeta=+)$ et $(\rho,A,B,\zeta=-)$ on prend $(\rho,A,B,\zeta=-)$. L'inclusion est, si $A>B$:
$$
\pi(\psi)\hookrightarrow <\rho\vert\,\vert^{\zeta B}, \cdots, \rho\vert\,\vert^{-\zeta A}> \times \pi(\psi',(\rho,A-1,B+1,\zeta) \times <\rho\vert\,\vert^{\zeta A}, \cdots, \rho\vert\,\vert^{-\zeta B}>;\eqno(1)
$$
et si $\psi$ est \'el\'ementaire et $B>0$
$$
\pi(\psi)\hookrightarrow \rho\vert\,\vert^{\zeta B}\times \pi(\psi',(\rho,A-1=B-1,B-1,\zeta))\times \rho\vert\,\vert^{-\zeta B}.\eqno(2)
$$Remarquons que dans (1), si $\zeta=+$, $<\rho\vert\,\vert^{\zeta B}, \cdots, \rho\vert\,\vert^{-\zeta A}>$ n'est autre que $St(a)\vert\,\vert^{-(b-1)/2}$ et si $\zeta=-$ cette repr\'esentation est isomorphe \`a $Sp(b)\vert\,\vert^{(a-1)/2}$. D'apr\`es  \ref{uneproprietecommuneacesactions}, si $\zeta=+$ (resp. $\zeta=-$), (1) entrelace $\theta_{W}(\psi)$ (resp. $\theta_{U}(\psi)$) et le prolongement canonique de $\theta_{W}(\psi',(\rho,A-1,B+1,\zeta))$ (resp. $\theta_{U}(\psi',(\rho,A-1,B+1,\zeta))$). Pour (2), si $\zeta=+$ on \'ecrit $\rho\vert\,\vert^{\zeta B}$ sous la forme $\rho\vert\,\vert^{(a-1)/2}$ et (2) entrelace $\theta_{U}(\psi)$ avec $\theta_{U}(\psi',(\rho,A-1,B-1,\zeta))$ et si $\zeta=-$, $\rho\vert\,\vert^{\zeta B}=\rho\vert\,\vert^{-(b-1)/2}$ et l'entrelacement est analogue en rempla\c{c}ant $U$ par $W$.
On aura donc d\'emontr\'e le th\'eor\`eme par r\'ecurrence pour les $\psi$ \'el\'ementaires et les $\psi$ de restriction discr\`ete \`a la diagonale si l'on peut d\'emontrer les propri\'et\'es suivantes:

(3) on suppose que $\psi$ est de restriction discr\`ete \`a la diagonale. Si $A>B$, $
{\cal Z}_{?}(\psi)_{(\rho,A,B,\zeta)}={\cal Z}_{?}(\psi',(\rho,A-1,B+1,\zeta))_{(\rho,A-1,B+1,\zeta)}$,   o\`u $?=W$ (resp $=U$) si $\zeta=+$ (resp. $-$) et en particulier est vide si $A=B+1$ (l'ensemble de droite est vide si $A=B+1$);

(4) on suppose que  $\psi$ est \'el\'ementaire et $B\geq 1/2$ (avec l'hypoth\`ese de minimalit\'e d\'ej\`a \'ecrite) ${\cal Z}_{?}(\psi)_{(\rho,A,B,\zeta)}={\cal Z}_{?}(\psi',(\rho,A-1,B-1,\zeta))_{(\rho,A-1,B-1,\zeta)}$,   o\`u $?=U$ (resp $=W$) si $\zeta=+$ (resp. $-$) (l'ensemble de droite est vide si $B=1/2$).

\noindent Plus g\'en\'eralement, car on en aura besoin pour le cas g\'en\'eral, on fixe $(\rho,A,B,\zeta)$ ou plut\^ot $(\rho,a,b)$; on suppose d'abord que $b\geq 2$; on remarque que $(\rho,a,b-2)$ correspond \`a $(\rho,A-1,B+1,\zeta)$ si $\zeta=+$ et \`a $(\rho,A-1,B-1,\zeta)$ si $\zeta=-$. Et on montre:

${\cal Z}_{W}(\psi)_{(\rho,a,b)}$ est vide si $b=2$ et  contient ${\cal Z}_{W}(\psi', (\rho,a,b-2))_{(\rho,a,b-2)}$ sinon, avec \'egalit\'e sauf s'il existe $(\rho,A',B',\zeta')\in Jord(\psi')$ tel que

soit $B'\neq 0$ et $\zeta B+\zeta' B'=-1$ et $\zeta' (A' -A+1)\geq 0$,

soit $B'=B\neq 0$, $\zeta\zeta'=-$ et $\zeta(A-A')<0$.

Supposons d'abord que $b=2$; comme ${\cal Z}_{W}(\psi)_{(\rho,a,b)}$ est inclus dans ${\cal Z}(\psi)_{(\rho,a,b)}$ les \'el\'ements $(\rho,a',b')$ de ${\cal Z}_{W}(\psi)_{(\rho,a,b)}$ v\'erifie certainement $inf(b,b')$ est impair et ne peut donc valoir $b$. Ainsi $b=sup(b,b')$ et $b'=inf(b,b')$ vaut n\'ecessairement 1. D'apr\`es les propri\'et\'es de ${\cal Z}_{W}$, on a $a+a'\leq b+b'=3$ et les conditions de parit\'e entra\^{\i}nent encore $sup(a,a')=2$, $inf(a,a')=1$. D'o\`u $a+a'=b+b'$ et il faut donc $\zeta\zeta'=-1$ $a-b=b'-a'\neq 0$. D'o\`u encore $a=1$ et $a'=2$ mais alors $\zeta (a+b)+\zeta'(a'+b')$ vaut 0 et n'est donc pas $<0$. D'o\`u la premi\`ere assertion. 

Supposons maintenant que $b>2$. Il est clair que tout \'el\'ement $(\rho,a',b')\in {\cal Z}_{W}(\psi',(\rho,a,b-2))_{(\rho,a,b-2)}$ v\'erifie $a+a'\leq b-2+b'<b+b'$ et est donc un \'el\'ement de ${\cal Z}_{W}(\psi)_{(\rho,a,b)}$. 
R\'eciproquement soit $(\rho,a',b')\in {\cal Z}_{W}(\psi)_{(\rho,a,b)}$. Si $a+a'<b+b' -2$, il est aussi dans ${\cal Z}_{W}(\psi', (\rho,a,b-2))_{(\rho,a,b-2)}$. Supposons maintenant que $a+a'=b+b'-2$; un tel \'el\'ement n'est pas dans ${\cal Z}_{W}(\psi',(\rho,a,b-2)_{(\rho,a,b-2)}$ si 

soit $B'=0$ d'o\`u n\'ecessairement $\zeta B=-1$, mais alors $a'=b'$, $b=a+1$ et les conditions de parit\'e sur $inf(a,a'),inf(b,b')$ et $sup(a,a'),sup(b,b')$ ne peuvent \^etre satisfaites car 3 de ces nombres ont m\^eme parit\'e;

soit $B'\neq 0$ d'o\`u n\'ecessairement $(a-(b-2))$ est du signe $-\zeta'$ et, pour que cette \'el\'ement ne soit pas dans ${\cal Z}_{W}(\psi',(\rho,a,b-2))_{(\rho,a,b-1)}$,  $-\zeta'(a+b-2)+\zeta' (a'+b')\geq 0$. Cette derni\`ere condition se r\'ecrit $-\zeta' (A-1)+\zeta' A'\geq 0$ ou encore $\zeta' A' \geq \zeta' (A-1)$ et la condition $a+a'=b+b'-2$ est $\zeta B +\zeta' B'=-1$.

Il reste \`a voir le cas o\`u $a+a'=b+b'$ et alors n\'ecessairement $B=B'\neq 0$, $\zeta\zeta'=-$ et $\zeta (A-A')<0$. Ce sont les conditions de l'assertion.

Utilisons les tout de suite pour d\'emontrer (3) sous l'hypoth\`ese $\zeta=+$: l'hypoth\`ese $b\geq 2$ est alors \'equivalente \`a $A\geq B+1$ les \'egalit\'es se correspondant. Ainsi, si $A=B+1$, on a vu que les 2 ensembles sont vides. Supposons  que $A>B+1$. Il n'existe certainement pas $(\rho,A',B',\zeta')\in Jord(\psi')$ avec  $B'=B$; supposons que l'on ait $\zeta B+\zeta' B'=-1$ et $\zeta' A'\geq \zeta'(A-1)$; comme $\zeta=+$ par hypoth\`ese, $B+1=-\zeta' B'$ d'o\`u $\zeta'=-$ et $B'=B+1$, $A\geq A'+1$, en particulier $[B', A']\subset [B,A]$ ce qui est aussi impossible puisque dans (3), on a suppos\'e que $\psi$ est de restriction discr\`ete \`a la diagonale. Cela prouve (3) sous l'hypoth\`ese $\zeta=+$.

On d\'emontre aussi (4) sous les hypoth\`eses $\zeta=-$;  l'hypoth\`ese $\zeta=-$ et le fait que $\psi$ est suppos\'e  \'el\'ementaire dans l'\'enonc\'e de (4), donnent  $B=(b-1)/2$. Ainsi le fait que  $b\geq 2$ est ici \'equivalente\`a $B\geq 1/2$ et $b>2$ \`a $B\geq 1$. On a donc imm\'ediatement l'\'egalit\'e cherch\'ee si $B=1/2$ puisque l'ensemble de gauche est vide lui aussi. Supposons donc que $B\geq 1$. Comme $\psi$ est \'el\'ementaire si $(\rho,A',B',\zeta')\in Jord(\psi')$ v\'erifie $B'=B$, il v\'erifie aussi $A'=A$; il ne g\^ene donc pas. Supposons que $\zeta B+\zeta' B'=-1$, $B'\neq 0$ et $\zeta' A' \geq \zeta' (A-1)$; cela se traduit par $\zeta' B'=B-1$ d'o\`u puisque $B'\neq 0$ n\'ecessairement $\zeta'=+$ et  $B'=B-1$ ce qui contredit la minimalit\'e de $B$. Il reste encore le cas o\`u $B'=0$, $B=1$, c'est-\`a-dire $(a',b')=(1,1)$ et $(a,b)=(1,3)$ mais un tel couple n'est pas dans ${\cal Z}{\psi}$ pour des conditions de parit\'e.

\

Supposons que $a=2$, alors  ${\cal Z}_{U}(\psi)_{(\rho,a,b)}=\emptyset$ sauf si $b=2$ et $(\rho,1,1)\in Jord(\psi')$ ou $b=1$ et $(\rho,1,2)\in Jord(\psi')$. 

En effet, soit $(\rho,a',b')\in {\cal Z}_{U}(\psi)_{(\rho,a,b)}$. Les conditions de parit\'e pour $sup(a,a')$ et $inf(a,a')$ donnent $a=sup(a,a')=2$ et $a'=inf(a,a')=1$. Comme $a+a'\geq b+b'$, les conditions de parit\'e pour $sup(b,b')$ et $inf(b,b')$ donnent $sup(b,b')=2$ et $inf(b,b')=1$. D'o\`u $(\rho,a',b')$ vaut soit $(\rho,1,1)$ mais alors $b=2$ soit vaut $(\rho,1,2)$ mais $b=1$. Ce sont les cas annonc\'es. 

Supposons maintenant que $a>2$, et montrons que ${\cal Z}_{U}(\psi)_{(\rho,a,b)}$ contient ${\cal Z}_{U}(\psi',(\rho,a-2,b))_{(\rho,a-2,b)}$ avec \'egalit\'e sauf s'il existe $(\rho,A',B',\zeta')\in Jord(\psi')$ tel que 
 
soit $B=B'=0$ 

soit $B=B'\neq 0$, $\zeta\zeta'=-1$ et $\zeta(A-A')\geq 0$ 

soit $\zeta B+\zeta' B'=1$, $B'\neq 0$ et $\zeta' (A'+1-A)<0$.

\noindent
Soit $(\rho,a',b')\in {\cal Z}_{U}(\psi',(\rho,a-2,b))_{(\rho,a-2,b)}$; n\'ecessairement $b+b'\leq a-2+a'<a+a'$ et $(\rho,a',b')\in {\cal Z}_{U}(\psi)_{(\rho,a,b)}$. Supposons maintenant que $(\rho,a',b')\in {\cal Z}_{U}(\psi)_{(\rho,a,b)}$; si $b+b' <a-2+a'$, cet \'el\'ement est dans ${\cal Z}_{U}(\psi',(\rho,a-2,b))_{(\rho,a-2,b)}$. Supposons que $b+b'=a-2+a'$; cela se retraduit en $\zeta B+\zeta' B'=1$. L'\'el\'ement $(\rho,a',b')$ est dans ${\cal Z}_{U}(\psi',(\rho,a-2,b))$ sauf si $B'\neq 0$ et $\zeta''(a-2+b)+\zeta'(a'+b')<0$, o\`u $\zeta''$ est le signe de $a-2-b$ c'est-\`a-dire $-\zeta'$. Cette derni\`ere in\'egalit\'e est donc $\zeta'(A'+1-A)<0$.  Reste le cas $a+a'=b+b'$, c'est-\`a-dire $\zeta B+\zeta' B'=0$; mais on a aussi soit $B=B'=0$ soit $B=B'\neq 0$ d'o\`u $\zeta=-\zeta'$ et $\zeta (A-A') \geq 0$. Ce sont les conditions annonc\'ees.

Utilisons-les pour compl\'eter la d\'emonstration de (3). Dans l'\'enonc\'e, on a comme hypoth\`ese que $inf(a,b)\geq 2$ puisque $A>B$ et que $\psi$ est de restriction discr\`ete \`a la diagonale. Comme on a d\'ej\`a r\'egl\'e le cas o\`u $\zeta=+$, on suppose que  $\zeta=-$ et .

Si $a=2,b=2$ on ne peut avoir aussi $(\rho,1,1)\in Jord(\psi)$ car $\vert 2-2\vert=\vert 1-1\vert$ et les segments $[B,A]$ et $[B',A']$ contiennent tous les deux 0.

Si $a=2,b=1$ on ne peut avoir $a'=1,b'=2$ pour le m\^eme argument.

Or $a=2$ et \'equivalent \`a $A=B+1$ et ainsi si $a=2$ les deux ensembles de (3) sont vides et on a l'assertion. On suppose donc que $a>2$ et de fa\c{c}on \'equivalente que $A>B+1$. L'hypoth\`ese de restriction discr\`ete \`a la diagonale assure que si $(\rho,A',B',\zeta')\in Jord(\psi')$, $B'\neq B$; il faut donc voir si l'on peut avoir $\zeta B+\zeta' B'=1$, $B'\neq 0$ et $\zeta' (A'+1-A)<0$. Comme $\zeta=-$, on aurait $\zeta' B'=B+1$, en particulier $\zeta'=+$, d'o\`u $A'+1<A$; d'o\`u encore $[B',A']\subset [B,A]$ ce qui est exclu. D'o\`u (3).

Compl\'etons maintenant la d\'emonstration de (4); ici $\psi$ est \'el\'ementaire et $\zeta=+$. D'o\`u $B=(a-1)/2$. Supposons d'abord que $a=2$; n\'ecessairement $b=1$. Si   $(\rho,1,2)\in $  $Jord(\psi')$
 cela contredirait le choix minimal fait pour $(\rho,a,b)$; pr\'ecis\'ement quand $(\rho,c,1)$ et $(\rho,1,c)$ sont dans $Jord(\psi)$ on choisit $(\rho,1,c)$ et non $(\rho,c,1)$. Supposons maintenant que $a>2$ d'o\`u $B\geq 1$. Soit $(\rho,A',B',\zeta')\in Jord(\psi')$ et il faut v\'erifier qu'il n'est pas dans les exceptions. Puisque $B\geq 1$, on n'a certainement pas $B=B'=0$. 
 
 Consid\'erons le cas o\`u $B=B'\neq 0$, $\zeta\zeta'=-1$ et $\zeta(A-A')\geq 0$. Comme ci-dessus, on aurait $a'=1$ et $b'=a$ ce qui contredit le choix de $(\rho,a,b)$ comme on vient de l'expliquer.  
 
 Supposons maintenant que $\zeta B+\zeta' B'=1$, $B'\neq 0$; d'o\`u $\zeta' B'=B-1$ et la seule possibilit\'e \'etant donn\'ee l'hypoth\`ese de minimalit\'e sur $B$ est $B'=0$ ce qui a \'et\'e exclu. Cela termine la preuve de (4).

\

Nous avons donc d\'emontr\'e le th\'eor\`eme pour tous les morphismes $\psi$ \'el\'ementaires, gr\^ace \`a (4)  et ceux de restriction discr\`ete \`a la diagonale, gr\^ace \`a (3). Il reste \`a traiter le cas g\'en\'eral. Ici la proc\'edure est inverse; on part d'un $\tilde{\psi}$ de restriction discr\`ete \`a la diagonale et on calcule $\pi(\psi)$ comme module de  Jacquet \`a partir de $\pi(\tilde{\psi})$ et on a mis l'action de $\theta$ qui s'en d\'eduit. Donc on a construit une inclusion:
$$
\pi(\tilde{\psi})\hookrightarrow \sigma \times \pi(\psi) \times \, ^\theta\sigma,\eqno(5)
$$
pour $\sigma$ convenable et $\theta(\psi)$ est tel que $\theta(\tilde{\psi})$ est le prolongement canonique de $\theta(\psi)$. Donc ici, on a le r\'esultat pour $\theta(\tilde{\psi})$ et on veut l'obtenir pour $\theta(\psi)$. On le fait en plusieurs \'etapes et on se ram\`ene \`a la situation suivante:

il existe $(\rho,A,B,\zeta)\in Jord(\tilde{\psi})$ avec $B\geq 1$ tel que en notant $\psi'$ le morphisme qui se d\'eduit de $\tilde{\psi}$ en enlevant ce bloc, $\psi$ s'obtient \`a partir de $\psi'$ en ajoutant le bloc $(\rho,A-1,B-1,\zeta)$. Et $(\rho,A,B,\zeta)$ a les propri\'et\'es suivantes: pour tout $(\rho,A',B',\zeta')\in Jord(\psi')$ soit $B'<B$ soit $B'>>A$ et si $B'=B-1$ alors $A-1\geq A'$ et si $B'=B-1, A'=A-1$ mais $\zeta \zeta'=-$ alors $\zeta=+$.

Ces propri\'et\'es viennent des choix sur l'ordre que nous avons mis. Pour pouvoir utiliser tranquillement les calculs d\'ej\`a faits sans s'embrouiller dans les notations, ou oublie la notation $\psi$ que l'on remplace par $(\psi',(\rho,A-1,B-1,\zeta))$.

Dans ces conditions dans (5) $\sigma=<\rho\vert\,\vert^{\zeta B}, \cdots, \rho\vert\,\vert^{\zeta A}>$ c'est-\`a-dire si $(\rho,a,b)=(\rho,A,B,\zeta)$ si $\zeta=+$, $\sigma=Sp(b)\vert\,\vert^{(a-1)/2}$ tandis que si $\zeta=-$, $\sigma=St(a)\vert\,\vert^{-(b-1)/2}$ et si $\zeta=+$ $(\rho,A-1,B-1,\zeta)=(\rho,a-2,b)$ tandis que si $\zeta=-$, $(\rho,A-1,B-1,\zeta)=(\rho,a,b-2)$. Ainsi avec \ref{uneproprietecommuneacesactions}, (5) entrelace $\theta_{U}(\tilde{\psi})$ (resp. $\theta_{W}$) et le prolongement canonique de $\theta_{U}(\psi',(\rho,a-2,b)$ (resp. $\theta_{W}$) si $\zeta=+$ (resp. $\zeta=-$).

Supposons que $\zeta=+$. Il faut donc d\'emontrer que ${\cal Z}_{U}(\tilde{\psi})_{(\rho,a,b)}={\cal Z}_{U}(\psi',(\rho,a-2,b))_{(\rho,a-2,b)}$. 
On a suppos\'e que $B\geq 1$ or $(a-1)/2\geq B \geq 1$ d'o\`u aussi $a>2$. Il faut donc montrer que les exceptions ne peuvent se produire. Fixons $(\rho,A',B',\zeta')\in Jord(\psi')$; on suppose d'abord que $B=B'$; or dans nos hypoth\`eses, soit $B'<B$ soit $B'>>A$ ce qui est contradictoire. 

On suppose maintenant que  $\zeta B+\zeta' B'=1$, $B'\neq 0$ et $\zeta' (A'+1-A)<0$. Ceci se retraduit en $B-1=-\zeta' B'$, $B'\neq 0$, d'o\`u d\'ej\`a $B>1$ et $\zeta'=-$ et encore $A'>A-1$. Ceci a \'et\'e exclu puisque $B'=B-1$ entra\^{\i}ne $A'\leq A-1$.

Supposons maintenant que $\zeta=-$. On reprend mot pour mot le d\'ebut de la d\'emonstration ci-dessus. Il faut ici d\'emontrer que ${\cal Z}_{W}(\tilde{\psi})_{(\rho,a,b)}={\cal Z}_{W}(\psi',(\rho,a-2,b))_{(\rho,a-2,b)}$. 
On a suppos\'e que $B\geq 1$ or $(b-1)/2\geq B \geq 1$ d'o\`u aussi $b>2$. Il faut donc montrer que les exceptions ne peuvent se produire. Fixons $(\rho,A',B',\zeta')\in Jord(\psi')$;  supposons d'abord que
$B'=B\neq 0$, $\zeta\zeta'=-$ et $\zeta(A-A')<0$, ici c'est $B=B'$ qui est exclu. Supposons maintenant que  $\zeta' B'\neq 0$ et $\zeta B+\zeta' B'=-1$ et $\zeta' (A' -A+1)\geq 0$. Comme $\zeta=-$ cela donne $\zeta' B'=B-1$, d'o\`u encore $B>1$ et $\zeta'=+$, $B'=B-1$, $A'\geq A-1$. A priori on peut avoir $A'=A-1$ mais alors, il faudrait $\zeta=+$; les hypoth\`eses excluent donc aussi cette possibilit\'e.

Cela termine la preuve du th\'eor\`eme.
\subsection{Signe et transfert\label{signeettransfert}}
Les diff\'erents choix pour l'action de $\theta$ entra\^{\i}nent un changement de signe dans le calcul des traces. 
Pour unifier les notations on pose $z(\psi)=z_{\emptyset}(\psi)$ et ${\cal Z}(\psi)={\cal Z}_{\emptyset}(\psi)$.
Aux diff\'erents signes $z_{\emptyset}(\psi)$, $z_{W}(\psi)$, $z_{U}(\psi)$, on associe une application de $Jord(\psi)$ dans $\{\pm 1\}$, en posant simplement
$$
\epsilon_{?}(\rho,a,b):=(-1)^{\vert {\cal Z}_{?}(\psi)_{(\rho,a,b)}\vert}.
$$

On suppose  que  $\psi$ est \`a image dans $^LH_{n}$ le groupe dual de $H_{n}$ (cf \ref{lesgroupes}, une telle application s'interpr\`ete alors essentiellement comme un caract\`ere du groupe $Cent_{^LH_{n}}\psi$. Le essentiellement vient du fait que si $H_{n}=Sp(2n,F)$ l'ensemble des applications de $Jord(\psi)$ dans $\{\pm 1\}$ s'identifie aux caract\`eres du centralisateur de $\psi$ dans $O(2n+1,F)$. Dans tous les cas, on note $z$ l'\'el\'ement non trivial soit du centre de $^LH_{n}$ (si $H_{n}\neq Sp(2n,F)$)soit du  centre de $O(2n+1)$. Pour toute application de $Jord(\psi)$ dans $\{\pm 1\}$ identifi\'ee comme ci-dessus en un caract\`ere, on a: $$\epsilon(z)=\prod_{(\rho,a,b)\in Jord(\psi)}\epsilon(\rho,a,b)$$
o\`u $Jord(\psi)$ est  vu comme un ensemble avec multiplicit\'es.

On pose aussi $$c_{2,\psi}:=\psi(1_{W_{F}},1_{SL(2,{\mathbb C})}\bigl(\begin{matrix} -1&0\\ 0&-1\end{matrix}\bigr))$$,
o\`u $1_{W_{F}}$ et $1_{SL(2,{\mathbb C})}$ sont les unit\'es des groupes en indice. Et on a:

\bf Th\'eor\`eme. \sl Pour $?$ valant soit $\emptyset$, soit $W$ soit $U$,  $\epsilon_{?}(z)=1$   et $\epsilon_{?}(c_{2,\psi})=z_{?}(\psi)$.\rm

\

D'apr\`es ce qui pr\'ec\`ede l'\'enonc\'e:
$
\epsilon_{?}(z)=\prod_{(\rho,a,b)\in Jord(\psi)}\epsilon_{?}(\rho,a,b)$. Cela vaut pr\'ecis\'ement $(-1)^T$ o\`u $T $ est le cardinal de $\cup_{(\rho,a,b)\in Jord(\psi)}{\cal Z}_{?}(\psi)_{(\rho,a,b)}$; mais  cet ensemble s'identifie exactement avec  ${\cal Z}_{?}(\psi)$; un tel ensemble est de cardinal pair  (cf. \ref{{comparaisonavecnotrenormalisation}}).
D'o\`u la premi\`ere assertion de l'\'enonc\'e. Pour la deuxi\`eme,  on a de fa\c{c}on tr\`es g\'en\'erale:
$$
\epsilon_{?}(\psi)(c_{2,\psi})=\prod_{(\rho,a,b)\in Jord(\psi);(-1)^b=1}\epsilon_{?}(\psi)(\rho,a,b)=(-1)^{T'},
$$o\`u ici $T'$ est le cardinal de $\cup_{(\rho,a,b)\in Jord(\psi);(-1)^b=1}{\cal Z}_{?}(\psi)_{(\rho,a,b)}$.
Or  tout couple $(\rho,a,b),(\rho,a',b')\in{\cal Z}(\psi)$ est tel que soit $b$ soit $b'$ est pair et l'autre est impair. Ainsi l'ensemble cherch\'e est exactement $1/2\vert {\cal Z}_{?}(\psi)\vert $ et l'assertion est exactement la d\'efinition de $z_{?}(\psi)$. D'o\`u le th\'eor\`eme.

\

\bf Remarque. \sl On garde l'hypoth\`ese sur l'image de $\psi$. Supposons que l'on ait construit une application de l'ensemble des caract\`eres du centralisateur de $\psi$ dans l'ensemble des repr\'esentations de $H_{n}$ (en acceptant la repr\'esentation $0$), $\epsilon \mapsto \pi(\psi,\epsilon)$ de telle sorte que pour $\pi(\psi)$ muni de l'action $\theta(\psi)$, la trace tordue de $\pi(\psi)$ soit un transfert du caract\`ere $\sum_{\epsilon\in Cent_{^LH_{n}}^{\hat{\empty}}}\epsilon(c_{2,\psi})tr(\pi(\psi,\epsilon))$. Alors la trace tordue de  $\pi(\psi)$ muni de l'action $\theta_{W}(\psi)$ est un transfert de la distribution
$$
\sum_{\epsilon\in Cent_{^LH_{n}}^{\hat{\empty}}}\epsilon_{W}(c_{2,\psi})\epsilon(c_{2,\psi})tr(\pi(\psi,\epsilon)).
$$
Et on a une assertion analogue pour $\theta_{U}(\psi)$ et $\epsilon_{U}(c_{2,\psi})$.\rm

\

Il est alors int\'eressant de changer la param\'etrisation des $\pi(\psi,\epsilon)$ c'est \`a dire de poser $\pi_{W}(\psi,\epsilon)=\pi(\psi,\epsilon_{W}(\psi)\epsilon)$ d'o\`u la trace tordue de $\pi(\psi)$ avec l'action $\theta_{W}(\psi)$ est un transfert de
$$\sum_{\epsilon\in Cent_{^LH_{n}}^{\hat{\empty}}}\epsilon(c_{2,\psi})tr(\pi_{W}(\psi,\epsilon)),
$$
ce qui est une jolie formule. Et la question est en plus de savoir si cette param\'etrisation est compatible avec l'inclusion du paquet de Langlands; pr\'ecis\'ement, on a d\'ej\`a consid\'er\'e en \ref{definitiondelanormalisationdewhittaker} $$\psi^2: W_{F}\times SL(2,{\mathbb C}) \rightarrow GL(n,{\mathbb C}),
$$
obtenu en composant $\psi$ avec le morphisme de $W_{F}$ dans la 2e copie de $SL(2,{\mathbb C})$, $w\mapsto \bigl(\begin{matrix} \vert w\vert^{1/2}&0\\ 0& \vert w\vert^{-1/2}\end{matrix}\bigr)$. A $\psi^2$ correspond un paquet de Langlands, en g\'en\'eral conjecturalement. Notons $\Pi_{L}(\psi)$ cet ensemble de repr\'esentations; il est param\'etr\'e par les caract\`eres du centralisateur de $\psi^2_{imp}$ (voir \ref{explicitationdelanormalisationdewhittaker}) dans un Levi de $^LH_{n}$. Il n'est pas a priori clair que $\Pi_{L}(\psi)$ soit inclus dans l'ensemble $\{\pi(\psi,\epsilon)\}$ o\`u $\epsilon$ parcourt l'ensemble des caract\`eres du centralisateur de $\psi$. Ceci r\'esulte tr\`es vraisemblablement du transfert  de $\pi(\psi)$ (quand celui-ci est \'etabli) admettons le ici pour poser le probl\`eme. On doit encore remarquer que $Cent_{^LH_{n}}(\psi)$ s'envoie naturellement dans le $Cent(\psi^2_{imp})$ et la question est alors de savoir si pour $\eta$ un caract\`ere du centralisateur de $\psi^2_{imp}$ et $\epsilon$ un caract\`ere de $Cent_{^LH_{n}}(\psi)$ tel que $\pi_{W}(\psi,\epsilon)$ contienne la repr\'esentation de $\Pi_{L}(\psi)$ associ\'ee \`a $\eta$, $\epsilon$ serait la restriction de $\eta$ \`a $Cent_{^LH_{n}}(\psi)$.  On peut travailler de la m\^eme fa\c{c}on, en rempla\c{c}ant $_{W}$ par $_{U}$ et poser la m\^eme question en utilisant $\psi^1_{imp}$ de \ref{definitiondelanormalisationunipotente}. Sur les exemples calcul\'es les r\'eponses sont partout oui et on reviendra ult\'erieurement sur ces questions.

\begin{verbatim}moeglin@math.jussieu.fr; waldspur@math.jussieu.fr
\end{verbatim}

\end{document}